%% file: main_v2.tex
\documentclass[hidelinks, 12pt]{amsart}
\usepackage{tikz}
\usetikzlibrary{calc,3d}
\allowdisplaybreaks
\usepackage[english]{babel}
\usepackage[pdftex,textwidth=400pt,marginratio=1:1]{geometry}
\usepackage{amsfonts}
\usepackage{graphicx}
\usepackage{caption}
\usepackage{subcaption}
\usepackage{amsmath}
\usepackage{amsthm}
\usepackage{amssymb}
\usepackage{bbm}
\usepackage{cancel}
\usepackage{color}
\usepackage[all,cmtip]{xy}
\usepackage{graphicx}
\usepackage{mathabx}
\usepackage{tikz-cd}
\usepackage{epigraph} %For a quote on the first page
\usepackage{hyperref}
\usepackage{stmaryrd}
\usepackage{enumitem}
\usepackage[mathscr]{eucal}
\newtheorem{theorem}{Theorem}[section]
\newtheorem{corollary}[theorem]{Corollary}
\newtheorem{proposition}[theorem]{Proposition}
\newtheorem{lemma}[theorem]{Lemma}
\newtheorem{conjecture}[theorem]{Conjecture}

\theoremstyle{definition}
\newtheorem{definition}[theorem]{Definition}
\newtheorem{example}[theorem]{Example}
\newtheorem{remark}[theorem]{Remark}
\newtheorem{question}[theorem]{Question}

\DeclareFontFamily{OT1}{rsfs}{}
\DeclareFontShape{OT1}{rsfs}{n}{it}{<-> rsfs10}{}
\DeclareMathAlphabet{\curly}{OT1}{rsfs}{n}{it}

\newcommand\LL{\mathbb L}

\renewcommand\O{\mathcal O}
\newcommand\PP{\mathbb P}

\newcommand\cE{\mathcal E}
\newcommand\EE{\mathbb E}
\newcommand\DD{\mathbb D}

\newcommand\E{\mathbb E}
\newcommand\C{\mathbb C}

\newcommand\FF{\mathbb F}
\newcommand\GG{\mathbb G}

\newcommand\II{\mathbb I}

\newcommand\Q{\mathbb Q}

\newcommand\R{\mathbb R}

\newcommand\Z{\mathbb Z}
\newcommand\cZ{\mathcal Z}

\newcommand\Coh{\mathrm{Coh}}

\newcommand\sfV{\mathsf V}

\newcommand\pure{\mathrm{pure}}

\newcommand\Exp{\mathrm{Exp}}

\newcommand\vd{\mathrm{vd}}
\newcommand\pt{\mathrm{pt}}
\newcommand\vir{\mathrm{vir}}

\newcommand\loc{\mathrm{loc}}

\newcommand\DT{\mathrm{DT}}
\newcommand\PT{\mathrm{PT}}

\newcommand\td{\mathrm{td}}
\newcommand\rk{\operatorname{rk}}

\newcommand\tr{\operatorname{tr}}
\newcommand\coker{\operatorname{coker}}

\newcommand\ch{\operatorname{ch}}

\newcommand\Hom{\operatorname{Hom}}
\renewcommand\hom{\mathcal{H}{\it{om}}}
\newcommand\Ext{\operatorname{Ext}}
\newcommand\ext{\curly Ext}
\newcommand\At{\operatorname{At}}

\newcommand\Spec{\operatorname{Spec}\,}
\newcommand\Hilb{\operatorname{Hilb}}

\newcommand\Sym{\operatorname{Sym}}

\newcommand\mdot{{\scriptscriptstyle\bullet}}

\newcommand\INTO{\ar@{^{(}->}[r]}

\DeclareRobustCommand{\SkipTocEntry}[4]{}

\def\dual{^{\vee}}
\def\fI{\mathfrak{I}}

\def\Supp{\mathrm{Supp}}

\def\udot{^{\mdot}}
\def\spl{\mathrm{spl}}

\def\Pic{\mathrm{Pic}}

\def\cPair{{\curly Pair}}

\def\cPerf{{\curly Perf}}

\def\PTqvX{{\curly P}^{(q)}_v(X)}

\def\curP{\curly P}
\def\curL{\curly L}

\def\Xbar{\overline{X}}

\def\and{\quad\mathrm{and}\quad}

\newcommand{\Ohat}{\widehat{\mathcal{O}}^\mathrm{vir}}

\begin{document}
\title[Counting surfaces on Calabi-Yau 4-folds II]{Counting surfaces on Calabi-Yau 4-folds II: $\DT$--$\PT_0$ correspondence}
\author[Y.~Bae, M.~Kool, H.~Park]{Y.~Bae, M.~Kool, and H.~Park}
\maketitle
\begin{abstract}
This is the second part in a series of papers on counting surfaces on Calabi-Yau 4-folds. In this paper, we introduce $K$-theoretic $\DT, \PT_0, \PT_1$ invariants and conjecture a $\DT$--$\PT_0$ correspondence. For certain tautological insertions, we derive Lefschetz principles in both the compact and toric case allowing reductions to 3-dimensional $\DT, \PT$ invariants. We also develop a topological vertex and conjecture a $\DT$--$\PT_0$ vertex correspondence. These methods enable us to verify our conjectures in several examples. 
\end{abstract}
\tableofcontents
\addtocontents{toc}{\protect\setcounter{tocdepth}{1}}
\vspace{-1.5cm}

\section{Introduction}

\subsection{Overview}

This is the second part of a series of papers on counting surfaces on Calabi-Yau 4-folds \cite{BKP}.
Let $X$ be a smooth quasi-projective Calabi-Yau 4-fold, and let
\[v= (0,0,\gamma,\beta,n-\gamma \cdot \td_2(X))\in H^*_c(X,\Q)\]
be a fixed class in cohomology with compact support. We denote by  $\curly I_v(X)$  the Hilbert scheme of compactly supported closed subschemes $Z \subset X$ of dimension $\leq 2$ satisfying $\ch(\O_Z) = v$ (and hence $\chi(\O_Z)=n$). 
By \cite{OT}, we have a virtual cycle in the Chow group with $\Q$-coefficients
\begin{equation} \label{eqn:vc}
    [\curly I_v(X)]^\vir \in A_{\vd}(\curly I_v(X))_\Q\,, \quad \vd := n-\frac{1}{2}\gamma^2\,
\end{equation}
depending on a choice of orientation on $\curly I_v(X)$.
Integration against $[\curly I_v(X)]^\vir$ defines {\em Donaldson-Thomas} (DT) {\em invariants} of $X$.

$\DT$ invariants contain contributions from free moving points on $X$. In order to get better behaved invariants, one should normalize the generating series of $\DT$ invariants by the degree zero contribution, i.e., the contribution of Hilbert schemes of points on $X$ \cite{MNOP, PT1,CKM1}.
The aim of this project is to give a moduli-theoretic approach to the normalized $\DT$ generating series in the case of surfaces on Calabi-Yau 4-folds.
\begin{question}\label{question:intro}
    Is there a moduli-theoretic interpretation of the normalized generating series of $\DT$ invariants for surfaces on Calabi-Yau 4-folds?
\end{question}
The analogous question has been answered for curve counting on 3-folds via {\em stable} (PT) {\em pairs} \cite{PT1}. A direct extension of their idea leads to a moduli space of pairs $(F,s)$ on $X$, which we call $\PT_1$ {\em pairs}, where $F$ is a pure 2-dimensional coherent sheaf on $X$ with proper support and $s\in H^0(F)$ is a section with at most 1-dimensional cokernel. The $\PT_1$ pairs were independently studied in \cite{GJL}. However, $\PT_1$ pairs cannot be used to answer Question \ref{question:intro}, because $\PT_1$ pairs and ideal sheaves differ by both point and curve contributions.

In \cite{BKP} we introduce a new type of stability condition which sits between the $\DT$ and $\PT_1$ stability. 
\begin{definition}
     Let $F$ be a 2-dimensional coherent sheaf with proper support on $X$ and $s\in H^0(F)$.  The pair $(F,s)$ is called $\PT_0$ {\em stable} if $F$ has no non-trivial 0-dimensional subsheaf and the cokernel of $s$ is $0$-dimensional.
\end{definition}
$\PT_0$ stability is not covered by Le Potier's theory of coherent systems \cite{Pot1} since the sheaf $F$ is generally not pure. In \cite{BKP} we prove that the moduli space $\curly P^{(0)}_v(X)$ of $\PT_0$ pairs is fine, quasi-projective, and it is related to the Hilbert scheme $\curly I_v(X)$ by GIT wall-crossing.

We use a specific $K$-theoretic insertion for the $\DT$ generating series, which we now discuss (for details, see Section \ref{sec:Ktheorinv}). By Oh-Thomas \cite{OT}, there exists a twisted virtual structure sheaf in the Grothendieck group of coherent sheaves
\[\Ohat_{\curly I} \in K_0(\curly I_v(X))_\Q\]
which is related to \eqref{eqn:vc} via the Riemann-Roch isomorphism.
Let $\FF$ be the structure sheaf of the universal subscheme in $X\times \curly I_v(X)$ and $\pi_{\curly I}$ the projection to $\curly I_v(X)$.
For any line bundle $L$ on $X$ and a formal variable $y$, we work with the {\em Nekrasov genus}  \cite{Nek,NO} (also considered for curves in \cite{CKM1})
\begin{equation} \label{eqn:Ktheorinvintro} 
\langle\!\langle L\rangle\!\rangle^{\DT}_{X,\gamma,\beta,n} :=\chi(\curly I_v(X), \Ohat_{\curly I} \otimes \widehat{\Lambda}^\mdot( R\pi_{\curly I *}(\FF \otimes L) \otimes y^{-1}))\,,
\end{equation}
where we suppress pull-back of $L$ from $X$ to $X \times \curly I_v(X)$. The rank of $R\pi_{\curly I *}(\FF \otimes L)$ equals $\rk:=(\tfrac{1}{2} L \gamma + \beta)L + n$. We define the analogous invariant for $\curly P_v^{(0)}(X)$, in which case we replace the superscript $\DT$ by $\PT_0$. In some parts, we consider the Nekrasov genera for $\PT_1$ pairs, in which case we use the superscript $\PT_1$.

We propose a $\DT$--$\PT_0$ correspondence for surface counting:
\begin{conjecture} \label{conj:KDTPT0intro}
Let $X$ be a projective Calabi-Yau 4-fold and $L$ a line bundle on $X$. For any $\gamma \in H^4(X,\Q)$, $\beta \in H^6(X,\Q)$, there exist orientations such that 
\begin{equation} \label{eq:DTPT0intro} 
\frac{\sum_n \langle \! \langle L \rangle \! \rangle_{X,\gamma,\beta,n}^{\DT} q^n}{\sum_n \langle \! \langle L \rangle \! \rangle_{X,0,0,n}^{\DT} q^n} = \sum_n \langle \! \langle L \rangle \! \rangle_{X,\gamma,\beta,n}^{\PT_0} q^n.
\end{equation}
\end{conjecture}
A formula for the denominator can be found in \cite{CK1, Par, Boj2}. 
We also provide a quasi-projective and equivariant analog of this conjecture for toric Calabi-Yau 4-folds.
The curve counting $\DT$--$\PT$ correspondence for Calabi-Yau 4-folds is studied in \cite{CK2,CKM1, Liu}.
Unlike the curve case, the rank $\rk$ of the tautological complex generally does not match the degree $\vd$ of the virtual class. Remarkably, we claim that the insertion \eqref{eqn:Ktheorinvintro} remains the \emph{correct} one.
We prove Conjecture~\ref{conj:KDTPT0intro} for various Calabi-Yau 4-fold examples. 
\begin{enumerate}
    \item For toric Calabi-Yau $4$-folds, we develop a vertex formalism  and propose a vertex $\DT$--$\PT_0$ correspondence, which we can check in examples (up to certain orders in $q$) using virtual localization.
    \item For compact Calabi-Yau $4$-folds and $L = \O_X(Y)$, for a smooth connected effective divisor $Y$ and $\gamma$ a class on $Y$, we relate the $\DT$--$\PT_0$-correspondence on $X$ to the (curve) $\DT$--$\PT$ correspondence on $Y$. 
\end{enumerate}

\subsection{\texorpdfstring{$\PT_0$}{PT0} pairs on Calabi-Yau 4-folds}

A geometric understanding of $\PT_0$ pairs serves as a first step towards establishing the $\DT$--$\PT_0$ correspondence. Let $(F,s)$ be a $\PT_0$ pair with proper scheme theoretic support $Z$ on a quasi-projective Calabi-Yau 4-fold $X$. We give a description of $\PT_0$ pairs with increasing level of complexity: 
\begin{enumerate}
    \item If $\coker(s)$ is trivial, then $F\cong\O_Z$ where $Z$ has no 0-dimensional connected components or 0-dimensional embedded subschemes. In particular, there are no free moving points.
    \item If the maximal pure 2-dimensional subscheme $Z_1 \subset Z$ is Cohen-Macaulay, then $\coker(s)$ is supported on the $1$-dimensional subscheme $Z \setminus Z_1$.
    \item If $Z$ is pure 2-dimensional with non-Cohen-Macaulay singularities, then $\coker(s)$ is supported on these non-Cohen-Macaulay singularities.
\end{enumerate}

These cases have more explicit descriptions on toric Calabi-Yau 4-folds as we discuss below (Theorem~\ref{thm:fixedintro}).
Let $X$ be a smooth quasi-projective toric Calabi-Yau 4-fold. Let $T_X\cong (\C^*)^4$ be the 4-dimensional torus acting on $X$ and let $T:= \{t_1t_2t_3t_4=1\} \leq T_X$ be the \emph{Calabi-Yau torus} preserving the holomorphic volume form.
The action of $T_X$ lifts to the moduli space $\curly P^{(0)}_v(X)$. 

To simplify the discussion, let us consider a single toric chart $U_\alpha\cong\C^4$ of $X$ (details can be found in Section \ref{sec:fixlocal}).  The $(\C^*)^4$-fixed $\PT_0$ pairs on $\C^4$ with scheme theoretic support $Z\subset \C^4$ correspond to the finitely generated $(\C^*)^4$-invariant submodules of 
\begin{equation} \label{eqn:limtocalcintro} 
\varinjlim \hom(\mathfrak{m}^r, \O_{Z}) / \O_Z,
\end{equation}
where $\mathfrak{m} \subset \O_{Z}$ denotes the ideal sheaf of the origin.
The pure part $Z_1 \subset Z$ of the support is described by six finite partitions (i.e., piles of boxes in $\Z_{\geq 0}^2$) $\boldsymbol{\lambda} = \{\lambda_{ab} \}_{1 \leq a < b \leq 4}$.
Consider the subscheme
\begin{equation*} 
    W = (Z_{\lambda_{12}} \cup Z_{\lambda_{13}} \cup Z_{\lambda_{23}}) \cap (Z_{\lambda_{14}} \cup Z_{\lambda_{24}} \cup Z_{\lambda_{34}}),
\end{equation*}
where $Z_{\lambda_{ab}}$ is the pure 2-dimensional subscheme set-theoretically supported on $\C^2 \cong Z(x_a,x_b) \subset \C^4$ with scheme structure determined by $\lambda_{ab}$ (Section \ref{sec:fixlocal}). Then $Z_1$ is Cohen-Macaulay if and only if $T_0(\O_W) \subset \O_W$ is trivial.\footnote{For any coherent sheaf $E$, we denote by $T_i(E) \subset E$ its maximal $i$-dimensional subsheaf.}
We give an explicit description of \eqref{eqn:limtocalcintro} which determines $(\C^*)^4$-fixed $\PT_0$ pairs. An analogous description of $(\C^*)^4$-fixed $\PT_1$ pairs is given in Section \ref{sec:fixlocal}.
\begin{theorem} \label{thm:fixedintro}
Let $(F,s)$ be a $(\C^*)^4$-fixed $\PT_0$ pair on $\C^4$ with scheme theoretic support $Z$. Let $Z_1 \subset \C^4$ be the pure part of $Z$.
\begin{enumerate}
    \item[$\mathrm{(1)}$] The induced sequence
$$
0 \to  \varinjlim \hom(\mathfrak{m}^r, T_{1}(\O_Z) ) \to \varinjlim \hom(\mathfrak{m}^r, \O_{Z}) \to \varinjlim \hom(\mathfrak{m}^r, \O_{Z_1}) \to 0
$$ 
is exact and there is a $(\C^*)^4$-equivariant isomorphism
\begin{equation*} 
\varinjlim \hom(\mathfrak{m}^r, \O_{Z_1}) / \O_{Z_1} \cong T_0(\O_W).
\end{equation*}
    \item[$\mathrm{(2)}$]  Suppose $Z = Z_1$. Then the $(\C^*)^4$-fixed $\PT_0$ pairs on $\C^4$ with scheme theoretic support $Z$ are in bijective correspondence with the (finitely many) closed points of
$$
\bigcup_{m \geq 0} {\curly Quot}_W(T_0(\O_W) ,m)^{(\C^*)^4}\,,
$$
where ${\curly Quot}_W(T_0(\O_W) ,m)$ denotes the Quot scheme of length $m$ 0-dimensional quotients of $T_0(\O_W)$.
\end{enumerate}
\end{theorem}

\begin{example}\label{ex:intro1}
Consider $Z=S_1\cup S_2 \subset \C^4$, where $S_1 = Z(x_1,x_2)$ and $S_2 = Z(x_3,x_4)$. This is the simplest example of a non-Cohen-Macaulay singularity. 
The $(\C^*)^4$-fixed $\PT_0$ pairs on $\C^4$ with scheme theoretic support $Z$ correspond to elements of the Hilbert schemes of points $\Hilb^m(S_1 \cap S_2)$. 
Therefore, we have $m=0,1$ and there are only two fixed points corresponding to the trivial pair $\O_Z \to 0$ and the extension $0 \to \O_{Z} \to \O_{S_1} \oplus \O_{S_2} \to \O_{S_1 \cap S_2} \to 0$. 
\end{example}

\subsection{Vertex formalism}

For toric Calabi-Yau 4-folds $X$, we use the virtual localization formula \cite{OT} to compute equivariant $\DT,\PT_0$, and $\PT_1$ invariants.

Denote by $V(X)$ the collection of $T_X$-fixed points of $X$, by $E(X)$ the collection of $T_X$-fixed projective lines in $X$, and by $F(X)$ the collection of $T_X$-fixed irreducible smooth projective toric surfaces in $X$. We refer to elements of $V(X)$, $E(X)$, $F(X)$ as vertices, edges, and faces. To a vertex $\alpha \in V(X)$ corresponds a $T_X$-fixed point $p_\alpha$ and to an edge $\alpha\beta \in E(X)$ corresponds a $T_X$-fixed $\PP^1 \subset X$ passing through $p_\alpha, p_\beta$. To a face $f \in F(X)$ corresponds a $T_X$-fixed surface $S_f \subset X$.
Each $p_\alpha$ lies in a toric chart $U_\alpha\cong\C^4$, which jointly cover $X$. 

Let $\curly I := \curly I_v(X)$.
Then $\curly I^{T_X} = \curly I^T$ and it is 0-dimensional and reduced. On each chart $U_\alpha$, a 2-dimensional $T_X$-fixed subscheme $Z_{\alpha} \subset U_\alpha$ corresponds to an infinite solid partition $\pi^{(\alpha)}$. The solid partition $\pi^{(\alpha)}$ can be visualized as a pile of boxes in $\Z_{\geq 0}^4$ with a limiting plane partition (i.e., pile of boxes in $\Z_{\geq 0}^3$ as in \cite{MNOP}) along each of the four coordinate directions $x_i \gg 0$. Therefore, $T_X$-fixed 2-dimensional subschemes $Z  \in \curly I^{T_X}$ correspond to certain collections of solid partitions $\boldsymbol{\pi} = \{\pi^{(\alpha)}\}_{\alpha \in V(X)}$ satisfying gluing conditions described in Sections \ref{sec:fixlocglob}. Consider the $K$-theory class of the virtual tangent bundle 
\begin{equation} \label{eqn:Ktheor1}
T^{\vir}_{\curly I}|_{Z} = R\Hom_X(I_{Z/X},I_{Z/X})_0[1] \in K^0_{T_X}(\pt) \cong \Z[t_1^{\pm 1}, t_2^{\pm 1}, t_3^{\pm 1},t_4^{\pm 1}],
\end{equation}
where $I_{Z/X} \subset \O_X$ denotes the ideal sheaf of $Z$ and $(\cdot)_0$ the trace-free part. %, and $K^0_{T_X}(\pt)$ is the  

After redistribution, the $K$-theory class \eqref{eqn:Ktheor1} can be separated into vertex $\mathsf{V}_\alpha$, edge $\mathsf{E}_{\alpha\beta}$, and face $\mathsf{F}_f$ contributions, defined in Section \ref{sec:redist} in terms of $\boldsymbol{\pi}$, which are Laurent \emph{polynomials} in $t_1,t_2,t_3,t_4$. 
\begin{theorem}
Let $X$ be a toric Calabi-Yau 4-fold and $Z \subset X$ a 2-dimensional $T_X$-fixed closed subscheme with proper support. Then we have
\begin{equation*} 
T^{\vir}_{\curly I}|_{Z} = R\Hom_X(I_{Z/X},I_{Z/X})_0[1]  = \sum_{\alpha \in V(X)} \mathsf{V}_{\alpha} + \sum_{\alpha\beta \in E(X)} \mathsf{E}_{\alpha\beta} + \sum_{f \in F(X)} \mathsf{F}_{f},
\end{equation*}
where $\mathsf{V}_{\alpha}, \mathsf{E}_{\alpha\beta}, \mathsf{F}_{f}$ are Laurent polynomials in $t_1,t_2,t_3,t_4$.
\end{theorem}
For curve counting theory on $3$-folds, the redistribution was developed in \cite{MNOP,PT2}. The surface asymptotics adds new features to the redistribution process. We establish the analog of this theorem for any $T_X$-fixed $\PT_0$ and $\PT_1$ pair on $X$ with proper support in Section \ref{sec:redist}.

Let $\curly P := \curly P_v^{(0)}(X)$. Unlike the $\DT$ case, in general we no longer have ``$\curly P^{T_X} = \curly P^{T}$, and it is 0-dimensional and reduced'' (and similarly in the $\PT_1$ case). Henceforth, let us assume we are in a case where $\curly P^{T_X} = \curly P^{T}$, and it is 0-dimensional and reduced. The same issue appears for $\PT$ curve theory on toric Calabi-Yau 4-folds. Then this assumption is rather restrictive and essentially rules out 3 or 4 legs coming together in any chart $U_\alpha$ \cite{PT2,CKM1}. In fact, the union of all $(\C^*)^4$-fixed loci in the 4-leg case satisfies Murphy's Law \cite{Sch}.
However, by the following proposition, for $\PT_0$ theory, our assumption is satisfied in many interesting ``truly 4-dimensional'' examples.

\begin{theorem} 
 Let $\curly P:= \curly P^{(0)}_v(X)$. For all $(F,s)\in \curly P^{T_X}$, suppose either (1)  the sheaf $F$ is pure, or (2) the scheme theoretic support of $F|_{U_\alpha}$ on each chart $U_\alpha$ has at most two embedded 1-dimensional components and its pure part is Cohen-Macaulay.
Then $\curly P^{T_X} = \curly P^{T}$ and it is 0-dimensional and reduced.  
\end{theorem}

The combinatorial description of edge (resp.~face) terms associated to $\DT$/$\PT_0$ pairs is related to vertex (resp.~edge) terms of $\DT$/$\PT$ pairs on toric 3-folds.
Fix plane partitions $\boldsymbol{\mu} = \{\mu_a\}_{a=1}^{4}$ corresponding to $(\C^*)^3$-fixed subschemes of $\C^3$ of dimension $\leq 1$, which are the ``asymptotic cross-sections'' that we fix along the $x_1,x_2,x_3,x_4$-axes. We twist each vertex contribution of \eqref{eq:redisteq} by the Nekrasov insertion \eqref{eqn:Ktheorinvintro} and define the $K$-theoretic $\DT$ and $\PT_0$ {\em topological vertex}
\begin{equation}\label{eqn:topvertex}
    \mathsf{V}^{\DT}_{\mu_1\mu_2\mu_3\mu_4}(q), \quad \mathsf{V}^{\PT_0}_{\mu_1\mu_2\mu_3\mu_4}(q)\,,
\end{equation}
where \eqref{eqn:topvertex} depend on a choice of $\pm \sqrt{\cdot}$ assigned to each fixed point.

We conjecture the following $\DT$--$\PT_0$ topological vertex correspondence.
\begin{conjecture} \label{conj:vertexDTPT0intro}
For plane partitions $\boldsymbol{\mu} = \{\mu_a\}_{a=1}^{4}$ as above, there exists a choice of $\pm \sqrt{\cdot}$ assigned to the $(\C^*)^4$-fixed points such that
$$
\frac{\mathsf{V}^{\DT}_{\mu_1\mu_2\mu_3\mu_4}(q)}{\mathsf{V}^{\DT}_{\varnothing\varnothing\varnothing\varnothing}(q)} = \mathsf{V}^{\PT_0}_{\mu_1\mu_2\mu_3\mu_4}(q).
$$
\end{conjecture}
The denominator 
is given by the Magnificent Four formula \cite{Nek,KR}. For curves counting on toric $3$-folds, the $\DT/\PT$ vertex and vertex correspondence has been studied in \cite{AKMV, MNOP,MNOP2,PT2,NO}. Its numerical version was proved combinatorially in \cite{JWY} and its full $K$-theoretic version was recently established in \cite{KLT}. For curve counting theories on toric Calabi-Yau 4-folds, the $\DT/\PT$ vertex correspondence has been studied in \cite{CKM1, Liu}. 

We check Conjecture \ref{conj:vertexDTPT0intro}, up to certain order in $q$, by direct calculation in several examples (Proposition \ref{prop:verif}).

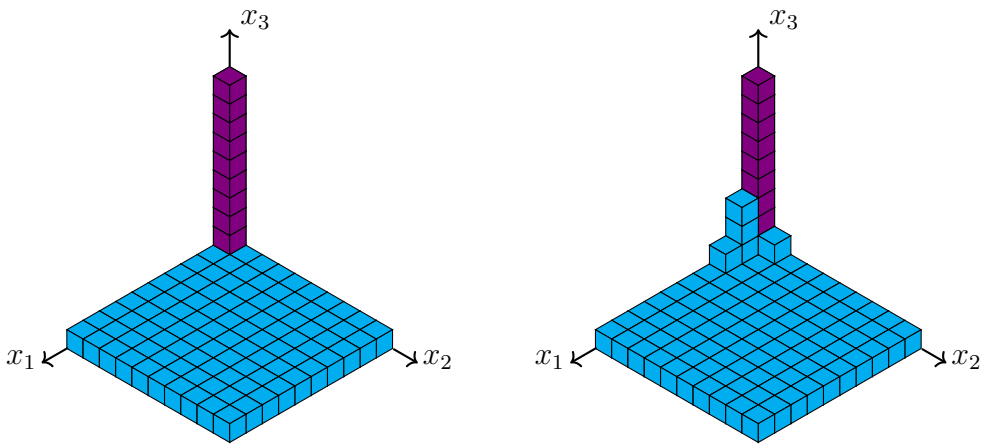
\begin{figure}[t]
\centering
\begin{subfigure}{.5\textwidth}
  \centering
  \input{fig0a.tex}
  %\caption{A subfigure}
  %\label{fig:sub1}
\end{subfigure}%
\begin{subfigure}{.5\textwidth}
  \centering
  \input{fig0b.tex}
  %\caption{A subfigure}
  %\label{fig:sub2}
\end{subfigure}
\caption{Adding boxes to $S_1 \cup S_2$.}
\label{fig0aand0b}
\end{figure}

\begin{example} 
Suppose we fix $\boldsymbol{\mu}$ such that the minimal $(\C^*)^4$-fixed 2-dimensional closed subscheme $Z \subset \C^4$ with asymptotics $\boldsymbol{\mu}$ is \emph{pure}. Then, by Theorem \ref{thm:fixedintro}(ii), $\mathsf{V}^{\PT_0}_{\mu_1\mu_2\mu_3\mu_4}(q)$ is a Laurent \emph{polynomial} in $q$. Moreover, its number of terms is at most one plus the length of the 0-dimensional sheaf $\ext^3(\O_Z,\O_{\C^4})$, which measures how far $Z$ is from being Cohen-Macaulay (Proposition \ref{Lem.PT0=PT1}). For $Z = S_1 \cup S_2$ as in Example \ref{ex:intro1}, $\ext^3(\O_Z,\O_{\C^4})$ has length 1 and 
$$
 \sfV^{\PT_0}_{\mu_1\mu_2\mu_3\mu_4}(q) = \frac{[t_1t_3][t_2t_3]}{[y]} q^{-1}  + [t_1t_2], 
$$
where $[x] := x^{\frac{1}{2}} - x^{-\frac{1}{2}}$. In Figure \ref{fig0aand0b} (left), we depict the $Z(x_4)$-slice of the solid partition corresponding to $S_1 \cup S_2$. The blue boxes have no boxes stacked on top of them in the positive $x_4$-direction, whereas the purple boxes have infinitely many boxes stacked on top in the positive $x_4$-direction. The box at $(0,0,0,0)$ is purple. On the $\DT$-side, the $(\C^*)^4$-fixed points are obtained by stacking boxes ``on top'' of this configuration as in Figure \ref{fig0aand0b} (right). On the $\PT_0$-side, there are only two $(\C^*)^4$-fixed points corresponding to putting zero or one box at the intersection $S_1 \cap S_2$ (Example \ref{ex:intro1}). 
\end{example}

We show that Conjecture~\ref{conj:vertexDTPT0intro} implies the \emph{global} $T$-equivariant $K$-theoretic $\DT$--$\PT_0$ correspondence.
\begin{proposition}
Let $L$ be a $T_X$-equivariant line bundle on $X$. Let $\gamma \in H_c^4(X,\Q)$ and $\beta \in H_c^6(X,\Q)$ such that $\curP_{v}^{(0)}(X)^T = \curP_{v}^{(0)}(X)^{T_X}$ is empty or 0-dimensional and reduced for $v = (0,0,\gamma,\beta,n-\gamma \cdot \td_2(X))$ and all $n$. Suppose Conjecture \ref{conj:vertexDTPT0intro} holds. Then there are choices of square roots such that the equivariant form of \eqref{eq:DTPT0intro} holds.
\end{proposition}

Our vertex formalism is useful for calculating $\DT, \PT_0, \PT_1$ invariants for \emph{global} toric geometries. There is a huge advantage of $\PT_0$ and $\PT_1$ theory over $\DT$ theory --- since there are no free moving points the moduli space can be compact (Proposition \ref{prop:localFanos}). 

\begin{example} \label{ex:intro2}
For $X = \mathrm{Tot}(\O_{\PP^2}(-2) \oplus \O_{\PP^2}(-1))$, the $\PT_0$ and $\PT_1$ moduli spaces are compact. Using the vertex formalism, we determined the following formulae\footnote{The $\PT_1$ calculations in Examples \ref{ex:intro2}, \ref{ex:intro3} are subject to Conjecture \ref{conj:0dimTXfix}: if $\curly P^{T_X}$ is 0-dimensional (e.g.~for local surfaces), then $\curly P^{T_X} = \curly P^T$ and it is 0-dimensional and reduced.} for certain $m,\delta,n$
\begin{align*}
\langle\!\langle \O_X \rangle\!\rangle_{1,\frac{5}{2},m}^{\PT_0} &= [y](y^{\frac{m-1}{2}}+y^{-\frac{m-1}{2}}), \\
\langle\!\langle \O_X \rangle\!\rangle_{1,\frac{3}{2}+\delta,n}^{\PT_1} &= \mathsf{CO}_{\chi(\O_{\PP^2}(\delta))-n}(\PP^2,\O(-1)) [y] (y^{\frac{n-1}{2}} + \cdots + y^{-\frac{n-1}{2}}),
\end{align*}
where $m > 3$, $\delta \geq 0$, and $1 \leq n \leq \chi(\O_{\PP^2}(\delta))$, and $\mathsf{CO}_{i}(\PP^2,\O(-1))$ are the Carlsson-Okounkov numbers associated to the $\O(-1)$-twisted tangent bundles of Hilbert schemes of points on $\PP^2$ (see Section \ref{sec:localP2P3}). 
\end{example}

The formulae in the previous example are for primitive $\gamma$ and are in fact proved using non-toric methods for \emph{all} $m,\delta,n$ in a sequel to this paper. However, the vertex formalism allows us to explore the non-primitive case:

\begin{example} \label{ex:intro3}
Continuing the notation of the previous example, we define
$$
\mathsf{G}^{\PT_1}_d(q,Q,y) := \sum_{m,n} \langle\!\langle \O_X \rangle\!\rangle_{d,m,n}^{\PT_1} Q^m q^n.
$$
One benefit of the $\PT_1$ generating series is that, for fixed $d,m$, it has only finitely many terms \cite[Prop.~2.14]{BKP}. Using the vertex formalism, for appropriate square roots we find
\begin{align*}
[y]^{-16} \mathsf{G}^{\PT_1}_4 = &\, q^{16}  Q^{11}+\Big(8q^{16}+4(y^{\frac{1}{2}}+y^{-\frac{1}{2}}) q^{17}+4(y+y^{-1}) q^{18}+ \\
&4(y^{\frac{3}{2}}+y^{-\frac{3}{2}}) q^{19} +10(y^2+8y+12+8y^{-1}+y^{-2})q^{20} \Big)Q^{12} + \cdots,
\end{align*}
where $\cdots$ means higher order terms in $Q$. For more examples, see Section \ref{sec:localP2P3}.
\end{example}

\subsection{Virtual Lefschetz principle} 

Let $X$ be a projective Calabi-Yau 4-fold and let $\iota : Y \hookrightarrow X$ be smooth connected effective divisor inducing an injection $\iota_* : H_*(Y,\Q) \hookrightarrow H_*(X,\Q)$. Let $L:=\O_X(Y)$, $v' \in H^*(Y,\Q)$, and 
$$
v := (0,0,\gamma,\beta,n- \gamma \cdot \td_2(X)) =  \iota_* v' \cdot \td(L)^{-1} \in H^*(X,\Q).
$$
Let $q \in \{-1,0,1\}$ and consider the moduli spaces $\curly P_X=\curly P^{(q)}_{v}(X)$ and $\curly P_Y := \curly P_{v'}^{(q)}(Y)$. Then we have a closed embedding $\imath = \iota_* : \curly P_Y \hookrightarrow \curly P_X$.

We relate surface counting on $X$ to surface counting on the divisor $Y$.
\begin{theorem} \label{thm:Lefschetzviapullbackintro}
Suppose $H^{>0}(X,F \otimes L) = 0$ for all $(F,s) \in \curly P_X$, $H^2(Y,F) = 0$ for all $(F,s) \in \curly P_Y$, and $\imath$ induces an injective map on the sets of connected components of $\curly P_Y$ and $\curly P_X$. Then there exists an orientation such that
\[\imath_* \Ohat_{\curly P_Y} = \widehat{\Lambda}\udot (R \pi_{\curly P_X *}(\FF \otimes L))^\vee \otimes  \Ohat_{\curly P_X}.\]
\end{theorem}
We use the third-named author’s virtual pullback formula \cite{Par} for the proof. Next, we identify surface counting invariants on $X$ with curve counting invariants on $Y$ by ``twisting away $\O_Y(S)$" when $S\subset Y$ is the unique pure 2-dimensional subscheme in class $[S]\in H^2(Y,\Q)$.
\begin{corollary} 
Continuing the setting of Theorem \ref{thm:Lefschetzviapullbackintro}, suppose $S\subset Y$ is the only pure 2-dimensional subscheme in class $[S] \in H^2(Y,\Q)$ and let $\gamma := \iota_* [S]$. Let $v'' := 1-(1-v') e^{[S]} \in H^*(Y,\Q)$ and let $\beta'' := v_2'' \in H^4(Y,\Q)$. Suppose that, for all $n \in \Z$ and $q \in \{-1,0\}$, the inclusion maps $\curly P^{(q)}_{v'}(Y) \hookrightarrow \curly P^{(q)}_{v}(X)$ induce injective maps on the sets of connected components. Then, for $y=1$, there exist orientations such that
\begin{align*} 
\sum_n \langle \! \langle L \rangle \! \rangle_{X,\gamma,\beta,n}^{\DT} q^n = q\udot \sum_n \langle \! \langle 1 \rangle \! \rangle_{Y,\beta'',n}^{\DT} q^n, \quad \sum_n \langle \! \langle L \rangle \! \rangle_{X,\gamma,\beta,n}^{\PT_0} q^n = q\udot \sum_n \langle \! \langle 1 \rangle \! \rangle_{Y,\beta'',n}^{\PT} q^n,
\end{align*}
where $q\udot$ is and overall factor with power depending only on $S,Y,X,\beta''$ and the right hand sides are the $K$-theoretic $\DT$/$\PT$ generating series of $Y$ \cite{NO}. 
\end{corollary}

In the setting of the corollary, we deduce that the $\DT$--$\PT_0$ correspondence (Conjecture \ref{conj:KDTPT0intro}) for $X$, $L$, $\gamma$, $\beta$, $y=1$ is equivalent to the (curve) $\DT$--$\PT$ correspondence for $Y$, $\beta''$. In particular, it holds in the following three cases (1) when $\beta''=0$, (2) when $Y$ is a Calabi-Yau 3-fold by \cite{Bri, Tod}, and (3) when $Y$ is a toric 3-fold by \cite{KLT}. We present two examples.\footnote{In the following two examples, we assume that the condition on connected components in the corollary is satisfied.}

\begin{example}
Let $Y \to \PP^2$ be a general Weierstra{\ss} fibration with unique section $\PP^2 \subset Y$. Let $X = Y \times E$, where $E$ is a smooth elliptic curve. Let $\gamma = [\PP^2]$ and $\beta = (\epsilon + \tfrac{3}{2}) [\PP^1]$, where $\PP^1 \subset \PP^2$ is the class of a line on the section and $\epsilon = 0,1$. Then for $L = \O_X(Y \times \{\mathrm{pt}\})$ and $y=1$, the $\DT$--$\PT_0$ correspondence of Conjecture \ref{conj:KDTPT0intro} holds. For $\epsilon = 1$, it reduces to the $\DT$--$\PT$ correspondence for curve class $[\PP^1]$ of the Calabi-Yau 3-fold $Y$ \cite{Bri, Tod}.
\end{example}

\begin{example}
Let $Y$ be the blow-up of $\PP^3$ in a point and let $\PP^2 \subset Y$ be the class of the exceptional divisor. Let $X \to Y$ be a general Weierstra{\ss} fibration with section, so we have $\PP^2 \subset Y \subset X$. Let $\gamma = [\PP^2]$ and $\beta = (\epsilon - \tfrac{1}{2})[\PP^1]$,  where $\PP^1 \subset \PP^2$ is the class of a line on the exceptional divisor and $\epsilon = 0,1$. Then for $L = \O_X(Y)$ and $y=1$, the $\DT$--$\PT_0$ correspondence of Conjecture \ref{conj:KDTPT0intro} holds. For $\epsilon = 1$, it reduces to the $K$-theoretic $\DT$--$\PT$ correspondence of the Fano 3-fold $Y$ for curve class $[\PP^1]$ \cite{KLT}.
\end{example}

We also provide a second virtual Lefschetz principle when $Y$ is a smooth quasi-projective toric 3-fold and $X = \mathrm{Tot}(K_Y)$ by combining the virtual localization formula with the spectral construction (Theorem \ref{thm:toricdimred}).

\subsection{Relations to physics} \label{sec:phys}

Our vertex formalism for $q=-1$ is related to recent work of Nekrasov-Piazzalunga \cite{NP2}. 
Their work originates from physics, namely maximally supersymmetric $(8+1)$-dimensional gauge theory on a toric Calabi-Yau 4-fold fibred over a circle or, string-theoretically, the counting of BPS states of a system of D0--D2--D4--D6--D8 branes on a toric Calabi-Yau 4-fold in the presence of a large Neveu-Schwarz $B$-field. We provide the mathematical explanation for why the vertex, edge, and face terms in \cite{NP2} calculate the virtual tangent bundle \eqref{eqn:Ktheor1}.
As opposed to loc.~cit., which only concerns the vertex formalism for the $q=-1$ case, we consider all three cases $q=-1,0,1$. For $q=0,1$, there are many examples of \emph{compact} moduli spaces, whereas this almost never happens for $q=-1$ due to free roaming points. This allows us to provide global computations of $\PT_q$ invariants (Proposition \ref{prop:localFanos}, Section \ref{sec:localP2P3}).

\subsection{Open questions}

We addressed the $\DT$--$\PT_0$ correspondence among the three stability conditions for pairs. It is interesting to find relations between $\PT_0$ and $\PT_1$ invariants. In a sequel, we will study this direction. 

In \cite{BKP} we introduced \emph{reduced} virtual cycles for surface counting on Calabi-Yau 4-folds. We expect Conjecture~\ref{conj:KDTPT0intro} to hold for reduced invariants. See Remark \ref{rem:nored} for why reducing does not play a significant role in this paper.
Furthermore, it is interesting to find the connection between Conjecture~\ref{conj:KDTPT0intro} and the new wall-crossing formalism introduced in \cite{GJT, Boj2}. 

It is desirable to remove the condition ``$\curly P^{T_X}=\curly P^T$, and it is 0-dimensional reduced''. For curve counting on Calabi-Yau 4-folds, this is pursued in \cite{Liu}. 

Our formulae depend on a choice of square root assigned to the fixed points. We expect that there exist \emph{canonical global orientations} on $\curly P$ which induce these signs. In all examples we calculated, our formulae hold for very special (often unique) choices of $\pm \sqrt{\cdot}$. For Hilbert schemes of points on $\C^4$, sign formulae were conjectured in \cite{NP1} and proved using global orientations in \cite{KR}. \\

\noindent \textbf{Acknowledgements.} We thank R.~Pandharipande for encouraging us to study surfaces on Calabi-Yau 4-folds, which prompted this project. The second-named author likes to thank Y.~Cao and S.~Monavari for their collaborations, which play an important role in this paper. We also thank A.~Bojko, N.~Kuhn, H.~Liu, N.~Piazzalunga, J.~Rennemo, R.~Schmiermann, F.~Thimm, R.P.~Thomas for crucial discussions. YB is supported by ERC Grant ERC-2017- AdG-786580-MACI and SNSF Postdoc.Mobility Grant. The project has received funding from the European Research Council (ERC) under the European Union Horizon 2020 research and innovation program (grant agreement No. 786580).  MK is supported by NWO Grant VI.Vidi.192.012 and ERC Consolidator Grant FourSurf 101087365. HP is supported by KIAS (SG089201).

\subsection{Notation and conventions}

We use the following notation and conventions throughout the paper.
\begin{itemize}
\item All schemes and algebraic stacks are assumed to be of finite type over the complex field $\C$, unless stated otherwise. 
A {\em variety} is a separated integral scheme.
A {\em point} $x\in X$ in a scheme or an algebraic stack $X$ means a $\C$-valued point. 
\item A {\em Calabi-Yau $n$-fold} is a smooth quasi-projective variety $X$ of dimension $n$ with $K_X \cong \O_X$. 
\item For a coherent sheaf $F$ on a smooth quasi-projective variety $X$, we use the following notation for the various duals:
\begin{itemize} 
\item $F^\vee := R\hom_X(F,\O_X)$;
\item $F^D:=\ext^c_X(F,K_X)$, where $c$ is the codimension of $F$ in $X$.
\end{itemize}
\item For a smooth quasi-projective variety $X$, we denote by 
\begin{itemize}
\item $\Coh(X)$ the abelian category of coherent sheaves on $X$;
\item $\Coh_{\leq q}(X) := \{F \in \Coh(X) : \dim(F)\leq q\}$; 
\item $\Coh_{\geq q+1}(X) :=\{ F \in \Coh(X) : T_q(F)=0\}$, where $T_q(F) \subset F$ is the maximal $q$-dimensional torsion subsheaf of $F$;
\item $D^b_{\mathrm{coh}}(X)$ the bounded derived category of coherent sheaves on $X$.
\end{itemize}
We also consider the full subcategories $\Coh(X)_c$, $\Coh_{\leq q}(X)_c$, $\Coh_{\geq q+1}(X)_c$ of coherent sheaves with proper (``compact'') support.
\item For a quasi-projective scheme $M$, let $K_0(M)$ be the Grothendieck group of coherent sheaves on $M$ and $K^0(M)$ the Grothendieck group of locally free sheaves on $M$. For an algebraic torus $T$ acting on $M$, we denote the $T$-equivariant versions of these $K$-groups by $K_0^T(M)$ and $K^0_T(M)$. We also work with $K$-groups of coherent sheaves with proper (``compact'') support, which we denote by $K_0(M)_c$, $K_0^T(M)_c$.
\item For any $V \in K^0(M)$ with $\rk(V)=1$, there is a well-defined square-root operation $V^\frac{1}{2} \in K^0(M,\mathbb{Z}[\tfrac{1}{2}])$ \cite[Lem.~5.1]{OT}.
\item For any rank $r$ vector bundle $V$ on $M$, we write $\Lambda_x V := \sum_{i=0}^{r} [\Lambda^i V] x^i \in K^0(M)[x]$. This operation can be extended to $K^0(M)[[x]]$ by setting $$\Lambda_x(-V) = 1/\Lambda_x V = \Sym_{-x} V \in K^0(M)[[x]].$$ Whenever it is well-defined, e.g.~in the case $M$ has trivial $T$-action and $V \in K^0_T(M)$ has no $T$-fixed part with negative coefficients, we set $\Lambda^\mdot V := \Lambda_{-1} V$. Finally, we define a normalized version
\begin{align*}
\widehat{\Lambda}^\mdot V &:= (\Lambda^\mdot V) \otimes \det(V)^{-\frac{1}{2}} \in K^0(M,\mathbb{Z}[\tfrac{1}{2}]).
\end{align*}
\item For any quasi-projective scheme $M$, we denote by $A_*(M)$ the Chow group with $\Q$-coefficients. For any $d \in \Q \setminus \Z$, we let $A_d(M):=0$.
\item For a morphism $f : M \to N$ of quasi-projective schemes, we denote by $\LL_{M/N}$ the full cotangent complex \cite{Ill} and by $\tau^{\geq -1}\LL_{M/N}$ the truncated cotangent complex. In the case $N = \Spec \C$, we write $\LL_{M} := \LL_{M / \Spec \C}$.
\item For any $d$-dimensional smooth toric variety $Y$ defined by a fan $\Sigma$ in a lattice $N \cong \Z^d$, we always assume that all cones of $\Sigma$ are contained in a $d$-dimensional cone. This implies that $Y$ can be covered by $(\C^*)^d$-invariant copies of $\C^d$ \cite{Ful}.
\end{itemize}

\section{Preliminaries}

\subsection{Basic properties of stable pairs}

Let $X$ be a smooth quasi-projective variety of dimension $n$ over the complex field $\C$. We first recall the following definition from \cite{BKP}.
\begin{definition}\label{def:stabilityconditionsofpairs'}
Let $F$ be a coherent sheaf on $X$ with 2-dimensional proper support and let $s: \O_X \to F$ be a section.
For $q \in \{-1,0,1\}$, we say that the pair $(F,s)$ is {\em $\PT_q$ stable} if
\[F \in \Coh_{\geq q+1}(X) \and Q:=\coker(\O_X \xrightarrow{s} F) \in \Coh_{\leq q}(X).\]
We also write $\DT := \PT_{-1}$.
\end{definition}

The following results were proved in \cite{BKP} when $X$ is projective. The proofs in loc.~cit.~also work in the quasi-projective case. 
\begin{lemma}
Given a $\PT_q$ pair $(F,s)$ on $X$, the scheme theoretic support of the coherent sheaf $F$ is the scheme theoretic support of the section $s$.
\end{lemma}

\begin{proposition} \label{lem:limPTq} Let $W \subseteq Z \subseteq X$ be closed subschemes such that 
$\O_Z \in \Coh_{\geq q+1}(X)$ and $\O_W \in \Coh_{\leq q}(X)$. Suppose $W$ is reduced 
and denote its ideal by $\fI \subseteq \O_Z$. 
Then $\PT_q$ pairs $(F,s)$ with support $Z$ and cokernel $Q$ satisfying $\Supp(Q)^{\mathrm{red}} \subseteq W$ are equivalent to coherent subsheaves of $\varinjlim \hom_X(\fI^r,\O_Z) / \O_Z$.
\end{proposition}

\begin{proposition}\label{Lem.PT0=PT1}
Let $S$ be a \emph{pure} 2-dimensional subscheme of $X$. Then the following three conditions are equivalent:
\begin{enumerate}
\item[$\mathrm{(i)}$] $S$ is not Cohen-Macaulay.
\item[$\mathrm{(ii)}$] $\ext^{n-1}_X(\O_S,K_X) \neq 0$.
\item[$\mathrm{(iii)}$] There exists a $\PT_0$ pair $(F,s)$ with support $S$ and non-zero cokernel $Q$.
\end{enumerate}
Furthermore, for any $m \geq 0$, the $\PT_0$ pairs $(F,s)$ with support $S$ and $\chi(Q) = m$ are in bijective correspondence with the closed points of the Quot scheme
$$
{\curly Quot}_X(\ext^{n-1}_X(\O_S,K_X),m)
$$
parametrizing length $m$ quotients of the 0-dimensional sheaf $\ext^{n-1}_X(\O_S,K_X)$.
\end{proposition}

\begin{proposition}  \label{prop:PT0PTcv}
For any $\PT_0$ pair $(F,s)$ on $X$ such that $S=\Supp(F)^{\pure}$ is Cohen-Macaulay, we have a natural quasi-isomorphism
\begin{equation*}
I\udot \cong [I_{S/X} \to T_1(F)],  
\end{equation*}
where $T_1(F)$ is zero or pure 1-dimensional and $I_{S/X}$ denotes the ideal sheaf of $S$ in $X$.

In particular, when $n = 3$, we have $I_{S/X} = \O_X(-S)$ and 
\[I\udot \cong [\O_X \to T_1(F)(S)](-S)\]
which, for $T_1(F) \neq 0$, is a 1-dimensional $\PT$ pair on $X$ as in \cite{PT1}. 

Moreover, for $n=3$ and an effective divisor $S \subset X$, this induces an equivalence between $\PT_0$ pairs $(F,s)$ on $X$ with $S = \Supp(F)^{\pure}$ and $T_1(F) \neq 0$, and $\PT$ pairs on $X$.
\end{proposition}
\begin{proof}
All but the last sentence were proved in \cite{BKP}. For the last part, given a $\PT$ pair $(G,t)$ on $X$, its support $S$ is a divisor, hence Cohen-Macaulay and we can form the push-out
\begin{displaymath}
\xymatrix
{
\O_X(-S) \ar[r] \ar[d] & G(-S) \ar[d] \\
\O_X \ar[r] & \O_X \sqcup_{\O_X(-S)} G(-S).
}
\end{displaymath}
Then the bottom row defines a $\PT_0$ pair with torsion subsheaf $G(-S)$, and pure part of its support equal to $S$. This construction is easily seen to be the inverse of the construction in loc.~cit.
\end{proof}

\begin{proposition} \label{lem:PT1PTcv}
For any $\PT_1$ pair $(F,s)$ on $X$ with Gorenstein support $S$ and non-zero cokernel $Q$, we have a natural quasi-isomorphism
\begin{equation*}
F\dual \otimes K_X[n-2] \cong [\O_S \to Q^D \otimes K_S\dual] \otimes K_S
\end{equation*}
where $Q^D:=\ext^{n-1}_X(Q,K_X)$ is a pure 1-dimensional sheaf scheme theoretically supported on $S$.
Moreover, this induces an equivalence between $\PT_1$ pairs supported on $S$ with non-zero cokernel and 1-dimensional $\PT$ pairs on $S$. 
\end{proposition}

\begin{proposition}\label{cor:PT1smsupp}
Let $\iota : S \hookrightarrow X$ be a smooth surface in a quasi-projective variety of dimension $n$. A $\PT_1$ pair $(F,s)$ on $X$ with support $S$ and non-zero cokernel $Q$ is equivalent to a pair $(C,Z)$ of a (non-zero) effective divisor $C \subseteq S$ and a 0-dimensional subscheme $Z \subseteq C$. Moreover, $Q \cong \jmath_* I_{Z/C}(C)$, where $\jmath : C \hookrightarrow X$ denotes the inclusion.
\end{proposition}
\begin{proof}
Except for the last sentence, this is proved in \cite[Sect.~2.3]{BKP}. By Proposition \ref{lem:PT1PTcv}, the stable pair $J\udot = [\O_S \to G]$ corresponding to a $\PT_1$ pair $(F,s)$ with non-zero cokernel $Q$ is obtained by $F^\vee \otimes K_X[n-2] \cong \iota_*( J\udot \otimes K_S)$. In particular, the cokernel $Q$ satisfies $Q^D \cong \iota_*(G \otimes K_S)$, i.e., $Q^\vee \otimes K_X [n-1] \cong \iota_*(G \otimes K_S)$ by purity of $Q$. Denote the scheme theoretic support of $G$ by $\jmath : C \hookrightarrow X$, then
\begin{align*}
(\iota_*(G \otimes K_S))^\vee \cong \jmath_* R\hom_C(G \otimes K_S|_C,\jmath^! \O_X) \cong (\jmath_*G^\vee(C))^\vee \otimes K_X^{-1}[1-n],
\end{align*}
where we used $K_C \cong K_S|_C(C)$ (in particular $C$ is Gorenstein). Hence $Q \cong \jmath_*G^\vee(C)$. Since $C$ is Gorenstein, we have $G^\vee \cong G^* \cong I_{Z/C}$ \cite[Lem.~B.2]{PT3} and the result follows.
\end{proof}

\subsection{Moduli spaces and virtual structures} \label{sec:moduli}

We now review the moduli spaces and virtual structures introduced in \cite{BKP} with a few adaptations required for the quasi-projective case.
For smooth quasi-projective varieties, it is convenient to work on a ``compactification'', so we can apply Serre duality. 
Let $X$ be a smooth quasi-projective variety. An open embedding $j: X\hookrightarrow \Xbar$ is a \textit{good compactification} of $X$ if $\Xbar$ is a smooth projective variety and $\Xbar\setminus X$ is a simple normal crossings divisor. Existence of good compactifications is guaranteed by Hironaka's theorem.

We denote by $H^*_c(X,\Q)$ the singular cohomology groups with compact support (with respect to the complex analytic topology on $X$). Then there is a natural Chern character map 
\[\ch : K_0(X)_c\to H^*_c(X,\Q)\,.\]

Let $X$ be a smooth quasi-projective 4-fold and denote by $\cPair(X)_c$ the moduli stack of non-zero pairs $(F,s)$ where $F\in \Coh(X)_c$ and $0\neq s\in H^0(X,F)$. As above, we fix a good compactification $j : X \hookrightarrow \Xbar$. For an $\O_X$-module $F$ on $X$, let $j_!F$ be the extension by zero. For compactly supported sheaves, $j_!$ preserves coherence, i.e., $j_!:\Coh(X)_c \to \Coh(\Xbar)$. This follows from the fact that for any object $F$ of $\Coh(X)_c$, we have $j_! F \cong j_* F$.
%Useful facts: $j_!$ left adjoint to $j^{-1} = j^*$ left adjoint to $j_*$ (on sheaves of modules). Moreover $j_!, j^*$ exact, and $j_! = j_*$ on coherent sheaves with proper support, and $j^* j_* = \mathrm{id}$.
We fix
\[\ch(F)=v := (0,0,\gamma,\beta,n-\gamma \cdot \td_2(X))\in H^*_c(X,\Q)\, .\]
Using a good compactification, the following proposition follows from \cite[Prop.~3.3, Thm.~3.10]{BKP}.

\begin{proposition}\label{prop:qproj}
For $q\in\{-1,0,1\}$ and $v\in H^*_c(X,\Q)$, there exists an open substack $\PTqvX$ of $\cPair(X)_c$ for which the $\C$-points are given by
\[\PTqvX (\C) = \{(F,s) : \ch(F)=v,\, (F,s) \text{ is a } \PT_q\text{ pair, } F \in \Coh(X)_c \}.\]
Then $\PTqvX  \to \curP^{(q)}_{j_*v}(\Xbar)$, $(F,s) \mapsto (j_* F, s')$, where $s'$ is the composition of the natural map $\O_{\overline{X}} \to j_* \O_X$ and $j_*s$, is an open embedding.
In particular, $\PTqvX$ is a quasi-projective scheme.\end{proposition}

Suppose $X$ is additionally Calabi-Yau and 
let $\curly P:=\PTqvX$. 
We consider the moduli stack \cite{Ina, Lie}
\[
\cPerf(\Xbar,1-j_*v)^\spl_{\O_{\Xbar}}
\]
of simple perfect complexes $I\udot$ on $\Xbar$ with $\ch(I\udot) = 1-j_*v$ and $\det(I\udot) \cong \O_{\Xbar}$. Then Proposition~\ref{prop:qproj} and \cite[Thm.~3.19]{BKP} imply:
\begin{proposition}
    Let $j:X\to\Xbar$ be a good compactification. The composition
    \[\PTqvX\to \curly P^{(q)}_{j_*v}(\Xbar)\to \cPerf(\Xbar,1-j_*v)^\spl_{\O_{\Xbar}}\,, \, (F,s)\to [\O_{\Xbar}\to j_*F]\]
    is an open embedding.
\end{proposition}

As a consequence the moduli space $\curly P$ has an obstruction theory 
\begin{equation}\label{eqn:obs}
    \phi: \EE = R\hom_{\pi_{\curly P}} (\II\udot,\II\udot)_0[3]\xrightarrow{\At(\II\udot)} \tau^{\geq-1}\LL_{\curly P}
\end{equation}
induced by the Atiyah class of the universal complex $\II\udot$ on $X\times\curly P$, where $\pi_{\curly P} : X \times \curly P \to \curly P$ denotes the projection, and $(\cdot)_0$ denote the trace-free part. We will show this obstruction theory is 3-term and symmetric. This requires using Serre duality for which we need the good compactification.

Suppose $\mathbb{F}$, $\mathbb{G}$ are coherent sheaves on $X \times B$ flat over a base scheme $B$ and with proper support. When $\FF$ has a section $s$, we have an induced complex $\II\udot = [\O_{X \times B} \to \FF]$. We write $$\II\udot_{\Xbar} := [\O_{\Xbar \times B} \to (j \times 1_B)_* \FF]$$ for the composition of the adjunction unit $\O_{\Xbar \times B} \to (j\times 1_B)_* \O_{X \times B}$ and $(j \times 1_B)_* s$. Then $R\hom_{\Xbar \times B}((j \times 1_B)_* \FF,(j \times 1_B)_*\GG) \cong (j \times 1_B)_* R\hom_{X \times B}(\FF,\GG)$, $R\hom_{\Xbar \times B}(\II\udot_{\Xbar},(j \times 1_B)_*\GG) \cong (j \times 1_B)_* R\hom_{X \times B}(\II\udot,\GG)$, $R\hom_{\Xbar \times B}(\II\udot_{\Xbar},\II\udot_{\Xbar})_0 \cong (j \times 1_B)_* R\hom_{X \times B}(\II\udot,\II\udot)_0$. One can now apply Grothendieck-Verdier duality along the projective morphism $\pi_B : \Xbar \times B \to B$ in order to deduce
\begin{align*}
R\hom_{\pi_B}(\FF,\GG) \cong (R\hom_{\pi_B}(\GG,\FF \otimes K_X[4])^\vee, \\
R\hom_{\pi_B}(\II\udot,\GG) \cong (R\hom_{\pi_B}(\GG,\II\udot \otimes K_X[4])^\vee, \\
R\hom_{\pi_B}(\II\udot,\II\udot)_0 \cong (R\hom_{\pi_B}(\II\udot,\II\udot \otimes K_X[4])_0^\vee,
\end{align*}
where we use $j^* K_{\Xbar} \cong K_X$ and for the last isomorphism we also use that we take the trace-free part.

Let $\overline{\curly P}:=\curly P^{(q)}_{j_*v}(\Xbar)$ be the moduli space of $\PT_q$ pairs on the good compactification $\Xbar$. Then it is straight-forward to show the following:
\begin{proposition} \label{prop:symmob}
    Let $X$ be a Calabi-Yau 4-fold. Let $j:X\to\Xbar$ be a good compactification. Denote by $\overline{\II}\udot$ the universal pair on $\Xbar \times \overline{\curly P}$ and consider the open embedding $\curly P \subset \overline{\curly P}$ of Proposition \ref{prop:qproj}.
    Then $\phi$ in \eqref{eqn:obs} is the restriction to $\curly P$ of the 3-term obstruction theory
    $$
    \overline{\phi} : \overline{\EE} = R\hom_{\pi_{\overline{\curly P}}}(\overline{\II}\udot,\overline{\II}\udot)_0[3]\xrightarrow{\At(\overline{\II}\udot)}\tau^{\geq -1}\LL_{\overline{\curly P}}.
    $$
    Therefore $\phi$ is 3-term and Grothendieck-Verdier duality for $\overline{\EE}$ induces a pairing $\theta : \EE^\vee[2] \cong \EE$ satisfying $\theta = \theta^\vee[2]$. Hence $\phi$ defines a 3-term \emph{symmetric} obstruction theory. 
\end{proposition}

Unlike the case of perfect (i.e., 2-term) obstruction theories, we require an additional property in order to construct the virtual cycle and virtual structure sheaf. Specifically, the intrinsic normal cone $\mathfrak{C}_{\curly P}$ should vanish under the quadratic function on $\mathfrak{C}(\mathbb{E}) := h^1 / h^0(\mathbb{E}^\vee)$ induced by the symmetric pairing \cite{OT,Par}. This property follows from the existence of a $(-2)$-shifted symplectic derived enhancement of $\curly P$ combined with \cite{BBJ}. The construction of $(-2)$-shifted symplectic structures was carried out in \cite{PTVV} in the projective case and in \cite{Pre} in the quasi-projective case.

Therefore, combining Proposition \ref{prop:symmob} with the work of Oh-Thomas \cite{OT}, for any choice of orientation, we obtain a virtual cycle and structure sheaf
\begin{align*}
&[\curly P]^{\vir} \in A_{\vd}(\curly P), \quad \vd := n - \frac{1}{2} (j_* \gamma)^2 \\
&\widehat{\O}_{\curly P}^{\vir} \in K_0(\curly P, \Z[\tfrac{1}{2}]).
\end{align*}
This virtual cycle is a lift of the Borisov-Joyce virtual cycle \cite{BJ} (constructed by Cao-Leung \cite{CL} in special cases) from homology groups to Chow groups.
The existence of orientations in the projective case was shown in \cite{CGJ} and in the quasi-projective case in \cite{Boj1}. 

The number $\vd$ is sometimes non-integer, in which case we define $A_{\vd}(\curly P)$ to be zero (notation and conventions). The following simple observation implies integrality of $\vd$:
\begin{lemma} \label{lem:intvd}
Suppose $\gamma = \sum_i m_i [S_i]$ where $m_i \in \mathbb{Z}$ and $S_i \subset X$ are irreducible smooth projective surfaces for which $N_{S_i / X}$ splits. Then $\vd \in \mathbb{Z}$. In particular, for $X$ a toric Calabi-Yau 4-fold and any $\gamma \in H^4_c(X,\Q)$ Poincar\'e dual to an integral algebraic cycle, we have $\vd \in \Z$.
\end{lemma}
\begin{proof}
It suffices to show that $(j_* [S_i])^2$ is even. By assumption, $N_{S_i / X} \cong L_1 \oplus L_2$, where $L_1 \otimes L_2 \cong K_S$ since $X$ is Calabi-Yau. Hence $(j_* [S_i])^2 = \int_{S_i} c_2(N_{S_i / X}) = L_1 \cdot (K_S - L_1) = 2(\chi(\O_{S_i})-\chi(L_1))$. The last statement of the lemma follows because any such $\gamma$ is Poincar\'e dual to an integral linear combination of $T_X$-fixed irreducible smooth projective surfaces $S_i \subset X$, and $N_{S_i / X}$ splits for any such $S_i$.
\end{proof}

\subsubsection{Equivariant case} \label{sec:equivcase} 

Finally, let $X$ be a Calabi-Yau 4-fold endowed with the action of an algebraic torus $T$ preserving the Calabi-Yau volume form. For example, if $X$ is a toric Calabi-Yau 4-fold, then it has an action of $T_X \cong (\mathbb{C}^*)^4$, which contains a 3-dimensional subtorus preserving the Calabi-Yau volume form. Note that a toric Calabi-Yau 4-fold is never projective. 

Assume there exists a $T$-equivariant good compactification, i.e., an irreducible smooth projective 4-fold $\Xbar$ with $T$-action and a $T$-equivariant open embedding $j : X \hookrightarrow \Xbar$ such that $\Xbar \setminus X$ is a simple normal crossings divisor. 
For example, a $T$-equivariant good compactification always exist for a toric Calabi-Yau 4-fold by turning its fan into a complete fan by adding sufficiently many 4-dimensional rational simplicial cones.\footnote{Any smooth quasi-projective toric variety has a toric completion to a projective toric variety \cite{Sum}, which in turn has a toric projective resolution \cite[Thm.~11.1.9]{CLS}.}

The results of this section hold $T$-equivariantly using the $T$-equivariant Atiyah class \cite{Ric}. 
The upshot is that we obtain classes 
\begin{align*}
&[\curly P]^{\vir} \in A_{\vd}^T(\curly P), \quad \vd := n - \frac{1}{2} (j_* \gamma)^2 \\
&\widehat{\O}_{\curly P}^{\vir} \in K_0^T(\curly P, \Z[\tfrac{1}{2}]).
\end{align*}

\subsection{Virtual localization} \label{sec:virtualloc}

In this section, we discuss the Oh-Thomas localization formula \cite{OT}, with special focus on the case when the fixed locus is isolated reduced.

Let $M$ be a quasi-projective scheme with 3-term symmetric obstruction theory $\phi : \EE \to \tau^{\geq -1} \LL_M$. 
Suppose $\phi$ is oriented and isotropic. Let $T$ be an algebraic torus acting on $M$ and suppose $\phi$ and the orientation are $T$-equivariant. Furthermore, suppose $\rk(\E) = 2\vd \in 2\Z$, then\footnote{If $\rk(\EE)$ is odd, then the virtual cycle and virtual structure sheaf are zero.}
$$
[M]^{\vir} \in A_{\vd}^T(M), \quad \widehat{\O}_M^{\vir}  \in K_0^T(M,\Z[\tfrac{1}{2}]).
$$
We will sometimes write $\widehat{\O}^{\vir}$ instead of $\widehat{\O}_M^{\vir}$. 

Let $\iota : M^T \hookrightarrow M$ denote the inclusion of the fixed point locus. The (general) localization formula of Thomason \cite[Thm.~2.1]{Tho}, see also \cite[Sect.~2.3.1]{Oko}, states$$
\iota_* : K_0^T(M^T) \otimes_{\Z[t,t^{-1}]} \Z(t) \stackrel{\cong}{\rightarrow} K_0^T(M) \otimes_{\Z[t,t^{-1}]} \Z(t),
$$
where $t$ denotes the equivariant parameters of $T$ and 
$$
K_0^T(M^T) \cong K_0(M^T) \otimes_{\Z} \Z[t,t^{-1}].
$$
Suppose $M^T$ is \emph{proper}. Then we can use this isomorphism to define
$$
\chi(M,\alpha), \quad \alpha \in K_0^T(M) \otimes_{\Z[t,t^{-1}]} \Z(t)
$$
by $\chi(M^T,\widetilde{\alpha})$, where $\widetilde{\alpha} \in K_0^T(M^T) \otimes_{\Z[t,t^{-1}]} \Z(t)$ is the unique class satisfying $\iota_* \widetilde{\alpha} = \alpha$. In general this procedure is not explicit. Fortunately, for $\alpha = \widehat{\O}_M^{\vir}$ we have the explicit virtual localization formula of Oh-Thomas, which we will now discuss. Their formula requires taking square roots of the equivariant parameters $t$ \cite[Sect.~7]{OT}, so we will work with 
$$
\iota_* : K_0^T(M^T) \otimes_{\Z[t,t^{-1}]} \Q(t^{\frac{1}{2}}) \stackrel{\cong}{\rightarrow} K_0^T(M) \otimes_{\Z[t,t^{-1}]} \Q(t^{\frac{1}{2}}) =: K_0^T(M)_{\loc}.
$$
\begin{theorem}[Oh-Thomas] \label{thm:OTloc}
In the above setup, and denoting the connected components of $M^T$ by $M_i$, we have
$$
\widehat{\O}_M^{\vir}  = \sum_i \iota_* \Bigg(\frac{\widehat{\O}_{M_i}^{\vir}}{\sqrt{\mathfrak{e}}(N^{\vir}|_{M_i})  }\Bigg).
$$
\end{theorem}

Let us explain this formula (see \cite[Sect.~7]{OT} for further details). On each component $M_i$, we have a decomposition into fixed and moving part
$$
\EE|_{M_i} = \EE|_{M_i}^f \oplus \EE|_{M_i}^m.
$$
We refer to $N^{\vir} = \EE^\vee|_{M^T}^m$ as the virtual normal bundle. The complex $\EE$ is $T$-equivariantly quasi-isomorphic to a $T$-equivariant self-dual complex
$$
\{F \stackrel{a}{\to} E \cong E^* \stackrel{a^*}{\to} F^*\},
$$
where $E$ is a $T$-equivariant vector bundle with \emph{$T$-invariant} non-degenerate quadratic form $q : E \cong E^*$. The induced map \cite{GP, CKL}
$$
\EE|_{M_i}^f \to (\tau^{\geq -1} \LL_M)|_{M_i}^f \to \tau^{\geq -1} \LL_{M_i}
$$
is a 3-term obtruction theory. 
Moreover, we see that  $(N^{\vir}|_{M_i})^{\vee}$, $\EE|_{M_i}^f$ have self-dual resolutions
\begin{align*}
(N^{\vir}|_{M_i})^{\vee} &\cong \{F|_{M_i}^m \to E|_{M_i}^m \to F^*|_{M_i}^m\}, \\
\EE|_{M_i}^f &\cong \{F|_{M_i}^f \to E|_{M_i}^f \to F^*|_{M_i}^f\}.
\end{align*}
Choosing a generic 1-parameter subgroup $\C^* \leq T$, the virtual normal bundle $N^{\vir}|_{M_i}$ can be oriented as follows. Serre duality provides a pairing between terms with weights $w>0$ and terms with weights $w<0$, which gives a splitting of $N^{\vir}|_{M_i}$. Since we oriented $\EE$ and $N^{\vir}|_{M_i}$, we obtain an induced orientation for $\EE|_{M_i}^f$ and we obtain objects $[M_i]^{\vir}$, $\widehat{\O}_{M_i}^{\vir}$. Finally, we have
$$
\sqrt{\mathfrak{e}}(N^{\vir}|_{M_i}) = \frac{\widehat{\Lambda}\udot (F^*|^m_{M_i})}{\sqrt{\mathfrak{e}}(E|_{M_i}^m)},
$$
where $\sqrt{\mathfrak{e}}(E|_{M_i}^m)$ is the $K$-theoretic square root Euler class defined in \cite{OT}. If $\Lambda$ is a positively oriented maximal isotropic subbundle of $E|_{M_i}^m$, then $\sqrt{\mathfrak{e}}(E|^m_{M_i}) = \widehat{\Lambda}\udot \Lambda^*$ and we have
\begin{equation*} 
\sqrt{\mathfrak{e}}(N^{\vir}|_{M_i}) = \widehat{\Lambda}\udot (F^* |_{M_i}^m - \Lambda^*).   
\end{equation*}

The following proposition was first proved in \cite{KR}.
\begin{proposition} \label{prop:OTisored}
Suppose $M^T$ consists of isolated reduced points. Then 
$$
2\mathrm{vd}(\{P\}) = \rk(\EE|_{P}^f) \leq 0, \quad \forall P \in M^T.
$$
Moreover
$$
\widehat{\O}_M^{\vir} = \sum_{\substack{P \in M^T \\ \vd(\{P\}) = 0}} (-1)^{\sigma_P} \iota_* \Bigg( \frac{1}{\sqrt{\mathfrak{e}}(\EE^\vee|_P)} \Bigg), 
$$
where $(-1)^{\sigma_P}$ is the orientation on $\{P\}$ induced by the globally chosen orientation on $M$ and
\begin{equation} \label{eqn:sqrtsqrd}
(\sqrt{\mathfrak{e}}(\EE^\vee|_P))^2 = (-1)^{\frac{\rk(\EE)}{2}} \widehat{\Lambda}\udot(\EE|_P).
\end{equation}
\end{proposition}
\begin{proof}
Consider the complex
\begin{equation} \label{eqn:EfixP}
\EE|_P^f \cong \{F|^f_P \stackrel{a}{\to} E|^f_P \cong E^*|^f_P \stackrel{a^*}{\to} F^*|^f_P\}.
\end{equation}
By reducedness, the second map is surjective and, consequently, the first map is injective. Moreover $F|^f_P \hookrightarrow E|^f_P$ is isotropic so $\dim_{\C}(F|^f_P) \leq \dim_{\C}(E|^f_P) / 2$ and hence $\mathrm{vd}(\{P\}) \leq 0$. 

When $\vd(\{P\}) < 0$, we have $[\{P\}]^{\vir} = 0$. Suppose $\vd(\{P\}) = 0$, then $F|^f_P \hookrightarrow E|^f_P$ is maximal isotropic and \eqref{eqn:EfixP} is acyclic, i.e., $\EE|_P^f \cong 0$. Hence $[\{P\}]^{\vir} = \pm \{P\}$, where the sign is determined by the choice of orientation.
Moreover, $N^{\vir}|_P = \EE^\vee|_P$ (as classes in $K_0^T(\pt)$) and \eqref{eqn:sqrtsqrd} follows from the definition of $K$-theoretic square root Euler class \cite[Prop.~5.4]{OT}.
\end{proof}
\begin{remark} \label{rem:KRsign}
In \cite{KR}, the following formula for the sign is derived when the orientation on $M$ is induced by a choice of maximal isotropic subbundle $\Lambda \subset E$. For $P \in M^T$ such that $\vd(\{P\})=0$, choose a splitting $E|_P = \Lambda|_P \oplus \Lambda^*|_P$ and denote projection by $p_\Lambda : E|_P \to \Lambda|_P$. Then
$$
\sigma_P \equiv \dim_{\C}\Big( \coker(p_\Lambda \circ a)^f \Big) \mod 2,
$$
where $a : F|^f_P \to E|^f_P$ as before. For the computations in this paper, we do not have an explicit self-dual resolution \eqref{eqn:EfixP}, so this formula is not enough to determine the sign.
\end{remark}

\section{Invariants and correspondences}

In this section, we define the virtual invariants which form the subject of this paper and we derive some of their structural properties. We also formulate conjectural $K$-theoretic $\DT$--$\PT_0$ correspondences. We first introduce all invariants in the projective case and then in the $T$-equivariant quasi-projective setup.

\subsection{Tautological invariants} \label{sec:tautinv}

Let $X$ be a projective Calabi-Yau 4-fold and fix a Chern character 
$$
v = (0,0,\gamma,\beta,n- \gamma \cdot \td_2(X)) \in H^*(X,\Q).
$$
We consider any of the $\PT_q$ moduli spaces $\curly P := \PTqvX$ and we denote the universal pair by
$$
\II\udot = [\O_{X \times \curly P} \to \FF].
$$
Consider the projections 
\begin{displaymath}
\xymatrix
{
& X \times \curly P \ar_{\pi_X}[dl] \ar^{\pi_{\curly P}}[dr] & \\
X & & \curly P.
}
\end{displaymath}
Fixing a line bundle $L \in \Pic(X)$, we are interested in the following tautological complex
\begin{equation} \label{def:tautrk}
R\pi_{\curly P *} (\FF \otimes L), \quad \rk := \rk(R\pi_{\curly P *} (\FF \otimes L)) = (\tfrac{1}{2} L \gamma + \beta)L + n,
\end{equation}
where we suppress the pull-back of $L$ from $X$ to $X \times \curly P$. We endow $\curly P$ with a trivial $\C^*$-action for which we denote a corresponding primitive character by $y$. We denote its equivariant Euler class by $m = c_1^{\C^*}(y)$. Writing $(R\pi_{\curly P *} (\FF \otimes L))^\vee = V_0 - V_1 \in K^0(\curly P)$, where $V_0, V_1$ are classes of vector bundles, we can therefore consider the equivariant Euler class
\[
e((R\pi_{\curly P *} (\FF \otimes L))^\vee \otimes y) := \frac{e(V_0 \otimes y)}{e(V_1 \otimes y)} \in A^*(\curly P)[m,m^{-1}],
\]
where 
\begin{equation} \label{eqn:equivEulerexp}
e(V_i \otimes y) = c_{\rk V_i}(V_i) + c_{\rk V_i - 1}(V_i) m + \cdots + m^{\rk(V_i)}.
\end{equation}

By \cite{OT} and \cite[Thm.~3.19]{BKP}, we have a 3-term symmetric obstruction theory \eqref{eqn:obs}
and we define its virtual tangent bundle by $T_{\curly P}^{\vir} := \EE^\vee$. Fixing an orientation, we obtain a virtual cycle 
\begin{align*}
[\curly P]^{\vir} \in A_{\vd}(\curly P), \quad \vd := n - \frac{1}{2}  \gamma^2.
\end{align*}

\begin{definition} \label{def:tautinv}
We define \emph{tautological $\PT_q$ invariants} of $X$ by
\begin{align*}
\langle L \rangle^{\PT_q}_{X,v} :=\int_{[\curly P]^{\vir}} e((R\pi_{\curly P *}(\FF \otimes L))^\vee \otimes y) \in \Z \cdot m^{\rk-\vd},
\end{align*}
where $\rk-\vd = \frac{1}{2}\gamma(\gamma+L^2) + \beta L$. We also write $\langle L\rangle^{\PT_q}_{X,\gamma,\beta,n}$ for this invariant.
\end{definition}

\begin{remark} \label{rem:EulervsChern}
By \eqref{eqn:equivEulerexp}, we have
\[
\langle L \rangle^{\PT_q}_{X,v} = m^{\rk - \vd} \int_{[\curly P]^{\vir}} c((R\pi_{\curly P *}(\FF \otimes L))^\vee),
\]
where $c(\cdot)$ denotes the total Chern class. 
\end{remark}

It is also reasonable to consider descendent insertions
\[
\tau_k(\sigma) := \pi_{\curly P *} \Big( \pi_X^* \sigma \cap \ch_{2+k}(\FF) \Big) \in H^{2(a+k-2)}(\curly P,\Q), \quad \sigma \in H^{2a}(X,\Q)
\]
and study \emph{descendent $\PT_q$ invariants} of $X$
\[
\langle \tau_{k_1}(\sigma_1) \cdots \tau_{k_\ell}(\sigma_{\ell}) \rangle_X^{\PT_q} := \int_{[\curly P]^{\vir}} \tau_{k_1}(\sigma_1) \cdots \tau_{k_\ell}(\sigma_{\ell}).
\]
We do not pursue descendent invariants in this paper.

\subsection{\texorpdfstring{$K$}{K}-theoretic invariants} \label{sec:Ktheorinv}

We continue using the notation from the previous section. We consider the virtual structure sheaf
\[
\widehat{\O}_{\curly P}^{\vir} \in K_0(\curly P,\Z[\tfrac{1}{2}]).
\]
Recall the notation for $\widehat{\Lambda}^\mdot(\cdot)$ from the introduction. We are interested in the \emph{symmetrized equivariant $K$-theoretic Euler class} of the tautological complex 
\begin{equation} \label{eqn:Kduality}
\widehat{\Lambda}^\mdot( R\pi_{ \curly P*}(\FF \otimes L) \otimes y^{-1})^\vee = (-1)^{\rk} \cdot \widehat{\Lambda}^\mdot( R\pi_{\curly P*}(\FF \otimes L) \otimes y^{-1}).
\end{equation}
More precisely, writing $R\pi_{\curly P*}(\FF \otimes L) = V_0 - V_1$ as before, we have 
\[
\widehat{\Lambda}^\mdot( R\pi_{\curly P*}(\FF \otimes L) \otimes y^{-1}) = y^{\frac{\rk}{2}} \frac{\det(V_1)^{\frac{1}{2}}}{\det(V_0)^{\frac{1}{2}}} \frac{\Lambda_{-y^{-1}}V_0}{\Lambda_{-y^{-1}}V_1} \in K^0(\curly P,\Z[\tfrac{1}{2}])(\!(y^{-\frac{1}{2}})\!).
\]
Note that identity \eqref{eqn:Kduality} only makes sense when the left hand side is expanded as formal power series in $y^{-1}$.

\begin{definition} \label{def:Ktheorinv}
We define \emph{$K$-theoretic invariants} of $X$ by
\begin{align*}
\langle\!\langle L\rangle\!\rangle^{\PT_q}_{X,v} &=\chi(\curly P, \Ohat_{\curly P} \otimes \widehat{\Lambda}^\mdot( R\pi_{\curly P*}(\FF \otimes L) \otimes y^{-1})) \in \Q(\!(y^{-\frac{1}{2}})\!).
\end{align*}
We also write $\langle\!\langle L\rangle\!\rangle^{\PT_q}_{X,\gamma,\beta,n}$ for this invariant.
\end{definition}

Recall that for any formal power series $f(x)$ over a commutative $\Q$-algebra $R$ such that $f(0) = 1$, we can consider the formal power series
\[
\sqrt{f}(x) := f^{\frac{1}{2}}(x) := e^{\frac{1}{2} \log f(x)} \in 1+ xR[[x]].
\]
It has the properties $(f^{\frac{1}{2}})^2  = (f^2)^{\frac{1}{2}}= f$ and $f^{\frac{1}{2}}g^{\frac{1}{2}} = (fg)^{\frac{1}{2}}$. An example relevant to us is $R = \Q$ and $f(x) = x / (1-e^{-x})$, which appears in the definition of the Todd class, i.e., for a line bundle $L$ we have
\[
\td(L) = \frac{c_1(L)}{1 - e^{-c_1(L)}}.
\]
We start with a general property which states that, at the level of virtual invariants, the moduli spaces $\curly P$ are themselves \emph{virtually Calabi--Yau}.\footnote{The results of this section hold for any moduli space endowed with 3-term symmetric isotropic obstruction theory.}
\begin{proposition} \label{prop:weakvirtualSD}
For $V \in K^0(\curly P)$, we have 
$$
\chi(\curly P,\widehat{\O}^{\vir}_{\curly P} \otimes V) = (-1)^{\vd} \chi(\curly P,\widehat{\O}^{\vir}_{\curly P} \otimes V^\vee).
$$
\end{proposition}
\begin{proof}
By \cite[Thm.~6.1]{OT} we have
$$
\chi(\curly P,\widehat{\O}^{\vir}_{\curly P} \otimes V) = \int_{[\curly P]^{\vir}} \ch(V) \sqrt{\td}(T_{\curly P}^{\vir})\,.
$$
On the other hand, using $T^{\vir}_{\curly P}[2] \cong (T_{\curly P}^{\vir})^\vee$, we have
$$
\chi(\curly P,\widehat{\O}^{\vir}_{\curly P} \otimes V^\vee) = \int_{[\curly P]^{\vir}} \ch(V^\vee) \sqrt{\td}((T_{\curly P}^{\vir})^{\vee})\,.
$$
Therefore, the second intersection number is obtained from the first by replacing all Chern roots $x$ by $-x$ and the result follows.
\end{proof}

\begin{remark}
This lemma provides an indication why $\chi(\curly P, \Ohat_{\curly P})$ is a bad invariant to consider: it is automatically zero when $\vd$ is odd. In fact, Bojko \cite{Boj2} showed (assuming a certain conjectural vertex algebra wall-crossing framework for virtual classes similar to \cite{GJT, Joy}) that for Hilbert schemes of points $$\chi(\Hilb^n(X), \Ohat_{\Hilb^n(X)}) = 0,$$ for all $n$.
\end{remark}

A priori, the $K$-theoretic invariants are Laurent series in $y^{-\frac{1}{2}}$. However, more can be said about their structure. We define
\[
[y]:=y^{\frac{1}{2}} - y^{-\frac{1}{2}}.
\]
A polynomial $f(y)$ of degree $n$ is called {\em palindromic} if $f(y)=y^nf(y^{-1})$.
\begin{proposition} \label{prop:palin}
There exists a Laurent polynomial $f(y) \in \Q[y^{\frac{1}{2}},y^{-\frac{1}{2}}]$ such that $y^{\frac{\vd}{2}} f(y)$ is a palindromic polynomial of degree $\leq \vd$ and
\[
\langle \! \langle L \rangle \! \rangle_{X,v}^{\PT_q} = [y]^{\rk - \vd} f(y).
\]
\end{proposition}
\begin{proof}
For a line bundle $L$, we have
\[
\ch(\widehat{\Lambda}\udot(L \otimes y^{-1})) 
= y^{\frac{1}{2}} e^{-\frac{c_1(L)}{2}} -  y^{-\frac{1}{2}} e^{\frac{c_1(L)}{2}}. 
\]
Let us write $R\pi_{\curly P *} (\FF \otimes L) = V_0 - V_1$ and $T_{\curly P}^{\vir} = E_0 - E_1$, where $V_0,V_1,E_0,E_1$ are vector bundles with Chern roots $\{v_i\}, \{w_j\}, \{x_k\}, \{ u_\ell\}$ respectively. Then
\begin{align} 
\begin{split} \label{eqn:Phi}
&\langle \! \langle L \rangle \! \rangle_{X,v}^{\PT_q} = \int_{[\curly P]^{\vir}} \Phi(y) \cdot \Bigg( \frac{\prod_k \frac{x_k}{1-e^{-x_k}}}{\prod_\ell \frac{u_\ell}{1-e^{-u_\ell}}} \Bigg)^{\frac{1}{2}}, \\
&\Phi(y) := \frac{\prod_i (y^{\frac{1}{2}}e^{-\frac{v_i}{2}} - y^{-\frac{1}{2}}e^{\frac{v_i}{2}})}{\prod_j (y^{\frac{1}{2}}e^{-\frac{w_j}{2}} - y^{-\frac{1}{2}}e^{\frac{w_j}{2}})}.
\end{split}
\end{align}
We rewrite $\Phi(y)$ as follows 
\[
\Phi(y) = [y]^{\rk}  \frac{\prod_i \Big(e^{\frac{v_i}{2}} - \frac{y^{\frac{1}{2}}}{[y]} (e^{\frac{v_i}{2}} - e^{-\frac{v_i}{2}}) \Big) }{\prod_j \Big(e^{\frac{w_j}{2}} - \frac{y^{\frac{1}{2}}}{[y]} (e^{\frac{w_j}{2}} - e^{-\frac{w_j}{2}}) \Big)}.
\]
Note that 
\[
e^{\frac{v_i}{2}} - \frac{y^{\frac{1}{2}}}{[y]} (e^{\frac{v_i}{2}} - e^{-\frac{v_i}{2}}) = 1+ \ldots - \frac{y^{\frac{1}{2}}}{[y]}(v_i+\ldots),
\]
where $\ldots$ are terms in $A^{\geq 1}(\curly P)$ and similarly for the terms involving $w_j$ in the denominator of $\Phi(y)$. Therefore, we can write
\[
\Phi(y) = [y]^{\rk} \sum_{p=0}^{\infty} [y]^{-p} y^{\frac{p}{2}} \Phi_p,
\]
where $\Phi_p \in A^{\geq p}(\curly P)$. Combined with \eqref{eqn:Phi}, where we integrate over a class of degree $\vd$, we obtain the proposition except for palindromicity. 

For the palindromic property, we consider 
\[
\int_{[\curly P]^{\vir}} \Phi(y^{-1}) \cdot \Bigg( \frac{\prod_k \frac{x_k}{1-e^{-x_k}}}{\prod_\ell \frac{u_\ell}{1-e^{-u_\ell}}} \Bigg)^{\frac{1}{2}}
\]
and use a similar strategy to the proof of Proposition \ref{prop:weakvirtualSD}. Indeed
\begin{align*}
\int_{[\curly P]^{\vir}} \Phi(y^{-1}) \Bigg( \frac{\prod_k \frac{x_k}{1-e^{-x_k}}}{\prod_\ell \frac{u_\ell}{1-e^{-u_\ell}}} \Bigg)^{\frac{1}{2}} = (-1)^{\rk} \int_{[\curly P]^{\vir}}  \frac{\prod_i (y^{\frac{1}{2}}e^{\frac{v_i}{2}} - y^{-\frac{1}{2}}e^{-\frac{v_i}{2}})}{\prod_j (y^{\frac{1}{2}}e^{\frac{w_j}{2}} - y^{-\frac{1}{2}}e^{-\frac{w_j}{2}})} \Bigg( \frac{\prod_k \frac{x_k}{1-e^{-x_k}}}{\prod_\ell \frac{u_\ell}{1-e^{-u_\ell}}} \Bigg)^{\frac{1}{2}} \\
= (-1)^{\rk} \int_{[\curly P]^{\vir}}  \frac{\prod_i (y^{\frac{1}{2}}e^{\frac{v_i}{2}} - y^{-\frac{1}{2}}e^{-\frac{v_i}{2}})}{\prod_j (y^{\frac{1}{2}}e^{\frac{w_j}{2}} - y^{-\frac{1}{2}}e^{-\frac{w_j}{2}})} \Bigg( \frac{\prod_k \frac{-x_k}{1-e^{x_k}}}{\prod_\ell \frac{-u_\ell}{1-e^{u_\ell}}} \Bigg)^{\frac{1}{2}},
\end{align*}
where we used $\det(T_{\curly P}^{\vir})^{\otimes 2} \cong \O_{\curly P}$, which implies $\prod_k e^{2 x_k} \prod_\ell e^{-2u_\ell} = 1$. Thus the integrand is obtained from \eqref{eqn:Phi} by simultaneously replacing all $v_i, w_j, x_k, u_\ell$ by $-v_i, -w_j, -x_k, -u_\ell$. Hence
\begin{align*}
\int_{[\curly P]^{\vir}} \Phi(y^{-1}) \cdot \Bigg( \frac{\prod_k \frac{x_k}{1-e^{-x_k}}}{\prod_\ell \frac{u_\ell}{1-e^{-u_\ell}}} \Bigg)^{\frac{1}{2}} = (-1)^{\rk+\vd} \int_{[\curly P]^{\vir}} \Phi(y) \cdot \Bigg( \frac{\prod_k \frac{x_k}{1-e^{-x_k}}}{\prod_\ell \frac{u_\ell}{1-e^{-u_\ell}}} \Bigg)^{\frac{1}{2}}, 
\end{align*}
which implies the palindromic property of $f$.
\end{proof}

The $K$-theoretic invariants are refinements of the tautological invariants.
\begin{proposition} \label{prop:limy=1}
We have
$$
\langle L \rangle_{X,v}^{\PT_q} = m^{\rk-\vd} \cdot \lim_{y \to 1} [y]^{\vd-\rk} \langle \! \langle L \rangle \! \rangle_{X,v}^{\PT_q}.
$$
\end{proposition}
\begin{proof}
By the definition of $\Phi(y), \Phi_p$ in the proof of Proposition \ref{prop:palin}, we have
\[
\lim_{y \to 1} [y]^{\vd-\rk} \langle \! \langle L \rangle \! \rangle_{X,v}^{\PT_q} = \int_{[\curly P]^{\vir}} \Phi_{\vd},
\]
where the square root Todd class does not contribute since $\deg[\curly P]^{\vir} = \vd$. The desired result is
\[
\int_{[\curly P]^{\vir}} \mathrm{Coeff}_{z^{\vd}} \Bigg( \frac{\prod_i \Big(e^{\frac{v_i}{2}} - z(e^{\frac{v_i}{2}} - e^{-\frac{v_i}{2}}) \Big) }{\prod_j \Big(e^{\frac{w_j}{2}} - z(e^{\frac{w_j}{2}} - e^{-\frac{w_j}{2}}) \Big)} \Bigg) = \int_{[\curly P]^{\vir}} \mathrm{Coeff}_{z^{\vd}} \Bigg( \frac{\prod_i (1 - z v_i) }{\prod_j (1 - z w_j)} \Bigg),
\]
where the equality follows from the fact that $\deg[\curly P]^{\vir} = \vd$. The right hand side is precisely
\[
\int_{[\curly P]^{\vir}} c((R\pi_{\curly P *} (\FF \otimes L))^\vee),
\]
which implies the result by Remark \ref{rem:EulervsChern}
\end{proof}

\begin{remark} 
Suppose we have a self-dual resolution $\EE \cong [F \stackrel{a}{\to} E \cong E^* \stackrel{a^*}{\to} F^*]$ such that $\rk(E)$ is even and $a(F) \subset \Lambda$ for some maximal isotropic subbundle $\Lambda$. Suppose furthermore that $\phi$ factors through the quotient $(\frac{1}{2} T_{\curly P}^{\vir})^\vee := [\Lambda^* \stackrel{a^*}{\to} F^*]$. 
Examples are when $\phi$ arises as a $(-2)$-shifted cotangent bundle of a scheme with perfect obstruction theory \cite[Sect.~8]{OT}. Then $\Ohat_{\curly P} \otimes (\det(\frac{1}{2} T_{\curly P}^{\vir}))^{\frac{1}{2}} \in K^0(\curly P)$, i.e., it is an \emph{integral} $K$-theory class by \cite[eqn.~(105)]{OT}. Thus, if furthermore there exists a line bundle $\curly L$ on $\curly P$ such that $$(\det(\tfrac{1}{2} T_{\curly P}^{\vir}) \otimes \det(R\pi_{\curly P *}(\FF \otimes L)))^{\frac{1}{2}} = \curly L \in K^0(\curly P),$$ then the coefficients of the Laurent polynomial in Proposition \ref{prop:palin} are integer (note that $[y]^{\rk - \vd}$ is an invertible element of $\Z(\!(y^{-\frac{1}{2}})\!)$ and its inverse has integer coefficients). This integrality is explored further in the sequels to this paper.
\end{remark}

\subsection{\texorpdfstring{$\DT$--$\PT_0$}{DT/PT0} correspondence}

We continue using the notation from the previous sections. We recall the three types of surface objects on $X$: $\DT:=\PT_{-1}$, $\PT_0$, $\PT_1$ (Definition \ref{def:stabilityconditionsofpairs'}). We note that $\Hilb^0(X) = \{\varnothing\}$, in which case $\rk = 0$ by \eqref{def:tautrk}, so $\langle \! \langle L \rangle \! \rangle_{X,0,0,0}^{\DT} = 1$ (for appropriate orientation).

Our main conjecture is as follows.
\begin{conjecture} \label{conj:KDTPT0inbodytext}
Let $X$ be a projective Calabi-Yau 4-fold and $L$ a line bundle on $X$. For any $\gamma \in H^4(X,\Q)$ and $\beta \in H^6(X,\Q)$, there exist orientations such that 
\begin{equation*} 
\frac{\sum_n \langle \! \langle L \rangle \! \rangle_{X,\gamma,\beta,n}^{\DT} q^n}{\sum_n \langle \! \langle L \rangle \! \rangle_{X,0,0,n}^{\DT} q^n} = \sum_n \langle \! \langle L \rangle \! \rangle_{X,\gamma,\beta,n}^{\PT_0} q^n.
\end{equation*}
\end{conjecture}

In fact, for fixed $\gamma, \beta$, the $\DT$ and $\PT_0$ moduli spaces are empty for $n \ll 0$ \cite[Prop.~2.14]{BKP}. Therefore, the series in the conjecture are all Laurent series in $q$. There are choices of $\gamma, \beta$ such that all $\PT_0$ pairs $(F,s)$ on $X$ with $\ch_2(F) = \gamma$ and $\ch_3(F) = \beta$ are also $\PT_1$ pairs, and vice versa. In this case, the $\PT_0 = \PT_1$ moduli spaces are empty for $n \gg 0$ \cite[Prop.~2.14]{BKP}. Then the right hand side of the formula is a Laurent \emph{polynomial} in $q$.

Using Proposition \ref{prop:limy=1}, and the fact that $\rk - \vd = \frac{1}{2} \gamma(\gamma+L^2) + \beta L$ is independent of $n$, we obtain the following corollary.  
\begin{corollary} \label{cor:cohoDTPT0}
Assume Conjecture \ref{conj:KDTPT0inbodytext} holds. Let $X$ be a projective Calabi-Yau 4-fold and $L$ a line bundle on $X$. For any $\gamma \in H^4(X,\Q)$ and $\beta \in H^6(X,\Q)$, there exist orientations such that 
$$
\frac{\sum_{n}  \langle L \rangle_{X,\gamma,\beta,n}^{\DT} q^n }{\sum_n \langle L \rangle_{X,0,0,n}^{\DT} q^n } = \sum_n \langle L \rangle_{X,\gamma,\beta,n}^{\PT_0} q^n.
$$
\end{corollary}

\begin{remark} \label{rem:nored}
In \cite{BKP}, we argued that for most ``interesting'' $(2,2)$ classes $\gamma \in H^4(X,\Q)$, one should not work with the Oh-Thomas virtual cycle but rather its \emph{reduced} version. We believe that Conjecture \ref{conj:KDTPT0inbodytext} holds with invariants defined via the reduced virtual cycle and reduced virtual structure sheaf as well. However, in this paper, the reduced class does not play a significant role for the following reasons:
\begin{itemize}
\item For the virtual Lefschetz principle in Section \ref{sec:Lefglob}, we consider a smooth connected effective divisor $Y \subset X$ and $\gamma$ pushed forward from a $(1,1)$-class on $Y$. Then reducing is often unnecessary by \cite[Prop.~4.37]{BKP}
\item When $X$ is a toric Calabi-Yau 4-fold, we work fully equivariantly with respect to the Calabi-Yau torus $T \cong (\C^*)^3$. The cosections of $\Ext_X^2(I\udot,I\udot)_0$ are typically $T$-equivariant and not $T$-invariant, so the toric theory is non-zero before reducing.
\end{itemize}
We refer to \cite{Sav} for reducing of the Borisov-Joyce cycle as an element of $H_*(\curly P, \Z)$.
\end{remark}

\subsection{Variations} \label{sec:vardef}

Let $X$ be a (quasi-projective) Calabi-Yau 4-fold. Fix $v = (0,0,\gamma,\beta,n-\gamma \cdot \td_2(X)) \in H_c^*(X,\Q)$. Then we have moduli spaces $\curly P := \curly P_v^{(q)}(X)$ and a virtual cycle and structure sheaf as discussed in Section \ref{sec:moduli} 
\[
[\curly P]^{\vir} \in A_*(\curly P), \quad \Ohat_{\curly P} \in K_0(\curly P,\Z[\tfrac{1}{2}])
\]
depending on a choice of orientation. Even when $X$ is not projective, it may happen that $\curly P$, which is a priori only quasi-projective (Proposition \ref{prop:qproj}), is in fact projective. The following proposition provides many such examples and is proved in essentially the same way as \cite[Prop.~3.1]{CKM2}:
\begin{proposition} \label{prop:localFanos}
\hfill
\begin{enumerate}
\item[$\mathrm{(1)}$] Let $Y$ be a smooth projective Fano 3-fold and $p : X = \mathrm{Tot}(K_Y) \to Y$ the total space of the canonical bundle. Then, for $q\in\{0,1\}$, all elements of $\curly P^{(q)}_v(X)$ are set theoretically supported on the zero section $Y \subset X$, hence $\curly P^{(q)}_v(X)$ is projective.
\item[$\mathrm{(2)}$] Let $S$ be a smooth projective surface and $p : X = \mathrm{Tot}(L_1 \oplus L_2) \to S$, where $L_1^{-1}$, $L_2^{-1}$ are ample line bundles on $S$ satisfying $L_1 \otimes L_2 \cong K_S$.  Then, for $q\in\{0,1\}$, all elements of $\curly P^{(q)}_v(X)$ are set theoretically supported on the zero section $S \subset X$, hence $\curly P^{(q)}_v(X)$ is projective.
\end{enumerate}
\end{proposition}

We assume the following: \\

\noindent \textbf{Assumption.} $\curly P$ is projective. \\

Then the tautological and $K$-theoretic invariants of Definitions \ref{def:tautinv}, \ref{def:Ktheorinv} are well-defined and all results from Sections \ref{sec:tautinv}, \ref{sec:Ktheorinv} hold. Moreover, we conjecture that the $\DT$--$\PT_0$ correspondence of Conjecture \ref{conj:KDTPT0inbodytext} holds in this setting. Finally, Corollary \ref{cor:cohoDTPT0} also holds in this setting. However, when $X$ is quasi-projective, assuming $\curly P$ to be projective is a very restrictive assumption. It almost always rules out the case $q=-1$ ($\DT$ case) due to points ``floating to infinity''.

\subsubsection{Equivariant case} 

Next, we suppose $X$ is a (quasi-projective) Calabi-Yau 4-fold with an action by an algebraic torus $T$ preserving the Calabi-Yau volume form and admitting a $T$-equivariant good compactification (e.g.~$X$ is a toric Calabi-Yau 4-fold). 
Fix $v = (0,0,\gamma,\beta,n-\gamma \cdot \td_2(X)) \in H_c^*(X,\Q)$ and consider $\curly P := \curly P_v^{(q)}(X)$ for $q\in\{-1,0,1\}$. Let $L$ be a $T$-equivariant line bundle on $X$.
We \emph{no longer} require $\curly P$ to be projective. As discussed in Section \ref{sec:equivcase}, we then have an equivariant virtual cycle and structure sheaf 
\[
[\curly P]^{\vir} \in A_*^T(\curly P), \quad \Ohat_{\curly P} \in K_0^T(\curly P,\Z[\tfrac{1}{2}]).
\]
The action of $T$ on $X$ lifts to $\curly P$. As we will see in Section \ref{sec:vertex}, there are many cases where $\curly P$ is not projective, but the the fixed locus $\curly P^T$ is projective. Let us suppose this is the case: \\

\noindent \textbf{Assumption.} $\curly P^T$ is projective. \\

As discussed in Section \ref{sec:virtualloc}, the $K$-theoretic invariants of Definition \ref{def:Ktheorinv} then make sense by the Thomason localization formula
\[
\langle\!\langle L \rangle\!\rangle_{X,v}^{\PT_q} \in K_0^T(\mathrm{pt})_{\mathrm{loc}}(\!(y^{-\frac{1}{2}})\!),
\]
where 
\[
K_0^T(\mathrm{pt})_{\mathrm{loc}}(\!(y^{-\frac{1}{2}})\!) = K_0^T(\mathrm{pt}) \otimes_{\Z[t,t^{-1}]} \Q(t^{\frac{1}{2}})(\!(y^{-\frac{1}{2}})\!) \cong \Q(t^{\frac{1}{2}})(\!(y^{-\frac{1}{2}})\!)
\]
is the localized equivariant $K$-group of a point, $t$ denote the equivariant parameters of $T$, and we added the formal variable $y^{-\frac{1}{2}}$. Similarly, using localization on Chow groups \cite{EG}, one can define tautological invariants
\[
\langle L \rangle_{X,v}^{\PT_q} \in A_*^T(\mathrm{pt})_{\mathrm{loc}}(\!(m^{-1})\!),
\]
where 
\[
A_0^T(\mathrm{pt})_{\mathrm{loc}}(\!(m^{-1})\!) = A_0^T(\mathrm{pt}) \otimes_{\Q[s]} \Q(s)(\!(m^{-1})\!) \cong \Q(s)(\!(m^{-1})\!),
\]
and $s,m$ are the equivariant first Chern classes of the equivariant parameters $t,y$ respectively.
However, one can no longer use the results of Sections \ref{sec:tautinv} and \ref{sec:Ktheorinv}, which involved intersection-theoretic arguments on $\curly P$. In the case $X$ is toric and $\curly P^T$ is isolated reduced, we show in a different way that $K$-theoretic invariants specialize to tautological invariants in Section \ref{sec:coholimittaut}.\footnote{Actually, we will derive stronger specializations at the level of vertex, edge, and face.}

As reviewed in Section \ref{sec:virtualloc}, fixing an orientation on $\curly P$, one may apply the Oh-Thomas localization formula (Theorem \ref{thm:OTloc}). 

In Section \ref{sec:vertex}, we will focus on the setting where $\curly P^T$ consists of isolated reduced points. Then Proposition \ref{prop:OTisored} applies. As indicated in Remark \ref{rem:KRsign}, we currently do not have a way to determine the signs $(-1)^{\sigma_P}$ for any given (global) choice of orientation (except in the case $\curly P = \Hilb^n(\C^4)$ studied by the second-named author and Rennemo \cite{KR}). Therefore, in Proposition \ref{prop:OTisored}, we allow \emph{any} choice of signs $(-1)^{\sigma_P}$.  This means we do \emph{not} require the signs $(-1)^{\sigma_P}$ to be induced from a global choice of orientation, and allow any sign choice.  Put differently: we define equivariant $K$-theoretic invariants by the Oh-Thomas localization formula:

\begin{definition} \label{def:equivKtheorinv}
We define \emph{equivariant $K$-theoretic invariants} by the Oh-Thomas localization formula
$$
\langle\!\langle L \rangle\!\rangle_{X,v}^{\PT_q} := \sum_{\substack{P \in \curly P^T \\ \vd(\{P\}) = 0}} (-1)^{\sigma_P} \frac{\widehat{\Lambda}\udot ( R\pi_{\curly P*}(\FF \otimes L)|_P \otimes y^{-1})}{\sqrt{(-1)^{\frac{\vd}{2}} \widehat{\Lambda}\udot( T_{\curly P}^{\vir}|_P^\vee)}} \in K_0^T(\pt)_{\loc}(y^{\frac{1}{2}}), 
$$
where the combination $(-1)^{\sigma_P} \sqrt{(\cdot)|_P}$ means that for each fixed point $P$ we make a choice of square root of $(\cdot)|_P$. This definition depends on a choice of square root attached to each fixed point $P$, which is suppressed from the notation of $\langle\!\langle L \rangle\!\rangle_{X,v}^{\PT_q}$.
\end{definition}

More concretely, this can be rewritten as
$$
\langle\!\langle L \rangle\!\rangle_{X,v}^{\PT_q} = \sum_{\substack{I\udot = [\O_X \to F] \in \curly P^T \\ \vd(\{I\udot\}) = 0}} (-1)^{\sigma_{I\udot}} \frac{\widehat{\Lambda}\udot (R\Gamma(X,F \otimes L) \otimes y^{-1})}{\sqrt{(-1)^{\frac{\vd}{2}} \widehat{\Lambda}\udot(R\Hom_X(I\udot,I\udot)_0[3] )}} \in \Q(t^{\frac{1}{2}},y^{\frac{1}{2}}).
$$
Similarly, we define \emph{equivariant tautological invariants}
\begin{align*}
\langle L \rangle_{X,v}^{\PT_q} &:= \sum_{\substack{P \in \curly P^T \\ \vd(\{P\}) = 0}} (-1)^{\sigma_P} \frac{e((R\pi_{\curly P*}(\FF \otimes L))|_P^\vee \otimes y)}{\sqrt{(-1)^{\frac{\vd}{2}} e( T_{\curly P}^{\vir}|_P})} \\
&= \sum_{\substack{I\udot = [\O_X \to F] \in \curly P^T \\ \vd(\{I\udot\}) = 0}} (-1)^{\sigma_{I\udot}} \frac{e(R\Gamma(X,F \otimes L)^\vee \otimes y)}{\sqrt{(-1)^{\frac{\vd}{2}} e(R\Hom_X(I\udot,I\udot)_0[1])}},
\end{align*}
which is an element of $A_*^T(\pt)_{\loc}(m) = \Q(s,m)$.

Although the above invariants are liberally defined, we of course expect that geometrically meaningful choices of square roots are induced by a global orientation. We then conjecture the following equivariant $K$-theoretic $\DT$--$\PT_0$ correspondence:
\begin{conjecture} \label{conj:KDTPT0inbodytext_equiv}
Let $X$ be a Calabi-Yau 4-fold with action by an algebraic torus $T$ preserving the Calabi-Yau volume form. Let $L$ be a $T$-equivariant line bundle on $X$. Let $\gamma \in H^4_c(X,\Q)$, $\beta \in H^6_c(X,\Q)$, and suppose $\curly P_{\gamma,\beta,n}^{(q)}(X)^T$ consists of isolated reduced points for all $n$ and $q \in \{-1,0\}$. Then there are choices of square roots associated to all $T$-fixed points such that
\begin{equation*} 
\frac{\sum_n \langle \! \langle L \rangle \! \rangle_{X,\gamma,\beta,n}^{\DT} q^n}{\sum_n \langle \! \langle L \rangle \! \rangle_{X,0,0,n}^{\DT} q^n} = \sum_n \langle \! \langle L \rangle \! \rangle_{X,\gamma,\beta,n}^{\PT_0} q^n
\end{equation*}
as an identity of Laurent series in $q$ with coefficients in $\Q(t^{\frac{1}{2}},y^{\frac{1}{2}})$.
\end{conjecture}

\section{Virtual Lefschetz principle} 

In this section, we present two virtual Lefschetz principles for a smooth connected 3-fold $Y$ in a Calabi-Yau 4-fold $X$: one in the case $X$ is projective, and one in the case $X = \mathrm{Tot}(K_Y)$ and $Y$ is toric.

\subsection{Via virtual pull-back} \label{sec:Lefglob}

For any smooth projective variety $X$ and smooth connected effective divisor $\iota : Y \hookrightarrow X$, we can push-forward a $\PT_q$ pair $(F,s)$ on $Y$ to a $\PT_q$ pair on $X$ as follows
\[
\O_X \to \iota_* \O_Y \stackrel{\iota_* s}{\to} \iota_* F.
\]
Fixing a Chern character $v \in H^*(X,\Q)$, we use the notation
\[
\curly P_v^{(q)}(Y) := \bigsqcup_{\iota_* v' = v \cdot \td(\O_X(Y))} \curly P_{v'}^{(q)}(Y),
\]
where the union is over all $v' \in H^*(Y,\Q)$ satisfying $\iota_* v' = v \cdot \td(\O_X(Y))$.
Then we have a closed embedding
\[
\imath = \iota_* : \curly P_v^{(q)}(Y) \hookrightarrow \curly P_{v}^{(q)}(X),
\]
where the Chern characters match by Grothendieck-Riemann-Roch. We denote the universal pairs of $\curly P_v^{(q)}(Y)$ and $\curly P_{\iota_* v}^{(q)}(X)$ by $\II\udot_Y = [\O \to \FF_Y]$ and $\II\udot = [\O \to \FF]$ respectively. 
The main result of this section is the following virtual Lefschetz principle.

\begin{theorem} \label{thm:Lefschetzviapullback}
Let $X$ be a projective Calabi-Yau 4-fold and $\iota : Y \hookrightarrow X$ a smooth connected effective divisor. Let $q \in \{-1,0,1\}$ and $v \in H^*(X,\Q)$. Let $\curly P_Y := \curly P_v^{(q)}(Y)$, $\curly P_X := \curly P_{v}^{(q)}(X)$, and $L:= \O_X(Y)$. Assume the following:
\begin{enumerate}
\item [$\mathrm{A1)}$] The tautological complex ${\curly L}:=R \pi_{\curly P_X *}(\FF \otimes L)$ is a vector bundle concentrated in degree 0.
\item [$\mathrm{A2)}$] The canonical map
\begin{equation*}
\phi_Y: \EE_Y:=R\hom_{\pi_{\curly P_Y}}(\II_Y^\mdot,\II_Y^\mdot \otimes K_Y)_0[2] \to \LL_{\curly P_Y}   
\end{equation*}
induced by the Atiyah class of $\II_Y^\mdot$ is a perfect, i.e.~2-term, obstruction theory on $\curly P_Y$.
\item [$\mathrm{A3)}$] The moduli space $\curly P_Y$ is connected.
\end{enumerate}
Then we have
\[\imath_* \Ohat_{\curly P_Y} = \widehat{\Lambda}\udot (R \pi_{\curly P_X *}(\FF \otimes L) )^\vee \otimes  \Ohat_{\curly P_X},\]
where $\imath : \curly P_Y \hookrightarrow \curly P_X$ is the inclusion map.
\end{theorem}

\begin{remark}
In the setting of the theorem, we also have the following cohomological virtual Lefschetz principle
\[\imath_*[\curly P_Y]^\vir = e(R \pi_{\curly P_X *}(\FF \otimes L)) \cap [\curly P_X]^\vir.\]
\end{remark}

We discuss the assumptions in the theorem.
\begin{itemize}
\item For $q=-1,0$ and $\gamma = 0$, condition A2) is automatic and $\curly P_Y$ are the familiar $\DT, \PT$ spaces. In this case, Theorem \ref{thm:Lefschetzviapullback} was proved by the third-named author \cite{Par}. 
\item For $q=-1$ and $\gamma = \beta = 0$, conditions A1), A2), A3) are all satisfied (Hilbert schemes of points on a connected quasi-projective scheme are connected). 
%E.g.~Bertin.
In this case Theorem \ref{thm:Lefschetzviapullback} was proved by the third-named author \cite{Par}, which establishes a conjecture of the second-named author and Cao \cite{CK1}.
\item Suppose for all closed points $(F,s)$ of $\curly P_X$, we have $H^{>0}(X,F \otimes L) = 0$. Then A1) is satisfied by cohomology and base change. 
\item The assumption A2) can be divided into two conditions: 
1.~$\phi_Y$ is an obstruction theory, 
and 2.~$\EE_Y$ is a $2$-term perfect complex (i.e.~has tor-amplitude $[-1,0]$).
The first part is not automatic for $\gamma \neq 0$ since the codimension of pairs in $Y$ is $1\not\geq 2$.
\item We show below that if $H^{>0}(X,F \otimes L) = 0$ for all closed points $(F,s)$ of $\curly P_X$ and $Y$ is Calabi-Yau, then A2) also follows. 
\item Assumption A3) is a technical assumption related to orientations. If we drop assumption A3), then we obtain the following weaker conclusion: for each connected component $\curly P_Y^i$ of $\curly P_Y$, there exists a sign $(-1)^{\sigma(i)}$ determined by the choice of orientation on $\curly P_X$, such that 
\begin{equation} \label{eqn:conncomppus}
\sum_{i}(-1)^{\sigma(i)} \imath_* \Ohat_{\curly P^i_Y}  =  \widehat{\Lambda}\udot (R \pi_{\curly P_X *}(\FF \otimes L))^\vee \otimes  \Ohat_{\curly P_X}.
\end{equation}
We expect that for any choice of orientation on $\curly P_X$, the signs $(-1)^{\sigma(i)}$ are independent of $i$.
\end{itemize} 

\begin{proof}[Proof of Theorem \ref{thm:Lefschetzviapullback}] This is a variant of the proof of \cite[Thm.~3.1]{Par}. (See also \cite[Cor.~B.4]{Par} for the $K$-theoretical statement.)
Here we give a sketch proof focussing on the differences with loc.~cit.
We use the following notation.
Both the projections $\curP_Y\times Y, \curP_Y \times X \to \curP_Y$ are denoted by $\pi$.
We let
\[\II_X\udot:=\II\udot|_{\curP_Y\times X},\qquad \EE_X:=\EE|_{\curP_Y}= R\hom_\pi(\II\udot_X,\II\udot_X)_0[3],\]
where $\mathbb{E}$ was defined in \eqref{eqn:obs}. Also note that $d_*\FF_Y\cong \FF|_{\curP_Y\times X}$ 
for the inclusion  $d: \curP_Y \times Y \hookrightarrow \curP_Y \times X$.

We first observe that  $\curly P_Y$ is the scheme-theoretical zero locus of the tautological section
$\tau\in\Gamma(\curP_X,\curL)$ induced by  $\O_{{\curly P_X}\times X}(-Y) \hookrightarrow \O_{{\curly P_X}\times X} \xrightarrow{s} \FF$
\[\xymatrix{
& \curL \ar[d]\\
\curP_Y \cong Z(\tau) \ar@{^{(}->}[r]^-{\imath} & \curP_X. \ar@/_1pc/[u]_{\tau}
}\] 
This is an analog of \cite[Prop.~3.6]{Par}.
The only difference is when we apply \cite[Lem.~3.7]{Par}, we use the following simple fact on a $\PT_q$ pair $(F,s)$ on $X$:
if $G$ is a $q$-dimensional sheaf on $X$, then $\Hom(G,F)=0$.
(An alternative proof is to use the derived moduli of pairs (cf.~\cite[Rem.~3.9]{BKP}) based on \cite[Lem.~3.7]{BKP}. Then we can show that the relative cotangent complex for the canonical map $R\curP_Y \to Z_{R\curP_X}(R\tau)$ is concentrated in degrees $\leq -2$ by the argument in \cite[Rem.~3.8]{Par} and thus its classical truncation is an isomorphism.)

We then note that there exists a $3$-term symmetric obstruction theory \[\overline{\phi}_Y :=\phi_Y \circ \epsilon : \overline{\EE}_Y \to \EE_Y \to \LL_Y\] given by some $3$-term symmetric complex $\overline{\EE}_Y$ and a self-dual exact triangle
\[\xymatrix{
\EE_Y\dual[2] \ar[r]^-{\epsilon\dual[2]} & \overline{\EE}_Y\dual[2]\cong\overline{\EE}_Y \ar[r]^-{\epsilon} & \EE_Y}\]
as in \cite[Prop.~3.10]{Par}.
Indeed, we have a canonical self-dual exact triangle
\[\xymatrix{
R\hom_{\pi}(\II\udot_Y,\II\udot_Y) [4]\ar[r] & R\hom_{\pi}(d_*\II\udot_Y,d_*\II\udot_Y)[3] \ar[r] & R\hom_{\pi}(\II\udot_Y,\II\udot_Y\otimes L) [2]
}\]
given by the Serre duality for the divisor $d: \curP_Y \times Y \hookrightarrow \curP_Y \times X$ (cf. \cite[Lem.~3.11]{Par}).
The desired $3$-term symmetric complex $\overline{\EE}_Y$ can be obtained from the above triangle after removing the traces in the left and right terms.
This is the technical part but can be done by the arguments in \cite[Lem.~3.12]{Par} after checking that
\begin{align*}
&\Hom_{{\curly P}_Y\times X} (d_*\O_{{\curly P}_Y\times Y},R\hom_{d}(\II\udot_Y,\II\udot_Y\otimes L)[-2])\\
&=\Hom_{{\curly P}_Y\times Y} (\II\udot_Y,\II\udot_Y\otimes L[-2]) \oplus \Hom_{{\curly P}_Y\times Y} (\II\udot_Y\otimes L\dual,\II\udot_Y\otimes L[-3]) =0.
\end{align*}
Consequently, the Behrend-Fantechi virtual cycle associated to the $2$-term perfect obstruction theory $\phi_Y$ is equal to the Oh-Thomas virtual cycle associated to the $3$-term symmetric obstruction theory $\overline{\phi}_Y$ (up to sign) by the reduction formula \cite[Prop.~1.18]{Par} (see \cite[B.2]{Par} for the $K$-theoretical version)
\[(\Ohat_{\curP_Y})^{\mathrm{OT}}_{\overline{\phi}_Y} = (\Ohat_{\curP_Y})^{\mathrm{BF}}_{\phi_Y} \in K_0(\curP_Y,\mathbb{Z}[\tfrac{1}{2}]),\]
where $(\Ohat_{\curP_Y})^{\mathrm{BF}}_{\phi_Y}:= (\O_{\curP_Y}^\vir)^{\mathrm{BF}}_{\phi_Y}\cdot (\det(\EE_Y))^{\frac{1}{2}}$ is the twisted virtual structure sheaf.

We show that the two $3$-term symmetric complexes $\EE_X$ and $\overline{\EE}_Y$ 
fit into a reduction diagram (in the sense of \cite[Prop.~1.5]{Par})
\[\xymatrix{
\DD\dual[2] \ar[r]^-{\alpha\dual[2]} \ar[d]^-{\beta\dual[2]} & \overline{\EE}_Y\dual[2]\cong\overline{\EE}_Y \ar[r]^-{\delta} \ar[d]^-{\alpha} & \curL|_{\curly P_Y}\dual[1] \ar@{=}[d] \\
\EE_X\dual[2]|_{\curly P_Y}\cong\EE_X|_{\curly P_Y} \ar[r]^-{\beta} & \DD \ar[r]^-{\gamma} & \curL|_{\curly P_Y}\dual[1],
}\]
for some perfect complex $\DD$ and maps $\alpha,\beta,\gamma,\delta$.
Indeed, it is straightforward to form a canonical diagram
\[\xymatrix{
R\hom_\pi(d_*\II\udot_Y,\II\udot_X)[3] \ar[r] \ar[d] & R\hom_{\pi}(d_*\II\udot_Y,d_*\II\udot_Y)[3] \ar[r] \ar[d] & R\hom_\pi(d_*\II\udot_Y,\O_{\curP_Y}\boxtimes L\dual)[4] \ar@{=}[d]\\
R\hom_\pi(\II\udot_X,\II\udot_X)_0[3] \ar[r] & R\hom_\pi(\II\udot_X,d_*\II\udot_Y)[3] \ar[r] &R\hom_\pi(d_*\II\udot_Y,\O_{\curP_Y}\boxtimes L\dual)[4]
}\]
where the rows are exact and the right square commutes as in \cite[3.14]{Par}.
Again the technical part is removing some traces,
but this can be done by the arguments in \cite[Prop.~3.13]{Par}.
A necessary additional ingredient is:
\begin{align*}
&\Hom_{\curP_Y \times X}(d_*\FF_Y,d_*\O_{\curP_Y \times Y})\\
&= \Hom_{\curP_Y \times Y}(\FF_Y,\O_{\curP_Y\times Y})\oplus \Hom_{\curP_Y \times Y}(\FF_Y,\O_{\curP_Y\times Y}\otimes L[-1])=0.
\end{align*}
This follows since the underlying sheaves of $\PT_q$ pairs on $Y$ are still torsion.

Finally, the $3$-term symmetric obstructions theories $\phi$ and $\overline{\phi}_Y$ are compatible
and thus we can apply the virtual pullback formula \cite[Thm.~2.2]{Par}.
Indeed, we have a commutative diagram
\[\xymatrix{
R\hom_\pi(d_*\II\udot_Y,\II\udot_X)[3] \ar[r] \ar[d]& R\hom_\pi(d_*\II\udot_Y,d_*\II\udot_Y)[3] \ar[r] \ar[rdd]^{\At_{\curP_Y\times X}(d_*\II\udot_Y)}  & R\hom_\pi(\II\udot_Y,\II\udot_Y\otimes L)[2] \ar[dd]^{\At_{\curP_Y\times Y}(\II\udot_Y)} 
\\
R\hom_\pi(\II\udot_X,\II\udot_X)[3] \ar[d]^{\At_{\curP_X\times X}(\II\udot)} \ar[rrd]^{\At_{\curP_Y\times X}(\II\udot_X)} && \\
\LL_{\curP_X}|_{\curP_Y} \ar[rr] && \LL_{\curP_Y}
}\]
by the basic properties of Atiyah classes (cf.~\cite[3.4.5]{Par}).
By removing the traces, we have the compatibility condition (in Definition \cite[Def.~2.1]{Par}) and thus we have the $K$-theoretic virtual pullback formula \cite[Thm.~B.3]{Par}
\[\Ohat_{\curP_Y} = \widehat{\imath}^!\Ohat_{\curP_X} \in K_0(\curP_Y,\Z[\tfrac{1}{2}]),\]
where $\imath_*\widehat{\imath}^!(-) = (-)\cdot \widehat{\Lambda}^{\bullet}(\curL)\dual$ by \cite[Lem.~4.11, Lem.~B.6]{Par}.
\end{proof}

The following cohomological conditions guarantee that A1), A2) are satisfied.
\begin{proposition} \label{prop:A123}
Let $X$ be a projective Calabi-Yau 4-fold and $\iota : Y \hookrightarrow X$ a smooth connected effective divisor. Let $q \in \{-1,0,1\}$ and $v \in H^*(X,\Q)$. Let $\curly P_Y := \curly P_v^{(q)}(Y)$, $\curly P_X := \curly P_{v}^{(q)}(X)$, and $L:= \O_X(Y)$. Assume the following:
\begin{enumerate}
\item[$\mathrm{(1)}$] $H^{>0}(X,F \otimes L) = 0$ for all closed points $(F,s) \in \curly P_X$.
\item[$\mathrm{(2)}$] $H^2(Y,F) = 0$ for all closed points $(F,s) \in \curly P_Y$. 
\end{enumerate}
Then A1) and A2) of Theorem \ref{thm:Lefschetzviapullback} are satisfied. 
\end{proposition}
\begin{proof}
A1) is satisfied by cohomology and base change. Denote by $\II_Y\udot = [\O \to \FF_Y]$ the universal complex on $Y \times \curly P_Y$ and consider the complexes
\begin{align*}
\EE_Y &:= (R\hom_{\pi_{\curly P_Y}}(\II_Y\udot, \II_Y\udot)_0[1])^\vee, \\
\EE'_Y &:= (R\hom_{\pi_{\curly P_Y}}(\II_Y\udot, \FF_Y) )^\vee.
\end{align*}
There is a natural exact triangle
\begin{equation} \label{eqn:triangleonY}
R\hom_{\pi_{\curly P_Y}}(\II_Y\udot, \FF_Y) \to R\hom_{\pi_{\curly P_Y}}(\II_Y\udot, \II_Y\udot)_0[1] \to R\hom_{\pi_{\curly P_Y}}(\FF_Y,\O)[2]
\end{equation}
which induces a map $\EE_Y \to \EE'_Y$. This map fits in the commutative diagram
\begin{displaymath}
\xymatrix
{
\EE_Y \ar[r] \ar_{\psi}[dr] & \EE'_Y \ar^{\psi'}[d] \\
 & \LL_{\curly P_Y},
}
\end{displaymath}
induced by the canonical map from the derived moduli of pairs on $Y$ to the derived moduli of perfect complexes on $Y$, see e.g.~\cite[Rem.~3.9]{BKP}.
In particular, $\psi'$ is an obstruction theory. Since $Y \times \curly P_Y$ is projective (Proposition \ref{prop:qproj}), the complex $\II_Y\udot$ is perfect. Therefore, the reasoning in the proof of \cite[Lem.~2.10]{PT1} shows that $\EE_Y$ is a perfect complex. Since $K_Y = L|_Y$, condition (1) implies $R^{>0} \pi_{\curly P_Y*}(\FF_Y \otimes K_Y) = 0$. The dual of \eqref{eqn:triangleonY} and the commutative diagram imply that $\psi$ is an obstruction theory.

Next, we have to show that $\EE_Y$ is 2-term, for which it suffices to prove that $\Ext^i_Y(I_Y\udot,I_Y\udot)_0 = 0$ for $i \neq 1,2$ for any $I_Y\udot = [\O_Y \to F] \in \curly P_Y$. Consider the absolute version of \eqref{eqn:triangleonY}, i.e.
\begin{align*}
R\Hom_Y(I_Y\udot,F) \to R\Hom_Y(I_Y\udot, I_Y\udot)_0[1] \to R\Hom_Y(F,\O_Y)[2].
\end{align*}
Using the fact that $I_Y\udot$ is a $\PT_q$ pair, it is easy to show that $\Ext^{<0}_Y(I_Y\udot, F) = 0$ (see \cite[Lem.~3.20]{BKP} and its proof). Moreover, $\Ext_Y^{<2}(F,\O_Y) \cong H^{>1}(Y,F \otimes K_Y)^* \cong H^{>1}(X,F \otimes L)^* \cong 0$, thus $\Ext^{<1}_Y(I_Y\udot,I_Y\udot)_0 = 0$.  Using the same exact triangle twisted by $K_Y$, and condition (2), we find $\Ext^{>2}_Y(I_Y\udot,I_Y\udot)_0 = 0$.
\end{proof}

We now relate the \emph{surface} counting on the smooth projective 3-fold $Y$ to \emph{curve} counting on $Y$. 
First, we recall curve counting on $Y$. For $q \in \{-1,0\}$, $v = (0,0,\beta,n - \beta \cdot \td_1(Y)) \in H^*(Y,\Q)$, we denote by $\curly P_{v}^{(q)}(Y)$ 
the moduli space parametrizing pairs $(F,s)$ such that $F$ is a coherent sheaf on $Y$ satisfying $\ch(F) = v$, and $s : \O_Y \to F$ is a section satisfying the following $\PT_q$ stability condition
\[
F \in \Coh_{\geq q+1}(Y) \and Q:=\coker(\O_Y \xrightarrow{s} F) \in \Coh_{\leq q}(Y).
\]
It is well-known that $\curly P := \curly P_v^{(q)}(Y)$ is a fine projective moduli space. For $q=-1$, these are the Hilbert schemes of curves on $Y$ \cite{Gro}, and for $q=0$, they are the moduli spaces of Pandharipande-Thomas pairs \cite{Pot1, Pot2, PT1}. Moreover, $\curly P$ has a perfect obstruction theory with accompanying virtual cycle and twisted virtual structure sheaf \cite{BF, Tho1, PT1}
\begin{align*}
&[\curly P]^{\vir} \in A_{\vd}(\curly P), \quad \vd = \beta \cdot c_1(Y) \\
&\Ohat_{\curly P} := \O^{\vir}_{\curly P} \otimes (\det(T_{\curly P}^{\vir}))^{-\frac{1}{2}} \in K_0(\curly P,\Z[\tfrac{1}{2}]),
\end{align*}
where $(\det(T_{\curly P}^{\vir}))^{-\frac{1}{2}}$ is the Nekrasov-Okounkov twist \cite{NO}.\footnote{Existence of square roots of $\det(T_{\curly P}^{\vir})$, as a line bundle, is discussed in \cite[Sect.~6]{NO}, but does not concern us here in this section.} The $K$-theoretic $\PT_q$ invariants are defined by
\[
\langle \! \langle 1 \rangle \! \rangle_{Y,v}^{\PT_q} := \chi\big(\curly P, \Ohat_{\curly P} \big) \in \Q.
\]
We also denote these invariants by $\langle \! \langle 1 \rangle \! \rangle_{Y,\beta,n}^{\PT_q}$ and we sometimes write $\DT := \PT_{-1}$ and $\PT := \PT_{0}$. For background, we refer to \cite{MNOP, PT2, Oko, NO}. Then the $K$-theoretic $\DT$--$\PT$ correspondence states \cite{NO}
\begin{equation} \label{eqn:3DKtheorDTPT}
\frac{\sum_n \langle \! \langle 1 \rangle \! \rangle_{Y,\beta,n}^{\DT} q^n}{\sum_n \langle \! \langle 1 \rangle \! \rangle_{Y,0,n}^{\DT} q^n} = \sum_n \langle \! \langle 1 \rangle \! \rangle_{Y,\beta,n}^{\PT} q^n.
\end{equation}
For $Y$ a Calabi-Yau 3-fold, the virtual dimension of $\curly P$ is zero, and we have $\langle \! \langle 1 \rangle \! \rangle_{Y,\beta,n}^{\PT_q} = \int_{[\curly P]^{\vir}} 1$. Then \eqref{eqn:3DKtheorDTPT} was proved by Bridgeland \cite{Bri}, and Toda \cite{Tod} in the Euler characteristic case. For $Y$ a toric 3-fold, the virtual dimension of $\curly P$ may be positive, and the correspondence was recently proved by Kuhn-Liu-Thimm \cite{KLT}. 

\begin{remark} \label{rem:zerocvclass}
It is a useful convention to declare $[\O_Y \to 0]$ to be the unique $\PT$ pair on $Y$ with $\beta = n = 0$. We then set $\curly P_v^{(0)}(Y) := \Spec \C$ and $\langle \! \langle 1 \rangle \! \rangle_{Y,0,0}^{\PT}:= 1 = 1 / \widehat{\Lambda}\udot(0)$.
This convention makes sense in view of \eqref{eqn:3DKtheorDTPT}. 
\end{remark}

\begin{corollary} \label{cor:virtualLefschetz}
Let $S \subset Y \subset X$, where $X$ is a projective Calabi-Yau 4-fold, $Y$ is a smooth connected 3-fold, and $S$ is a smooth connected surface. Suppose $h^0(N_{S/Y}) = h^1(N_{S/Y}) = h^2(\O_S) = 0$, and $S$ is the only pure 2-dimensional subscheme in class $[S] \in H^2(Y,\Q)$. Suppose furthermore that the inclusion map $\iota : Y \hookrightarrow X$ induces an injection $H_*(Y,\Q) \hookrightarrow H_*(X,\Q)$. Let $v' \in H^*(Y,\Q)$ and $v = (0,0,\gamma,\beta,n-\td_2(X) \cdot \gamma) = \iota_* v' \cdot \td(L)^{-1}$ where $L:=\O_X(Y)$. Let $v'' := 1-(1-v') e^{[S]} \in H^*(Y,\Q)$, and set $\beta'' := v_2'' \in H^4(Y,\Q)$. Suppose for any 2-dimensional closed subscheme $Z \subset X$ satisfying $\ch_2(\O_Z) = \gamma$ and $\ch_3(\O_Z) = \beta$, we have 
\[
H^1(X, (T_1(\O_Z) / T_0(\O_Z)) \otimes L) = 0,
\]
where $T_0(\O_Z) \subset T_1(\O_Z) \subset \O_Z$ is the torsion filtration of $\O_Z$. Finally suppose that, for all $n \in \Z$, the inclusion maps $\curly P^{(q)}_{v'}(Y) \hookrightarrow \curly P^{(q)}_{v}(X)$ for $q \in \{-1,0\}$ induce injective maps on the sets of connected components. 

Then, for $y=1$, we have
\begin{align*} 
\sum_n \langle \! \langle L \rangle \! \rangle_{X,\gamma,\beta,n}^{\DT} q^n &= q\udot \sum_n \langle \! \langle 1 \rangle \! \rangle_{Y,\beta'',n}^{\DT} q^n, \\
\sum_n \langle \! \langle L \rangle \! \rangle_{X,\gamma,\beta,n}^{\PT_0} q^n &= q\udot \sum_n \langle \! \langle 1 \rangle \! \rangle_{Y,\beta'',n}^{\PT} q^n,
\end{align*}
where $q\udot$ is and overall factor with power depending only on $S,Y,X,\beta''$. Moreover, the $\DT$--$\PT_0$ correspondence for $X$, $L$, $\gamma$, $\beta$, $y=1$ is equivalent to the $\DT$--$\PT$ correspondence for $Y$, $\beta''$. Thus it holds in the following three cases:
\begin{itemize}
\item when $\beta''=0$ (cf.~Remark \ref{rem:zerocvclass});
\item when $Y$ is a Calabi-Yau 3-fold by \cite{Bri, Tod};
\item when $Y$ is a toric 3-fold by \cite{KLT}.
\end{itemize}
\end{corollary}

\begin{proof}
Let $q \in \{-1,0\}$. For any $(F,s) \in \PTqvX$ with cokernel $Q$ and scheme theoretic support $Z$, we have short exact sequences
\begin{align*}
&0 \to \O_Z \to F \to Q \to 0 \\
&0 \to T_1(\O_Z) \to \O_Z \to \O_S \to 0,
\end{align*}
from which (1) and (2) of Proposition \ref{prop:A123} are easily deduced. For the chosen cohomology classes, we have morphisms
\[
\curly P_{v''}^{(q)}(Y) \to \curly P_{v'}^{(q)}(Y), \quad J\udot \mapsto J\udot(-S) = I_Y\udot.
\]
Since $S$ is the only pure 2-dimensional subscheme in $\gamma$, these are isomorphisms by Proposition \ref{prop:PT0PTcv}. Moreover, the universal sheaves are related by $\mathbb{J}\udot(-S) \cong \mathbb{I}_Y\udot$. Hence the isomorphisms preserve the virtual cycles, from which the two equalities in the proposition follow.

The assumption on connected components implies that we can choose orientations so that all signs $(-1)^{\sigma(i)}$ in \eqref{eqn:conncomppus} are $+1$. Thus the statement on the $\DT$--$\PT_0$ correspondence follows from the fact that for $y=1$, we have
$$
\sum_n \langle \! \langle L \rangle \! \rangle_{X,0,0,n}^{\DT} q^n = \sum_n \langle \! \langle 1 \rangle \! \rangle_{Y,0,n}^{\DT} q^n, 
$$
which was proved by third-named author \cite{Par}.
\end{proof}

In the following two examples, we use the notation of Corollary \ref{cor:virtualLefschetz} and we assume that the condition on connected components in the corollary is satisfied.

\begin{example}
Let $Y \to \PP^2$ be a general Weierstra{\ss} fibration with unique section $\PP^2 \subset Y$. Let $X = Y \times E$, where $E$ is a smooth elliptic curve. Let $\gamma = [\PP^2]$ and $\beta'' = \epsilon [\PP^1]$, where $\PP^1 \subset \PP^2$ and $\epsilon = 0,1$. Then the assumptions of Corollary \ref{cor:virtualLefschetz} are satisfied.
\end{example}

\begin{example}
Let $Y$ be the blow-up of $\PP^3$ in a point and let $\PP^2 \subset Y$ be the class of the exceptional divisor. Let $X \to Y$ be a general Weierstra{\ss} fibration with section, so we have $\PP^2 \subset Y \subset X$. Let $\gamma = [\PP^2]$ and $\beta'' = \epsilon [\PP^1]$,  where $\PP^1 \subset \PP^2$ and $\epsilon = 0,1$. Then the assumptions of Corollary \ref{cor:virtualLefschetz} are satisfied.
\end{example}

\subsection{Via virtual localization} \label{sec:Leflocal}

Let $Y$ be a smooth quasi-projective toric 3-fold. Denote the dense open torus of $Y$ by $T_Y \cong (\C^*)^3$. For $q \in \{-1,0\}$, $v = (0,0,\beta,n - \beta \cdot \td_1(Y)) \in H^*_c(Y,\Q)$, let $\curly P := \curly P_{v}^{(q)}(Y)$. Similar to the previous section, $\curly P$ has a perfect obstruction theory with accompanying \emph{equivariant} virtual class and twisted virtual structure sheaf 
\begin{align*}
&[\curly P]^{\vir} \in A_{\vd}^{T_Y}(\curly P), \quad \vd = \beta \cdot c_1(Y) \\
&\Ohat_{\curly P} := \O^{\vir}_{\curly P} \otimes (\det(T_{\curly P}^{\vir}))^{-\frac{1}{2}} \in K_0^{T_Y}(\curly P,\Z[\tfrac{1}{2}]).
\end{align*}
The action of $T_Y$ on $Y$ lifts to $\curly P$. Then the fixed locus $\curly P^{T_Y}$ is proper.\footnote{For $q = -1$ the fixed locus is 0-dimensional reduced \cite{MNOP}, and for $q=0$ it is a disjoint union of products of $\PP^1$ \cite{PT2}.} The $K$-theoretic $\PT_q$ invariants are defined by
\[
\langle \! \langle 1 \rangle \! \rangle_{Y,v}^{\PT_q} := \chi\big(\curly P, \Ohat\big) \in \Q(t^{\frac{1}{2}}_1, t^{\frac{1}{2}}_2, t^{\frac{1}{2}}_3),
\]
where $t_1,t_2,t_3$ are the equivariant parameters of $T_Y$. We also denote these invariants by $\langle \! \langle 1 \rangle \! \rangle_{Y,\beta,n}^{\PT_q}$. Abstractly, these invariants are defined by Thomason localization as in Section \ref{sec:virtualloc}, and concretely they are given by the Graber-Pandharipande localization formula \cite{GP} (or its $K$-theoretic version due to Qu \cite{Qu}). Then the equivariant $K$-theoretic $\DT$--$\PT$ correspondence is given by the same formula \eqref{eqn:3DKtheorDTPT} and was recently established in \cite{KLT}. 

We now consider the toric Calabi-Yau 4-fold 
$$
p : X = \mathrm{Tot}(K_Y) \to Y,
$$
i.e., the total space of the canonical line bundle on $Y$. Let $T \leq T_X := T_Y \times \C^*$ be the algebraic torus preserving the Calabi-Yau volume form, then $T \cong (\C^*)^3 \cong T_Y$. Let $v = (0,0,\gamma,\beta,n-\gamma \cdot \td_2(X)) \in H_c^*(X,\Q)$ and consider $\curly P := \curly P_v^{(q)}(X)$. We suppose the following: \\

\noindent \textbf{Assumption.} $\curly P^T$ is 0-dimensional and reduced. \\

As we will see in Section \ref{sec:vertex}, there are numerous examples where this assumption is satisfied, e.g., it is always satisfied for $q=-1$. We consider the $T$-equivariant line bundle $L := \O_X(Y)$. In Definition \ref{def:equivKtheorinv}, we defined the $T$-equivariant $K$-theoretic invariants
$$
\langle\!\langle \O_X(Y) \rangle\!\rangle_{X,v}^{\PT_q}, \quad q \in \{-1,0,1\}.
$$

Let $\cE$ be a $T$-equivariant coherent sheaf with proper support on $X$ and define $E  := p_* \cE$. Then there exists a distinguished triangle \cite[Prop.~2.14]{TT}
$$
R\Hom_X(\cE,\cE) \to R\Hom_Y(E,E) \to R\Hom_Y(E,E \otimes K_Y).
$$
Therefore, in $K$-theory, we have a splitting
\begin{equation} \label{eqn:Kthyeqn}
R\Hom_X(\cE,\cE) = R\Hom_Y(E,E) + R\Hom_Y(E,E)^\vee \in K_0^{T}(\pt).
\end{equation}
Serre duality in this setting is justified by Section \ref{sec:moduli}. 

\begin{lemma} 
Let $I\udot = [\O_X \to F] \in \curly P^T$. Then
$$
(\sqrt{\mathfrak{e}}(T_{\curly P}^{\vir}|_{I\udot}))^2 = (\widehat{\Lambda}\udot ( R\Gamma(Y, p_* F) - R\Hom_Y(p_* F, p_* F))^\vee)^2.
$$
\end{lemma}
\begin{proof}
Applying \eqref{eqn:Kthyeqn} to $\mathcal{E} = F$, we obtain the following identity in $K_0^T(\pt)$
\begin{align}
\begin{split} \label{eqn:Tvirspectral}
T_{\curly P}^{\vir}|_{I\udot} &= -R\Hom_X(I^\mdot, I^\mdot)_0 \\
&= R\Gamma(X,F) + R\Gamma(X,F)^\vee - R\Hom_Y(p_* F,p_* F) - R\Hom_Y(p_* F,p_* F)^\vee.
\end{split}
\end{align}
The result follows by combining with Proposition \ref{prop:OTisored} and using the fact that for any complex we have $\widehat{\Lambda}\udot V = (-1)^{\rk V} \widehat{\Lambda}\udot V^\vee$.
\end{proof}

We therefore have
\begin{equation} \label{eqn:invKY}
\langle\!\langle L \rangle\!\rangle_{X,v}^{\PT_q} = \sum_{\substack{I\udot = [\O_X \to F] \in \curly P^T \\ \vd(\{I\udot\}) = 0}} (-1)^{\sigma_{I\udot}} \frac{\widehat{\Lambda}\udot (R\Gamma(X,F(Y)) \otimes y^{-1})}{\widehat{\Lambda}\udot (R\Gamma(Y,p_*F) - R\Hom_Y(p_* F,p_*F))^\vee},
\end{equation}
which depends on a choice of sign attached to each $T$-fixed point. Henceforth, we make any choice satisfying the following: \\

\noindent \textbf{Choice.} For any $I\udot \in \curly P^T$ corresponding to a $\PT_q$ pair which is scheme theoretically supported on the zero section $Y \subset X$, we take $(-1)^{\sigma_{I\udot}} = 1$. \\

\begin{remark}
We expect that for any choice of (global) orientation on $\curly P$, we have $(-1)^{\sigma_{I\udot}} = (-1)^{\sigma_{J\udot}}$ for any two $I\udot, J\udot \in \curly P^T$ corresponding to $\PT_q$ pairs supported on the zero section $Y \subset X$. This can be seen as a toric analog of the expectation discussed below \eqref{eqn:conncomppus}. In other words, we expect that for signs coming from (global) orientations on $\curly P$, the choice made above is automatic (up to overall sign). 
\end{remark}

\begin{theorem} \label{thm:toricdimred}
Let $Y$ be a smooth quasi-projective toric 3-fold and consider $p : X = \mathrm{Tot}(K_Y) \to Y$. Denote by $\iota : Y \hookrightarrow X$ the inclusion of the zero section and by $T \leq T_X = T_Y \times \C^*$ the subtorus preserving the Calabi-Yau volume form. Fix $v \in H^*_c(X,\Q)$, and let $v' \in H^*_c(Y,\Q)$ be the unique class determined by the equation $\iota_* v' = v \cdot \td(\O_X(Y))$. Let $q \in \{-1,0,1\}$ and assume the following conditions hold:
\begin{enumerate}
\item [$\mathrm{B1)}$] $\curly P_{v}^{(q)}(X)^T$ is 0-dimensional and reduced.
\item [$\mathrm{B2)}$] For any $(F,s) \in \curly P_{v}^{(q)}(X)^T$ \emph{not} scheme-theoretically supported on $Y \subset X$, the class
$R\Gamma(X,F(Y)) \in K^T_0(\mathrm{pt})$ has a \emph{positive} $T$-fixed term. 
\item [$\mathrm{B3)}$] There exists an effective $T_Y$-invariant divisor\footnote{Recall that any effective divisor $S = \sum_i n_i D_i$ on $Y$, where 
$D_i$ are prime divisors, gives rise to a scheme $S \subset Y$ defined by the ideal $\prod_i I_{D_i / Y}^{n_i}$.} $\jmath : S \hookrightarrow Y$ with proper support, such that for all $(F,s)  \in \curly P_{v'}^{(q)}(Y)^{T_Y}$ we have $S = \Supp(F)^{\pure}$. 
\end{enumerate}
Then, for $q \in \{-1,0\}$, $\lim_{y \to 1} \langle\!\langle \O_X(Y) \rangle\!\rangle_{X,v}^{\PT_q}$ is well-defined and 
\begin{align*}
\lim_{y \to 1} \langle\!\langle \O_X(Y) \rangle\!\rangle_{X,v}^{\DT} &=   \langle\!\langle 1 \rangle\!\rangle_{Y,v''}^{\DT}, \\ 
\lim_{y \to 1} \langle\!\langle \O_X(Y) \rangle\!\rangle_{X,v}^{\PT_0} &=   \langle\!\langle 1 \rangle\!\rangle_{Y,v''}^{\PT},
\end{align*}
where $v'' = 1-(1-v')\cdot \exp(v'_1)$ and $v'_1$ denotes the component of $v'$ in $H^2_c(Y,\Q)$.

For $q=1$, $\lim_{y \to 1} \langle\!\langle \O_X(Y) \rangle\!\rangle_{X,v}^{\PT_1}$ is well-defined and for $S$ sufficiently ample\footnote{The required ``sufficient ampleness'' is defined in the proof.}, we have
\begin{align*}
\lim_{y \to 1} \langle\!\langle \O_X(Y) \rangle\!\rangle_{X,v}^{\PT_1} =  \langle\!\langle 1 \rangle\!\rangle_{Y,v'''}^{\PT},
\end{align*}
where $v''' = (v')^{\vee} \cdot (\ch(\O_Y(S)) / \td(\O_Y(S)))^{-1} + \jmath_* \ch(\O_S)$. 
\end{theorem}
\begin{proof}
We are concerned with the calculation of \eqref{eqn:invKY}. Note that we may drop the requirement $\vd(\{ I\udot \}) = 0$. Indeed, if $\vd(\{I\udot\}) < 0$, the complexes $T_{\curly P}^{\vir}|_{I\udot}$ and
$$
R\Gamma(Y, p_* F) - R\Hom_Y(p_* F, p_* F)
$$
have a negative $T$-fixed term (by \eqref{eqn:Tvirspectral}), thus contributing zero to \eqref{eqn:invKY}.

By assumption B2), we can also restrict attention to $T$-fixed pairs which are scheme theoretically supported on $Y$, i.e., they are in the image of
\[
\curly P_{v'}^{(q)}(Y)^{T_Y} \hookrightarrow \curly P_{v}^{(q)}(X)^T.
\]
Such a pair is of the form $[\O_X \to \iota_* F]$ for some $\PT_q$ pair $I\udot = [\O_Y \to F]$, where the map $\O_X \to \iota_* F$ is the composition $\O_X \to \iota_* \O_Y \to \iota_* F$. Since $p_* ((\iota_* F)(Y)) \cong F \otimes K_Y$, we are interested in 
$$
R\Gamma(Y, F) - R\Gamma(Y,F \otimes K_Y)^\vee y - R\Hom_Y(F, F),
$$
which, for $y=1$, equals $-R\Hom_Y(I^\mdot, I^\mdot)_0$ by Serre duality on $Y$. In each of the cases $q \in \{-1,0,1\}$, we show below that $-R\Hom_Y(I^\mdot, I^\mdot)_0$ has no positive $T_Y$ fixed part. Then the limit in the theorem is well-defined and, by Choice, we have
\begin{equation} \label{eqn:invKY2}
\lim_{y \to 1} \langle\!\langle \O_X(Y) \rangle\!\rangle_{X,v}^{\PT_q} = \sum_{I\udot \in \curly P_{v'}(Y)^{T_Y}} \frac{1}{\widehat{\Lambda}\udot (R\Hom_Y(I\udot,I\udot)_0[1])^\vee}.
\end{equation}
Our remaining task is to rewrite the complex $R\Hom_Y(I\udot,I\udot)_0[1]$ in each of the cases $q \in \{-1,0,1\}$ and show in the process that it has no positive $T_Y$-fixed part. Take $I\udot \in \curly P_{v'}(Y)^{T_Y}$ and denote the scheme-theoretic support of $I^\mdot$ by $Z \subset Y$ and let $Z^{\mathrm{pure}}$ be the underlying effective divisor. By B3), we have $Z^\pure = S$. \\

\noindent \textbf{$\DT$ case.} For $q = -1$, we have 
$$
I^\mdot \cong I_{Z/Y} \cong I_{C/Y}(-S),
$$
where $C \subset Y$ has dimension $\leq 1$. We deduce the following equality in $K_0^T(\pt) \cong K_0^{T_Y}(\pt)$
\begin{equation} \label{eqn:DTcv}
-R\Hom_Y(I^\mdot, I^\mdot)_0 = -R\Hom_Y(I_{C/Y},I_{C/Y})_0.
\end{equation}
For $v''$ as defined in the theorem, consider the following bijection (at the level of $\C$-valued points)
$$
\curP_{v'}^{(-1)}(Y)^{T_Y} \to \curP_{v''}^{(-1)}(Y)^{T_Y}, \quad I_{Z/Y} \mapsto I_{Z/Y}(S).
$$
Moreover, by \cite[Lem.~6 \& 8]{MNOP}, $\curP_{v''}^{(-1)}(Y)^{T_Y}$ is 0-dimensional and reduced, and the complex $R \Hom_Y(I_{C/Y}, I_{C/Y})_0[1]$ has no $T_Y$-fixed part. This establishes the $q=-1$ case by \eqref{eqn:invKY2} and \eqref{eqn:DTcv}. \\

\noindent \textbf{$\PT_0$ case.} For $q=0$, by Proposition \ref{prop:PT0PTcv}, we have
$$
I^\mdot \cong J^\mdot(-S), \quad  J^\mdot:=[\O_Y \to T_1(F)(S)],
$$
where $J^\mdot$ is a $\PT$ pair on $Y$. 
We obtain the following equality in $K_0^{T_Y}(\pt)$
\begin{equation} \label{eqn:PTcv}
-R\Hom_Y(I^\mdot, I^\mdot)_0 = -R\Hom_Y(J^\mdot, J^\mdot)_0.
\end{equation}
For $v''$ as defined in the theorem, consider the map (at the level of $\C$-valued points)
$$
\curP_{v'}^{(0)}(Y)^{T_Y} \to \curP_{v''}^{(0)}(Y)^{T_Y}, \quad I^\mdot \mapsto I^\mdot(S).
$$
This map is a bijection by Proposition \ref{prop:PT0PTcv} and B3).\footnote{Note that, by Remark \ref{rem:zerocvclass}, we allow the case $T_1(F) = 0$.}
By B1), $\curP_{v''}^{(0)}(Y)^{T_Y}$ is 0-dimensional. Moreover, by \cite[Thm.~1]{PT2} it is reduced and therefore $R\Hom_Y(J^\mdot, J^\mdot)_0[1]$ has no positive $T_Y$-fixed part. The $q=0$ case then follows by \eqref{eqn:invKY2} and \eqref{eqn:PTcv}. \\

\noindent \textbf{$\PT_1$ case.} For $q=1$, by Proposition \ref{lem:PT1PTcv}, we have
$$
F\dual \otimes K_Y [1] \cong \jmath_* (J\udot \otimes K_S),
$$
where $J^\mdot = [\O_S \to G]$ is a $\PT$ pair on $S$ and $\jmath_* (G \otimes K_S) = Q^D$. We first write
$$
-R\Hom_Y(I^\mdot, I^\mdot)_0 = R \Hom_Y(F^\vee,\O_Y) + R \Gamma(Y,F^\vee) - R\Hom_Y(F^\vee,F^\vee). 
$$
Using the substitutions
\begin{align*}
F^\vee &= -\jmath_*(K_S) \otimes K_Y^{-1} + \jmath_*(G \otimes K_S) \otimes K_Y^{-1}, \\ 
K_S &= K_Y(S)|_S, \quad \O_S = \O_Y - \O_Y(-S),
\end{align*}
we find 
\begin{align}
\begin{split} \label{eqn:PTcv2}
-R\Hom_Y(I^\mdot, I^\mdot)_0 &= -R\Hom_Y(J^\mdot, J^\mdot)_0 - R\Gamma(Y,\O_S(S)) - R\Hom_Y(\O_S,\O_Y) \\
&= -R\Hom_Y(J^\mdot, J^\mdot)_0.
\end{split}
\end{align}
For the second equality, we used Serre duality on $Y$ and $S$ ($S$ is Gorenstein with proper support):
\[
-R\Hom_Y(\O_S,\O_Y) = R\Gamma(S,K_S(-S))^\vee = R\Gamma(S,\O_S(S)).
\]
For $v'''$ as defined in the theorem, consider the map (at the level of $\C$-valued points)
$$
\curP_{v'}^{(1)}(Y)^{T_Y} \to \curP_{v'''}^{(0)}(Y)^{T_Y}, \quad I^\mdot \mapsto [\O_Y \to \jmath_* G],
$$
where $G$ was defined above.\footnote{Note that, by Remark \ref{rem:zerocvclass}, we allow the case $Q = 0$, i.e., $G = 0$.} By Proposition \ref{lem:PT1PTcv}, this is a bijection onto the locus of $\PT$ pairs scheme theoretically supported on $S$. By definition, we say that $S$ is sufficiently ample when this map is surjective. Then the proof follows from \eqref{eqn:invKY2} and \eqref{eqn:PTcv2}
\end{proof}

Let $L_1, L_2$ be line bundles on a smooth projective toric surface $S$, such that $L_1 \otimes L_2 \cong K_S$, and $Y = \mathrm{Tot}(L_1)$. Then $X = \mathrm{Tot}(L_1 \oplus L_2)$ and, taking $\gamma = d[S]$ for any $d>0$, condition B3) is automatic. We establish in Theorem \ref{prop:fixlocmain1} that condition B1) is satisfied for $q \in \{-1,0\}$ and it conjecturally holds for $q=1$ too (Conjecture \ref{conj:0dimTXfix}). Hence, for local surfaces, the main condition is B2). 

Suppose $S$ is del Pezzo and $Z \subset X$ is a $T$-fixed 2-dimensional proper closed subscheme with class $\gamma$ and \emph{no embedded curves}, i.e., its torsion filtration $T_0(\O_Z) \subset T_1(\O_Z) \subset \O_Z$ satisfies $T_0(\O_Z) = T_1(\O_Z)$. Using push-forward along the projection $p : X \to S$, it is not hard to show that if $Z$ is not scheme-theoretically supported on the zero section $Y \subset X$, then $R\Gamma(X,L_2|_Z)$ has positive $T$-fixed part so B2) is satisfied.

Suppose $S = \PP^2$ and $L_1 = \O_{\PP^2}(-a)$, $L_2 = \O_{\PP^2}(-3+a)$ for any $a \geq 1$. We claim that for $q=1$, condition B2) is always satisfied. Indeed, for a $T$-fixed $\PT_1$ pair on $X$ with proper support and class $\gamma$, we have a short exact sequence $0 \to \O_{Z} \to F \to Q \to 0$, where $Z \subset X$ is pure 2-dimensional. Thus, by the previous paragraph, it suffices to show that $R\Gamma(X,Q \otimes L_2)$ does not have negative $T$-fixed part. This follows from the toric description of $T$-fixed points in Section \ref{sec:fixlocal}. Taking $a=3$ and $d \gg 0$, Theorem \ref{thm:toricdimred} implies that the $\PT_1$ invariants recover the refined $\PT$ invariants of $Y=\mathrm{Tot}(\O_{\PP^2}(-3))$ studied by Choi-Katz-Klemm \cite{CKK}. 

Condition B2) is not always satisfied. Continuing with $S = \PP^2$ and $X = \mathrm{Tot}(\O_{\PP^2}(-a) \oplus \O_{\PP^2}(-3+a))$, let $C := Z(x_0x_1x_2) \subset \PP^2$ and consider the $T$-fixed 2-dimensional subscheme $Z \subset X$ with reduced support $\PP^2$ and determined by the short exact sequence
\[
0 \to \O_C \otimes L_2^{-1} \to \O_Z \to \O_{\PP^2} \to 0,
\]
where $\O_Z \to \O_{\PP^2}$ denotes restriction to the reduced support. Since $R\Gamma(\PP^2,\O_C) = 0$, the $K$-theory class $R\Gamma(X,L_2|_{Z})$ does not have positive $T$-fixed part. However, for a $T$-fixed $\DT$ or $\PT_0$ pair $(F,s)$ on $X$ with proper support and class $\gamma = d[S]$ (for any $d>0$) and scheme theoretic support $Z$ satisfying $\ch_3(T_1(\O_Z)) = [\PP^1]$ or $2[\PP^1]$, the toric description of $T$-fixed points in Section \ref{sec:fixlocal} can be used to show that $R\Gamma(X,F \otimes L_2)$ has a positive $T$-fixed part. 

\begin{example}
Let $X = \mathrm{Tot}(\O_{\PP^2}(-a) \oplus \O_{\PP^2}(-3+a))$ and $Y = \mathrm{Tot}(\O_{\PP^2}(-a))$ for any $a \leq 1$. Fix $\gamma = d [\PP^2]$ for any $d > 0$ and $\beta = (\epsilon + \tfrac{3}{2} + \tfrac{1}{2}d(d-1)a) [\PP^1]$ with $\epsilon\in\{0,1,2\}$. Then for any $v=(0,0,\gamma,\beta,*) \in H^*_c(X,\Q)$, we are in the setting of the previous paragraph and the conditions of Theorem \ref{thm:toricdimred} are satisfied for $q \in \{-1,0\}$. Hence, for $L = \O_X(Y)$ and $y=1$, the $\DT$--$\PT_0$ correspondence reduces to the $K$-theoretic $\DT$--$\PT$ correspondence for $Y$ which was proved in \cite{KLT}.\footnote{In this example, we use that for $y=1$, $\langle\!\langle \O_X(Y) \rangle\!\rangle_{X,0,0,n}^{\DT} = \langle\!\langle 1 \rangle\!\rangle_{Y,0,n}^{\DT}$ which is known to hold by the 0-dimensional analog of Proposition \ref{prop:dimredvertex} (cf.~\cite{KR}). See also Remark \ref{rem:otherL}.} In particular, Conjecture \ref{conj:KDTPT0inbodytext_equiv} holds in this setting.
\end{example}

\section{Fixed loci} \label{sec:fixloc}

In this section, for $X$ a toric Calabi-Yau 4-fold, we discuss the fixed loci of the moduli spaces $\PTqvX$ and the representations of tangent spaces at the fixed points.

\subsection{Toric geometry} \label{sec:toricgeomnot}

For  any smooth quasi-projective toric 4-fold $X$, we denote its dense algebraic torus by $T_X \cong (\C^*)^4$. 
We denote by $V(X)$ its collection of $T_X$-fixed points, by $E(X)$ its collection of \emph{proper} $T_X$-fixed irreducible 1-dimensional subvarieties (which are ``$\PP^1$'s''), and by $F(X)$ its collection of \emph{proper} $T_X$-fixed irreducible 2-dimensional subvarieties (which are irreducible smooth projective toric surfaces). We refer to elements of $V(X)$, $E(X)$, $F(X)$ as vertices, edges, and faces respectively.

The elements $\alpha \in V(X)$ correspond to $T_X$-fixed points $p_\alpha \in X$. Since these fixed points correspond to 4-dimensional cones in the fan, they are also in bijection with $T_X$-invariant affine open subsets $U_\alpha \subset X$. For each $\alpha \in V(X)$, there exist coordinate functions $x_1,x_2,x_3,x_4$ on $U_\alpha$ such that the action of $T_X$ is given by
$$
(t_1,t_2,t_3,t_4) \cdot (x_1,x_2,x_3,x_4) = (t_1x_1,t_2x_2,t_3x_3,t_4x_4).
$$
If $X$ is Calabi-Yau, then the subtorus 
\begin{equation} \label{eqn:CYtorus}
T = Z(t_1t_2t_3t_4 - 1) \leq T_X
\end{equation} 
preserves the Calabi-Yau volume form, and we refer to $T$ as the \emph{Calabi-Yau torus}. 
Using labels $\alpha$ for the vertices, we can label edges by $\alpha\beta$ according to the vertices which they connect. Any edge $\alpha\beta \in E(X)$ corresponds to a $T_X$-fixed line $\PP^1 \cong L_{\alpha\beta} \subset X$ with normal bundle
\begin{align} \label{weightsN}
N_{L_{\alpha\beta}/X} \cong \O(m_{\alpha\beta}) \oplus \O(m_{\alpha\beta}') \oplus \O(m_{\alpha\beta}'').
\end{align}
When $X$ is Calabi-Yau, we have
$$
m_{\alpha\beta} + m_{\alpha\beta}' + m_{\alpha\beta}'' = -2.
$$
For coordinates $x_1,x_2,x_3,x_4$ on $U_\alpha$ and $x_1',x_2',x_3',x_4'$ on $U_\beta$ such that $L_{\alpha\beta} \cap U_\alpha = Z(x_2,x_3,x_4)$ and $L_{\alpha\beta} \cap U_\beta = Z(x_2',x_3',x_4')$, the coordinate transformation from $U_\alpha$ to $U_{\beta}$ is given by
\begin{equation} \label{eqn:torictransf}
(x_1,x_2,x_3,x_4) \mapsto (x_1^{-1},x_1^{-m_{\alpha\beta}}x_2,x_1^{-m_{\alpha\beta}'}x_3,x_1^{-m_{\alpha\beta}''} x_4). 
\end{equation}
Faces $f \in F(X)$ correspond to irreducible smooth projective toric surfaces $S_f$. The normal bundle $N_{S_f/X}$ is a (split!) rank 2 equivariant bundle on $S_f$.
%Split: follow the gluing description along a maximal $T_S$-invariant chain of $\PP^1$'s in $S_f$. This determines $N_{S_f/X}$ outside a single point, and and on this open it splits $T_X$-equivariantly. Then use reflexivity.
 Note that any smooth projective toric surface can arise in this way. 

In this section, we use partitions in various dimensions. 
\begin{definition} \label{def:partitions}
Let $\mu \subset \Z_{\geq 0}^d$. We say that $\mu$ is a $d$-dimensional partition when
$$
(i_1, \ldots, i_d) \in \mu \Rightarrow (j_1, \ldots, j_d) \in \mu \quad \forall (j_1, \ldots, j_d) \in \Z_{\geq 0}^d \, \textrm{with} \, j_k \leq i_k \, \textrm{for all} \, k.
$$
For $d=2,3,4$, we refer to $\mu$ as a \emph{partition}, \emph{plane partition}, \emph{solid partition} respectively. A partition is called finite when $\mu$ is a finite set, in which case its \emph{size} $|\mu|$ is the cardinality of this set.
\end{definition}

A partition $\mu$ can also be written as a sequence
$$
\mu = \{\mu(i_1, \ldots, i_{d-1}) \}_{i_1, \ldots, i_{d-1} \geq 0}
$$
determined by
$$
(i_1, \ldots, i_d) \in \mu \Leftrightarrow i_d \leq \mu(i_1, \ldots, i_{d-1}).
$$
In this way, solid partitions can be visualized as plane partitions in $\mathbb{R}^3$, where each box located at $(i,j,k) \in \mathbb{Z}_{\geq 0}^3$ carries a number $\mu(i,j,k) \geq 0$ indicating how many boxes are stacked ``on top of it'' in the positive $x_4$-direction. Note that the number $\mu(i,j,k)$ can be infinite.

\begin{figure}
\input{fig1.tex}
\caption{The solid partition of Example \ref{ex:solidpart}.} \label{fig1}
\end{figure}
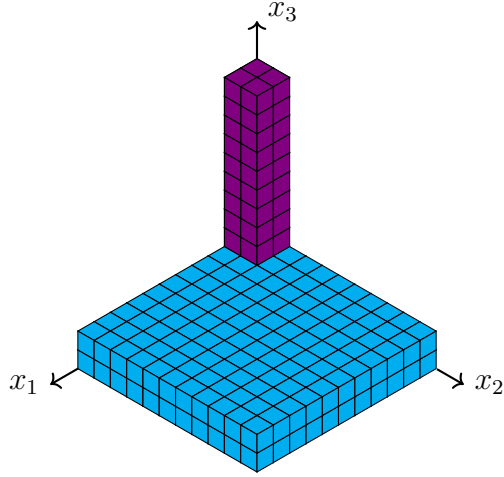

\begin{example} \label{ex:solidpart}
In Figure \ref{fig1}, we give a graphical representation of the solid partition corresponding to the scheme theoretic union of the schemes cut out by $(x_3^2,x_4)$ and $(x_1^2,x_2^2)$.
The blue boxes extend infinitely far into the first quadrant of the $x_1x_2$-plane. The blue boxes have height $0$ in the positive $x_4$-direction. The purple boxes extend infinitely far along the positive $x_3$-axis. The purple boxes have height $\infty$ in the positive $x_4$-direction. The eight ``hidden'' boxes at $(i,j,k)$ with $0 \leq i,j,k \leq 1$ are also coloured purple.  
\end{example}

\subsection{Fixed loci --- Local} \label{sec:fixlocal}

In this section, we study $(\C^*)^4$-fixed $\PT_q$ pairs on $\C^4$ with respect to the standard action, where we do not require the pair to have proper support. These form the building blocks for all further investigations. \\

\noindent \textbf{$\DT$ case.} As usual, $(\C^*)^4$-invariant subschemes of $\C^4$ correspond to monomial ideals of $\C[x_1,x_2,x_3,x_4]$, which in turn correspond to solid partitions. We denote the subscheme corresponding to a solid partition $\pi$ by $Z_\pi$.\footnote{We use the same correspondence in dimensions 2 and 3 where solid partitions are replaced by partitions and plane partitions respectively.} More precisely, for all $(i_1, i_2,i_3,i_4) \in \Z_{\geq 0}^{4}$
$$
x_1^{i_1}x_2^{i_2}x_3^{i_3}x_4^{i_4} \notin I_{Z_\pi} \Leftrightarrow (i_1,i_2,i_3,i_4) \in \pi.
$$

Let $Z \subset \C^4$ be a $(\C^*)^4$-invariant closed subscheme of dimension $\leq 2$. Consider the torsion filtration
\begin{equation*} 
T_0(\O_Z) \subset T_1(\O_Z) \subset \O_{Z}.
\end{equation*}
We have corresponding schemes
$$
Z \supset Z_0 \supset Z_1,
$$
where $Z_1$ has no embedded components of dimension $\leq 1$, i.e.~$\O_{Z_1}$ is pure, and $Z_0$ has no embedded components of dimension 0. 

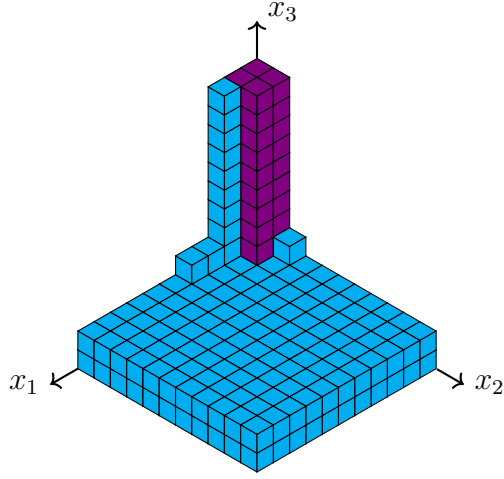
\begin{figure}[t]
\input{fig1a.tex}
\caption{The solid partition of Example \ref{ex:solidpart2}.} 
\label{fig1a}
\end{figure}

\begin{example} \label{ex:solidpart2}
This is a continuation of Example \ref{ex:solidpart}. Using the notation of that example, we place an embedded curve at weights $(2,0,i,0)$, $i \geq 2$ and we place three more embedded boxes (all with height 0 in direction $x_4$). We can view this as stacking boxes \emph{on top} of the solid partition of Example \ref{ex:solidpart} with gravity pulling into direction $(-1,-1,-1,-1)$ (Figure \ref{fig1a}). In this example, all inclusions of $\varnothing \subset Z_0 \subset Z_1 \subset Z$ are strict.
\end{example}

The subscheme $Z_1$ is described by six finite partitions
$$
\boldsymbol{\lambda} = \{\lambda_{ab} \}_{1 \leq a < b \leq 4}.
$$
More precisely, a monomial $x_1^{i_1}x_2^{i_2}x_3^{i_3}x_4^{i_4}$ lies in $\O_{Z_{1}}$ if and only if 
$$
i_a,i_b \geq 0, \quad \textrm{and } (i_c,i_d) \in \lambda_{ab}, \quad \exists a,b,c,d \in \{1,2,3,4\} \textrm{ \, distinct \, }  a<b, c<d.
$$

Next we consider the scheme $Z_0$. Since it has no embedded points, it is determined by the restrictions
$$
Z_0|_{\C^* \times \C^3}, \quad Z_0|_{\C \times \C^* \times \C^2}, \quad Z_0|_{\C^2 \times \C^* \times \C}, \quad Z_0|_{\C^3 \times \C^*}.
$$
There are corresponding $(\C^*)^3$-invariant subschemes $W_1$, $W_2$, $W_3$, $W_4 \subset \C^3$ such that
\begin{align*}
&\O_{Z_0|_{\C^* \times \C^3}} = \O_{W_1}[x_1,x_1^{-1}], \quad \O_{Z_0|_{\C \times \C^* \times \C^2}} = \O_{W_2}[x_2,x_2^{-1}], \\
&\O_{Z_0|_{\C^2 \times \C^* \times \C}} = \O_{W_3}[x_3,x_3^{-1}], \quad \O_{Z_0|_{\C^3 \times \C^*}} = \O_{W_4}[x_4,x_4^{-1}].
\end{align*}
As in \cite{MNOP}, the subschemes $W_1, W_2, W_3, W_4$ have dimension $\leq 1$ and correspond to plane partitions 
$$
\boldsymbol{\mu} = \{\mu_a\}_{a=1}^{4}.
$$ 
The plane partitions $\boldsymbol{\mu}$ are of course compatible with the finite partitions $\boldsymbol{\lambda}$. Specifically
\begin{align}
\begin{split} \label{eqn:planepartcompat}
\forall j \gg 0: \quad (j,k,l) \in \mu_1 \Leftrightarrow (k,l) \in \lambda_{12}, \\
\forall k \gg 0: \quad (j,k,l) \in \mu_1 \Leftrightarrow (j,l) \in \lambda_{13},\\
\forall l \gg 0: \quad (j,k,l) \in \mu_1 \Leftrightarrow (j,k) \in \lambda_{14},
\end{split}
\end{align}
and similarly for $\mu_2, \mu_3, \mu_4$. Note that $\mu_a$ has infinite size unless $\lambda_{ab} = \varnothing$ for all $b \neq a$ (Definition \ref{def:partitions}). 

We summarize our discussion as follows:
\begin{lemma} \label{lem:fixlocaffDT}
Let $Z \subset \C^4$ be a 2-dimensional $(\C^*)^4$-invariant closed subscheme. Let $Z \supset Z_0 \supset Z_1$ be the subschemes corresponding to the torsion filtration $T_0(\O_Z) \subset T_1(\O_Z) \subset \O_{Z}$. Then 
\begin{itemize}
\item $Z$ is uniquely determined by a solid partition $\pi$,
\item $Z_1$ is uniquely determined by finite partitions $\boldsymbol{\lambda} = \{\lambda_{ab} \}_{1 \leq a < b \leq 4}$,
\item $Z_0$ is uniquely determined by plane partitions $\boldsymbol{\mu} = \{\mu_a \}_{a=1}^{4}$ which are compatible with finite partitions $\boldsymbol{\lambda} = \{\lambda_{ab} \}_{1 \leq a < b \leq 4}$.
\end{itemize}
\end{lemma} 

This allows us to make the following definition.
\begin{definition} \label{def:partitionnotation}
Any collection of finite partitions $\boldsymbol{\lambda} = \{\lambda_{ab} \}_{1 \leq a < b \leq 4}$ (not all empty), determines a unique corresponding $(\C^*)^4$-invariant pure 2-dimensional subscheme which we denote by $Z_{\boldsymbol{\lambda}} \subset \C^4$. 

If $\boldsymbol{\mu} = \{\mu_a \}_{a=1}^{4}$ is a choice of plane partitions (not all finite), each corresponding to a subscheme of dimension $\leq 1$ in $\C^3$, then it uniquely determines a $(\C^*)^4$-invariant 2-dimensional subscheme without embedded 0-dimensional components, which we denote by $Z_{\boldsymbol{\mu}} \subset \C^4$.
\end{definition} 

It is often convenient to represent the partitions in Lemma \ref{lem:fixlocaffDT} by formal power series in $t_1, t_2, t_3, t_4$. Suppose $Z \subset \C^4$ is a 2-dimensional $(\C^*)^4$-invariant closed subschheme with corresponding solid partition $\pi$. Then we write
$$
\pi = \sum_{i_1,i_2,i_3 \geq 0} \sum_{i_4=0}^{\pi(i_1,i_2,i_3)} t_1^{i_1}t_2^{i_2}t_3^{i_3}t_4^{i_4}.
$$
Denoting the partitions corresponding to $Z_0, Z_1$ by $\boldsymbol{\mu}$, $\boldsymbol{\lambda}$ respectively, we write
\begin{align*}
\mu_1 = \sum_{i_2,i_3 \geq 0} \sum_{i_4=0}^{\mu_1(i_2,i_3)} t_2^{i_2}t_3^{i_3}t_4^{i_4}, \ldots \quad \lambda_{12} =  \sum_{i_3 \geq 0} \sum_{i_4=0}^{\lambda_{12}(i_3)} t_3^{i_3}t_4^{i_4}, \ldots
\end{align*}

\begin{example}
Let $I_Z = (x_2^2, x_2x_3, x_3^2, x_2x_4,x_3x_4)$. Then $T_0(\O_Z) = 0$ and
$$
T_1(\O_Z)  = x_2 \C[x_1] \oplus x_3 \C[x_1]. 
$$
In particular, $Z_1 \cong \C^2$. Therefore $\lambda_{14} = 1$, $\mu_1 = t_2+t_3+ \frac{1}{1-t_4}$, $\mu_4 = \frac{1}{1-t_1}$, and all other $\lambda_{ab}$, $\mu_a$ are empty. Here $1/(1-t_i)$ is the formal power series obtained by expanding in ascending powers of $t_i$.
\end{example}

\noindent \textbf{$\PT_0$ case.} Before we give a combinatorial description of $(\C^*)^4$-fixed $\PT_0$ pairs on $\C^4$, we study pure 2-dimensional $(\C^*)^4$-fixed subschemes $Z \subset \C^4$. By Lemma \ref{lem:fixlocaffDT}, $Z = Z_{\boldsymbol{\lambda}}$ for a collection of finite partitions $\boldsymbol{\lambda} = \{\lambda_{ab}\}_{1 \leq a < b \leq 4}$. Clearly, each partition $\lambda_{ab}$ also determines a pure 2-dimensional scheme $Z_{\lambda_{ab}}$, with irreducible support, and we have
$$
Z_{\boldsymbol{\lambda}} = \bigcup_{1 \leq a<b \leq 4} Z_{\lambda_{ab}}.
$$
The following proposition classifies when $Z_{\boldsymbol{\lambda}}$ is Cohen-Macaulay:
\begin{proposition} \label{lem:casesCM}
Consider the scheme-theoretic intersection
$$
W = (Z_{\lambda_{12}} \cup Z_{\lambda_{13}} \cup Z_{\lambda_{23}}) \cap (Z_{\lambda_{14}} \cup Z_{\lambda_{24}} \cup Z_{\lambda_{34}}).
$$
There are four cases:
\begin{enumerate}
\item[$\mathrm{(i)}$] $W$ is empty, then $Z_{\boldsymbol{\lambda}}$ is Cohen-Macaulay.
\item[$\mathrm{(ii)}$] $W$ is 0-dimensional, then  $Z_{\boldsymbol{\lambda}}$ is not Cohen-Macaulay.
\item[$\mathrm{(iii)}$] $W$ is 1-dimensional and impure, then  $Z_{\boldsymbol{\lambda}}$ is not Cohen-Macaulay.
\item[$\mathrm{(iv)}$] $W$ is 1-dimensional and pure, then  $Z_{\boldsymbol{\lambda}}$ is Cohen-Macaulay.
\end{enumerate}
\end{proposition}
\begin{proof}
We use the following characterization: a pure 2-dimensional subscheme $S$ in a smooth quasi-projective 4-fold $X$ is Cohen-Macaulay if and only if $\ext^3_X(\O_S,\O_X) = 0$. We also note that each $Z_{\lambda_{ab}}$ is Cohen-Macaulay (being the base change of an element of $\Hilb^n(\C^2)^{(\C^*)^2}$). Let $X = \C^4$.
%Irreducible 0-dimensional schemes of finite type over $\C$ are CM + properties CM morphisms Stacks Project.

For (i), there are two cases. Suppose $Z_{\lambda_{14}} \cup Z_{\lambda_{24}} \cup Z_{\lambda_{34}} = \varnothing$. We want to conclude that $Z_{\lambda_{12}} \cup Z_{\lambda_{13}} \cup Z_{\lambda_{23}}$ is Cohen-Macaulay. We first observe that $Z_{\lambda_{12}} \cup Z_{\lambda_{13}}$ is Cohen-Macaulay. This follows from dualizing the short exact sequence
$$
0 \to \O_{Z_{\lambda_{12}} \cup Z_{\lambda_{13}}} \to \O_{Z_{\lambda_{12}}} \oplus \O_{Z_{\lambda_{13}}} \to \O_{Z_{\lambda_{12}} \cap Z_{\lambda_{13}}} \to 0
$$
combined with the facts that $Z_{\lambda_{12}}$, $Z_{\lambda_{13}}$ are Cohen-Macaulay and $Z_{\lambda_{12}} \cap Z_{\lambda_{13}}$ is pure (so $\ext^3_X(\O_{Z_{\lambda_{12}}},\O_X) = \ext^3_X(\O_{Z_{\lambda_{13}}},\O_X)  = \ext^4_X(\O_{Z_{\lambda_{12}} \cap Z_{\lambda_{13}}},\O_X) = 0$). Similarly, dualizing
$$
0 \to \O_{Z_{\lambda_{12}} \cup Z_{\lambda_{13}} \cup Z_{\lambda_{23}}} \to \O_{Z_{\lambda_{12}} \cup Z_{\lambda_{13}}} \oplus \O_{Z_{\lambda_{23}}} \to \O_{(Z_{\lambda_{12}} \cup Z_{\lambda_{13}}) \cap Z_{\lambda_{23}}} \to 0
$$
and using purity of $(Z_{\lambda_{12}} \cup Z_{\lambda_{13}}) \cap Z_{\lambda_{23}}$, we find that $Z_{\lambda_{12}} \cup Z_{\lambda_{13}} \cup Z_{\lambda_{23}}$ is Cohen-Macaulay. The second case is $Z_{\lambda_{12}} \cup Z_{\lambda_{13}} \cup Z_{\lambda_{23}} = \varnothing$. Then one can similarly deduce that $Z_{\lambda_{14}} \cup Z_{\lambda_{24}} \cup Z_{\lambda_{34}}$ is Cohen-Macaulay.

For (ii), we have ---after suitable relabelling--- $Z_{\lambda_{12}}, Z_{\lambda_{34}} \neq \varnothing$ and $Z_{\lambda_{ab}} = \varnothing$ otherwise. Then we have a short exact sequence
$$
0 \to \O_{Z_{\lambda_{12}} \cup Z_{\lambda_{34}}} \to \O_{Z_{\lambda_{12}}} \oplus \O_{Z_{\lambda_{34}}} \to \O_{Z_{\lambda_{12}} \cap Z_{\lambda_{34}}} \to 0,
$$
which is a $\PT_0$ pair with non-zero cokernel. By Proposition \ref{Lem.PT0=PT1}, we deduce that $Z_{\lambda_{12}} \cup Z_{\lambda_{34}}$ is not Cohen-Macaulay.

For (iii) and (iv), we consider the short exact sequence
$$
0 \to \O_{Z_{\boldsymbol{\lambda}}} \to \O_{Z_{\lambda_{12}} \cup Z_{\lambda_{13}} \cup Z_{\lambda_{23}}} \oplus \O_{Z_{\lambda_{14}} \cup Z_{\lambda_{24}} \cup Z_{\lambda_{34}}} \to \O_W \to 0,
$$
where $Z_{\lambda_{12}} \cup Z_{\lambda_{13}} \cup Z_{\lambda_{23}}$ and $Z_{\lambda_{14}} \cup Z_{\lambda_{24}} \cup Z_{\lambda_{34}}$ are pure 2-dimensional. If there exists a 0-dimensional subsheaf $0 \neq Q \subset \O_W$, we obtain a $\PT_0$ pair with nonzero cokernel
$$
0 \to \O_{Z_{\boldsymbol{\lambda}}} \to F \to Q \to 0
$$
and $Z_{\boldsymbol{\lambda}}$ is not Cohen-Macaulay by Proposition \ref{Lem.PT0=PT1}. Otherwise
$$
\ext_X^3(\O_{Z_{\lambda_{12}} \cup Z_{\lambda_{13}} \cup Z_{\lambda_{23}}},\O_X) = \ext_X^3(\O_{Z_{\lambda_{14}} \cup Z_{\lambda_{24}} \cup Z_{\lambda_{34}}},\O_X) = \ext_X^4(\O_W,\O_X) = 0,
$$
since $W$ is pure 1-dimensional and $Z_{\lambda_{12}} \cup Z_{\lambda_{13}} \cup Z_{\lambda_{23}}$, $Z_{\lambda_{14}} \cup Z_{\lambda_{24}} \cup Z_{\lambda_{34}}$ are Cohen-Macaulay.
\end{proof}

The subscheme $W$ can be easily read off from the data $\boldsymbol{\lambda} = \{\lambda_{ab}\}_{1 \leq a < b \leq 4}$. Here are some examples:
\begin{example} \label{ex:CMclass}
\hfill
\begin{enumerate}
\item[(i)] If $Z_{\boldsymbol{\lambda}}$ is set theoretically, but not necessarily scheme theoretically, supported in a hyperplane $Z(x_i)$, then $Z_{\boldsymbol{\lambda}}$ is Cohen-Macaulay.
\item[(ii)] For $\lambda_{12} = \lambda_{34} = 1$ and $\lambda_{ab} = \varnothing$ otherwise, we have $W = \{\pt\}$ and $Z_{\boldsymbol{\lambda}}$ is not Cohen-Macaulay.
\item[(iii)] For $\lambda_{12} = \lambda_{13} = \lambda_{34} =  1$ and $\lambda_{ab} = \varnothing$ otherwise, we see that $W = Z(x_1,x_2,x_4) \cong \mathbb{A}^1$ is pure. Therefore $Z_{\boldsymbol{\lambda}}$ is Cohen-Macaulay. However, for $\lambda_{12} = 1+t_4$,  $\lambda_{13} = \lambda_{34} =  1$ and $\lambda_{ab} = \varnothing$ otherwise, we have $W = Z(x_1,x_2,x_3x_4,x_4^2)$ which has an embedded point. Therefore $Z_{\boldsymbol{\lambda}}$ is not Cohen-Macaulay
\end{enumerate}
\end{example}

\begin{figure}[t]
\input{fig1b.tex}
\caption{The scheme $W$ of Example \ref{ex:PT0=PT1ex1}.} \label{fig1b}
\end{figure}
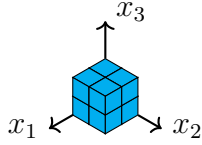

\begin{example} \label{ex:PT0=PT1ex1}
This is a continuation of Example \ref{ex:solidpart}. Similar to Example \ref{ex:CMclass}(ii), $W$ is 0-dimensional for the scheme of Example \ref{ex:solidpart}. It has length 4 and its weight spaces are depicted in Figure \ref{fig1b}, where we recall that blue boxes have height 0 in direction $x_4$. 
\end{example}

\begin{figure}[t]
\centering
\begin{subfigure}{.5\textwidth}
  \centering
  \input{fig3.tex}
  %\caption{A subfigure}
  %\label{fig:sub1}
\end{subfigure}%
\begin{subfigure}{.5\textwidth}
  \centering
  \input{fig3a.tex}
  %\caption{A subfigure}
  %\label{fig:sub2}
\end{subfigure}
\caption{The solid partition and scheme $W$ of Example \ref{ex:PT0=PT1ex2}.}
\label{fig3and3a}
\end{figure}
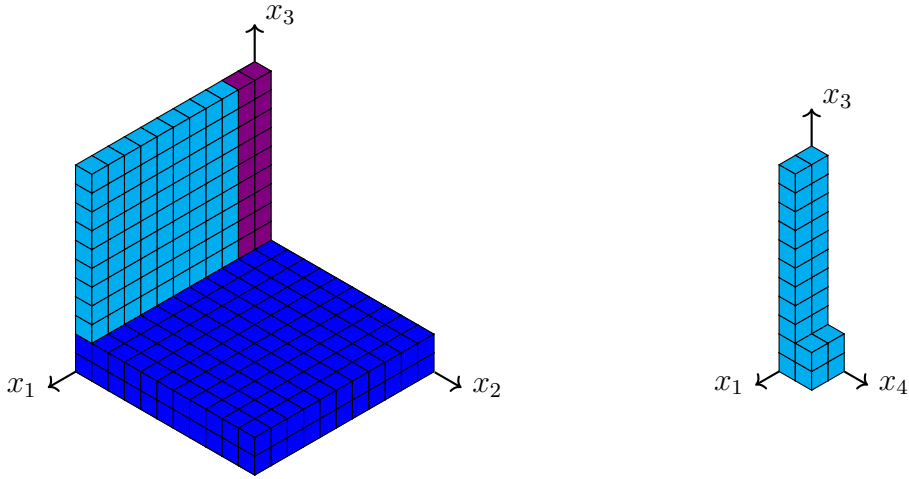

\begin{example} \label{ex:PT0=PT1ex2}
We provide another example similar to Example \ref{ex:CMclass}(iii). As in Example \ref{ex:solidpart}, we use (light) blue for boxes with height 0 in direction $x_4$ and purple for boxes with height $\infty$ in direction $x_4$. We also use dark blue for boxes with height 1 in direction $x_4$. In Figure \ref{fig3and3a} (left), we have two layers of dark blue boxes stacked on top of each other (along the $x_3$-axis) and extending infinitely far into the first quadrant of the $x_1x_2$-plane. We also have a layer of (light) blue boxes stacked on top and extending infinitely far into the first quadrant of the $x_1x_3$-plane. Finally, we have two columns of purple boxes extending infinitely far along the positive $x_3$-axis. There are 4 ``hidden'' boxes at weights $(i,0,j)$ with $0 \leq i,j \leq 1$ with colour purple. Note that the reduced support of the scheme corresponding to this solid partition is the union of the $x_1x_2$-, $x_1x_3$-, $x_3x_4$-planes. In this example, $W$ is a 1-dimensional scheme with 4 embedded points depicted in Figure \ref{fig3and3a} (right). It is scheme-theoretically supported in the hyperplane $(x_2 = 0)$.
\end{example}

By Proposition \ref{lem:limPTq},  the $(\C^*)^4$-fixed $\PT_0$ pairs on $\C^4$ with support $Z$ correspond to the finitely generated $(\C^*)^4$-invariant submodules of 
\begin{equation} \label{eqn:limtocalc}
\varinjlim \hom(\mathfrak{m}^r, \O_{Z}) / \O_Z,
\end{equation}
where $\mathfrak{m} \subset \O_{Z}$ denotes the ideal sheaf of the origin.  The support $Z$ may have embedded components of dimension 1, and we denote by $Z_1 \subset Z$ the underlying pure 2-dimensional subscheme. In particular, $Z_1 = Z_{\boldsymbol{\lambda}}$ for some finite partitions $\boldsymbol{\lambda}$ (Lemma \ref{lem:fixlocaffDT}) and we have
$$
Z_{\boldsymbol{\lambda}} = \bigcup_{1 \leq a<b \leq 4} Z_{\lambda_{ab}}.
$$
The following key proposition helps to determine \eqref{eqn:limtocalc} in many cases.
\begin{proposition} \label{lem:PT0limses}
Let $Z \subset \C^4$ be a 2-dimensional $(\C^*)^4$-fixed subscheme without embedded 0-dimensional components. Let $Z_1 = Z_{\boldsymbol{\lambda}} \subset Z$ be the underlying pure subscheme defined by $T_1(\O_Z) \subset \O_Z$. Define the scheme
$$
W = (Z_{\lambda_{12}} \cup Z_{\lambda_{13}} \cup Z_{\lambda_{23}}) \cap (Z_{\lambda_{14}} \cup Z_{\lambda_{24}} \cup Z_{\lambda_{34}})
$$
and consider the torsion subsheaf $T_0(\O_W) \subset \O_W$. Then the induced sequence
$$
0 \to  \varinjlim \hom(\mathfrak{m}^r, T_{1}(\O_Z) ) \to \varinjlim \hom(\mathfrak{m}^r, \O_{Z}) \to \varinjlim \hom(\mathfrak{m}^r, \O_{Z_1}) \to 0
$$ 
is exact and there is a $(\C^*)^4$-equivariant isomorphism
\begin{equation} \label{eqn:limT0}
\varinjlim \hom(\mathfrak{m}^r, \O_{Z_1}) / \O_{Z_1} \cong T_0(\O_W).
\end{equation}
\end{proposition}
\begin{proof}
We start with \eqref{eqn:limT0}. We first observe that 
$$
\varinjlim \hom(\mathfrak{m}^r, \O_{Z_{\lambda_{ab}}}) \cong \O_{Z_{\lambda_{ab}}}
$$
%For these calculations, it is helpful to calculate weight spaces for each fixed $t_4^{w_4}$ separately.
and more generally
\begin{align*}
&\varinjlim \hom(\mathfrak{m}^r,  \O_{Z_{\lambda_{12}} \cup Z_{\lambda_{13}} \cup Z_{\lambda_{23}}}) \cong  \O_{Z_{\lambda_{12}} \cup Z_{\lambda_{13}} \cup Z_{\lambda_{23}}}, \\
&\varinjlim \hom(\mathfrak{m}^r, \O_{Z_{\lambda_{14}} \cup Z_{\lambda_{24}} \cup Z_{\lambda_{34}}}) \cong \O_{Z_{\lambda_{14}} \cup Z_{\lambda_{24}} \cup Z_{\lambda_{34}}}.
\end{align*}
In particular
$$
\varinjlim \hom(\mathfrak{m}^r, \O_{Z_1}) \subset \O_{Z_{\lambda_{12}} \cup Z_{\lambda_{13}} \cup Z_{\lambda_{23}}} \oplus \O_{Z_{\lambda_{14}} \cup Z_{\lambda_{24}} \cup Z_{\lambda_{34}}},
$$
and $\varinjlim \hom(\mathfrak{m}^r, \O_{Z_1}) / \O_{Z_1}$ is supported on $W$. More precisely, we have
\begin{equation*}
\varinjlim \hom(\mathfrak{m}^r, \O_{Z_1}) / \O_{Z_1} \cong T_0(\O_W).
\end{equation*}
%This last isomorphism is true but non-trivial.

The first part follows by applying $\varinjlim \hom(\mathfrak{m}^r, \cdot ) $ to the short exact sequence
$$
0 \to T_{1}(\O_Z) \to \O_Z \to \O_{Z_1} \to 0
$$
and noting that $\lim \hom(\mathfrak{m}^r, \O_{Z}) \to \lim \hom(\mathfrak{m}^r, \O_{Z_1})$ is surjective.
\end{proof}

The point of Proposition \ref{lem:PT0limses} is that we understand $\varinjlim \hom(\mathfrak{m}^r, \O_{Z_1}) / \O_{Z_1}$. Furthermore, the quotient $\varinjlim \hom(\mathfrak{m}^r, T_{1}(\O_Z) ) / T_{1}(\O_Z)$ is 1-dimensional and is described similarly to the original work of Pandharipande-Thomas \cite{PT2} (extended to the 4-fold case in \cite{CK2} and subsequently studied further in \cite{Liu, CZ, Pia}). 

\begin{example} \label{ex:PT0=PT1lims}
For Example \ref{ex:solidpart}, we have $Z = Z_1$ and $\varinjlim \hom(\mathfrak{m}^r, \O_{Z}) / \O_{Z} \cong T_0(\O_W) = \mathcal{O}_W$ is depicted in Figure \ref{fig1b}. For Example \ref{ex:PT0=PT1ex2}, we also have $Z = Z_1$ but $T_0(\O_W) \neq \O_W$. In particular, the non-zero weight spaces of $\varinjlim \hom(\mathfrak{m}^r, \O_{Z}) / \O_{Z} \cong T_0(\O_W)$ correspond to the boxes at weights $(i,0,j,1)$ with $0 \leq i,j \leq 1$ in Figure \ref{fig3and3a}(right). 
\end{example}

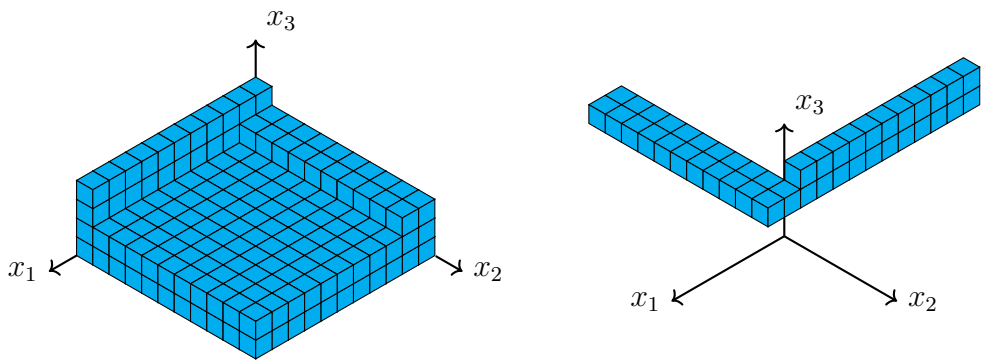
\begin{figure}[t]
\centering
\begin{subfigure}{.5\textwidth}
  \centering
  \input{fig2.tex}
  %\caption{A subfigure}
  %\label{fig:sub1}
\end{subfigure}%
\begin{subfigure}{.5\textwidth}
  \centering
  \input{fig2a.tex}
  %\caption{A subfigure}
  %\label{fig:sub2}
\end{subfigure}
\caption{The solid partition and direct limit of Example \ref{ex:PT0lims}.}
\label{fig2and2a}
\end{figure}

\begin{example} \label{ex:PT0lims}
In Example \ref{ex:PT0=PT1lims}, we always have $\varinjlim \hom(\mathfrak{m}^r, \O_{Z_1}) / \O_{Z_1} \cong T_0(\O_W) \neq 0$ and $\varinjlim \hom(\mathfrak{m}^r, T_{1}(\O_Z) ) / T_1(\O_Z) = 0$. We now discuss an example for which the former direct limit is zero and the latter is non-zero. Consider the solid partition of Figure \ref{fig2and2a}(left), which is defined by the ideal $(x_1^2 x_2 x_3^2,x_2x_3^3,x_3^4,x_4)$. The corresponding scheme is scheme-theoretically supported in $Z(x_4)$. Then $\varinjlim \hom(\mathfrak{m}^r, T_{1}(\O_Z) ) / T_1(\O_Z)$ is very similar to the type of modules studied in \cite[Sect.~2]{PT2}: its weight spaces are depicted in Figure \ref{fig2and2a}(right), where the only boxes in $\Z_{\geq 0}^4$ are the the ones at weights $(0,0,2,0)$ and $(1,0,2,0)$, and the two ``legs'' extend to $-\infty$ along the $x_1,x_2$ axes as indicated. 
\end{example}

The following corollary is as an equivariant analog of Proposition \ref{Lem.PT0=PT1}
\begin{corollary} \label{cor:PT0=PT1Quot}
Let $Z = Z_{\boldsymbol{\lambda}} \subset \C^4$ be a \emph{pure} 2-dimensional $(\C^*)^4$-fixed subscheme. Then the $(\C^*)^4$-fixed $\PT_0$ pairs on $\C^4$ with support $Z$ are in bijective correspondence with the closed points of
\footnote{${\curly Quot}_W(T_0(\O_W) ,m)$ is the Quot scheme of length $m$ 0-dimensional quotients of $T_0(\O_W)$.}
$$
\bigcup_{m \geq 0} {\curly Quot}_W(T_0(\O_W) ,m)^{(\C^*)^4},
$$
where $T_0(\O_W) \subset \O_W$ is the torsion subsheaf of 
$$
W = (Z_{\lambda_{12}} \cup Z_{\lambda_{13}} \cup Z_{\lambda_{23}}) \cap (Z_{\lambda_{14}} \cup Z_{\lambda_{24}} \cup Z_{\lambda_{34}}).
$$
More precisely, a $(\C^*)^4$-equivariant quotient $[T_0(\O_W) \twoheadrightarrow \zeta]$ with kernel $Q \subset T_0(W)$ corresponds to a $(\C^*)^4$-equivariant $\PT_0$ pair $(F,s)$ via the short exact sequences
\begin{displaymath}
\xymatrix
{
0 \ar[r] & \O_Z \ar@{=}[d] \ar[r] & F \ar[r] \ar@{^(->}[d] & Q \ar@^{(->}[d] \ar[r] & 0 \\
0 \ar[r] & \O_Z \ar[r] & \O_{Z_{\lambda_{12}} \cup Z_{\lambda_{13}} \cup Z_{\lambda_{23}}} \oplus \O_{Z_{\lambda_{14}} \cup Z_{\lambda_{24}} \cup Z_{\lambda_{34}}} \ar[r] & \O_{W} \ar[r] & 0, \\
}
\end{displaymath}
where the second square is Cartesian.
\end{corollary}

\begin{definition} \label{def:PT0boxconfig}
Fix plane partitions $\boldsymbol{\mu}$ with corresponding 2-dimensional subscheme $Z = Z_{\boldsymbol{\mu}} \subset \C^4$ (Definition \ref{def:partitionnotation}) and let $M:= \varinjlim \hom(\mathfrak{m}^r,\O_{Z})$, where $\mathfrak{m} \subset \O_{Z}$ denotes the ideal sheaf of the origin. We refer to the finitely generated $(\C^*)^4$-invariant submodules $B \subset M / \O_Z$ as \emph{$\PT_0$ box configurations} with asymptotics $\boldsymbol{\mu}$. We denote the collection of $\PT_0$ box configurations with asymptotics $\boldsymbol{\mu}$ by $\mathcal{B}^{(0)}_{\boldsymbol{\mu}}$. We say that $B \in \mathcal{B}^{(0)}_{\boldsymbol{\mu}}$ has \emph{no $\PT_0$ moduli} when there exists no $B' \in \mathcal{B}^{(0)}_{\boldsymbol{\mu}}$ with $B' \neq B$ and
$$
\dim_{\C}(B_w) = \dim_{\C}(B'_w), \quad \forall w \in \Z^4,
$$
where $B_w, B'_w$ denote the weight spaces of $B,B'$ for the character $w \in \Z^4$ respectively. We say that the choice of partitions $\boldsymbol{\mu}$ has no $\PT_0$ moduli when $\mathcal{B}^{(0)}_{\boldsymbol{\mu}}$ contains no elements with $\PT_0$ moduli. We refer to $|B| = \sum_w \dim_{\C}(B_w)$ as the \emph{size} of $B$.
\end{definition}

\begin{remark} \label{rem:PT0boxstacking}
Fix plane partitions $\boldsymbol{\mu}$ with corresponding 2-dimensional subscheme $Z = Z_{\boldsymbol{\mu}} \subset \C^4$ (Definition \ref{def:partitionnotation}). Evidently, the choice of partitions $\boldsymbol{\mu}$ has no $\PT_0$ moduli if and only if\footnote{Similar to the previous definition, $(M / \O_Z)_w$ denotes the weight space of $M / \O_Z$ for the character $w \in \Z^4$.}
$$
\dim_{\C} (M / \O_Z)_w \leq 1, \quad \forall w \in \Z^4.
$$
When this is the case, the elements of $\mathcal{B}^{(0)}_{\boldsymbol{\mu}}$ are particularly easy to describe. They are the finite subsets $B \subset \mathcal{R} := \{w \in \Z^4 \, : \, \dim_{\C} (M / \O_Z)_w = 1 \}$ such that \\

\noindent if $w = (w_1,w_2,w_3,w_4) \in   \mathcal{R}$ and one of $(w_1-1,w_2,w_3,w_4)$, $(w_1,w_2-1,w_3,w_4)$, $(w_1,w_2,w_3-1,w_4)$, $(w_1,w_2,w_3,w_4-1)$ lies in $B$, then $w \in B$.\footnote{Informally, the elements of $\mathcal{B}^{(0)}_{\boldsymbol{\mu}}$ are obtained by piling finitely many boxes in the region $ \mathcal{R}$ subject to the condition that gravity pulls in direction $(1,1,1,1)$.} 
\end{remark}

Proposition \ref{lem:PT0limses} is helpful for finding choices of plane partitions $\boldsymbol{\mu}$ with no $\PT_0$ moduli. For finite partitions $\boldsymbol{\lambda}$, there is a \emph{minimal} plane partition $\mu_1^{\boldsymbol{\lambda}}$ compatible with $\{\lambda_{ab}\}_{b \neq 1}$, i.e.,~satisfying \eqref{eqn:planepartcompat}.\footnote{Note that the scheme $Z_{\mu_1^{\boldsymbol{\lambda}}} \subset \C^3$ has no embedded 0-dimensional components.} Defining $\mu_a^{\boldsymbol{\lambda}}$ analogously for all $a$, we deduce that the choice of partitions 
\begin{equation} \label{eqn:lambdamu}
\boldsymbol{\mu}^{\boldsymbol{\lambda}} := (\mu_1^{\boldsymbol{\lambda}},\mu_2^{\boldsymbol{\lambda}},\mu_3^{\boldsymbol{\lambda}},\mu_4^{\boldsymbol{\lambda}})
\end{equation} 
has no $\PT_0$ moduli. In other words, $(\C^*)^4$-fixed $\PT_0$ pairs on $\C^4$, which are also $\PT_1$ pairs, have no moduli. 

\begin{example}
In Example \ref{ex:CMclass}(ii), we have $\O_{W_0} = \C$ (with weight $(0,0,0,0)$). In the second case of Example \ref{ex:CMclass}(iii) we have $\O_{W_0} = \C \cdot t_4$. In these cases there are no $\PT_0$ moduli.
\end{example}

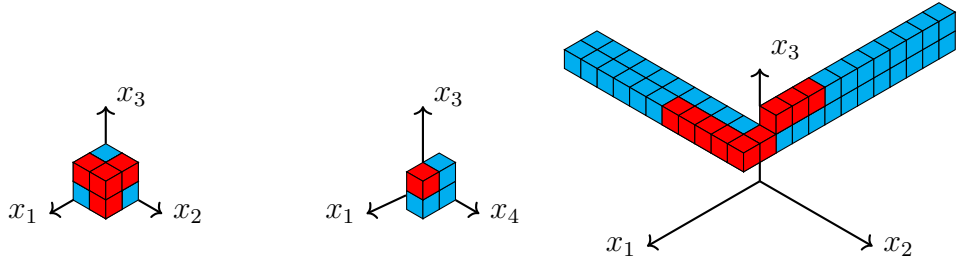
\begin{figure}[t]
\centering
\begin{subfigure}{.3\textwidth}
  \centering
  \input{fig1c.tex}
  %\caption{A subfigure}
  %\label{fig:sub1}
\end{subfigure}%
\begin{subfigure}{.3\textwidth}
  \centering
  \input{fig3b.tex}
  %\caption{A subfigure}
  %\label{fig:sub2}
\end{subfigure}
\begin{subfigure}{.3\textwidth}
  \centering
  \input{fig2b.tex}
  %\caption{A subfigure}
  %\label{fig:sub2}
\end{subfigure}
\caption{Examples of $\PT_0$ box configurations.}
\label{fig1cand3band2b}
\end{figure}

\begin{example} \label{ex:PT0boxconfs}
For Examples \ref{ex:PT0=PT1ex1}, \ref{ex:PT0=PT1ex2}, \ref{ex:PT0lims} we have drawn a blue box for each element of the region $\mathcal{R} \subset \mathbb{Z}^4$ defined in Remark \ref{rem:PT0boxstacking} --- see Figure \ref{fig1cand3band2b}. Recall that these boxes have no boxes stacked on top of them in the direction of the 4th axis not drawn in the diagram (in these examples: directions $x_4$, $x_2$, $x_4$ respectively). These blue boxes indicate the ``container'' inside of which we stack $\PT_0$ box configurations with gravity pulling into the $(1,1,1,1)$ direction, i.e., according to the rule of Remark \ref{rem:PT0boxstacking}. The boxes of the $\PT_0$ box configurations are drawn in red. From left to right, we have depicted a $\PT_0$ box configuration of size 4, 1, 9 respectively. Note that $\mathcal{R}$ has finite size in the first two examples and infinite size in the third example.
\end{example}

The following corollary gives many examples of choices of plane partitions with no $\PT_0$ moduli.
\begin{corollary} \label{lem:noPT0moduli}
Let $\boldsymbol{\lambda}$ be finite partitions and let $\boldsymbol{\mu}$ be compatible plane partitions, i.e., satisfying \eqref{eqn:planepartcompat}. Let $a,b,c,d \in \{1,2,3,4\}$ be mutually distinct, suppose $\mu_a$, $\mu_b$ have no further restrictions, and $\mu_c$, $\mu_d$ are minimal, i.e., $\mu_c = \mu_c^{\boldsymbol{\lambda}}$, $\mu_d = \mu_d^{\boldsymbol{\lambda}}$.\footnote{I.e., the scheme $Z_{\boldsymbol{\mu}}$ has embedded 1-dimensional components in at most two directions.} Then $\boldsymbol{\mu}$ has no $\PT_0$ moduli. 
\end{corollary}
\begin{proof}
We use Proposition \ref{lem:PT0limses}. Consider the weight spaces
$$
\Big( \varinjlim \hom(\mathfrak{m}^r, T_{1}(\O_Z) ) / T_{1}(\O_Z) \Big)_w, \quad \Big( \varinjlim \hom(\mathfrak{m}^r, \O_{Z_1}) / \O_{Z_1} \Big)_w
$$
These weight spaces are at most 1-dimensional and they are not 1-dimensional for the same $w$.
\end{proof}

\begin{example}
\hfill
\begin{enumerate}
\item[(i)] In general, when there are $\geq 3$ embedded legs, there will be $\PT_0$ moduli (as in the curve case). For example, let $I_Z = (x_1x_2x_3, x_1x_3^2, x_2 x_3^2, x_4)$. This corresponds to $\lambda_{12} = 1$, $\mu_1 = t_3 + \frac{1}{1-t_2}$, $\mu_2 = t_3 + \frac{1}{1-t_1}$, $\mu_3 = 1$, and $\lambda_{ab} = \mu_4 = \varnothing$ otherwise. Then $(\varinjlim \hom(\mathfrak{m}^r, \O_{Z}) / \O_{Z})_{(0,0,1,0)} = \C^2$ and there are $\PT_0$ moduli. 
\item[(ii)] However, there are examples with $\geq 3$ embedded legs for which there are no $\PT_0$ moduli such as $I_Z = (x_2x_3, x_1 x_3^2, x_1x_4, x_3x_4, x_4^2)$. This corresponds to $\lambda_{12} = 1$, $\mu_1 = t_3 + \frac{1}{1-t_2}$, $\mu_2 = t_4 + \frac{1}{1-t_1}$, $\mu_3 = 1$, and $\lambda_{ab} = \mu_4 = \varnothing$ otherwise.
\end{enumerate}
\end{example}

\noindent \textbf{$\PT_1$ case.} By Proposition \ref{lem:limPTq},  the $(\C^*)^4$-fixed $\PT_1$ pairs on $\C^4$ with support $Z$ correspond to the finitely generated $(\C^*)^4$-invariant submodules of 
\begin{equation*}
\varinjlim \hom(\mathfrak{I}^r, \O_{Z}) / \O_Z,
\end{equation*}
where $\mathfrak{I} \subset \O_{Z}$ denotes the ideal sheaf of the union of coordinate axes contained in $Z$. Recall that $Z = Z_{\boldsymbol{\lambda}}$ for some finite partitions $\boldsymbol{\lambda}$ (Lemma \ref{lem:fixlocaffDT}) and 
$$
Z_{\boldsymbol{\lambda}} = \bigcup_{1 \leq a<b \leq 4} Z_{\lambda_{ab}}.
$$
\begin{proposition} \label{lem:PT1lim}
Let $Z = Z_{\boldsymbol{\lambda}} \subset \C^4$ and denote by $\mathfrak{I} \subset \O_Z$ the ideal sheaf of the union of the coordinate axes contained in $Z$. Then
$$
\varinjlim \hom(\mathfrak{I}^r,\O_Z) = \bigoplus_{1 \leq a<b \leq 4} \O_{Z_{\lambda_{ab}}}[x_a^{-1}, x_b^{-1}].
$$
\end{proposition}
\begin{proof}
From the restrictions $ \O_Z \to \O_{Z_{\lambda_{ab}}}$, we obtain a short exact sequence
\begin{equation*} 
0 \to \O_Z \to \bigoplus_{1 \leq a<b \leq 4} \O_{Z_{\lambda_{ab}}} \to C \to 0,
\end{equation*}
where the cokernel $C$ is supported on the union of coordinate axes contained in $Z$. Then
$$
\varinjlim \hom(\mathfrak{I}^r,C) = 0
$$
and therefore
$$
\varinjlim \hom(\mathfrak{I}^r,\O_Z) = \bigoplus_{1 \leq a<b \leq 4} \varinjlim \hom(\mathfrak{I}^r,\O_{Z_{\lambda_{ab}}}).
$$
Finally, we have
\begin{equation*}
\varinjlim \hom(\mathfrak{I}^r,\O_{Z_{\lambda_{ab}}}) \cong \O_{Z_{\lambda_{ab}}}[x_a^{-1},x_b^{-1}]. \qedhere
\end{equation*}
\end{proof}

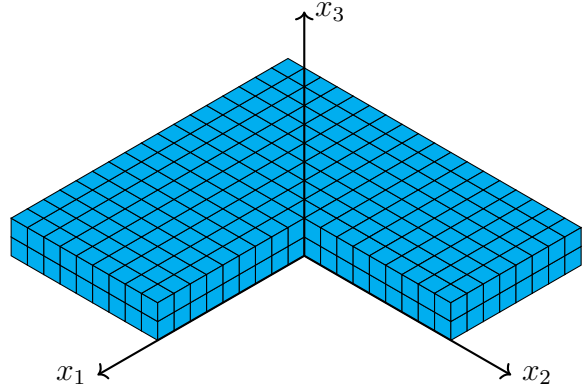
\begin{figure}[t]
\input{fig4.tex}
\caption{The direct limit of Example \ref{ex:PT1lim}.} \label{fig4}
\end{figure}

\begin{example} \label{ex:PT1lim}
Consider the solid partition corresponding to the scheme $Z$ defined by the ideal $(x_3^2,x_4)$. In Figure \ref{fig4} we have drawn the weight spaces of $\varinjlim \hom(\mathfrak{I}^r, \O_{Z}) / \O_Z$, where we recall that the blue colour indicates that there are no boxes stacked on top into the positive $x_4$-direction. The double layer of boxes fill all quadrants of the $x_1x_2$-plane except for the first quadrant (since we take the quotient by $\O_Z$).
\end{example}

 Similar to the $\PT_0$ case, we define:
 \begin{definition} \label{def:PT1boxconfig}
 Fix finite partitions $\boldsymbol{\lambda}$ with corresponding 2-dimensional subscheme $Z = Z_{\boldsymbol{\lambda}} \subset \C^4$ (Definition \ref{def:partitionnotation}). Let $M:= \varinjlim \hom(\mathfrak{I}^r, \O_{Z})$, where $\mathfrak{I} \subset \O_{Z}$ denotes the ideal sheaf of the union of the coordinate axes contained in $Z$. We refer to finitely generated $(\C^*)^4$-invariant submodules $B \subset M / \O_Z$ as \emph{$\PT_1$ box configurations} with asymptotics $\boldsymbol{\lambda}$. We denote the collection of $\PT_1$ box configurations with asymptotics $\boldsymbol{\lambda}$ by $\mathcal{B}^{(1)}_{\boldsymbol{\lambda}}$. We say that $B \in \mathcal{B}^{(1)}_{\boldsymbol{\lambda}}$ has \emph{no $\PT_1$ moduli} when there exists no $B' \in \mathcal{B}^{(1)}_{\boldsymbol{\lambda}}$ such that $B' \neq B$ and
$$
\dim_{\C}(B_w) = \dim_{\C}(B'_w), \quad \forall w \in \Z^4.
$$
We say that the choice of partitions $\boldsymbol{\lambda}$ has no $\PT_1$ moduli when $\mathcal{B}^{(1)}_{\boldsymbol{\lambda}}$ contains no elements with $\PT_1$ moduli. We refer to $|B| = \sum_w \dim_{\C}(B_w)$ as the \emph{size} of $B$. Unlike the $\PT_0$ case, $|B|$ may be infinite.
\end{definition}

\begin{remark} \label{rem:PT1boxstacking}
Fix finite partitions $\boldsymbol{\lambda}$ with corresponding 2-dimensional subscheme $Z = Z_{\boldsymbol{\lambda}} \subset \C^4$ (Definition \ref{def:partitionnotation}).
Then the choice of partitions has no $\PT_1$ moduli if and only if
$$
\dim_{\C} (M / \O_Z)_w \leq 1, \quad \forall w \in \Z^4.
$$
When this is the case, the elements of $\mathcal{B}^{(1)}_{\boldsymbol{\lambda}}$ are again easy to describe. They are the (not necessarily finite!) subsets $B \subset \mathcal{R} := \{w \in \Z^4 \, : \, \dim_{\C} (M / \O_Z)_w = 1 \}$ such that \\

\noindent if $w = (w_1,w_2,w_3,w_4) \in   \mathcal{R}$ and one of $(w_1-1,w_2,w_3,w_4)$, $(w_1,w_2-1,w_3,w_4)$, $(w_1,w_2,w_3-1,w_4)$, $(w_1,w_2,w_3,w_4-1)$ lies in $B$, then $w \in B$. Moreover, $B$ is required to support a \emph{finitely generated} module.
\end{remark}

\begin{figure}[t]
\input{fig4a.tex}
\caption{Example of a $\PT_1$ box configuration.} \label{fig4a}
\end{figure}
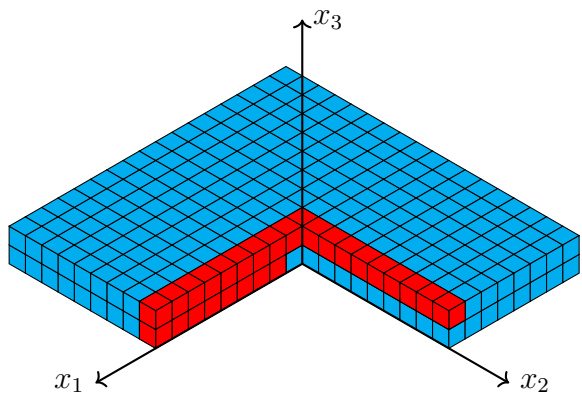

\begin{example}
For Example \ref{ex:PT1lim} we have drawn a blue box for each element of the region $\mathcal{R} \subset \mathbb{Z}^4$ defined in Remark \ref{rem:PT0boxstacking} --- see Figure \ref{fig4}. Similar to Example \ref{ex:PT0boxconfs}, these blue boxes indicate the ``container'' \emph{inside} of which we stack $\PT_1$ box configurations with gravity pulling into the $(1,1,1,1)$ direction, i.e., according to the rule of Remark \ref{rem:PT1boxstacking}. In Figure \ref{fig4a}, we have drawn in red an example of a $\PT_1$ box configuration. In this case, there are no finite $\PT_1$ box configurations. Note that placing a red box at $(-1,-1,1,0)$ forces us to put red boxes at $(i,-1,1,0)$ and $(-1,j,1,0)$ for all $i,j \geq 0$. Similarly, placing a red box at $(1,-1,0,0)$ forces us to put a red box at $(i,-1,0,0)$ for all $i \geq 2$.
\end{example}

The following is an immediate consequence of Proposition \ref{lem:PT1lim}.
\begin{corollary} \label{lem:noPT1moduli}
A collection of finite partitions $\boldsymbol{\lambda}$ has no $\PT_1$ moduli if and only if $\lambda_{ab} = \varnothing$ unless $(a,b) = (1,2)$ or $(3,4)$ (up to relabelling).
\end{corollary}

\begin{remark}
By Example \ref{ex:CMclass}(iii) it can happen that $\boldsymbol{\lambda}$ has $\PT_1$ moduli but that $\boldsymbol{\mu}^{\boldsymbol{\lambda}}$ (the minimal plane partitions compatible with $\boldsymbol{\lambda}$) has no $\PT_0$ moduli. This simply means that $(\C^*)^4$-equivariant $\PT_1$ pairs on supported on $Z_{\boldsymbol{\lambda}}$ with 1-dimensional cokernel may have moduli but those with 0-dimensional cokernel do not.
\end{remark}

We will see shortly that Corollaries \ref{lem:noPT0moduli} and \ref{lem:noPT1moduli} together imply that for toric Calabi-Yau 4-folds of the form $X = \mathrm{Tot}(L_1 \oplus L_2) \to S$, where $S$ is a smooth projective toric surface and $L_1, L_2$ line bundles on $S$ satisfying $L_1 \otimes L_2 \cong K_S$, the $T_X \cong (\C^*)^4$ fixed loci of the $\DT,\PT_0,\PT_1$ moduli spaces are always 0-dimensional (Remark \ref{rem:localsurf}).

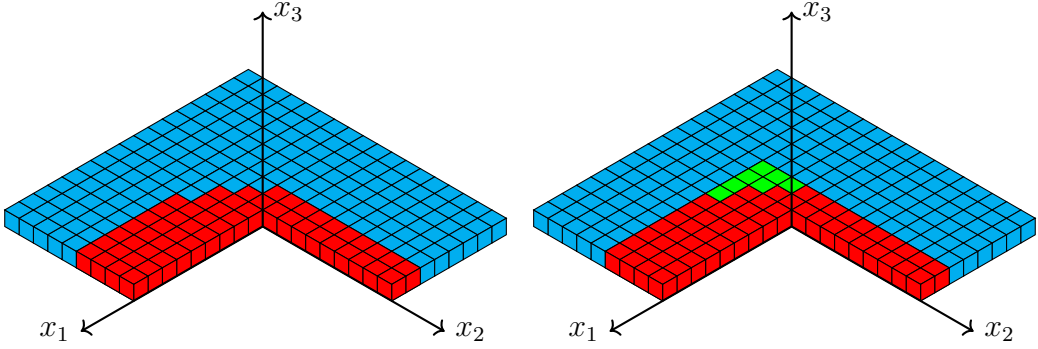
\begin{figure}[t]
\centering
\begin{subfigure}{.5\textwidth}
  \centering
  \input{fig5a.tex}
  %\caption{A subfigure}
  %\label{fig:sub1}
\end{subfigure}%
\begin{subfigure}{.5\textwidth}
  \centering
  \input{fig5b.tex}
  %\caption{A subfigure}
  %\label{fig:sub2}
\end{subfigure}
\caption{$\PT_1$ pair of Example \ref{ex:PT0lims}.}
\label{fig5aand5b}
\end{figure}

\begin{example}
Consider the solid partition corresponding to $S \cong \C^2$ defined by the ideal $(x_3,x_4)$. Then $\varinjlim \hom(\mathfrak{I}^r, \O_{S}) / \O_S = \C[x_1^{\pm 1},x_2^{\pm 1}] / \C[x_1,x_2]$. In Figure \ref{fig5aand5b}(left) we depict an example of a $\PT_1$ box configuration. Since $S$ is smooth, by Proposition \ref{cor:PT1smsupp}, $\PT_1$ pairs $(F,s)$ with support $S$ and non-zero cokernel $Q$ correspond to nestings $Z \subset C \subset S$ where $Z$ is 0-dimensional and $C$ is an effective divisor. Furthermore, the cokernel $Q$ satisfies $Q \cong \jmath_*I_{Z/C}(C)$, and $(F,s)$ is $(\C^*)^4$-equivariant if and only if $Z \subset C \subset S$ is a $(\C^*)^2$-invariant nesting. In Figure \ref{fig5aand5b}(right), we have drawn in red the weight spaces corresponding to $I_{Z/C}(C) \subset \O_C(C)$ and in green the weight spaces corresponding to $\O_Z(C)$.
\end{example}

\subsection{Fixed loci --- Global} \label{sec:fixlocglob}

We turn our attention to the characterization of $T_X \cong (\C^*)^4$ fixed $\PT_q$ stable pairs on an arbitrary smooth quasi-projective toric 4-fold $X$. Given a $T_X$-fixed $\PT_q$ stable pair $(F,s)$ on $X$, we can consider the restrictions
$$
(F^{(\alpha)}, s^{(\alpha)}) := (F,s)|_{U_\alpha}, \quad \alpha \in V(X).
$$
Since $U_\alpha \cong \C^4$, we can use the previous section to describe $(F^{(\alpha)}, s^{(\alpha)})$ combinatorially. Denote the scheme theoretic support of $(F,s)$ by $Z$ and recall the closed subschemes
$$
Z \supset Z_0 \supset Z_1
$$
defined by the torsion filtration of $\O_Z$. Recall that for $q=0$, we have $Z = Z_0$ and for $q=1$ we have $Z = Z_1$. We define $Z^{(\alpha)} := Z \cap U_\alpha$, $Z^{(\alpha)}_0 := Z_0 \cap U_\alpha$, and $Z^{(\alpha)}_1 := Z_1 \cap U_\alpha$. Then there exist collections of finite partitions and compatible \eqref{eqn:planepartcompat} plane partitions
$$
\boldsymbol{\lambda}^{(\alpha)} = \{\lambda_{ab}^{(\alpha)}\}_{1 \leq a < b \leq 4 \atop \alpha \in V(X)}, \quad \boldsymbol{\mu}^{(\alpha)} = \{\mu_{a}^{(\alpha)} \}_{1 \leq a \leq 4 \atop \alpha \in V(X)}
$$
such that (Definition \ref{def:partitionnotation})
$$
Z_0^{(\alpha)} = Z_{\boldsymbol{\mu}^{(\alpha)} }, \quad Z_1^{(\alpha)} = Z_{\boldsymbol{\lambda}^{(\alpha)}}.
$$ 
The partitions $\boldsymbol{\lambda}^{(\alpha)}, \boldsymbol{\mu}^{(\alpha)}$ satisfy gluing conditions. For two charts $U_\alpha,U_\beta$ intersecting in $\C^* \times \C^3$, let us label the coordinates in each chart such that the first coordinate corresponds to $\C^*$. Then
\begin{equation} \label{eqn:glucon1}
\mu_{1}^{(\alpha)}  = \mu_{1}^{(\beta)}, \quad \lambda^{(\alpha)}_{1b} = \lambda^{(\beta)}_{1b}, \quad \forall b = 2,3,4. 
\end{equation}
By the compatibility \eqref{eqn:planepartcompat}, the equations on plane partitions imply the ones on partitions.
For two charts $U_\alpha,U_\beta$ intersecting in $(\C^*)^2 \times \C^2$, let us label the coordinates in each chart such that the first two coordinates correspond to $(\C^*)^2$. Then
\begin{equation} \label{eqn:glucon2}
\lambda^{(\alpha)}_{12} = \lambda^{(\beta)}_{12}. 
\end{equation}
For $q \in \{-1,0\}$, the gluing conditions \eqref{eqn:glucon1}, \eqref{eqn:glucon2} completely determine when a collection of $(\C^*)^4$-equivariant $\PT_q$ pairs $\{(F^{(\alpha)}, s^{(\alpha)})\}_{\alpha \in V(X)}$ on the charts $U_\alpha$ glues to a $T_X$-equivariant $\PT_q$ pair on $X$. For $q=1$, one additionally requires that the $\PT_1$ box configurations glue.

In Section \ref{sec:moduli}, we discussed the moduli spaces 
\begin{equation} \label{eq:toricsecmoduli}
\PTqvX, \quad q \in \{-1,0,1\}, \quad v = (0,0,\gamma,\beta,n -  \gamma \cdot \td_2(X)) \in H_c^*(X,\Q).
\end{equation}
Let $X$ be a toric Calabi-Yau 4-fold and denote the (3-dimensional) subtorus preserving the Calabi-Yau volume form by $T \leq T_X$. The following theorem provides a nice collection of cases where Oh-Thomas localization for 0-dimensional reduced $T$-fixed loci as in Proposition \ref{prop:OTisored} applies.

\begin{theorem} \label{prop:fixlocmain1}
Let $\curly P:= \PTqvX$. 
\begin{enumerate}
\item[$\mathrm{(1)}$] For $q=-1$, we have $\curly P^{T_X} = \curly P^{T}$ and it is 0-dimensional and reduced. 
\item[$\mathrm{(2)}$] For $q=0$ and $\PT_0 = \PT_1$, i.e.~$F$ is pure for all $(F,s) \in \curly P^{T_X}$, we have $\curly P^{T_X} = \curly P^{T}$ and it is 0-dimensional and reduced.  
\item[$\mathrm{(3)}$] For $q=0$, suppose for any $(F,s) \in \curly P^{T_X}$ and any $\alpha \in V(X)$ the scheme theoretic support of $F|_{U_\alpha}$ has at most two embedded 1-dimensional components and its pure part is Cohen-Macaulay. Then $\curly P^{T_X} = \curly P^{T}$ and it is 0-dimensional and reduced.
\end{enumerate}
\end{theorem}

\begin{proof}
\noindent \textbf{$\DT$ case.} Clearly, we have
$$
\curP_v^{(-1)}(X)^{T_X} \subset \curP_v^{(-1)}(X)^T.
$$
By comparing $T$-fixed ideals and $(\C^*)^4$-fixed ideals on each chart $U_\alpha \cong \C^4$, it is easy to see that this is an equality at the level of closed points. The argument is written out in \cite[Lem.~3.1]{CK1} for $\gamma=\beta = 0$, but it works for 2-dimensional subschemes as well. Next, we claim that for any $Z \in \curP_v^{(-1)}(X)^{T_X}$, we have $\Hom_X(I_Z, \O_Z)^T = 0$, which establishes (1). 

In fact, we claim the stronger statement 
$$
\Hom_{U_\alpha}(I_{Z_{\alpha}}, \O_{Z_{\alpha}})^T = 0
$$
on each $U_\alpha \cong \C^4$. Since $I_{Z_{\alpha}}$ is a monomial ideal, $\Hom_{U_\alpha}(I_{Z_{\alpha}}, \O_{Z_{\alpha}})$ has an explicit combinatorial basis in terms of so-called Haiman arrows \cite{MS}. From this description, it is clear that there are no non-zero $T$-fixed Haiman arrows. This is worked out in \cite[Lem.~3.6]{CK1} for $\gamma=\beta = 0$, but the argument holds for 2-dimensional subschemes as well (but it is crucial that $Z_\alpha$ has codimension $\geq 2$). \\

\noindent \textbf{$\PT_0$ case.}
Since the fixed locus $\curly P^T$ is proper and $\curly P^{T_X} \subset \curly P^T$, in order to show that $\curly P^{T} = \curly P^{T_X}$ is 0-dimensional reduced, it suffices to show that $\Hom_X(I^\mdot, F)^T = 0$ for any $I^\mdot = [\O_X \to F] \in \curly P^{T_X}$. 

We start by following \cite[Sect.~3.1, 3.2]{PT2}. Take $I^\mdot = [\O_X \to F] \in \curly P^{T_X}$, denote its cokernel by $Q$, and its scheme theoretic support by $Z$. We write $I^\mdot_\alpha := I^\mdot|_{U_\alpha}$, $F_\alpha := F|_{U_\alpha}$, $Q_\alpha := Q|_{U_\alpha}$, and $Z_\alpha = Z|_{U_\alpha}$. Consider the restriction morphism
$$
\Hom_X(I^\mdot, F) \to \bigoplus_{\alpha \in V(X)} \Hom_{U_\alpha}(I_\alpha^\mdot, F_\alpha).
$$
From the exact triangle $I_{Z_\alpha} \to I^\mdot_\alpha \to Q_\alpha[-1]$ we obtain the left exact sequence
$$
0 \to \Ext^1(Q_\alpha, F_\alpha) \to \Hom(I_\alpha^\mdot, F_\alpha) \to \Hom(I_{Z_{\alpha}}, F_\alpha),
$$
which remains left exact after taking $T$-fixed part. Suppose $\Hom(I_{Z_\alpha},F_\alpha)^T = 0$ for all $\alpha \in V(X)$. Then $\Ext^1(Q_\alpha, F_\alpha)^T \cong \Hom(I_\alpha^\mdot, F_\alpha)^T$. Moreover, since $Q_\alpha$ is 0-dimensional, we deduce (by \cite{PT2} paragraph below equation (3.2)) that the restriction morphism induces an isomorphism
\begin{equation*} 
\Hom_X(I^\mdot, F)^T \cong \bigoplus_{\alpha \in V(X)} \Hom_{U_\alpha}(I_\alpha^\mdot, F_\alpha)^T.
\end{equation*}

Next, we define $M_\alpha:=\varinjlim \hom(\mathfrak{m}_\alpha^r, \O_{Z_\alpha})$, where $\mathfrak{m}_\alpha \subset \O_{Z_\alpha}$ is the ideal sheaf of the origin. In fact, it is easier to work with a sufficient approximation $M_\alpha^r:=\hom(\mathfrak{m}_\alpha^r, \O_{Z_\alpha})$, with $r \gg 0$ such that $F_\alpha \subset M_\alpha^r$ (so $M_\alpha^r$ is coherent).
We then obtain an exact sequence
$$
\Hom(Q_\alpha,M_\alpha^r) \to \Hom(Q_\alpha,M_\alpha^r / F_\alpha ) \to \Ext^1(Q_\alpha,F_\alpha) \to \Ext^1(Q_\alpha,M_\alpha^r).
$$

Summarizing: in order to prove that $\curly P^T$ is 0-dimensional reduced, it suffices to show that for any $I^\mdot = [\O_X \to F] \in \curly P^{T_X}$ we have $\Hom(I_{Z_{\alpha}}, F_\alpha)^T = \Ext^1(Q_\alpha, F_\alpha)^T = 0$ for all $\alpha \in V(X)$. In turn, for this it suffices to show
\begin{align*}
\Hom(I_{Z_{\alpha}}, F_\alpha)^T &= 0, \\
\Ext^1(Q_\alpha,M^r_\alpha)^T &= 0, \\
\Hom(Q_\alpha,M_\alpha^r / F_\alpha )^T &=0,
\end{align*}
for all $I^\mdot = [\O_X \to F] \in \curly P^{T_X}$ and $\alpha \in V(X)$.

We now prove (2). We first study $\Hom(I_{Z_{\alpha}}, F_\alpha)^{T}$. Since we assume $\PT_0 = \PT_1$, $Z_\alpha$ is pure 2-dimensional. 
We first observe that, by the $\DT$ case, we have $\Hom(I_{Z_{\alpha}}, \O_{Z_\alpha})^{T} = 0$. We now consider $\Hom(I_{Z_{\alpha}}, Q_\alpha)^{T}$. Take a homogeneous element of $\Hom(I_{Z_{\alpha}}, Q_\alpha)^{T}$. By Proposition \ref{lem:PT0limses}, it will map a homogeneous element $g$ of $I_{Z_{\alpha}}$ to a weight space contained in $T_0(\O_W)$. Then $g$ maps to zero by an argument similar to the $\DT$ case. 

It suffices to show that $\Ext^1(Q_\alpha, F_\alpha)^T = 0$. Writing $Z_\alpha = Z_{\boldsymbol{\lambda}}$, we set $S_1 := Z_{\lambda_{12}} \cup Z_{\lambda_{13}} \cup Z_{\lambda_{23}}$ and $S_2 := Z_{\lambda_{14}} \cup Z_{\lambda_{24}} \cup Z_{\lambda_{34}}$. Then Corollary \ref{cor:PT0=PT1Quot} gives a short exact sequence
$$
\cdots \rightarrow \Hom(Q_\alpha, \O_{S_1 \cap S_2} / Q_\alpha) \rightarrow \Ext^1(Q_\alpha, F_\alpha) \rightarrow \Ext^1(Q_\alpha,\O_{S_1} \oplus \O_{S_2}) \rightarrow \cdots,
$$
where the scheme theoretic support of $Q_\alpha$ is contained in $S_1 \cap S_2$. 
Since $Q_\alpha$ is 0-dimensional and $\O_{S_i}$ is pure 2-dimensional and Cohen-Macaulay (Proposition \ref{lem:casesCM}), we have $\Ext^1(Q_\alpha, \O_{S_i}) = 0$. Indeed, since $Q_\alpha$ is 0-dimensional $\Ext^1(Q_\alpha, \O_{S_i})$ is finite-dimensional and dual to $\Ext^3(\O_{S_i},Q_\alpha)$. At the stalk at 0, the module $\O_{S_i,0}$ has depth 2 (Cohen-Macaulay property) and hence projective dimension 2 (Auslander-Buchsbaum), so 
$$
\Ext_{\O_{\C^4,0}}^3(\O_{S_i,0},Q_{\alpha,0}) = 0.
$$
Furthermore, similar to the $\DT$ case, one can show
$$
\Hom(Q_\alpha, \O_{S_1 \cap S_2} / Q_\alpha)^T = 0.
$$

Next, we prove (3). By assumption, the support of $T_1(\O_{Z_\alpha})$ has at most two irreducible components for all $\alpha \in V(X)$. We first study $\Hom(I_{Z_{\alpha}}, F_\alpha)^{T}$ and note that, by the $\DT$ case, $\Hom(I_{Z_{\alpha}}, \O_{Z_\alpha})^{T} = 0$. We consider $\Hom(I_{Z_{\alpha}}, Q_\alpha)^{T}$. Take a homogeneous element of $\Hom(I_{Z_{\alpha}}, Q_\alpha)^{T}$. It maps a homogeneous element $g$ of $I_{Z_{\alpha}}$ to a weight space of $Q_\alpha$ supported on $$\varinjlim \hom(\mathfrak{m}^r_\alpha, T_1(\O_{Z_\alpha})) / T_1(\O_{Z_\alpha})$$ by Proposition \ref{lem:PT0limses}. It can then be seen that $g$ maps to zero by a reasoning similar to \cite[Lem.~2]{PT2}. Thus $\Hom(I_{Z_{\alpha}}, F_\alpha)^{T} = 0$.

By assumption, the pure part $Z_1 \subset Z_\alpha$ is Cohen-Macaulay. Propositions \ref{lem:PT0limses}, \ref{lem:casesCM} give a short exact sequence
$$
0 \to \varinjlim \hom(\mathfrak{m}_\alpha^r, T_{1}(\O_{Z_\alpha}) ) \to M_\alpha \to \O_{Z_1} \to 0.
$$ 
Furthermore, by assumption, the support of $T_{1}(\O_{Z_\alpha})$ has at most two irreducible components. In particular, $(F_\alpha,s_\alpha)$ has no $\PT_0$ moduli (Corollary \ref{lem:noPT0moduli}). Without loss of generality, we take the reduced support of $T_{1}(\O_{Z_\alpha})$ to be contained in the union of the $x_1$- and $x_2$-axis. We first show $\Ext^1(Q_\alpha,M^r_\alpha)^T = 0$.
There exists a $(\C^*)^3$-equivariant 0-dimensional $\C[x_2,x_3,x_4]$-module $N_1$ and $(\C^*)^3$-equivariant 0-dimensional $\C[x_1,x_3,x_4]$-module $N_2$ such that 
$$
\varinjlim \hom(\mathfrak{m}^r_\alpha, T_{1}(\O_{Z_\alpha}) ) = M_1 \oplus M_2, \quad M_i  = N_i[x_i,x_i^{-1}], \quad i=1,2.
$$
Using the truncated versions $M_\alpha^r, M_1^r, M_2^r$ of $M_\alpha, M_1, M_2$, we have a short exact sequence
\begin{equation} \label{eqn:truncses}
0 \to M_1^r \oplus M_2^r \to M_\alpha^r \to \O_{Z_1} \to 0.
\end{equation}
Note that
$$
M_i^{r+N} = M_i^r \otimes t_i^{-N},
$$
for all $N > 0$. Since $\Ext^1(Q_\alpha,M_i^r)$ are finite $T$-representations ($Q_\alpha$ is 0-dimensional), we can take $N$ sufficiently large such that
$$
\Ext^1(Q_\alpha,M_i^{r+N}) = \Ext^1(Q_\alpha,M_i^r) \otimes t_i^{-N}
$$
has no $T$-fixed part (this trick is also used in \cite[Lem.~1]{PT2}). Furthermore, since $Q_\alpha$ is 0-dimensional and $Z_1$ is 2-dimensional Cohen-Macaulay, we have $\Ext^1(Q_\alpha, \O_{Z_1}) = 0$ (as we saw above). Consequently, for $r \gg 0$, \eqref{eqn:truncses} implies $\Ext^1(Q_\alpha,M^r_\alpha)^T = 0$. 

Finally, we show that $\Hom(Q_\alpha,M_\alpha^r / F_\alpha )^T = 0$. Consider 
$$
\Hom(Q_\alpha,M_\alpha^r / F_\alpha ) \cong  \Hom(Q_\alpha, (\hom(\mathfrak{m}^r_\alpha, T_{1}(\O_{Z_\alpha}) ) \cap F_\alpha) / T_{1}(\O_{Z_\alpha})).
$$
The vanishing of the $T$-fixed part on the right hand side follows by an argument which is almost the same as \cite[Prop.~2.6, Step 3]{CK2}.
\end{proof}

We expect that there are many more cases where $\curly P^{T_X} = \curly P^{T}$ is 0-dimensional and reduced, but the combinatorics involved in the tangent space analysis becomes increasingly complicated:
\begin{conjecture} \label{conj:0dimTXfix}
Let $X$ be a toric Calabi-Yau 4-fold and $\curly P := \PTqvX$ for $q \in \{-1,0,1\}$. If $\curly P^{T_X}$ is 0-dimensional, then $\curly P^{T_X} = \curly P^{T}$ and it is reduced.
\end{conjecture}

Evidence for the conjecture is provided by the cases proved in Theorem \ref{prop:fixlocmain1}. Further indications that this conjecture is true come from the fact that in all examples we calculated where $\curly P^{T_X}$ is 0-dimensional (including, e.g., for $q=1$ the examples in Section \ref{sec:localP2P3}), we never found a $T$-fixed term in $T_{\curly P}^{\vir}|_{I\udot}$. This does not imply that $\Hom_X(I\udot, F)^{T}$ is zero, because it could cancel against a $T$-fixed term of the same dimension in the obstructions. See also \cite[Sect.~3.3]{PT2} for an analogous result in the 3-fold case.

\begin{remark} \label{rem:localsurf}
Suppose $X = \mathrm{Tot}(L_1 \oplus L_2) \to S$, where $S$ is a smooth projective toric surface and $L_1, L_2$ line bundles on $S$ satisfying $L_1 \otimes L_2 \cong K_S$. Let $\curly P := \PTqvX$ for $q \in \{-1,0,1\}$ and $v = (0,0,\gamma,\beta,n - \gamma \cdot \td_2(X)) \in H_c^*(X,\Q)$. By Corollaries \ref{lem:noPT0moduli}, \ref{lem:noPT1moduli}, and the gluing description at the beginning of this section, it is clear that $\curly P^{T_X}$ is always 0-dimensional. For $q \in \{-1,0\}$, we have $\curly P^{T_X} = \curly P^{T}$ is reduced by Theorem \ref{prop:fixlocmain1}(1) and (3). For $q=1$ and assuming Conjecture \ref{conj:0dimTXfix}, we have $\curly P^{T_X} = \curly P^{T}$ is reduced.
\end{remark}

\section{Vertex formalism} \label{sec:vertex}

In this section, we investigate the calculation of $K$-theoretic $\PT_q$ invariants of toric Calabi-Yau 4-folds for $q \in \{-1,0,1\}$. We also formulate a $K$-theoretic $\DT$--$\PT_0$ vertex correspondence, which we verify in some examples.

\subsection{Redistribution} \label{sec:redist}

Let $X$ be a toric Calabi-Yau 4-fold. Let $\curly P:= \PTqvX$ with $q \in \{-1,0,1\}$ with $v$ as in \eqref{eq:toricsecmoduli}. We now present a method to determine the virtual tangent representation\footnote{Under this isomorphism, $[\O_{\pt}]$ corresponds to $1$.}
$$
R\Hom_X(I^\mdot,I^\mdot)_0[1] \in K_0^{T_X}(\pt) \cong \Z[t_1^{\pm 1},t_2^{\pm 1},t_3^{\pm 1},t_4^{\pm 1}],
$$
for any $T_X$-fixed $I^\mdot \in \curly P$. Note that the spaces $\Ext_X^i(I\udot,I\udot)_0$ are finite-dimensional (Section \ref{sec:moduli}). We observe that the rank of this complex is even because $N_{S_f/X}$ is split for all $f \in F(X)$ (Lemma \ref{lem:intvd}). There exists a smooth projective toric 4-fold $\overline{X}$ and a $(\C^*)^4$-equivariant good compactification $j : X \hookrightarrow \overline{X}$. We then have
\[
\vd := n - \frac{1}{2} (j_* \gamma)^2 \in \Z.
\]
We write out the case $q=-1$ and comment on the differences for $q \in \{0,1\}$ in Remarks \ref{rem:modPT0}, \ref{rem:modPT1}.
Then $I^\mdot  = I_Z$ for some $T_X$-fixed subscheme $Z \subset X$. The $q=-1$ was recently also developed by Nekrasov-Piazzalunga \cite{NP2} (see Section \ref{sec:phys}).

Consider the open affine invariant cover $\{U_\alpha\}_{\alpha \in V(X)}$ of $X$. We write $U_{\alpha\beta} = U_\alpha \cap U_{\beta}$, $U_{\alpha\beta\gamma} = U_\alpha \cap U_\beta \cap U_\gamma$, etc.~and we choose a total order on the indices $\alpha$. We also write $Z_\alpha = Z \cap U_\alpha$,  $I_\alpha = I_{Z_{\alpha}}$, $Z_{\alpha\beta} = Z \cap U_{\alpha\beta}$,  $I_{\alpha\beta} = I_{Z_{\alpha\beta}}$, $Z_{\alpha\beta\gamma} = Z \cap U_{\alpha\beta\gamma}$,  $I_{\alpha\beta\gamma} = I_{Z_{\alpha\beta\gamma}}$ etc. Using the local-to-global spectral sequence and \v{C}ech cohomology, as in \cite[Sect.~4.6]{MNOP}, we obtain
\begin{align*}
R \mathrm{Hom}_X(I_Z,I_Z)_0[1] = &\sum_{\alpha \in V(X)} R\mathrm{Hom}_{U_\alpha}(I_\alpha,I_\alpha)_0[1] - \sum_{\alpha<\beta \in V(X)} R\mathrm{Hom}_{U_{\alpha\beta}}(I_{\alpha\beta},I_{\alpha\beta})_0[1] \\
&+ \sum_{\alpha<\beta<\gamma \in V(X)} R\mathrm{Hom}_{U_{\alpha\beta\gamma}}(I_{\alpha\beta\gamma},I_{\alpha\beta\gamma})_0[1] - \cdots.
\end{align*}
We always label edges $\alpha\beta \in E(X)$ so that $\alpha < \beta$. The previous equation should be read as an identity in the additive group $\mathbb{Q}(\!(t_1,t_1^{-1}, \ldots, t_4,t_4^{-1})\!)$: the individual summands on the right hand side are generally not finite but, when adding them, there are infinitely many cancellations resulting in an element in $\mathbb{Q}[t_1,t_1^{-1}, \ldots, t_4,t_4^{-1}] \subset \mathbb{Q}(\!(t_1,t_1^{-1}, \ldots, t_4,t_4^{-1})\!)$.

\begin{remark} \label{rem:alphabetagamma}
Since $Z$ is $T_X$-invariant, 2-dimensional, and proper, the only contributions in the above sum are:
\begin{itemize}
\item $\alpha \in V(X)$, which corresponds to $U_\alpha \cong \C^4$,
\item $\alpha < \beta \in V(X)$ such that $\alpha\beta \in E(X)$, which corresponds to $U_{\alpha\beta} \cong \C^* \times \C^3$,
\item $\alpha < \beta \in V(X)$ such that $\alpha\beta \notin E(X)$ but $\alpha, \beta \in f$ for some $f \in F(X)$, which corresponds to $U_{\alpha\beta} \cong (\C^*)^2 \times \C^2$,
\item $\alpha_1 < \cdots < \alpha_N \in V(X)$ with $N \geq 3$ such that $\alpha_1, \ldots, \alpha_N \in f$ for some $f \in F(X)$, which corresponds to $U_{\alpha_1 \cdots \alpha_N} \cong (\C^*)^2 \times \C^2$.
\end{itemize}
%The isomorphisms in 3 and 4 are not true for any smooth quasi-projective toric 4-fold. But they are true for toric Calabi-Yau 4-folds (always with the convention on cones in ``conventions''). 
\end{remark}

In order to deduce these statements, it is useful to recall the description of fans of toric Calabi-Yau 4-folds. Let $N = \mathbb{Z}^4$ and $N_{\R} = N \otimes_{\Z} \R$. Consider the hyperplane $\{x_4 = 1\} \subset N_{\R}$ and any collection of (not necessarily regular) rational tetrahedra in $\{x_4=1\}$ such that any two are disjoint or intersect in a common vertex, edge, or facet. Then we can consider the induced fan $\Delta$ of strongly convex rational polyhedral cones in $N_{\R}$ spanned by the origin and the tetrahedra in $\{x_4=1\}$. Fans obtained in this way are precisely the fans of toric Calabi-Yau 4-folds.

By Lemma \ref{lem:fixlocaffDT}, $I_\alpha \subset \O_{U_\alpha} \cong  \O_{\C^4}$ corresponds to a solid partition $\pi^{(\alpha)} \subset \Z_{\geq 0}^4$ with asymptotic finite partitions and (compatible) plane partitions
$$
\boldsymbol{\lambda}^{(\alpha)} =  \{\lambda_{ab}^{(\alpha)}\}_{1 \leq a<b \leq 4}, \quad \boldsymbol{\mu}^{(\alpha)} = \{ \mu_{a}^{(\alpha)}\}_{a=1}^{4}.
$$
Suppose we choose coordinate functions $(x_1,x_2,x_3,x_4)$ on $U_\alpha \cong \C^4$ such that $t:=(t_1,t_2,t_3,t_4) \in T_X \cong (\C^*)^4$ acts by
\begin{equation} \label{eqn:stndact}
t \cdot (x_1,x_2,x_3,x_4) = (t_1x_1,t_2x_2,t_3x_3,t_4x_4).
\end{equation}
Then we denote the $T_X$-character of $\O_{Z_\alpha}$ by\footnote{For any $M \in K_0^{(\C^*)^4}(\C^4)$, we denote by $\tr_M$ the virtual $(\C^*)^4$-representation of $M$. It is viewed as a rational function in $t_1,t_2,t_3,t_4$ and the dimensions of the weight spaces of $M$ are obtained by expanding this rational function in \emph{ascending} powers of $t_1,t_2,t_3,t_4$. For example, for $M = \O_Z$ and $Z = Z(x_1,x_2) \cup Z(x_3,x_4)$, we have $\tr_{\O_Z} = (1-t_1)^{-1}(1-t_2)^{-1}+(1-t_3)^{-1}(1-t_4)^{-1} -1$.}\begin{equation} \label{Zalpha}
Z_\alpha := \tr_{\O_{Z_\alpha}} = \sum_{(i_1,i_2,i_3,i_4) \in \pi^{(\alpha)}} t_1^{i_1} t_2^{i_2} t_3^{i_3} t_4^{i_4}.
\end{equation}
Let $\alpha\beta \in E(X)$. Then $U_{\alpha\beta} \cong \C^* \times \C^3$, and we label coordinates $(x_1,x_2,x_3,x_4)$ such that $x_1$ corresponds to the first $\C^*$ factor (and $T_X \cong (\C^*)^4$ acts by the standard torus action \eqref{eqn:stndact}).
Let $\alpha_1, \ldots, \alpha_N \in V(X)$ such that $\alpha_1, \ldots, \alpha_N \in f$ for some $f \in F(X)$, where for $N=2$, we assume $\alpha_1\alpha_2 \not\in E(X)$. 
Then $U_{\alpha_1 \cdots \alpha_N} \cong (\C^*)^2 \times \C^2$ (Remark \ref{rem:alphabetagamma}).
We label coordinates $(x_1,x_2,x_3,x_4)$ on $U_{\alpha_1 \cdots \alpha_N}$ such that $(x_1,x_2)$ corresponds to the $(\C^*)^2$ factor (and $T_X \cong (\C^*)^4$ acts by the standard torus action \eqref{eqn:stndact}).
Similar to \eqref{Zalpha}, we define 
$$
Z_{\alpha\beta} := \tr_{\O_{Z_{\alpha\beta}}}, \quad Z_{\alpha_1\cdots\alpha_N} := \tr_{\O_{Z_{\alpha_1\cdots\alpha_N}}}.
$$

Using the Taylor resolution \cite{Tay}
$$
0 \to F_s \to \cdots \to F_1 \to F_0 \to I_\alpha \to 0,
$$
a calculation similar to \cite[Sect.~4.7]{MNOP} gives
\begin{align} \label{prevertex}
\tr_{R \mathrm{Hom}_{U_\alpha}(I_\alpha,I_\alpha)_0[1]} =  Z_\alpha + \frac{\overline{Z}_{\alpha}}{t_1t_2t_3t_4} - \frac{P_{1234}}{t_1t_2t_3t_4} Z_\alpha \overline{Z}_\alpha,
\end{align}
where $\overline{(\cdot)}$ denotes the dual representation, which amounts to replacing each $t_i$ by $t_i^{-1}$, and
\begin{align*}
P_{1234} :=(1-t_1)(1-t_2)(1-t_3)(1-t_4).
\end{align*}
The rational function \eqref{prevertex} in $t_1,t_2,t_3,t_4$ should be expanded in \emph{ascending} powers of $t_i$.

Furthermore, we have
\begin{align} \label{preedge} 
\tr_{R \mathrm{Hom}_{U_{\alpha\beta}}(I_{\alpha\beta},I_{\alpha\beta})_0[1]} =  \delta(t_1)\Bigg(Z_{\alpha\beta} - \frac{\overline{Z}_{\alpha\beta}}{t_2t_3t_4} + \frac{P_{234}}{t_2t_3t_4} Z_{\alpha\beta} \overline{Z}_{\alpha\beta}\Bigg),
\end{align}
where 
\begin{align*}
P_{ijk}:=(1-t_i)(1-t_j)(1-t_k), 
 \quad\delta(t_i) := \sum_{i \in \Z} t_i.
\end{align*}
The rational function \eqref{preedge} in $t_2,t_3,t_4$ between brackets should be expanded in \emph{ascending} powers of $t_2, t_3, t_4$.

Finally, we have
\begin{align*} 
\tr_{R \mathrm{Hom}_{U_{\alpha_1\cdots\alpha_N}}(I_{\alpha_1\cdots\alpha_N},I_{\alpha_1\cdots\alpha_N})_0[1]} =  \delta(t_1)\delta(t_2)\Bigg(Z_{\alpha_1\cdots\alpha_N} + \frac{\overline{Z}_{\alpha_1\cdots\alpha_N}}{t_3t_4} - \frac{P_{34}}{t_3t_4} Z_{\alpha_1\cdots\alpha_N} \overline{Z}_{\alpha_1\cdots\alpha_N}\Bigg),
\end{align*}
where
\begin{align*}
\quad P_{ij}:=(1-t_i)(1-t_j). 
\end{align*}
Note that 
$$
Z_{\alpha_1\cdots\alpha_N} + \frac{\overline{Z}_{\alpha_1\cdots\alpha_N}}{t_3t_4} - \frac{P_{34}}{t_3t_4} Z_{\alpha_1\cdots\alpha_N} \overline{Z}_{\alpha_1\cdots\alpha_N}
$$
is a Laurent \emph{polynomial} in $t_3, t_4$. However, \eqref{prevertex} and the term between brackets in \eqref{preedge} are in general infinite Laurent \emph{series} (in ascending powers of the $t_i$). We therefore redistribute the terms. 
This means we want to attach Laurent polynomials to all elements of $V(X)$, $E(X)$, $F(X)$ which sum up to the $T_X$-representation $R\Hom_X(I_Z,I_Z)_0[1]$.

In this section, we denote the vertices neighbouring a vertex $\alpha$ by $\beta_1,\beta_2,\beta_3,\beta_4$.\footnote{Since $X$ is non-compact, $\alpha$ may lie on a non-compact edge in which case there are $<4$ neighbouring vertices.} Then we label the coordinates $(x_1,x_2,x_3,x_4)$ on $U_\alpha \cong \C^4$ such that
$$
L_{\alpha\beta_i} = Z(x_j,x_k,x_l),
$$
for all $i,j,k,l \in \{1,2,3,4\}$ distinct (and $T_X \cong (\C^*)^4$ acts by the standard torus action \eqref{eqn:stndact}). 

In the following four definitions, we fix $Z \subset X$ a proper 2-dimensional $T_X$-fixed subscheme. Recall that, for any $\alpha$, $Z_\alpha$ determines partitions in various dimensions (Lemma \ref{lem:fixlocaffDT})
$$
\{\pi^{(\alpha)}\}_{\alpha \in V(X)}, \quad \{\boldsymbol{\mu}^{(\alpha)}\}_{\alpha \in V(X)}, \quad \{\boldsymbol{\lambda}^{(\alpha)}\}_{\alpha \in V(X)}.
$$
\begin{definition} \label{def:renorm}
By the inclusion-exclusion principle, there exist unique Laurent \emph{polynomials} $W_{\alpha\beta_i}$, $W_\alpha$ such that $W_{\alpha\beta_i}$ does not depend on $t_i$ and 
\begin{align*}
Z_\alpha =  \sum_{1 \leq i<j \leq 4} \frac{Z_{\alpha\beta_i\beta_j}}{(1-t_i)(1-t_j)} + \sum_{i=1}^{4} \frac{W_{\alpha\beta_i}}{1-t_i} + W_\alpha.
\end{align*}
Recall that $Z_{\alpha\beta_i\beta_j}$ are also Laurent polynomials. Then $|Z_\alpha| := \rk(W_\alpha)$ is called the \emph{renormalized volume} of $Z_\alpha$. We also write 
$$
|\pi^{(\alpha)}| := \rk(W_\alpha).
$$ 
Furthermore, we have
\begin{align*}
Z_{\alpha\beta_i} =  \sum_{a \in \{j,k,l\}} \frac{Z_{\alpha\beta_i\beta_a}}{1-t_a} + W_{\alpha\beta_i}.
\end{align*}
Then $|Z_{\alpha\beta_i}| := \rk(W_{\alpha\beta_i})$ is called the renormalized volume of $Z_{\alpha\beta_i}$. Moreover
$$
|\mu_i^{(\alpha)}| := \rk(W_{\alpha\beta_i})
$$
coincides with the renormalized volume of plane partitions defined in \cite[Sect.~4.4]{MNOP}.
\end{definition}

\begin{definition} \label{def:vertex}
We define the \emph{vertex term}
\begin{align*}
\mathsf{V}_{\alpha} &:= \tr_{R \mathrm{Hom}_{U_\alpha}(I_\alpha,I_\alpha)_0[1]} + \sum_{i=1}^{4} \frac{\mathsf{A}_{\alpha\beta_i}}{1-t_i} + \sum_{1 \leq i<j \leq 4} \frac{\mathsf{A}_{\alpha\beta_i\beta_j}}{(1-t_i)(1-t_j)},
\end{align*}
where
\begin{align*}
\mathsf{A}_{\alpha\beta_i} &:= - Z_{\alpha\beta_i} + \frac{\overline{Z}_{\alpha\beta_i}}{t_{j} t_{k} t_{l}} - \frac{P_{jkl}}{t_{j}t_{k}t_{l}}  Z_{\alpha\beta_i}   \overline{Z}_{\alpha\beta_i}, \\
\mathsf{A}_{\alpha\beta_i\beta_{j}} &:= Z_{\alpha\beta_i\beta_j} + \frac{\overline{Z}_{\alpha\beta_i\beta_j}}{t_{k} t_{l}} - \frac{P_{kl}}{t_{k}t_{l}}  Z_{\alpha\beta_i\beta_j}   \overline{Z}_{\alpha\beta_i\beta_j},
\end{align*}
and $\tr_{R \mathrm{Hom}_{U_\alpha}(I_\alpha,I_\alpha)_0[1]}$ is given by \eqref{prevertex}.
Here $i,j,k,l \in \{1,2,3,4\}$ are mutually distinct.
As before, the expressions in this definition should be expanded in ascending powers of $t_i$ (though we will see shortly that $\mathsf{V}_\alpha$ is a Laurent polynomial in the $t_i$).
When we want to stress the dependence of these expressions on $Z_\alpha, Z_{\alpha\beta}, Z_{\alpha\beta\gamma}$, we write $\mathsf{V}_{Z_\alpha}, \mathsf{A}_{Z_{\alpha\beta}}, \mathsf{A}_{Z_{\alpha\beta\gamma}}$. Referring to their corresponding ideals, we also denote these expressions by $\mathsf{V}_{I_\alpha}, \mathsf{A}_{I_{\alpha\beta}}, \mathsf{A}_{I_{\alpha\beta\gamma}}$ respectively. 
\end{definition}

Recall that the normal bundle to $L_{\alpha\beta_i} \cong \PP^1$ is of the form (Section \ref{sec:toricgeomnot})
\begin{align*}
&N_{L_{\alpha\beta_i}/X} \cong \O(m_{\alpha\beta_i}) \oplus \O(m_{\alpha\beta_i}') \oplus \O(m_{\alpha\beta_i}''), \\
&m_{\alpha\beta_i} + m'_{\alpha\beta_i} + m''_{\alpha\beta_i} = -2.
\end{align*}
\begin{definition} \label{def:edge}
We define the \emph{edge term}
\begin{align*}
\mathsf{E}_{\alpha\beta_i} &:= \frac{t_i^{-1}  \mathsf{B}_{\alpha\beta_i}(t_j,t_k,t_l) }{1-t_i^{-1}} -  \frac{\mathsf{B}_{\alpha\beta_i}(t_i^{-m_{\alpha\beta_i}}t_j,t_i^{-m_{\alpha\beta_i}'} t_k,t_i^{-m_{\alpha\beta_i}''}t_l) }{1-t_i^{-1}}, 
\end{align*}
where 
\begin{align*}
\mathsf{B}_{\alpha\beta_i} &:= \mathsf{A}_{\alpha\beta_i}  + \sum_{a \in \{j,k,l\}} \frac{\mathsf{A}_{\alpha\beta_i\beta_a} }{1-t_a},
\end{align*}
where $i,j,k,l \in \{1,2,3,4\}$ are mutually distinct.
%with $j<k<l$.
The expressions in this definition should be expanded in ascending powers of $t_i^{-1}$ (!), $t_j, t_k, t_l$ (though we will see shortly that $\mathsf{E}_{\alpha\beta_i}$ is a Laurent polynomial in $t_1,t_2,t_3,t_4$).
When we want to stress the dependence of these expressions on $Z_\alpha, Z_{\alpha\beta}$, we write $\mathsf{E}_{Z_{\alpha\beta}}, \mathsf{B}_{Z_{\alpha\beta}}$. Referring to their corresponding ideals, we also denote these expressions by $\mathsf{E}_{I_{\alpha\beta}}, \mathsf{B}_{I_{\alpha\beta}}$ respectively.
\end{definition}

One can also consider the expressions of Definitions \ref{def:vertex}, \ref{def:edge} for $Z \subset X$ of dimension $\leq 1$. Then they reduce to the redistributions studied in \cite{CK2} (in turn based on \cite{MNOP}). 
\begin{definition} \label{def:face}
Fix $f \in F(X)$. Denote by $Z_f \subset X$ the unique $T_X$-fixed 2-dimensional closed subscheme set-theoretically supported on $S_f$ such that 
$$
Z_f \cap U_{\alpha_1\alpha_2\alpha_3} = Z \cap U_{\alpha_1\alpha_2\alpha_3}
$$
for any $\alpha_1, \alpha_2, \alpha_3 \in f$ distinct. Then $Z_f$ is determined by a single finite partition. Let $I_f := I_{Z_f} \subset \O_X$ be the corresponding ideal, then we define the \emph{face term} by
$$
\mathsf{F}_f := \tr_{R\Hom_X(I_{f},I_{f})_0[1]}.
$$
Note that, since $Z_f$ is proper, $\mathsf{F}_f$ is a Laurent polynomial in the $t_i$.
When we want to stress the dependence of $\mathsf{F}_f $ on $Z_f$ we write $\mathsf{F}_{Z_f}$. Referring to the corresponding ideal, we also denote this expression by $\mathsf{F}_{I_f}$.
\end{definition}

The main result of this section is the following:
\begin{theorem} \label{prop:redist}
Let $X$ be a toric Calabi-Yau 4-fold and $Z \subset X$ a 2-dimensional $T_X$-fixed closed subscheme with proper support. Then
\begin{equation} \label{eq:redisteq}
R \mathrm{Hom}(I_Z,I_Z)_0[1] = \sum_{\alpha \in V(X)} \mathsf{V}_{\alpha} + \sum_{\alpha\beta \in E(X)} \mathsf{E}_{\alpha\beta} + \sum_{f \in F(X)} \mathsf{F}_{f},
\end{equation}
where $\mathsf{V}_{\alpha}, \mathsf{E}_{\alpha\beta}, \mathsf{F}_{f}$ are Laurent polynomials.
\end{theorem}
\begin{proof}
We first show the equality \eqref{eq:redisteq} and then establish polynomiality of $\mathsf{V}_{\alpha}$, $\mathsf{E}_{\alpha\beta}$, $\mathsf{F}_{f}$. Note that
\begin{align*}
&\sum_{\alpha \in V(X)} \tr_{R \mathrm{Hom}_{U_\alpha}(I_\alpha,I_\alpha)_0[1]} - \sum_{\alpha\beta \in E(X)} \tr_{R \mathrm{Hom}_{U_{\alpha\beta}}(I_{\alpha\beta},I_{\alpha\beta})_0[1]} \\
&= \sum_{\alpha \in V(X)} \Bigg( -R\mathrm{Hom}_{U_\alpha}(I_\alpha,I_\alpha)_0 + \sum_{i=1}^{4} \frac{\mathsf{A}_{\alpha\beta_i}}{1-t_i} \Bigg) \\
&+ \sum_{\alpha\beta \in E(X)} \Bigg( \frac{t_1^{-1}  \mathsf{A}_{\alpha\beta}(t_2,t_3,t_4) }{1-t_1^{-1}} -  \frac{\mathsf{A}_{\alpha\beta}(t_1^{-m_{\alpha\beta}}t_2 ,t_1^{-m'_{\alpha\beta}}t_3 ,t_1^{-m''_{\alpha\beta}}t_4 ) }{1-t_1^{-1}} \Bigg),
\end{align*}
where the coordinates of summand $\alpha\beta$ are labelled such that $L_{\alpha\beta} \cap U_\alpha = Z(x_2,x_3,x_4)$ and $L_{\alpha\beta} \cap U_\beta = Z(x_2',x_3',x_4')$. 

By Remark \ref{rem:alphabetagamma}, in order to calculate $R \mathrm{Hom}(I_Z,I_Z)_0[1]$, we are left with the contributions of (1) $\alpha < \beta \in V(X)$ such that $\alpha\beta \notin E(X)$ and $\alpha, \beta \in f$ for some $f \in F(X)$; and (2) $\alpha_1 < \cdots <\alpha_N \in V(X)$, with $N \geq 3$, such that $\alpha_1, \ldots, \alpha_N \in f$ for some $f \in F(X)$. For any face $f \in F(X)$, we now label the vertices of $f$ as $\alpha_{1} < \cdots  < \alpha_{v_f}$, where $v_f = e(S_f)$ is the topological Euler characteristic of $S_f$. On $U_{\alpha_i}$, we use coordinates $(x_1^{(i)},x_2^{(i)},x_3^{(i)},x_4^{(i)})$ such that $(t_1^{(i)},t_2^{(i)},t_3^{(i)},t_4^{(i)}) \in (\mathbb{C}^*)^4 \cong T_X$ acts by the standard torus action. We label the coordinates such that for any  distinct $\alpha_i, \alpha_j,\alpha_k$, the first two coordinates of $U_{\alpha_i \alpha_j \alpha_k} \cong (\C^*)^2 \times \C^2$ are $x_1^{(i)}, x_2^{(i)}$, then
\begin{align*}
\tr_{R \mathrm{Hom}_{U_{\alpha_i \alpha_j \alpha_k}}(I_{\alpha_{i}\alpha_{j}\alpha_{k}},I_{\alpha_{i}\alpha_{j}\alpha_{k}})_0[1]} &= \delta(t_1^{(i)})\delta(t_2^{(i)}) \mathsf{A}_{f}(t_3^{(i)},t_4^{(i)}) \\
&= \delta(t_1^{(j)})\delta(t_2^{(j)}) \mathsf{A}_{f}(t_3^{(j)},t_4^{(j)})\\
&= \delta(t_1^{(k)})\delta(t_2^{(k)}) \mathsf{A}_{f}(t_3^{(k)},t_4^{(k)}),
\end{align*}
where 
\[
\mathsf{A}_f(t_3,t_4) := Z_{\lambda_f} + \frac{\overline{Z}_{\lambda_f}}{t_3t_4} - \frac{P_{34}}{t_3 t_4} Z_{\lambda_f} \overline{Z}_{\lambda_f} \in \mathbb{Z}[t_3^{\pm 1}, t_4^{\pm 1}]
\]
and $\lambda_f$ is the finite partition determining the scheme $Z_f$ in Definition \ref{def:face} and $Z_{\lambda_f} := \sum_{(i,j) \in \lambda_f} t_3^{i} t_4^{j}$. More generally, for any $\alpha_{i_1} < \cdots < \alpha_{i_l}$ with $l \geq 2$, where for $l=2$ we assume $\alpha_{i_1}, \alpha_{i_2}$ are not adjacent, we have
$$
\tr_{R \mathrm{Hom}_{U_{\alpha_{i_1} \ldots \alpha_{i_l}}} (I_{\alpha_{i_1} \ldots \alpha_{i_l}},I_{\alpha_{i_1} \ldots \alpha_{i_l}})_0[1]} = \delta(t_1^{(1)}) \delta(t_2^{(1)}) \mathsf{A}_f(t_3^{(1)}, t_4^{(1)}).
$$
We define $N_f := -(-v_f + \sum_{i=2}^{v_f} (-1)^{i} \binom{v_f}{i}) =1$ and 
\begin{align}
\begin{split} \label{def:Fprime}
&\mathsf{F}_f'  := N_f \cdot \delta(t_1^{(1)}) \delta(t_2^{(1)}) \mathsf{A}_f(t_3^{(1)}, t_4^{(1)}) - \sum_{i} \frac{\mathsf{A}_f(t_3^{(i)}, t_4^{(i)})}{(1-t_1^{(i)})(1-t_2^{(i)})} \\
&- \sum_{i} \Bigg( \frac{(t_1^{(i)})^{-1} \mathsf{A}_{f}(t_3^{(i)},t_4^{(i)})}{(1 - (t_1^{(i)})^{-1})(1-t_2^{(i)})} - \frac{\mathsf{A}_{f} ( (t_1^{(i)})^{-m_{\alpha_i \alpha_{i+1}}'} t_3^{(i)}, (t_1^{(i)})^{-m_{\alpha_i \alpha_{i+1}}''} t^{(i)}_4) }{(1 - (t_1^{(i)})^{-1})(1- (t_1^{(i)})^{-m_{\alpha_i\alpha_{i+1}}}  t_2^{(i)} )}\Bigg).
\end{split}
\end{align}
Then, by Definitions \ref{def:vertex} and \ref{def:edge}, we have
\begin{align*}
R \mathrm{Hom}(I_Z,I_Z)_0[1] = &\sum_{\alpha \in V(X)} \Bigg( -R\mathrm{Hom}_{U_\alpha}(I_\alpha,I_\alpha)_0 + \sum_{i=1}^{4} \frac{\mathsf{A}_{\alpha\beta_i}}{1-t_i} \Bigg) \\
&+ \sum_{\alpha\beta \in E(X)} \Bigg( \frac{t_1^{-1}  \mathsf{A}_{\alpha\beta}(t_2,t_3,t_4) }{1-t_1^{-1}} -  \frac{\mathsf{A}_{\alpha\beta}(t_1^{-m_{\alpha\beta}}t_2 ,t_1^{-m'_{\alpha\beta}}t_3 ,t_1^{-m''_{\alpha\beta}}t_4 ) }{1-t_1^{-1}} \Bigg) \\
&+ \sum_{f \in F(X)} N_f \cdot \delta(t_1^{(1)}) \delta(t_2^{(1)}) \mathsf{A}_f(t_3^{(1)}, t_4^{(1)}) \\
= &\sum_{\alpha \in V(X)} \mathsf{V}_{\alpha} + \sum_{\alpha\beta \in E(X)} \mathsf{E}_{\alpha\beta} + \sum_{f \in F(X)} \mathsf{F}_{f}'.
\end{align*}
In order to complete the proof, we need to show $\mathsf{F}_f' = \mathsf{F}_f$ (Definition \ref{def:face}), and $\mathsf{V}_{\alpha}, \mathsf{E}_{\alpha\beta}, \mathsf{F}_{f}$ are Laurent polynomials.

Consider a vertex $\alpha \in V(X)$ with neighbouring vertices $\beta_1,\beta_2,\beta_3, \beta_4 \in V(X)$ (possibly fewer if $\alpha$ lies on a non-compact edge). Recall that $Z_{\alpha\beta_i\beta_j}$ are Laurent polynomials. Consider the Laurent \emph{polynomials} $W_{\alpha\beta_i}$, $W_\alpha$ of Definition \ref{def:renorm}, which satisfy
\begin{align}
\begin{split} \label{eq:Ztosubst}
Z_\alpha &=  \sum_{1 \leq i<j \leq 4} \frac{Z_{\alpha\beta_i\beta_j}}{(1-t_i)(1-t_j)} + \sum_{i=1}^{4} \frac{W_{\alpha\beta_i}}{1-t_i} + W_\alpha, \\
Z_{\alpha\beta_i} &=  \sum_{a \in \{j,k,l\}} \frac{Z_{\alpha\beta_i\beta_a}}{1-t_a} + W_{\alpha\beta_i}.
\end{split}
\end{align}
Substituting into $\mathsf{V}_{\alpha}$ gives an expression in terms of $W_{\alpha\beta_i}$, $W_\alpha$ (and their duals). An explicit (lengthy) calculation shows that all denominators cancel and $\mathsf{V}_{\alpha}$ is a Laurent polynomial. 
Similarly, substituting \eqref{eq:Ztosubst} into $\mathsf{B}_{\alpha\beta_i}$, all denominators cancel and we find that $\mathsf{B}_{\alpha\beta_i}$ is a Laurent polynomial. Since the numerator and denominator of
\begin{align*}
\mathsf{E}_{\alpha\beta_i} = \frac{t_i^{-1}  \mathsf{B}_{\alpha\beta_i}(t_j,t_k,t_l) }{1-t_i^{-1}} -  \frac{\mathsf{B}_{\alpha\beta_i}(t_j t_i^{-m_{\alpha\beta_i}},t_k t_i^{-m_{\alpha\beta_i}'},t_l t_i^{-m_{\alpha\beta_i}''}) }{1-t_i^{-1}}
\end{align*}
are zero for $t_i^{-1} = 1$, we deduce that $\mathsf{E}_{\alpha\beta_i} $ is a Laurent polynomial in $t_1, t_2, t_3, t_4$ as well. 

Recall the definition of $Z_f \subset X$ of Definition \ref{def:face}, and let $Z_{f,\alpha} = Z_f \cap U_\alpha$ and $Z_{f,\alpha\beta} = Z_{f} \cap U_{\alpha\beta}$. Then
\begin{align*}
\mathsf{F}_{Z_f} := R\Hom_X(I_{Z_f},I_{Z_f})_0[1] = \sum_{\alpha \in V(X)} \mathsf{V}_{Z_{f,\alpha}} + \sum_{\alpha\beta \in E(X)} \mathsf{E}_{Z_{f,\alpha\beta}} + \mathsf{F}_{Z_f}'.
\end{align*}
It is easy to calculate $\mathsf{V}_{Z_{f,\alpha}} = \mathsf{E}_{Z_{f,\alpha\beta}} = 0$ for all $\alpha \in V(X)$ and $\alpha\beta \in E(X)$. Hence $\mathsf{F}_{Z_f}' = \mathsf{F}_{Z_f} = R\Hom_X(I_{Z_f},I_{Z_f})_0[1]$ and it is a Laurent polynomial because $Z_f \subset X$ is proper.
\end{proof}

\begin{remark} \label{rem:expFf}
The proof of Theorem \ref{prop:redist} leads to a nice expression for $\mathsf{F}_f$. For a $T_X$-fixed 2-dimensional subscheme $Z \subset X$ with proper support, denote by $\lambda_f$ the finite partition corresponding to the subscheme $Z_f$ from Definition \ref{def:face}. Let
\[
\mathsf{A}_f(t_3,t_4) := Z_{\lambda_f} + \frac{\overline{Z}_{\lambda_f}}{t_3t_4} - \frac{P_{34}}{t_3 t_4} Z_{\lambda_f} \overline{Z}_{\lambda_f} \in \mathbb{Z}[t_3^{\pm 1}, t_4^{\pm 1}],
\]
where $Z_{\lambda_f} := \sum_{(i,j) \in \lambda_f} t_3^{i} t_4^{j}$. Denote by $\alpha_1< \cdots < \alpha_{e(S_f)}$ the vertices of $f$, where $e(S_f)$ denotes the topological Euler characteristic of $S_f$. As before, for any $\alpha_i \in f$, we use coordinates $(x_1^{(i)}, x_2^{(i)}, x_3^{(i)}, x_4^{(i)})$ on $U_{\alpha_i}$ such that $(t_1^{(i)},t_2^{(i)},t_3^{(i)},t_4^{(i)}) \in (\mathbb{C}^*)^4 \cong T_X$ acts by the standard torus action. Then 
\begin{equation} \label{eqn:exprFf}
\mathsf{F}_f = \sum_{i=1}^{e(S_f)} \frac{\mathsf{A}_f(t_3^{(i)},t_4^{(i)})}{(1-t_1^{(i)})(1-t_2^{(i)})}.
\end{equation}
This equality should be interpreted carefully and can be proved as follows.\footnote{We thank Schmiermann for pointing out the correct interpretation of \eqref{eqn:exprFf}.} Firstly, writing $t_1 := t_1^{(1)}$, $t_2 := t_2^{(1)}$, we view each $t_1^{(i)}$, $t_2^{(i)}$ as a function of $t_1,t_2$ determined by the coordinate change between $U_{\alpha_{1}}$ and $U_{\alpha_i}$ \eqref{eqn:torictransf}. Secondly, the left hand side of \eqref{eqn:exprFf} is the expression in \eqref{def:Fprime}, which lives in $\mathbb{Q}[t_1^{\pm 1}, \ldots, t_4^{\pm 1}] \subset \mathbb{Q}(\!(t_1,t_1^{-1}, \ldots ,t_4, t_4^{-1})\!)$, whereas the right hand side of \eqref{eqn:exprFf} (a priori) is an element of $\mathbb{Q}(t_1,\ldots,t_4)$. We consider the field extensions
\[
\mathbb{Q}(t_1, \ldots, t_4) \hookrightarrow \mathbb{Q}(t_1,t_2,t_3)(\!(t_4)\!) \hookrightarrow \cdots \hookrightarrow \mathbb{Q}(\!(t_1)\!)\cdots(\!(t_4)\!),
\]
obtained by expanding denominators in ascending powers of the $t_i$ (e.g.~for the first extension we expand $1/(1-f(t_4)) = 1 + f(t_4) + f(t_4)^2 + \ldots$ for any $f(t_4) \in t_4 \Q(t_1,t_2,t_3)[t_4]$).
Note that $\mathbb{Q}(\!(t_1)\!)\cdots(\!(t_4)\!) \subset \mathbb{Q}(\!(t_1,t_1^{-1}, \ldots, t_4,t_4^{-1})\!)$. This embeds the right hand side of \eqref{eqn:exprFf} into $\mathbb{Q}(\!(t_1,t_1^{-1}, \ldots ,t_4, t_4^{-1})\!)$ and after multiplying left and right hand side of \eqref{eqn:exprFf} by $\prod_i(1-t_1^{(i)})(1-t_2^{(i)})$, and using $\delta(t_a^{(i)})(1-t_a^{(i)}) = 0$, they are obviously equal.\footnote{Although $\mathbb{Q}(\!(t_1,t_1^{-1}, \ldots ,t_4, t_4^{-1})\!)$ is not a ring, it is a $\mathbb{Q}[t_1^{\pm 1}, \ldots, t_4^{\pm 1}]$-module.} This shows that the right hand side of \eqref{eqn:exprFf} is in fact an element $\mathbb{Q}[t_1^{\pm 1}, \ldots, t_4^{\pm 1}] \subset \mathbb{Q}(t_1,\ldots,t_4)$ and it is equal to the left hand side of \eqref{eqn:exprFf}.
\end{remark}

We use the following notation, borrowed from Nekrasov-Piazzalunga \cite{NP1}, for the symmetrized $T$-equivariant $K$-theoretic Euler class. As before, $T \leq T_X$ denotes the Calabi-Yau torus \eqref{eqn:CYtorus}.
\begin{definition} 
Consider the action of $T \times \C^*$ on $\pt$, where $T = Z(t_1t_2t_3t_4- 1) \leq T_X$ as before, and we denote the equivariant parameter of the second torus by $y$. We define
$$
[E] := \widehat{\Lambda}\udot E^\vee,
$$
for any class $E \in K_0^{T \times \C^*}(\mathrm{pt},\Z[\tfrac{1}{2}]) \cong \Z[t_1^{\pm \frac{1}{2}},t_2^{\pm \frac{1}{2}},t_3^{\pm \frac{1}{2}},t_4^{\pm \frac{1}{2}},y^{\pm \frac{1}{2}}] / (t_1t_2t_3t_4-1)$. In particular
$$
[L] = L^{\frac{1}{2}} - L^{-\frac{1}{2}},
$$
for the class of a $T \times \C^*$-equivariant line bundle $L$.
\end{definition}

Next, we incorporate the $K$-theoretic insertion (Sections \ref{sec:Ktheorinv}, \ref{sec:vardef})
$$
\widehat{\Lambda}\udot (R\pi_{\curly P*} \mathcal{O}_{\cZ} \otimes y^{-1})
$$
into the vertex formalism. At a fixed point $Z \in \curP_v^{(-1)}(X)^T$, we have
$$
\widehat{\Lambda}\udot ( R\pi_{\curly P *} \mathcal{O}_{\cZ} \otimes y^{-1})|_Z = \widehat{\Lambda}\udot (R\Gamma(X,\O_Z) \otimes y^{-1}).
$$
Using \v{C}ech cohomology, we deduce
\begin{align*}
R\Gamma(X,\O_Z) = &\sum_{\alpha \in V(X)} \Gamma(U_\alpha,\O_{Z_\alpha}) - \sum_{\alpha < \beta \in V(X)} \Gamma(U_{\alpha\beta}, \O_{Z_{\alpha\beta}}) \\
&+ \sum_{\alpha<\beta<\gamma \in V(X)} \Gamma(U_{\alpha\beta\gamma},\O_{Z_{\alpha\beta\gamma}}) - \cdots.
\end{align*}

These terms can be redistributed similar to Definitions \ref{def:vertex}, \ref{def:edge}, \ref{def:face}
\begin{align*}
\widetilde{\mathsf{V}}_{\alpha} := &\mathsf{V}_{\alpha} - y \overline{Z}_\alpha + \sum_{i=1}^{4} \frac{y \overline{Z}_{\alpha\beta_i}}{1-t_i^{-1}} - \sum_{1 \leq i<j \leq 4} \frac{y \overline{Z}_{\alpha\beta_i\beta_j}}{(1-t_i^{-1})(1-t_j^{-1})} \\
&\quad \ - y^{-1} Z_\alpha + \sum_{i=1}^{4} \frac{y^{-1} Z_{\alpha\beta_i}}{1-t_i} - \sum_{1 \leq i<j \leq 4} \frac{y^{-1} Z_{\alpha\beta_i\beta_j}}{(1-t_i)(1-t_j)} \\
\widetilde{\mathsf{E}}_{\alpha\beta_i} := &\mathsf{E}_{\alpha\beta_i} + \frac{t_i}{1-t_i} \Bigg\{ y \overline{Z}_{\alpha\beta_i} - \sum_{a \in \{j,k,l\}} \frac{y \overline{Z}_{\alpha\beta_i\beta_a}}{1-t_a^{-1}} \Bigg\} \\
&\quad\quad -\frac{1}{1-t_i} \Bigg\{ y \overline{Z}_{\alpha\beta_i} - \sum_{a \in \{j,k,l\}} \frac{y \overline{Z}_{\alpha\beta_i\beta_a}}{1-t_a^{-1}} \Bigg\}\Bigg|_{(t_j t_i^{-m_{\alpha\beta_i}},t_k t_i^{-m_{\alpha\beta_i}'},t_l t_i^{-m_{\alpha\beta_i}''})} \\
&\quad\quad + \frac{t_i^{-1}}{1-t_i^{-1}} \Bigg\{ y^{-1} Z_{\alpha\beta_i} - \sum_{a \in \{j,k,l\}} \frac{y^{-1} Z_{\alpha\beta_i\beta_a}}{1-t_a} \Bigg\} \\
&\quad\quad- \frac{1}{1-t_i^{-1}} \Bigg\{ y^{-1} Z_{\alpha\beta_i} - \sum_{a \in \{j,k,l\}} \frac{y^{-1} Z_{\alpha\beta_i\beta_a}}{1-t_a} \Bigg\}\Bigg|_{(t_j t_i^{-m_{\alpha\beta_i}},t_k t_i^{-m_{\alpha\beta_i}'},t_l t_i^{-m_{\alpha\beta_i}''})} \\
\widetilde{\mathsf{F}}_f  := &\mathsf{F}_f - y R\Gamma(X, \O_{Z_f})^\vee - y^{-1} R\Gamma(X, \O_{Z_f}).
\end{align*}

\begin{remark}
Using the notation of Remark \ref{rem:expFf}, we have
$$
\widetilde{\mathsf{F}}_f = \mathsf{F}_f - \sum_{i=1}^{e(S_f)} \frac{y \overline{Z}_{\lambda_f}(t_3^{(i)},t_4^{(i)})}{(1-(t_1^{(i)})^{-1})(1-(t_2^{(i)})^{-1})} - \sum_{i=1}^{e(S_f)} \frac{y^{-1} Z_{\lambda_f}(t_3^{(i)},t_4^{(i)})}{(1-t_1^{(i)})(1-t_2^{(i)})},
$$
where the two sums are elements of $\Q[t_1^{\pm 1}, \ldots, t_4^{\pm 1}, y^{\pm 1}] \subset \Q(t_1,\cdots,t_4)[y^{\pm 1}]$.
\end{remark}

Recall that, in the toric case, $\vd \in \Z$ (Lemma \ref{lem:intvd}) and we define invariants by the virtual localization formula (Definition \ref{def:equivKtheorinv}). Therefore, combining Proposition \ref{prop:OTisored} and Theorems \ref{prop:redist}, \ref{prop:fixlocmain1}(i), we obtain:
\begin{proposition} \label{prop:DTsumsqrt}
  Let $X$ be a toric Calabi-Yau 4-fold and fix $v = (0,0,\gamma,\beta,n - \gamma \cdot \td_2(X)) \in H_c^*(X,\Q)$. Let $\curly P := \curly P_v^{(-1)}(X)^T$ and denote by $\vd$ the virtual dimension and $\rk$ the rank of the tautological complex. Then 
\[
\langle\!\langle \O_X \rangle\!\rangle_{X,v}^{\DT} = \sum_{Z \in \curly P^T} \pm \sqrt{(-1)^{\vd - \rk}\Big[-\sum_{\alpha \in V(X)} \widetilde{\mathsf{V}}_{Z_\alpha} - \sum_{\alpha\beta \in E(X)} \widetilde{\mathsf{E}}_{Z_{\alpha\beta}} - \sum_{f \in F(X)} \widetilde{\mathsf{F}}_{Z_f} \Big]},
\]
where the right hand side depends on a choice of $\pm \sqrt{\cdot}$ for each fixed point.
\end{proposition}

\begin{remark} \label{rem:otherL}
For notational convenience, we considered the insertion 
$$
\widehat{\Lambda}\udot ( R\pi_{\curly P*} \O_{\cZ} \otimes y^{-1}).
$$ 
The insertion $\widehat{\Lambda}\udot ( R\pi_{\curly P *} (\O_{\cZ} \otimes L)  \otimes y^{-1})$, for an arbitrary $T_X$-equivariant line bundle $L$ on $X$, can be easily incorporated as follows. Define $L_\alpha = L|_{U_\alpha}$, $L_{\alpha\beta} = L|_{U_{\alpha\beta}}$ etc. Then, in the definition of $\mathsf{V}_\alpha$, we replace $y$ by $L_\alpha^* y$ and we denote the resulting expression by $\mathsf{V}_\alpha(L_\alpha)$ (and similarly for the edge and face). 
\end{remark}

We will now show that the vertex, edge, and face terms in Proposition \ref{prop:DTsumsqrt} contain no positive $T$-fixed terms, which justifies writing\footnote{It follows from \eqref{eqn:rkVtilde} in Section \ref{sec:coholimittaut} that $\rk(\widetilde{\mathsf{V}}_{Z_\alpha})$ is even. Similarly, it can be seen that $\rk(\widetilde{\mathsf{E}}_{Z_{\alpha\beta}})$ is even. It is also clear from Definition \ref{def:face} and Lemma \ref{lem:intvd} that $\rk(\widetilde{\mathsf{F}}_{Z_f})$ is even.} 
\begin{align*}
\langle\!\langle \O_X \rangle\!\rangle_{X,v}^{\DT} = \sum_{Z \in \curly P^T} \pm &\prod_{\alpha \in V(X)} \sqrt{(-1)^{\frac{\rk(\widetilde{\mathsf{V}}_{Z_\alpha})}{2}} [-\widetilde{\mathsf{V}}_{Z_\alpha}]} \cdot \prod_{\alpha\beta \in E(X)} \sqrt{(-1)^{\frac{\rk \widetilde{\mathsf{E}}_{Z_{\alpha\beta}}}{2}}  [-\widetilde{\mathsf{E}}_{Z_{\alpha\beta}}]} \cdot \\
&\prod_{f \in F(X)} \sqrt{(-1)^{\frac{\rk \widetilde{\mathsf{F}}_{Z_f}}{2}} [-\widetilde{\mathsf{F}}_{Z_f}]}.
\end{align*}

\begin{lemma} \label{lem:noposTfix}
$\mathsf{V}_\alpha$, $\mathsf{E}_{\alpha\beta}$, $\mathsf{F}_f$ have no positive $T$-fixed terms.
\end{lemma}
\begin{proof}
For $\mathsf{E}_{\alpha\beta}$, $\mathsf{F}_f$ we in fact show that they have no $T$-fixed terms at all. First note that $\mathsf{A}_f(t_3^{(i)}, t_4^{(i)})$ in \eqref{eqn:exprFf} has no terms of the form $(t_3^{(i)} t_4^{(i)})^n$ because tangent spaces to Hilbert schemes of points of $\C^2$ at monomial ideals have no such terms (by a simple version of the argument with Haiman arrows in the proof of Theorem \ref{prop:fixlocmain1}(i)). From this it follows that $\mathsf{F}_f$ has no $T$-fixed terms.

Next, we note that $\mathsf{B}_{\alpha\beta_1}$ in Definition \ref{def:edge}, which only depends on variables $t_2,t_3,t_4$, is precisely \emph{minus} the vertex term of Donaldson-Thomas theory on toric 3-folds \cite[Sect.~4.9]{MNOP}. This is known to have no terms of the form $(t_2t_3t_4)^n$ \cite[p.~1279]{MNOP}. Therefore $\mathsf{E}_{\alpha\beta_1}$ has no $T$-fixed terms. 

Finally, we turn to the vertex term. We fix a $(\C^*)^4$-invariant 2-dimensional subscheme $Z_{\alpha_0} \subset U_{\alpha_0} = \C^4$ and we want to show that $\mathsf{V}_{Z_{\alpha_0}}$ has no positive $T$-fixed term (we use a distinguished label $\alpha_0$). We pick a toric Calabi-Yau 4-fold $X$ which is a partial toric compactification of $U_{\alpha_0}$ in the following sense: $\alpha_0$ has four neighbouring vertices $\beta_1, \ldots, \beta_4$ and each $\alpha_0,\beta_i,\beta_j$, with $i \neq j$, are vertices of a \emph{proper} smooth toric surface in $X$.\footnote{In fact, it suffices to add ten copies of $\C^4$ in such a way that $X$ contains precisely six $\PP^2$'s which compactify the six coordinate 2-planes in $U_{\alpha_0} = \C^4$.} We consider the \emph{minimal} proper $T_X$-fixed 2-dimensional subscheme $Z \subset X$ such that $Z \cap U_{\alpha_0} = Z_{\alpha_0}$. It is then easy to see that for each $\alpha \neq \alpha_0$, we have $\mathsf{V}_{Z_\alpha} = 0$. More precisely, in Definitions \ref{def:renorm}, \ref{def:vertex}: putting $Z_{\alpha\beta_i\beta_j} = 0$ for $(i,j) \neq (1,4), (2,4), (3,4)$, $W_{\alpha\beta_i} = 0$ for $i \neq 4$, and $W_\alpha = 0$ gives $\mathsf{V}_{Z_\alpha} = 0$ by explicit calculation. Therefore
\begin{equation} \label{eqn:completionarg}
R\Hom_X(I_Z,I_Z)_0[1] = \mathsf{V}_{\alpha_0} + \sum_{\alpha\beta \in E(X)} \mathsf{E}_{\alpha\beta} + \sum_{f \in F(X)} \mathsf{F}_f.
\end{equation}
The left hand side of \eqref{eqn:completionarg} has no positive $T$-fixed terms since $\Ext_X^3(I_Z,I_Z)_0^* \cong \Ext_X^1(I_Z,I_Z)_0$ and $\Ext_X^1(I_Z,I_Z)_0 \cong \Hom_X(I_Z,\O_Z)$ (\cite[Thm.~3.19]{BKP}), which has no $T$-fixed term by Theorem \ref{prop:fixlocmain1}(i). Since the edge and face terms have no $T$-fixed terms, the result follows.
\end{proof}

It is interesting to provide an explicit ``half in $K_0^T(\pt)$'' of the vertex, edge, and face terms and hence the term under the square root symbol in Proposition \ref{prop:DTsumsqrt}. Using the spectral construction, we give such expressions when $X = \mathrm{Tot}(K_Y)$ (Section \ref{sec:local3fold}). In general, such an expression is exhibited in \cite{NP2}. However, it should be noted that even after such a ``halving'', one can still make any choice of signs assigned to the $T$-fixed points (Definition \ref{def:equivKtheorinv}). As we will see, in order to get interesting answers, one typically has to assign non-trivial signs even after ``halving''.

\begin{definition} \label{def:DTvertex}
Let $\boldsymbol{\mu} = \{\mu_{a} \}_{a = 1}^{4}$ be a collection of plane partitions corresponding to subschemes of dimension $\leq 1$ of $\C^3$ (and not all finite). We denote by $\mathcal{P}_{\boldsymbol{\mu}}$ the collection of solid partitions $\pi$ such that $Z_\pi$ has the property that its maximal subscheme without embedded points is $Z_{\boldsymbol{\mu}}$ (Definition \ref{def:partitionnotation}). We define the $(K$-theoretic) \emph{$\DT$ topological vertex} by
$$
\mathsf{V}^{\DT}_{\mu_1\mu_2\mu_3\mu_4}(q) := \sum_{\pi \in \mathcal{P}_{\boldsymbol{\mu}}} \pm  \sqrt{(-1)^{\frac{\rk(\widetilde{\mathsf{V}}_{Z_\pi})}{2}} [-\widetilde{\mathsf{V}}_{Z_\pi}]} \, q^{|\pi|},
$$
where $|\pi|$ denotes renormalized volume (Definition \ref{def:renorm}) and $\pm\sqrt{\cdot}$ is a choice of square root with sign for each $\pi \in \mathcal{P}_{\boldsymbol{\mu}}$. In particular, $\mathsf{V}^{\DT}_{\mu_1\mu_2\mu_3\mu_4}(q)$ depends on a choice of sign attached to each $\pi \in \mathcal{P}_{\boldsymbol{\mu}}$. 
\end{definition}

The previous definition extends the $\DT$ topological vertex for curves studied in \cite{CK2, CKM1, Mon} and for points in \cite{NP1}. It also matches the DT topological vertex in \cite{NP2}

\begin{remark} \label{rem:modPT0}
Let $I^\mdot \in \curP_v^{(0)}(X)^{T_X}$. Its scheme theoretic support $Z \subset X$ is a $T_X$-fixed proper 2-dimensional subscheme with no embedded 0-dimensional components.
Then $\mathsf{V}_\alpha$, $\mathsf{E}_{\alpha\beta}$, $\mathsf{F}_f$, $\widetilde{\mathsf{V}}_\alpha$, $\widetilde{\mathsf{E}}_{\alpha\beta}$, $\widetilde{\mathsf{F}}_f$ can be determined in exactly the same way as in the $\DT$ case with the following replacement. Instead of \eqref{Zalpha} we now take
$$
Z_\alpha := \tr_{\O_{Z_\alpha}} + \tr_{B_\alpha}
$$
where $\tr_{B_\alpha}$ is the $T$-representation of the $\PT_0$ box configuration corresponding to $I^\mdot|_{U_\alpha}$ (Definition \ref{def:PT0boxconfig}). Since the local description of $I^\mdot|_{U_{\alpha_1\cdots \alpha_l}}$ for $l \geq 2$ is identical to the $\DT$ case, no modifications are needed for the edge and face terms. 
Under the assumption $\curP_v^{(0)}(X)^{T_X} = \curP_v^{(0)}(X)^{T}$ is 0-dimensional reduced (e.g.~in the cases of Theorem \ref{prop:fixlocmain1}(ii),(iii)), we conclude that Proposition \ref{prop:DTsumsqrt} holds in the $\PT_0$ case. 

Let $\boldsymbol{\mu}$ be a choice of plane partitions with no $\PT_0$ moduli (Definition \ref{def:PT0boxconfig}). Recall that $\mathcal{B}_{\boldsymbol{\mu}}^{(0)}$ denotes the collection of $\PT_0$ box configurations with asymptotics $\boldsymbol{\mu}$ (Definition \ref{def:PT0boxconfig}). We denote by $\pi^{\boldsymbol{\mu} }$ the solid partition corresponding to $Z_{\boldsymbol{\mu}}$, i.e., the 2-dimensional subscheme without embedded 0-dimensional components determined by $\boldsymbol{\mu}$ (Definition \ref{def:partitionnotation}). Each element $B \in \mathcal{B}^{(0)}_{\boldsymbol{\mu}}$ then determines a $\PT_0$ pair $I_B\udot$ on $\C^4$. We define the ($K$-theoretic) \emph{$\PT_0$ topological vertex} by
$$
\mathsf{V}^{\PT_0}_{\mu_1\mu_2\mu_3\mu_4}(q) := \sum_{B \in \mathcal{B}^{(0)}_{\boldsymbol{\mu}}} \pm  \sqrt{(-1)^{\frac{\rk(\widetilde{\mathsf{V}}_{I_B\udot})}{2}} [-\widetilde{\mathsf{V}}_{I_B\udot}]} \, q^{|B|+|\pi^{\boldsymbol{\mu} }|}.
$$
This definition is subject to the following restriction. We must assume that for all $B \in \mathcal{B}^{(0)}_{\boldsymbol{\mu}}$, the vertex term $\mathsf{V}_{I_B\udot}$ has no positive $T$-fixed part (so we do not divide by zero). We expect that this is always the case when $\boldsymbol{\mu}$ is a choice of plane partitions with no $\PT_0$ moduli. It is certainly the case when $\boldsymbol{\mu}$ is chosen such that $\PT_0 = \PT_1$; or such that there are embedded curves in at most 2 directions and the underlying pure 2-dimensional subscheme is Cohen-Macaulay. Then the proof of Lemma \ref{lem:noposTfix} goes through.\footnote{The only modification is that this time the analog of the left hand side of \eqref{eqn:completionarg} has no positive $T$-fixed part by the proof of Theorem \ref{prop:fixlocmain1}(ii), (iii).} 
\end{remark}

The previous definition extends the $\PT$ topological vertex for curves studied in \cite{CK2, CKM1, Liu, Pia}. 

\begin{remark} \label{rem:modPT1}
Let $I^\mdot \in \curP_v^{(1)}(X)^{T_X}$. Its scheme theoretic support $Z \subset X$ is a $T_X$-fixed proper pure 2-dimensional subscheme.
Then $\mathsf{V}_\alpha$, $\mathsf{E}_{\alpha\beta}$, $\mathsf{F}_f$, $\widetilde{\mathsf{V}}_\alpha$, $\widetilde{\mathsf{E}}_{\alpha\beta}$, $\widetilde{\mathsf{F}}_f$ can be determined in exactly the same way as in the $\DT$ case with the following replacements. Instead of \eqref{Zalpha} we now take
$$
Z_\alpha := \tr_{\O_{Z_\alpha}} + \tr_{B_\alpha}
$$
where $\tr_{B_\alpha}$ is the $T$-representation of the $\PT_1$ box configuration corresponding to $I^\mdot|_{U_\alpha}$ (Definition \ref{def:PT1boxconfig}). In contrast to the $\PT_0$ case, the $T$-representation $\tr_{B_\alpha}$ can now be infinite. Moreover, the edge terms need to be modified as well. Let $\alpha\beta \in E(X)$. Suppose we label our coordinates such that $U_{\alpha\beta} = \C^* \times \C^3$. From Proposition \ref{lem:PT1lim}, it follows that there exists a $(\C^*)^3$-equivariant $\PT$ pair $J^\mdot$ on $\C^3$ such that
$$
I^\mdot|_{U_{\alpha\beta}} \cong J^\mdot [x_1^{-1},x_1].
$$
Therefore it determines a box configuration $B_{\alpha\beta}$ as defined in \cite[Sect.~2.5]{PT2}. Then instead of $Z_{\alpha\beta} = \tr_{\O_{Z_{\alpha\beta}}}$, we now take   
$$
Z_{\alpha\beta} := \tr_{\O_{Z_{\alpha\beta}}} + \tr_{B_{\alpha\beta}}.
$$
Recall that a box configuration $B_{\alpha\beta}$ may also have moduli (labelled by $\PP^1$'s, \cite[Sect.~2.5]{PT2}). 
The face term needs no modification. With the above replacements in mind, we have now defined $\mathsf{V}_\alpha$, $\mathsf{E}_{\alpha\beta}$, $\mathsf{F}_f$, $\widetilde{\mathsf{V}}_\alpha$, $\widetilde{\mathsf{E}}_{\alpha\beta}$, $\widetilde{\mathsf{F}}_f$ and they are again Laurent polynomials.\footnote{It is also instructive to note that $\mathsf{B}_{\alpha\beta}$ in Definition \ref{def:edge} is identical to minus the PT vertex for toric 3-folds studied in \cite{PT2}.} 
Under the assumption $\curP_v^{(1)}(X)^{T_X} = \curP_v^{(1)}(X)^{T}$ is 0-dimensional reduced, we conclude that Proposition \ref{prop:DTsumsqrt} holds in the $\PT_1$ case.
\end{remark}

\subsection{Vertex \texorpdfstring{$\DT$--$\PT_0$}{dtpt0} correspondence}

The easiest $\DT$ topological vertex is for $\boldsymbol{\mu}$ all empty. In this case, the $\DT$ topological vertex is a sum over all finite solid partitions viewed as the fixed points of the Hilbert schemes
$$
\bigcup_{n \geq 0} \Hilb^n(\C^4)^{T_X} = \bigcup_{n \geq 0} \Hilb^n(\C^4)^{T}.
$$
In \cite{KR}, it is shown that $\Hilb^n(\C^4)$ is \emph{globally} the the zero locus of an isotropic section of an orthogonal bundle on the smooth non-commutative Hilbert scheme $\mathrm{ncHilb}^n(\C^4)$. Moreover, from this description, one can construct an orientation on $\Hilb^n(\C^4)$ such that the following holds.
\begin{theorem} \cite{Nek, NP1, KR} \label{thm:KR}
There exist orientations on $\Hilb^n(\C^4)$ for which the corresponding signs, induced by the Oh-Thomas localization formula, give
\begin{align*}
\sum_{n=0}^{\infty} \langle \! \langle \O_{\mathbb{C}^4} \rangle \! \rangle_{\C^4,0,0,n}^{\DT} q^n = &\mathsf{V}^{\DT}_{\varnothing\varnothing\varnothing\varnothing}(q) = \Exp\Bigg( \frac{[t_1t_2][t_1t_3][t_2t_3][y]}{[t_1][t_2][t_3][t_4][y^{\frac{1}{2}}q][y^{\frac{1}{2}}q^{-1}]}\Bigg). 
\end{align*}
\end{theorem}
%Warning. The second equality is the same in \cite{KR}, \cite{CKM1} and this paper. The first equality is the same in \cite{CKM1} and this paper, but in \cite{KR} we have
%\begin{align*}
%\sum_{n=0}^{\infty} \chi(\Hilb^n(\C^4), \widehat{\O}^{\vir} \otimes \widehat{\Lambda}_{-1}(\O^{[n]} \otimes y^{-1})^\vee ) q^n = &\mathsf{V}^{\DT}_{\varnothing\varnothing\varnothing\varnothing\varnothing\varnothing\varnothing\varnothing\varnothing\varnothing}(t,y,q) = \Exp\Bigg( \frac{[t_1t_2][t_1t_3][t_2t_3][y]}{[t_1][t_2][t_3][t_4][y^{\frac{1}{2}}q][y^{\frac{1}{2}}q^{-1}]}\Bigg),
%\end{align*}
%where the orientation on each $\Hilb^n(\C^4)$ is induced by a specific choice of $\Lambda$ (as given in the paper). Switching the orientation in \cite{KR} to the one induced by $\Lambda^*$, we get the formula in the above theorem.

The second equality was conjectured in physics by Nekrasov \cite{Nek} and conjecturally extended to punctual Quot schemes\footnote{In physics parlance: ``the coloured case''} of $\O_{\C^4}^{\oplus r}$ by Nekrasov-Piazzalunga \cite{NP1}, who also proposed the formula for the signs (which we discuss in Section \ref{sec:local3fold}). We refer to the second equality as the \emph{Magnificent Four formula}. Let us explain the notation on the right hand side. For any formal variable $x$, we set 
$$
[x] := x^{\frac{1}{2}} - x^{-\frac{1}{2}}.
$$ 
Furthermore, for any formal power series $f(p_1, \ldots, p_s; q)$ in $\Q(p_1, \ldots, p_r)[\![q]\!]$ with $f(p_1, \ldots, p_s,0) = 0$,
its plethystic exponential is defined by
\begin{align*} 
\Exp(f(p_1, \ldots, p_s;q)) &:= \exp\Big( \sum_{n=1}^{\infty} \frac{1}{n} f(p_1^n, \ldots, p_s^n;q^n) \Big).
\end{align*}
In the theorem, before taking the plethystic exponential, the Exponent is written as a rational function in $q$ and expanded as a formal power series in $q$. For $\mathsf{V}^{\DT}_{\varnothing\varnothing\varnothing\varnothing}(q)$, we always fix the signs as in Theorem \ref{thm:KR}. 

We now formulate the main conjecture of this section.
\begin{conjecture} \label{conj:vertexDTPT0}
For any plane partitions $\boldsymbol{\mu} = \{\mu_a\}_{a=1}^{4}$ with no $\PT_0$ moduli (Definition \ref{def:PT0boxconfig}), there exists a choice of $\pm \sqrt{\cdot}$ assigned to the $(\C^*)^4$-fixed points such that\footnote{Note that the expectation at the end of Remark \ref{rem:modPT0} is part of the conjecture: for $\boldsymbol{\mu}$ with no $\PT_0$ moduli and for all $B \in \mathcal{B}_{\boldsymbol{\mu}}^{(0)}$, the $\PT_0$ vertex $\mathsf{V}_{I_B\udot}$ has no positive $T$-fixed term.}
$$
\frac{\mathsf{V}^{\DT}_{\mu_1\mu_2\mu_3\mu_4}(q)}{\mathsf{V}^{\DT}_{\varnothing\varnothing\varnothing\varnothing}(q)} = \mathsf{V}^{\PT_0}_{\mu_1\mu_2\mu_3\mu_4}(q).
$$
\end{conjecture}

For $\boldsymbol{\mu}$ all finite (the curve case), the right hand side should be viewed as the $\PT$ vertex \cite{CKM1, Liu}. In this case, the above was first conjectured in \cite{CKM1}. Before we provide computational evidence for this conjecture, we discuss some consequences. 

Let $\boldsymbol{\lambda} = \{\lambda_{ab}\}_{1 \leq a < b \leq 4}$ be finite partitions and $\boldsymbol{\mu}^{\boldsymbol{\lambda}} = \{\boldsymbol{\mu}^{\boldsymbol{\lambda}}_a\}_{a=1}^{4}$ be the minimal plane partitions compatible with $\boldsymbol{\lambda}$ \eqref{eqn:lambdamu}. Then there are no $\PT_0$ moduli (Corollary \ref{lem:noPT0moduli}). Furthermore, by Proposition \ref{Lem.PT0=PT1}, any $\PT_0$ pair is a $\PT_1$ pair with cokernel $Q$ of length
$$
\ell(Q) \leq \ell(\ext_{\C^4}^3(\O_{Z_{\boldsymbol{\lambda}}}, \O_{\C^4} )).
$$
Therefore
$$
\mathsf{V}^{\PT_0}_{\mu_1^{\boldsymbol{\lambda}}\mu_2^{\boldsymbol{\lambda}}\mu_3^{\boldsymbol{\lambda}}\mu_4^{\boldsymbol{\lambda}}}(q)
$$
is a Laurent \emph{polynomial} in $q$ with at most $\ell(\ext_{\C^4}^3(\O_{Z_{\boldsymbol{\lambda}}}, \O_{\C^4} ))+1$ terms. Conjecture \ref{conj:vertexDTPT0} implies that $\mathsf{V}^{\DT}_{\mu_1^{\boldsymbol{\lambda}}\mu_2^{\boldsymbol{\lambda}}\mu_3^{\boldsymbol{\lambda}}\mu_4^{\boldsymbol{\lambda}}}(q)/\mathsf{V}^{\DT}_{\varnothing\varnothing\varnothing\varnothing}(q)$ is equal to this Laurent polynomial.

\begin{example} \label{ex:PT0=PT1vertexex}
In the following examples, $Z_{\boldsymbol{\lambda}} $ is not Cohen-Macaulay and we used Section \ref{sec:redist} to determine the $\PT_0$ topological vertex.
\begin{enumerate}
\item For $\lambda_{12} = \lambda_{34} = 1$, and $\lambda_{ab} = \varnothing$ otherwise, we have
$$
 \sfV^{\PT_0}_{\mu_1^{\boldsymbol{\lambda}}\mu_2^{\boldsymbol{\lambda}}\mu_3^{\boldsymbol{\lambda}}\mu_4^{\boldsymbol{\lambda}}}(q) = \frac{[t_1t_3][t_2t_3]}{[y]} q^{-1}  + [t_1t_2]. 
 $$
\item For $\lambda_{12} = 1+t_4$, $\lambda_{13} = \lambda_{34} = 1$, and $\lambda_{ab} = \varnothing$ otherwise, we have
$$ \sfV^{\PT_0}_{\mu_1^{\boldsymbol{\lambda}}\mu_2^{\boldsymbol{\lambda}}\mu_3^{\boldsymbol{\lambda}}\mu_4^{\boldsymbol{\lambda}}}(q) = \frac{[t_1t_3][t_2t_3][t_1 t_2^2 t_3^2]}{[t_1t_2t_3 y]} q^{-1} + [t_1t_2][t_2t_3].
$$
\item For $\lambda_{12} = 1+t_3$, $\lambda_{34} = 1$, and $\lambda_{ab} = \varnothing$ otherwise, we have
$$
 \sfV^{\PT_0}_{\mu_1^{\boldsymbol{\lambda}}\mu_2^{\boldsymbol{\lambda}}\mu_3^{\boldsymbol{\lambda}}\mu_4^{\boldsymbol{\lambda}}}(q) = \frac{[t_1t_3][t_2t_3][t_1t_3^2][t_2t_3^2]}{[y][y t_3^{-1}]} q^{-2} + \frac{[t_1t_2][t_1t_3][t_2t_3][t_3^2]}{[y][t_3]} q^{-1} + [t_1t_2][t_1t_2t_3^{-1}].
 $$
 \end{enumerate}
The fact that $y$ can appear in the denominator is caused by the redistribution. The negative powers in $q$ come from the renormalized volumes (Definition \ref{def:renorm}).
\end{example}

Using the vertex formalism of Section \ref{sec:redist}, we verified Conjecture \ref{conj:vertexDTPT0} in several cases. In the following proposition, the ``bar'' means that we normalized the generating series by its leading term (which is the same for $\DT/\PT_0$), so the generating series starts with 1.
\begin{proposition} \label{prop:verif}
There exist choices of $\pm \sqrt{\cdot}$ such that
$$
\overline{\mathsf{V}}^{\DT}_{\mu_1\mu_2\mu_3\mu_4}(q) = \overline{\mathsf{V}}^{\PT_0}_{\mu_1\mu_2\mu_3\mu_4}(q) \, \mathsf{V}^{\DT}_{\varnothing\varnothing\varnothing\varnothing}(q)
$$
holds modulo $q^N$ in the following cases:
\begin{enumerate}
\item[$\mathrm{(1)}$] $\lambda_{12} = 1$, $\lambda_{ab} = \varnothing$ otherwise, $\boldsymbol{\mu} = \boldsymbol{\mu}^{\boldsymbol{\lambda}}$, and $N=4$. 
\item[$\mathrm{(2)}$] $\lambda_{12} = 1+t_3$, $\lambda_{ab} = \varnothing$ otherwise, $\boldsymbol{\mu} = \boldsymbol{\mu}^{\boldsymbol{\lambda}}$, and $N=3$.
\item[$\mathrm{(3)}$] $\lambda_{12} = 1+t_3+t_4$, $\lambda_{ab} = \varnothing$ otherwise, $\boldsymbol{\mu} = \boldsymbol{\mu}^{\boldsymbol{\lambda}}$, and $N=3$.
\item[$\mathrm{(4)}$] $\lambda_{12} = \lambda_{13} = 1$, $\lambda_{ab} = \varnothing$ otherwise, $\boldsymbol{\mu} = \boldsymbol{\mu}^{\boldsymbol{\lambda}}$, and $N=3$. 
\item[$\mathrm{(5)}$] $\lambda_{12} = \lambda_{34} = 1$, $\lambda_{ab} = \varnothing$ otherwise,  $\boldsymbol{\mu} = \boldsymbol{\mu}^{\boldsymbol{\lambda}}$, and $N=3$.
\item[$\mathrm{(6)}$] $\lambda_{12} = \lambda_{13} = \lambda_{23} = 1$, $\lambda_{ab} = \varnothing$ otherwise, $\boldsymbol{\mu} = \boldsymbol{\mu}^{\boldsymbol{\lambda}}$, and $N=3$. 
\item[$\mathrm{(7)}$] $\lambda_{ab} = 1$ for all $a<b$, $\boldsymbol{\mu} = \boldsymbol{\mu}^{\boldsymbol{\lambda}}$, and $N=3$.  
\item[$\mathrm{(8)}$] $\lambda_{12} = \lambda_{34} =\lambda_{23}  = 1$, $\lambda_{ab} = \varnothing$ otherwise, $\boldsymbol{\mu} = \boldsymbol{\mu}^{\boldsymbol{\lambda}}$, and $N=3$.
\item[$\mathrm{(9)}$] $\lambda_{12} = 1+t_4$, $\lambda_{13} = \lambda_{34} = 1, \lambda_{ab} = \varnothing$ otherwise,  $\boldsymbol{\mu} = \boldsymbol{\mu}^{\boldsymbol{\lambda}}$, and $N=3$.
\item[$\mathrm{(10)}$] $\lambda_{12} = 1+t_3$, $\lambda_{34} = 1, \lambda_{ab} = \varnothing$ otherwise,  $\boldsymbol{\mu} = \boldsymbol{\mu}^{\boldsymbol{\lambda}}$, and $N=3$.

\item[$\mathrm{(11)}$] $\mu_1 = \frac{1}{1-t_2}+t_3$, $\mu_2 = \frac{1}{1-t_1}$, $\mu_a = \varnothing$ otherwise, and $N=4$.
\item[$\mathrm{(12)}$] $\mu_1 = \frac{1}{1-t_2}+t_3+t_2 t_3$, $\mu_2 = \frac{1}{1-t_1}$, $\mu_a = \varnothing$ otherwise, and $N=3$.
\item[$\mathrm{(13)}$] $\mu_1 = \frac{1}{1-t_2}+ t_3+t_3^2$, $\mu_2 = \frac{1}{1-t_1}$, $\mu_a = \varnothing$ otherwise, and $N=3$.
\item[$\mathrm{(14)}$] $\mu_1 = \frac{1}{1-t_2}$, $\mu_2 = \frac{1}{1-t_1}$, $\mu_3 = 1$, $\mu_4 = \varnothing$, and $N=3$.
\item[$\mathrm{(15)}$] $\mu_1 = \frac{1}{1-t_2} + t_3$, $\mu_2 = \frac{1}{1-t_1} + t_4$, $\mu_3=1$, $\mu_4 = \varnothing$, and $N=2$. 
\end{enumerate}
\end{proposition}
In this list, the number of fixed points on the $\DT$ side can be rather large. In case (1), for 3 embedded points we have 19 torus fixed points. In case (11), for 3 embedded points we have 33 torus fixed points.

Since face and edge terms coincide for $\DT$ and $\PT_0$ theory, Proposition \ref{prop:DTsumsqrt}, Remark \ref{rem:modPT0}, and Conjecture \ref{conj:vertexDTPT0} imply the \emph{global} $T$-equivariant $K$-theoretic $\DT$--$\PT_0$ correspondence (Conjecture \ref{conj:KDTPT0inbodytext_equiv}):
\begin{corollary} \label{cor:globaltoricDTPT0}
Let $X$ be a toric Calabi-Yau 4-fold and let $L$ be a $T_X$-equivariant line bundle on $X$. Let $\gamma \in H_c^4(X,\Q)$ and $\beta \in H_c^6(X,\Q)$ such that $\curP_{v}^{(0)}(X)^{T_X} = \curP_{v}^{(0)}(X)^{T}$ is empty or 0-dimensional and reduced for $v = (0,0,\gamma,\beta,n-\gamma \cdot \td_2(X))$ and all $n \in \Z$. Suppose Conjecture \ref{conj:vertexDTPT0} holds. Then there are choices of $\pm \sqrt{\cdot}$ such that
$$
\frac{\sum_{n} \langle \! \langle L  \rangle \! \rangle_{X,\gamma,\beta,n}^{\DT} q^n}{\sum_{n} \langle \! \langle L  \rangle \! \rangle_{X,0,0,n}^{\DT} q^n} = \sum_{n} \langle \! \langle L  \rangle \! \rangle_{X,\gamma,\beta,n}^{\PT_0} q^n.
$$
\end{corollary}

\subsection{Cohomological limits} \label{sec:coholimittaut}

Let $X$ be a toric Calabi-Yau 4-fold and $L$ a $T_X$-equivariant line bundle on $X$. Let $\curly P := \curly P_v^{(q)}(X)$ for any $v$ as in \eqref{eq:toricsecmoduli} and $q \in \{-1,0,1\}$. Besides the $K$-theoretic invariants $\langle\!\langle L \rangle\!\rangle_{X,v}^{\PT_q}$,
we are also interested in the following ``tautological'' and ``cohomological'' invariants 
\begin{align*}
\langle L \rangle_{X,v}^{\PT_q} \in \frac{\Q(s_1,s_2,s_3,s_4,m)}{(s_1+s_2+s_3+s_4)}, \quad \int_{[\curly P]^{\vir}} 1 \in \frac{\Q(s_1,s_2,s_3,s_4)}{(s_1+s_2+s_3+s_4)},
\end{align*}
where the former were defined in Section \ref{sec:vardef}, and we recall that $c_1^T(t_i) = s_i$ and $c_1^{\C^*}(y) = m$ are the cohomological equivariant parameters. 
Similar to the $K$-theoretic case, these invariants can be calculated using the vertex, edge, face terms $\widetilde{\mathsf{V}}_{\alpha}$, $\widetilde{\mathsf{E}}_{\alpha\beta}$, $\widetilde{\mathsf{F}}_{f}$ of Section \ref{sec:redist}. The only differences are the following:
\begin{itemize}
\item tautological invariants: replace 
$$
[t_1^{i_1} t_2^{i_2} t_3^{i_3} t_4^{i_4} y^{i_5}] := t_1^{\frac{i_1}{2}} t_2^{\frac{i_2}{2}} t_3^{\frac{i_3}{2}} t_4^{\frac{i_4}{2}} y^{\frac{i_5}{2}} - t_1^{-\frac{i_1}{2}} t_2^{-\frac{i_2}{2}} t_3^{-\frac{i_3}{2}} t_4^{-\frac{i_4}{2}} y^{-\frac{i_5}{2}}$$ 
by 
$$e(t_1^{i_1} t_2^{i_2} t_3^{i_3} t_4^{i_4} y^{i_5}) := i_1 s_1+ i_2 s_2+ i_3 s_3+i_4 s_4+ i_5 m,$$
\item cohomological invariants: the same as the previous case but omit all terms involving $y$ from the vertex, edge, and face terms of Section \ref{sec:redist}, in other words, use $\mathsf{V}_{\alpha}$, $\mathsf{E}_{\alpha\beta}$, $\mathsf{F}_{f}$ instead of $\widetilde{\mathsf{V}}_{\alpha}$, $\widetilde{\mathsf{E}}_{\alpha\beta}$, $\widetilde{\mathsf{F}}_{f}$.
\end{itemize}

For the vertex formalism in the curve case we have \cite[Prop.~2.3]{CKM1}
$$
\rk(\widetilde{\mathsf{V}}_{\alpha}) = \rk(\widetilde{\mathsf{E}}_{\alpha\beta}) = 0.
$$
This is no longer the case for surfaces. Then $\widetilde{\mathsf{V}}_{\alpha}$, $\widetilde{\mathsf{E}}_{\alpha\beta}$, $\widetilde{\mathsf{F}}_{f}$ can have non-zero rank (e.g.~see Example \ref{ex:PT0=PT1vertexex}). 
However, we have the following weaker result, which we formulate for the case $q=-1$, but it also holds for $q=0,1$.
\begin{lemma} \label{lem:samerk}
Let $Z, Z' \subset \C^4$ be 2-dimensional $(\C^*)^4$-invariant closed subschemes with the same underlying pure 2-dimensional sheaf, i.e., $Z_1 = Z_1'$.  Then
$$
\rk(\widetilde{\mathsf{V}}_{Z}) = \rk(\widetilde{\mathsf{V}}_{Z'}).
$$
\end{lemma}
\begin{proof}
Consider the definition of $\widetilde{\mathsf{V}}_Z$ (Definition \ref{def:vertex}) and make the substitutions of Definition \ref{def:renorm} (where we write $Z = Z_\alpha$, $W = W_\alpha$, $W_i := W_{\alpha\beta_i}$, $W_{ij} = Z_{\alpha\beta_i\beta_j}$). Taking $t_1=\cdots=t_4=1$, an explicit calculation gives
\begin{equation} \label{eqn:rkVtilde}
\widetilde{\mathsf{V}}_Z \Big|_{t_1=\cdots=t_4=1} = \Big( W + \overline{W} - y^{-1} W - y \overline{W} - \sum_{i,j,k,l \in \{1,2,3,4\} \atop \textrm{distinct } i<j, k<l} W_{ij} \overline{W}_{kl} \Big)\Big|_{t_1=\cdots=t_4=1}, 
\end{equation}
which is independent of $W_i$. We see that for fixed 2-dimensional asymptotics $W_{ij}$, adding or removing embedded points or curves keeps the rank of $\widetilde{\mathsf{V}}_Z$ unchanged as claimed. 
\end{proof}

We formulate the following result for $q=-1$, but it holds analogously for $q=0,1$.
\begin{proposition} 
Let $Z, Z' \subset \C^4$ be 2-dimensional $(\C^*)^4$-invariant closed subschemes with the same underlying pure 2-dimensional sheaf, i.e., $Z_1 = Z_1'$. Suppose $\widetilde{V}_{Z'}$ has no $T$-fixed term. Then
\begin{align*}
\lim_{b \mapsto 0} \frac{[-\widetilde{\mathsf{V}}_{Z}]}{[-\widetilde{\mathsf{V}}_{Z'}]} \Big|_{t_i = e^{b s_i}, y = e^{b m}} &= \frac{e(-\widetilde{\mathsf{V}}_{Z})}{e(-\widetilde{\mathsf{V}}_{Z'})}, \\
\lim_{b \mapsto 0 \atop m \to \infty}  \frac{[-\widetilde{\mathsf{V}}_{Z}] / (-m^2)^{|Z|}}{[-\widetilde{\mathsf{V}}_{Z'}]/ (-m^2)^{|Z'|}} \Big|_{t_i = e^{b s_i}, y = e^{b m}} &= \frac{e(-\mathsf{V}_{Z})}{e(-\mathsf{V}_{Z'})},
\end{align*}
where $|Z|$, $|Z'|$ denote the renormalized volumes of $Z,Z'$.
\end{proposition}
%Strictly speaking, we assume that $\widetilde{\mathsf{V}}_{Z}$ has no positive $T$-fixed term and $\widetilde{\mathsf{V}}_{Z'}$ has no $T$-fixed term.
\begin{proof}
Using multi-index notation $\eta = (t_1,t_2,t_3,t_4,y)$, for any \emph{rank 0} finite complex $E = \sum_v \eta^v - \sum_w \eta^w$, we have
$$
[E] = \frac{\prod_v [\tau^v]}{\prod_w [\tau^w]} \Big|_{t_i = e^{b s_i}, y = e^{b m}} = \frac{\prod_v (v_1 s_1+v_2s_2+v_3s_3+v_4s_4 + v_5 m + O(b))}{\prod_w (w_1s_1+w_2s_2+w_3s_3+w_4s_4 + w_5 m + O(b))},
$$
and therefore
$$
\lim_{b \mapsto 0} \widehat{\Lambda}\udot E^* = \lim_{b \mapsto 0} [E] = e(E).
$$
The first statement then follows from Lemma \ref{lem:samerk}. 

Now we turn to the second statement. Let $Z \subset \C^4$ be the 2-dimensional subscheme corresponding to $\pi$ and use the notation of Definition \ref{def:renorm} (with $Z = Z_\alpha$, $W = W_\alpha$). The only term involving $y$ in $[-\widetilde{\mathsf{V}}_{Z}]/(-m^2)^{|Z| }$ is $[y^{-1} W + y \overline{W}]/(-m^2)^{|Z| }$, where $|Z| = \rk(W)$ is the renormalized volume (Definition \ref{def:renorm}). We use multi-index notation $\tau = (t_1,t_2,t_3,t_4)$. Writing $W = \sum_v \tau^v - \sum_w \tau^w$ and taking the limit $b \to 0$, it contributes
\begin{align*}
\frac{[y^{-1} W + y \overline{W}]}{(-m^2)^{|Z| }} &=  (-m^2)^{-(\sum_v 1 -\sum_w 1) } \frac{\prod_v (-1)(v_1s_1+v_2s_2+v_3s_3+v_4s_4 - m)^2}{\prod_w (-1)(w_1 s_1+w_2s_2+w_3s_3+w_4s_4 - m)^2} \\
&=  \frac{\prod_v (\frac{v_1}{m}s_1+\frac{v_2}{m}s_2+\frac{v_3}{m}s_3+\frac{v_4}{m}s_4 -1)^2}{\prod_w (\frac{w_1}{m} s_1+\frac{w_2}{m}s_2+\frac{w_3}{m}s_3+\frac{w_4}{m}s_4 - 1)^2}.
\end{align*}
This term goes to 1 as $m \to \infty$. 
\end{proof}

For $q=-1,0$, we accordingly define the tautological and cohomological topological vertex
$$
\mathsf{V}^{\PT_q, \textrm{taut}}_{\mu_1\mu_2\mu_3\mu_4}(q), \quad \mathsf{V}^{\PT_q, \textrm{coho}}_{\mu_1\mu_2\mu_3\mu_4}(q).
$$
In the following corollary, the ``bar'' means that we normalized the generating series by its leading term, so the generating series starts with 1.

\begin{corollary} 
Let ${\boldsymbol{\mu}}$ be plane partitions corresponding to subschemes of dimension $\leq 1$ of $\C^3$. Then for $q=-1,0$ we have
\begin{align*}
\lim_{b \mapsto 0} \overline{\mathsf{V}}^{\PT_q, \mathrm{taut}}_{\mu_1\mu_2\mu_3\mu_4}(q) \Big|_{t_i = e^{b s_i}, y = e^{b m}} &= \overline{\mathsf{V}}^{\PT_q}_{\mu_1\mu_2\mu_3\mu_4}(q), \\
\lim_{b \mapsto 0 \atop m \to \infty} \overline{\mathsf{V}}^{\PT_q, \mathrm{coho}}_{\mu_1\mu_2\mu_3\mu_4}(q') \Big|_{t_i = e^{b s_i}, y = e^{b m}, q' = - \frac{q}{m}} &= \overline{\mathsf{V}}^{\PT_q}_{\mu_1\mu_2\mu_3\mu_4}(q).
\end{align*}
\end{corollary}
%Strictly speaking, for this corollary and the paragraph below: we assume that $\widetilde{V}_\pi$, for $\pi \in \mathcal{P}_{\boldsymbol{\mu}}$ minimal has no negative $T$-fixed term.

In particular, the $K$-theoretic vertex $\DT$--$\PT_0$ correspondence (Conjecture \ref{conj:vertexDTPT0}) implies a tautological and cohomological vertex $\DT$--$\PT_0$ correspondence (for the same choice of signs). Consequently, we also obtain global analogs of Corollary \ref{cor:globaltoricDTPT0} for tautological and numerical insertions.

\subsection{Local 3-folds} \label{sec:local3fold}

For $Y$ a smooth quasi-projective toric 3-fold, we consider the toric Calabi-Yau 4-fold 
$$
p : X = \mathrm{Tot}(K_Y) \to Y,
$$
i.e., the total space of the canonical line bundle on $Y$. Denote by $T_Y$ the torus acting on $Y$ and by $\C^*$ the torus acting on the fibres, then $T_X = T_Y \times \C^*$. The torus $T \leq T_X$ preserving the Calabi-Yau volume form is isomorphic to $T_Y$. Let $v = (0,0,\gamma,\beta,n-\gamma \cdot \td_2(X)) \in H_c^*(X,\Q)$ and $\curly P := \curly P_v^{(q)}(X)$. In Section \ref{sec:Leflocal}, assuming  $\curly P^T$ is 0-dimensional and reduced, we deduced
\begin{equation} \label{eqn:invKYinitial} 
\langle\!\langle L \rangle\!\rangle_{X,v}^{\PT_q} = \sum_{I\udot = [\O_X \to F] \in \curly P^T} (-1)^{\sigma_{I\udot}} \frac{\widehat{\Lambda}\udot (R\Gamma(X,F \otimes L) \otimes y^{-1})}{\widehat{\Lambda}\udot (R\Gamma(Y,p_*F) - R\Hom_Y(p_* F,p_*F))^\vee},
\end{equation}
for any $T_X$-equivariant line bundle $L$ on $X$ and choice of sign attached to each $T$-fixed point. Recall that in Section \ref{sec:Leflocal}, we focused on the case $L = \O_X(Y)$, but now $L$ is arbitrary.

Let $I\udot = [\O_X \to F] \in \curly P^{T_X}$. In this section, we develop a vertex formalism for the calculation of
$$
R\Gamma(Y, p_* F) - R\Hom_Y(p_* F, p_* F))
$$
and comment on the signs $(-1)^{\sigma_{I\udot}}$.

We start with the following definitions. Note that for $X = \mathrm{Tot}(K_Y)$, we have $V(X) = V(Y)$, $E(X) = E(Y)$, $F(X) = F(Y)$, because we only consider compact edges and faces. For notational simplicity, we consider the case $q = -1$, but the cases $q=0,1$ hold similarly (subject to the modifications of Remarks \ref{rem:modPT0} and \ref{rem:modPT1}). Let $Z \subset X$ be a closed subscheme of dimension $\leq 2$ with proper support. For $\alpha_1, \ldots, \alpha_N \in V(X)$, such that $U_{\alpha_1 \cdots \alpha_N} \cong (\C^*)^k \times \C^{4-k}$ with $k=0,1,2$, we always label coordinates on $(\C^*)^k \times \C^{4-k}$ in such a way that $Z \cap U_{\alpha_1 \cdots \alpha_N}$ is set theoretically (but not necessarily scheme theoretically!) supported on $Z(x_4)$. We also choose the coordinates such that $T_X$ acts on $U_{\alpha_1 \cdots \alpha_N}$ by the standard torus action \eqref{eqn:stndact}. In other words, locally, the fibre coordinate of $X = \mathrm{Tot}(K_Y)$ always corresponds to $x_4$. We use the notation of Section \ref{sec:redist}.
\begin{definition} 
We define the \emph{halved vertex term}
\begin{align*}
\mathsf{v}_{\alpha} &:= Z_\alpha + \frac{P_{123}}{t_1t_2t_3} Z_\alpha \overline{Z}_{\alpha} + \sum_{i=1}^{3} \frac{\mathsf{a}_{\alpha\beta_i}}{1-t_i} + \sum_{1 \leq i<j \leq 3} \frac{\mathsf{a}_{\alpha\beta_i\beta_j}}{(1-t_i)(1-t_j)},
\end{align*}
where, for any distinct $i,j,k \in \{1,2,3\}$, we set
\begin{align*}
\mathsf{a}_{\alpha\beta_i} &:= - Z_{\alpha\beta_i}  + \frac{P_{jk}}{t_{j}t_{k}}  Z_{\alpha\beta_i}   \overline{Z}_{\alpha\beta_i}, \\
\mathsf{a}_{\alpha\beta_i\beta_{j}} &:= Z_{\alpha\beta_i\beta_j} + \frac{P_{k}}{t_{k}}  Z_{\alpha\beta_i\beta_j}   \overline{Z}_{\alpha\beta_i\beta_j}.
\end{align*}
When we want to stress the dependence of these expressions on $Z_\alpha, Z_{\alpha\beta}, Z_{\alpha\beta\gamma}$, we write $\mathsf{v}_{Z_\alpha}, \mathsf{a}_{Z_{\alpha\beta}}, \mathsf{a}_{Z_{\alpha\beta\gamma}}$. 
\end{definition}

\begin{definition} \label{def:halfedge}
For any distinct $i,j,k \in \{1,2,3\}$, we define the \emph{halved edge term} by
\begin{align*} 
\mathsf{e}_{\alpha\beta_i} &:= \frac{t_i^{-1}  \mathsf{b}_{\alpha\beta_i}(t_j,t_k,t_4) }{1-t_i^{-1}} -  \frac{\mathsf{b}_{\alpha\beta_i}(t_i^{-m_{\alpha\beta_i}} t_j,t_i^{-m_{\alpha\beta_i}'} t_k,t_i^{-m_{\alpha\beta_i}''}t_4) }{1-t_i^{-1}}, 
\end{align*}
where 
\begin{align*}
\mathsf{b}_{\alpha\beta_i} &:= \mathsf{a}_{\alpha\beta_i}  + \sum_{a \in \{i,j,k\} \setminus \{i\}} \frac{\mathsf{a}_{\alpha\beta_i\beta_a} }{1-t_a} .
\end{align*}
When we want to indicate the dependence of these expressions on $Z_\alpha, Z_{\alpha\beta}$, we write $\mathsf{e}_{Z_{\alpha\beta}}, \mathsf{b}_{Z_{\alpha\beta}}$.
\end{definition}

\begin{definition} 
Fix $f \in F(X) = F(Y)$. Denote by $Z_f \subset X$ the unique $T_X$-fixed 2-dimensional closed subscheme set-theoretically supported on $S_f \subset Y \subset X$ such that 
$$
Z_f \cap U_{\alpha_1\alpha_2\alpha_3} = Z \cap U_{\alpha_1\alpha_2\alpha_3}
$$
for any distinct $\alpha_1, \alpha_2, \alpha_3 \in f$.  We define the \emph{halved face term} by
$$
\mathsf{f}_f := \tr_{R\Gamma(Y,p_* \O_{Z_f})-R\Hom_Y(p_* \O_{Z_f}, p_*\O_{Z_f})}.
$$
When we want to stress the dependence of $\mathsf{f}_f $ on $Z_f$ we write $\mathsf{f}_{Z_f}$. 
\end{definition}

For $\dim Z = 0$, $\mathsf{v}_{\alpha}$ reduces to the halved vertex term first introduced by Nekrasov-Piazzalunga \cite{Nek, NP1}. For $\dim Z = 1$, $\mathsf{v}_{\alpha}$ and $\mathsf{e}_{\alpha\beta}$ reduce to the halved vertex and edge terms studied in \cite{CKM1}. The following lemma justifies the use of the adjective ``halved''.
\begin{lemma} \label{lem:KYhalf}
For all $\alpha \in V(X)$, $\alpha\beta \in E(X)$, $f \in F(X)$, the following equations hold modulo the relation $t_1t_2t_3t_4=1$
\begin{align*}
\mathsf{v}_{\alpha} + \overline{\mathsf{v}}_{\alpha} = \mathsf{V}_\alpha, \quad \mathsf{e}_{\alpha\beta} + \overline{\mathsf{e}}_{\alpha\beta} = \mathsf{E}_{\alpha\beta}, \quad \mathsf{f}_{f} + \overline{\mathsf{f}}_{f} = \mathsf{F}_f.
\end{align*}
\end{lemma}
\begin{proof}
The first two equations follow from a straight-forward calculation involving identities such as 
$$
\frac{P_{123}}{t_1t_2t_3} + \frac{\overline{P}_{123}}{t_1^{-1}t_2^{-1}t_3^{-1}} = - \frac{P_{1234}}{t_1t_2t_3t_4},
$$
which uses the relation  $t_1t_2t_3t_4=1$. The third identity follows from \eqref{eqn:Tvirspectral}.
\end{proof}

The following proposition is proved similarly to Theorem \ref{prop:redist}.
\begin{proposition} \label{prop:KYvef}
Consider $p  : X = \mathrm{Tot}(K_Y) \to Y$, where $Y$ is a smooth quasi-projective toric 3-fold. Then for any $T_X$-fixed closed subscheme $Z \subset X$ with proper support and dimension $\leq 2$, we have
$$
R\Gamma(Y, p_* \O_Z) - R\Hom_Y(p_* \O_Z, p_* \O_Z) = \sum_{\alpha \in V(X)} \mathsf{v}_{\alpha} + \sum_{\alpha\beta \in E(X)} \mathsf{e}_{\alpha\beta} + \sum_{f \in F(X)} \mathsf{f}_{f},
$$
where $\mathsf{v}_{\alpha}, \mathsf{e}_{\alpha\beta}, \mathsf{f}_{f}$ are Laurent polynomials.
\end{proposition}

\begin{proof}
Denote by $\{U_\alpha\}_{\alpha \in V(Y)}$ the open cover of $Y$ by maximal $T_Y$-invariant affine open subsets. Choose a total order on the indices $\alpha \in V(Y)$. Using \v{C}ech cohomology and the local-to-global spectral sequence as in the beginning of Section \ref{sec:redist}, we obtain
\begin{align*}
&R\Gamma(Y, p_* \O_Z) - R\Hom_Y(p_* \O_Z,p_* \O_Z) = \\
&\quad\quad\quad\quad \sum_{\alpha \in V(Y)} \Big( R\Gamma(U_\alpha, p_* \O_Z) - R\Hom_{U_{\alpha}}(p_* \O_Z,p_* \O_Z) \Big) \\
&\quad\quad\quad\quad- \sum_{\alpha < \beta \in V(Y)} \Big( R\Gamma(U_{\alpha\beta}, p_* \O_Z) - R\Hom_{U_{\alpha\beta}}(p_* \O_Z,p_* \O_Z) \Big) \\
&\quad\quad\quad\quad+ \sum_{\alpha < \beta < \gamma \in V(Y)} \Big( R\Gamma(U_{\alpha\beta\gamma}, p_* \O_Z) - R\Hom_{U_{\alpha\beta\gamma}}(p_* \O_Z,p_* \O_Z) \Big) + \cdots,
\end{align*}
where we suppress restrictions to $U_\alpha$, $U_{\alpha\beta}$, $U_{\alpha\beta\gamma}$, $\ldots$

There exist finitely many $T_Y$-invariant closed subschemes 
$$
Y \supset Z(0) \supset Z(1) \supset \cdots 
$$ 
such that
$$
p_* \O_Z = \bigoplus_{i \geq 0} \O_{Z(i)} \otimes t_4^i.
$$
We denote their restrictions by $Z_\alpha(i) = Z(i) \cap U_\alpha$, $Z_{\alpha\beta}(i) = Z(i) \cap U_{\alpha\beta}$, $Z_{\alpha\beta\gamma}(i) = Z(i) \cap U_{\alpha\beta\gamma}$, etc. Using Taylor resolutions \cite{Tay}, it is easy to see that 
\begin{align*}
&R\Gamma(U_\alpha, p_* \O_Z) - R\Hom_{U_{\alpha}}(p_* \O_Z,p_* \O_Z) \\
&= \sum_{i \geq 0} R\Gamma(U_\alpha,\O_{Z_\alpha(i)})  t_4^i - \sum_{i,j \geq 0} R\Hom_{U_\alpha}(\O_{Z_{\alpha}(i)}, \O_{Z_{\alpha}(j)}) t_4^{j-i} \\
&=\sum_{i \geq 0} Z_{\alpha}(i) t_4^i - \sum_{i,j \geq 0} \frac{P_{123}}{t_1t_2t_3} \overline{Z_\alpha(i)} Z_\alpha(j) t_4^{j-i} = Z_\alpha + \frac{P_{123}}{t_1t_2t_3} Z_\alpha \overline{Z}_\alpha.
\end{align*}
Let $\alpha\beta \in E(Y)$. Let $\alpha_1,\ldots, \alpha_N \in V(Y)$ such that $\alpha_1, \ldots, \alpha_N \in f$ for some $f \in F(Y)$, where for $N=2$ we assume $\alpha_1\alpha_2$ is not an edge. Then $U_{\alpha\beta} \cong \C^* \times \C^2$ and $U_{\alpha_1\cdots\alpha_N} \cong (\C^*)^2 \times \C$, and we similarly obtain 
\begin{align*}
&R\Gamma(U_{\alpha\beta}, p_* \O_Z) - R\Hom_{U_{\alpha\beta}}(p_* \O_Z,p_* \O_Z) = \delta(t_1)\Big(Z_{\alpha\beta}  - \frac{P_{23}}{t_{2}t_{3}}  Z_{\alpha\beta}   \overline{Z}_{\alpha\beta}\Big), \\
&R\Gamma(U_{\alpha_1\cdots \alpha_N}, p_* \O_Z) - R\Hom_{U_{\alpha_1\cdots \alpha_N}}(p_* \O_Z,p_* \O_Z) = \\ &\quad\quad\quad\quad\quad\quad\quad\quad\quad\quad\quad\quad\quad\quad\quad\quad \delta(t_1)\delta(t_2)\Big(  Z_{\alpha_1\cdots \alpha_N} + \frac{P_{3}}{t_{3}}  Z_{\alpha_1\cdots \alpha_N}   \overline{Z}_{\alpha_1\cdots \alpha_N}  \Big).
\end{align*}
For the remainder of the proof, we can follow the proof of Theorem \ref{prop:redist}.
\end{proof}

\begin{remark}
Similar to Remark \ref{rem:expFf}, we can obtain a nice expression for $\mathsf{f}_f$ and $f \in F(Y)$. For a 2-dimensional $T_X$-fixed closed subscheme $Z \subset X$ with proper support, denote by $\lambda_f$ the finite partition corresponding to the subscheme $Z_f$ from Definition \ref{def:face}. Let
\[
\mathsf{a}_f(t_3,t_4) := Z_{\lambda_f} + \frac{P_{3}}{t_3} Z_{\lambda_f} \overline{Z}_{\lambda_f} \in \mathbb{Z}[t_3^{\pm 1}, t_4^{\pm 1}],
\]
where $Z_{\lambda_f} = \sum_{(i,j) \in \lambda_f} t_3^{i} t_4^{j}$. Then
\begin{equation} \label{eqn:exprff}
\mathsf{f}_f = \sum_{i=1}^{e(S_f)} \frac{\mathsf{a}_f(t_3^{(i)},t_4^{(i)})}{(1-t_1^{(i)})(1-t_2^{(i)})},
\end{equation}
where the right hand side is an element of $\Q[t_1^{\pm 1}, \ldots, t_4^{\pm 1}] \subset \Q(t_1, \ldots, t_4)$.
We observe that, modulo the relation $t_1t_2t_3t_4=1$, it is also clear by direct computation that the sum of \eqref{eqn:exprff} and its conjugate equals \eqref{eqn:exprFf}. In fact, \eqref{eqn:exprff} provides a ``half'' in $K$-theory of \eqref{eqn:exprFf} for any toric Calabi-Yau 4-fold (not necessarily of the form $X = \mathrm{Tot}(K_Y)$).
\end{remark}

Next, we incorporate the tautological insertion 
$$
\widehat{\Lambda}\udot ( R\pi_{\curly P*}\O_{\cZ} \otimes y^{-1})
$$
into the vertex formalism, which we already carried out  in Section \ref{sec:redist}. We define for all $i,j,k \in \{1,2,3\}$ distinct
\begin{align*}
\widetilde{\mathsf{v}}_{\alpha} := &\mathsf{v}_{\alpha} - y \overline{Z}_\alpha + \sum_{i=1}^{3} \frac{y \overline{Z}_{\alpha\beta_i}}{1-t_i^{-1}} - \sum_{1 \leq i<j \leq 3} \frac{y \overline{Z}_{\alpha\beta_i\beta_j}}{(1-t_i^{-1})(1-t_j^{-1})}, \\
\widetilde{\mathsf{e}}_{\alpha\beta_i} := &\mathsf{e}_{\alpha\beta_i} + \frac{t_i}{1-t_i} \Bigg\{ y \overline{Z}_{\alpha\beta_i} - \sum_{a \in \{j,k\}} \frac{y \overline{Z}_{\alpha\beta_i\beta_a}}{1-t_a^{-1}} \Bigg\} \\
&\quad\quad- \frac{1}{1-t_i} \Bigg\{ y \overline{Z}_{\alpha\beta_i} - \sum_{a \in \{j,k\}} \frac{y \overline{Z}_{\alpha\beta_i\beta_a}}{1-t_a^{-1}} \Bigg\}\Bigg|_{(t_j t_i^{-m_{\alpha\beta_i}},t_k t_i^{-m_{\alpha\beta_i}'},t_4 t_i^{-m_{\alpha\beta_i}''})}, \\
\widetilde{\mathsf{f}}_f  := &\mathsf{f}_f -y R\Gamma(X, \O_{Z_f})^\vee.
\end{align*}

Combining \eqref{eqn:invKYinitial} and Proposition \ref{prop:KYvef}, we obtain the following result.
\begin{proposition} \label{prop:DTsumsqrtKY}
  Let $X = \mathrm{Tot}(K_Y)$ for a smooth quasi-projective toric 3-fold $Y$ and fix $v = (0,0,\gamma,\beta,n - \gamma \cdot \td_2(X)) \in H_c^*(X,\Q)$. Let $\curly P^T := \curly P_v^{(-1)}(X)^T$. Then 
\[
\langle\!\langle \O_X \rangle\!\rangle_{X,v}^{\DT} = \sum_{Z \in \curly P^T} \pm \Big[-\sum_{\alpha \in V(Y)} \widetilde{\mathsf{v}}_{Z_\alpha} - \sum_{\alpha\beta \in E(Y)} \widetilde{\mathsf{e}}_{Z_{\alpha\beta}} - \sum_{f \in F(Y)} \widetilde{\mathsf{f}}_{Z_f} \Big],
\]
where the right hand side depends on a choice of $\pm$ for each fixed point.
\end{proposition}

By Lemmas \ref{lem:noposTfix} and \ref{lem:KYhalf}, $\mathsf{v}_\alpha$, $\mathsf{e}_{\alpha\beta}$, $\mathsf{f}_f$ have no positive $T$-fixed part, so we may also write
\begin{align*}
\langle\!\langle \O_X \rangle\!\rangle_{X,v}^{\DT} = \sum_{Z \in \curly P^T} \pm \prod_{\alpha \in V(Y)} [-\widetilde{\mathsf{v}}_{Z_\alpha}] \cdot \prod_{\alpha\beta \in E(Y)} [-\widetilde{\mathsf{e}}_{Z_{\alpha\beta}}] \cdot \prod_{f \in F(Y)}  [-\widetilde{\mathsf{f}}_{Z_f}].
\end{align*}
We stress that the elements of $\curly P^T$ are set theoretically, but not necessarily scheme theoretically, supported on the zero section $Y \subset X$. The invariants $\langle\!\langle L \rangle\!\rangle_{X,v}^{\DT} $, for other $T_X$-equivariant line bundles $L$ on $X$, can also easily be incorporated into the vertex formalism by Remark \ref{rem:otherL}.

Let $\boldsymbol{\mu} = \{\mu_{a} \}_{a = 1}^{4}$ be a collection of plane partitions corresponding to subschemes of dimension $\leq 1$ of $\C^3$ such that $Z_{\boldsymbol{\mu}}$ (Definition \ref{def:partitionnotation}) is \emph{set theoretically} supported on $Z(x_4)$. In particular $\mu_4 = \varnothing$. Then the $\DT$ topological vertex equals (Definition \ref{def:DTvertex})
\begin{equation} \label{eqn:DTvertexKY}
\mathsf{V}^{\DT}_{\mu_1\mu_2\mu_3\varnothing}(q) = \sum_{\pi \in \mathcal{P}_{\boldsymbol{\mu}}} \pm  [-\widetilde{\mathsf{v}}_{\pi}] \, q^{|\pi|},
\end{equation}
where $\pm$ indicates that this expression involves a choice of sign for each $\pi \in \mathcal{P}_{\boldsymbol{\mu}}$. 

\begin{remark}
Just like in Remark \ref{rem:modPT0} and \ref{rem:modPT1}, the results of this section so far hold similarly for $q=0,1$. 
\end{remark}

In the remainder of this section, we take $q=-1$ and comment on the signs in \eqref{eqn:DTvertexKY}. It is tempting to think that the ``natural'' sign of every solid partition in \eqref{eqn:DTvertexKY} is $+1$. However, neither the Magnificent Four formula (Theorem \ref{thm:KR}) nor the $\DT$--$\PT_0$ vertex correspondence (Conjecture \ref{conj:vertexDTPT0}) hold for this choice of signs. 

Suppose $\boldsymbol{\mu}$ are all empty. Then it was first found in physics by Nekrasov-Piazzalunga \cite{Nek, NP1} and proved mathematically by the second-named author and Rennemo \cite{KR}, that the Magnificent Four formula holds when we assign the sign
\begin{equation} \label{eqn:sign0dim}
(-1)^{|\pi| + |\{(i,i,i,j) \in \pi  \, : \, i<j \}|}
\end{equation}
to each finite solid partition $\pi$. Therefore, \eqref{eqn:sign0dim} provides the correct sign choice when all $\boldsymbol{\mu}$ are empty. 

Suppose $\boldsymbol{\mu}$ are \emph{finite} plane partitions and $\mu_4 = \varnothing$. For all verifications of the $\DT$--$\PT$ (curve) vertex correspondence performed in \cite{CKM1}, it was found that, on the $\DT$ side, the following choice of sign for each $\pi \in \mathcal{P}_{{\boldsymbol{\mu}}}$ works
\begin{equation} \label{eqn:sign1dim}
(-1)^{|\pi|} \prod_{w = (i,i,i,j) \in \pi \atop i< j} (-1)^{1  -  |\{\textrm{legs containing } w\}|}.
\end{equation}
When $ \boldsymbol{\mu}$ are all empty, \eqref{eqn:sign1dim} reduces to \eqref{eqn:sign0dim}. Conjecturally, on the $\DT$ side, this provides the sign choice for which the $\DT$--$\PT$ (curve) vertex correspondence holds when all $\boldsymbol{\mu}$ are finite plane partitions and $\mu_4 = \varnothing$.

Now let $\boldsymbol{\mu} = \{\mu_{a} \}_{a = 1}^{4}$ be a collection of plane partitions corresponding to subschemes of dimension $\leq 1$ of $\C^3$ such that $Z_{\boldsymbol{\mu}}$ is set theoretically supported on $Z(x_4)$ (in particular $\mu_4 = \varnothing$). Let $\pi \in \mathcal{P}_{\boldsymbol{\mu}}$ and denote the corresponding subscheme by $Z_\pi \subset \C^4$. Denote by $W:=W_\alpha$ the Laurent polynomial in Definition \ref{def:renorm} (it is in fact a polynomial since $q=-1$). Then $|\pi| = \rk(W)$ is the renormalized volume. We associate to $\pi$ the sign
\begin{equation} \label{eqn:sign2dim}
(-1)^{|\pi| + \mu_\pi}, \quad \mu_\pi := \sum_{(i,i,i,j) \in \pi \atop i < j} \frac{1}{(i!)^3 j!} \frac{\partial^{3i+j} W}{\partial t_1^i \partial t_2^i \partial t_3^i \partial t_4^j} \Bigg|_{t_1=t_2=t_3=t_4=0},
\end{equation}
where $\mu_\pi$ is simply the sum of all coefficients of all terms of the form $t_1^i t_2^i t_3^i t_4^j$ with $(i,i,i,j) \in \pi$ and $i<j$. When all $\boldsymbol{\mu}$ are finite plane partitions and $\mu_4 = \varnothing$, we have $\mu_{\pi} = \sum_{w = (i,i,i,j) \in \pi : i<j} (1  -  |\{\textrm{legs containing } w\}|)$ so \eqref{eqn:sign2dim} is a generalization of \eqref{eqn:sign1dim}.

\begin{remark}
For all the verifications of the $\DT$--$\PT_0$ vertex correspondence in Proposition \ref{prop:verif} for which $Z_{\boldsymbol{\mu}}$ is set theoretically supported on $Z(x_4)$, on the $\DT$ side the correct choice of $\pm \sqrt{\cdot}$ is given by \eqref{eqn:DTvertexKY} with sign \eqref{eqn:sign2dim}. We therefore expect that, when $Z_{\boldsymbol{\mu}}$ is \emph{set theoretically} supported on $Z(x_4)$, the $\DT$ topological vertex with sign choice for which Conjecture \ref{conj:vertexDTPT0} holds is given by the following formula
$$
\mathsf{V}^{\DT}_{\mu_1\mu_2\mu_3\varnothing}(q) = \sum_{\pi} (-1)^{|\pi| + \mu_\pi} [-\widetilde{\mathsf{v}}_{\pi}] \, q^{|\pi|}.
$$
\end{remark}

\subsection{Virtual Lefschetz principle -- Local}

In this section, we study the halved vertex and edge terms $\widetilde{\mathsf{v}}_\alpha$, $\widetilde{\mathsf{e}}_{\alpha\beta}$ of the previous section for the specialization $y = t_4 = (t_1t_2t_3)^{-1}$. For a 2-dimensional scheme $Z$ with torsion filtration $T_0(\O_Z) \subset T_1(Z) \subset \O_Z$ we recall that we denote the corresponding subschemes by $Z \supset Z_0 \supset Z_1$ (Section \ref{sec:fixlocal}). 
\begin{proposition} \label{prop:dimredvertex}
Let $X = \C^4$, $Y = Z(x_4) \cong \C^3$, $S = Z(x_3,x_4) \cong \C^2$. Let $dS \subset Y$ be the $d$ times thickening of $S$ in $Y$.
Let $I^\mdot$ be a $T_X$-equivariant $\PT_q$ pair on $X$ for $q \in \{-1,0\}$ with scheme theoretic support $Z$. We assume $Z_1 = dS$ and $Z_0 \subset Y$. In particular, only the 0-dimensional embedded components of $Z$ may not lie scheme-theoretically inside $Y$.
\begin{itemize}
\item Suppose $Z$ does \emph{not} lie scheme theoretically in $Y$, which can only happen for $q=-1$. Then $\widetilde{\mathsf{v}}_{I^\mdot} |_{y=t_4}$ has negative $T$-fixed term, in particular  
$$
[-\widetilde{\mathsf{v}}_{I^\mdot}] \big|_{y = t_4} = 0.
$$
\item Suppose $Z \subset Y$. Denote by $I_Y^\mdot = [ \O_Y \to F ]$ the corresponding $\PT_q$ pair on $Y$. Then there exists a unique 1-dimensional $T_Y$-equivariant pair $J^\mdot = [\O_Y \to G]$ on $Y$ such that 
$$
I_Y\udot \cong J\udot(-d S),
$$
where $J^\mdot$ is a 1-dimensional $\DT$ pair on $Y$ for $q=-1$ and a 1-dimensional $\PT$ pair on $Y$ for $q=0$. Moreover
$$
[-\widetilde{\mathsf{v}}_{I^\mdot}] \big|_{y = t_4} = [-\mathsf{V}^{\mathrm{3D}}_{J^\mdot}],
$$
where $\mathsf{V}^{\mathrm{3D}}_{J^\mdot}$ denotes the equivariant vertex term of 3-dimensional $\DT$ theory \cite[Sect.~4.9]{MNOP} resp.~$\PT$ theory \cite[Sect.~4.6]{PT2}.
\end{itemize}
\end{proposition}
\begin{proof}
Let $q=-1$ and suppose $Z \subsetneq Y$. Before specializing $y = t_4$, $\widetilde{\mathsf{v}}_{I^\mdot}$ has no positive $T$-fixed term (Lemma \ref{lem:noposTfix}), so we only need to consider terms containing $y$. Using the notation of Definition \ref{def:renorm} (with $Z := Z_\alpha$, $W := W_\alpha$, $W_i := W_{\alpha\beta_i}$, $W_{ij} = Z_{\alpha\beta_i\beta_j}$), these terms are 
$ -y \overline{W}.$
Note that $W$ is a Laurent polynomial in $t_1,t_2,t_3,t_4$ but can have terms with negative coefficients. 
However, since only the 0-dimensional embedded components of $Z$ may lie outside $Y$, any term in $W$ containing $t_4$ has positive coefficient. Moreover, since $Z \subsetneq Y$, one of these terms is $+t_4$. Therefore, after specializing $y=t_4$, we conclude that $\widetilde{\mathsf{v}}_{I^\mdot}$ has a negative $T$-fixed term and $[-\widetilde{\mathsf{v}}_{I\udot}] |_{y = t_4} = 0$. 

Next, we assume $Z \subset Y$. From Proposition \ref{prop:PT0PTcv}, it is clear that there exists a unique 1-dimensional $(\C^*)^3$-equivariant pair $J^\mdot = [\O_Y \to G]$ on $Y$ such that $I\udot \cong J\udot(-dS)$ (see also the proof of Theorem \ref{thm:toricdimred}). Continuing the notation of Definition \ref{def:renorm}, we have
\begin{align*}
\tr_F &= \frac{1+t_3 + \cdots + t_3^{d-1}}{(1-t_1)(1-t_2)}  + \frac{W_1}{1-t_1} + \frac{W_2}{1-t_2} + W, \\
\tr_G &= \frac{W_1 t_3^{-d}}{1-t_1} + \frac{W_2 t_3^{-d}}{1-t_2}  + W t_3^{-d}.
\end{align*}
An explicit calculation shows
$$
\widetilde{\mathsf{v}}_{I\udot} - \mathsf{V}_{J\udot}^{\mathrm{3D}} = (t_4 - y) \overline{W}.
$$
The statement of the proposition follows by setting $y = t_4$.
\end{proof}

The condition that $Z \subset Y$, except for 0-dimensional embedded components, is necessary as shown by the following example.
\begin{example}
Let $Z \subset \C^4$ be defined by the ideal 
$$
I_Z = (x_3,x_3x_4,x_1x_2x_4,x_4^2),
$$
which is a copy of $\C^2$ with an embedded nodal curve. Then $Z$ has an embedded 1-dimensional component which is not contained in $Y = Z(x_4)$. Moreover, $\widetilde{\mathsf{v}}_{Z}|_{y=t_4}$ has positive $T$-fixed term $+1$ and the specialization $[-\widetilde{\mathsf{v}}_{Z}]|_{y=t_4}$ is not defined.
\end{example}

As an immediate consequence of Proposition \ref{prop:dimredvertex}, we obtain a dimensional reduction of the $K$-theoretic $\DT$--$\PT_0$ vertex correspondence of Conjecture \ref{conj:vertexDTPT0} to the 3-dimensional $K$-theoretic $\DT$--$\PT$ (curve) vertex correspondence, which was formulated by Nekrasov-Okounkov \cite[Sect.~2.4.8]{NO} and which generalizes the cohomological version of Pandharipande-Thomas \cite[Conj.~4]{PT2}. The 3-dimensional $K$-theoretic $\DT$--$\PT$ vertex correspondence was recently proved in \cite{KLT}. The 3-dimensional numerical DT/PT vertex correspondence (i.e., for the specialization $t_1t_2t_3=1$) was proved combinatorially by Jenne--Webb--Young \cite{JWY}. More precisely, we obtain the following.
\begin{corollary}
Fix $\mu_1,\mu_2$ plane partitions and let $\mu_3 = \mu_4 = \varnothing$. Suppose $Z_{\boldsymbol{\mu}}$ (Definition \ref{def:partitionnotation}) is scheme theoretically supported in $\C^3 = Z(x_4) \subset \C^4$. Let $\lambda_1$, $\lambda_2$ be the unique finite partitions with the property that $I_{Z_{\boldsymbol{\mu}}/ \C^3} = I_{Z_{\boldsymbol{\lambda}} / \C^3}(-dS)$, for some $d \geq 0$, where $S = Z(x_3) \subset \C^3$ and $Z_{\boldsymbol{\lambda}} \subset \C^3$ is the Cohen-Macaulay curve determined by $\boldsymbol{\lambda} = (\lambda_1,\lambda_2,\varnothing)$.
Then there exist choices of $\pm \sqrt{\cdot}$ such that
$$
\frac{\mathsf{V}^{\DT}_{\mu_1\mu_2\varnothing\varnothing}(q)}{\mathsf{V}^{\DT}_{\varnothing\varnothing\varnothing\varnothing}(q)} \Big|_{y = t_4} = \mathsf{V}^{\PT_0}_{\mu_1\mu_2\varnothing\varnothing}(q)|_{y = t_4}
$$ 
is equivalent to 2-leg 3-dimensional $K$-theoretic $\DT$--$\PT$ vertex correspondence
$$
\frac{\mathsf{V}^{\DT}_{\lambda_1\lambda_2\varnothing}(q)}{\mathsf{V}^{\DT}_{\varnothing\varnothing\varnothing}(q)} = \mathsf{V}^{\PT}_{\lambda_1\lambda_2\varnothing}(q).
$$ 
\end{corollary}

The final result of this section concerns a \emph{$\PT_0$--$\PT_1$ edge correspondence} for the specialization $y = t_4$. Recall that the halved edge term $\widetilde{\mathsf{e}}_{I^\mdot}$ in Definition \ref{def:halfedge} depends on integers $m,m',m''$, i.e., the weights of the three normal directions \eqref{weightsN}. 
\begin{proposition}
Let $X = \C^* \times \C^3$, $Y = Z(x_4) \cong \C^* \times  \C^2$, $S = Z(x_3,x_4) \cong \C^* \times \C$. Let $dS \subset Y$ be the $d$ times thickening of $S$ in $Y$.
Let $I^\mdot$ be a $T_X$-equivariant $\PT_q$ pair on $X$ for $q \in \{0,1\}$ and denote its scheme theoretic support by $Z$.\footnote{Recall that on $X = \C^* \times \C^3$, $T_X$-equivariant $\DT$ and $\PT_0$ pairs are the same.} Assume $Z_1 = dS$ and $m''=0$.
\begin{itemize}
\item Suppose $Z$ does \emph{not} lie scheme theoretically in $Y$, which can only happen for $q=0$, then $\widetilde{\mathsf{e}}_{I^\mdot} |_{y=t_4}$ has negative $T$-fixed term, in particular
$$
[-\widetilde{\mathsf{e}}_{I^\mdot}] \big|_{y = t_4} = 0.
$$
\item Suppose $Z \subset Y$. Denote by $I_Y^\mdot = [ \O_Y \to F ]$ the corresponding $\PT_q$ pair on $Y$ with cokernel $Q$. Then there exists a unique 1-dimensional $T_Y$-equivariant $\DT = \PT$ pair\footnote{Note that on $Y = \C^* \times \C^2$, 1-dimensional $T_Y$-equivariant $\DT$ and $\PT$ pairs coincide.} $J^\mdot = [\O_Y \to \O_C]$ on $Y$ such that 
\begin{align*}
I_Y\udot &\cong I_{C/Y}(-dS), \quad \textrm{for } q=0 \\
\O_C \cdot t_2 t_3  &= \overline{Q(dS)}, \quad\quad \ \ \ \textrm{for } q=1.
\end{align*}
Moreover
$$
[-\widetilde{\mathsf{e}}_{I^\mdot}] \big|_{y = t_4} = [-\mathsf{E}^{\mathrm{3D}}_{J^\mdot}],
$$
where $\mathsf{E}^{\mathrm{3D}}_{J^\mdot}$ denotes the equivariant edge term of 3-dimensional $\DT$ theory from \cite[Sect.~4.9]{MNOP} or, equivalently\footnote{Recall that the equivariant edge terms of 3-dimensional $\DT$ and $\PT$ theory coincide.}, 3-dimensional $\PT$ theory from \cite[Sect.~4.6]{PT2}. In particular, the $\PT_0$ and $\PT_1$ edge terms are equal for $y=t_4$.
\end{itemize}
\end{proposition}
\begin{proof}
Let $q=0$ and $Z \subsetneq Y$. We use the notation of Definition \ref{def:renorm} with $Z_i := Z_{\alpha\beta_i}$, $W_i := W_{\alpha\beta_i}$, $W_{ij} := Z_{\alpha\beta_i\beta_j}$.
Since $\mathsf{e}_{I^\mdot}$ has no positive $T$-fixed terms (Lemma \ref{lem:noposTfix}), any $T$-fixed term resulting from the specialization $y = t_4$ comes from 
\begin{equation} \label{eqn:yWedge}
\frac{t_1}{1-t_1} \Big( y \overline{W}_1 \Big) - \frac{1}{1-t_1} \Big( y \overline{W}_1 |_{(t_2 t_1^{-m}, t_3 t_1^{-m'}, t_4 t_1^{-m''} )} \Big),
\end{equation}
where $m'' = 0$. Note that all terms of $W_1$ are of the form $+t_2^j t_3^k t_4^{l}$ for some $j,k,l \geq 0$ (by our conditions on $Z$). Such a term can only produce a $T$-fixed term under the specialization $y=t_4$ when $j=k=l-1$. Using the fact that $m''=0$, the contribution to \eqref{eqn:yWedge} of such a term is
\[
\frac{1}{1-t_1} y t_2^{-j} t_3^{-j} t_4^{-j-1} (1 - t_1^{jm +jm'-1}) = y t_2^{-j} t_3^{-j} t_4^{-j-1} (-1-t_1-\cdots -t_1^{-2j}),
\]
where we used that $m''=0$ and $m+m'=-2$ by the Calabi-Yau condition. It follows that such a term produces the $T$-fixed term $-1$ under the specialization $y=t_4$. Since we assume $Z \subsetneq Y$, such terms occur.

Next, we assume $Z \subset Y$. It is easy to see from the the local description of fixed loci (Lemma \ref{lem:fixlocaffDT}) that there exists  a unique 1-dimensional $(\C^*)^3$-equivariant $\DT = \PT$ pair $J^\mdot = [\O_Y \to \O_C]$ on $Y$ such that 
\begin{align*}
I_Y\udot &\cong I_{C/Y}(-dS), \quad \textrm{for } q=0 \\
\O_C \cdot t_2 t_3  &= \overline{Q(dS)}, \quad\quad \ \ \ \textrm{for } q=1.
\end{align*}
Let $C := \tr_{\O_C}$.
Then, for $q=0$, from the definitions of the halved edge term and 3-dimensional edge term, we obtain
\begin{align*}
\widetilde{e}_{I\udot} - \mathsf{E}_{J\udot}^{\mathrm{3D}} = &\frac{t_1^{-1}}{1-t_1^{-1}} \Big((1 - yt_1t_2t_3) t_2^{-1}t_3^{-d-1} \overline{C} \Big) - \\
&\frac{1}{1-t_1^{-1}} \Big( \big((1 - yt_1t_2t_3) t_2^{-1}t_3^{-d-1} \overline{C} \big) \Big|_{(t_1^{-1},t_2t_1^{-m},t_3t_1^{-m'})} \Big).
\end{align*}
On the other hand, for $q=1$, we obtain
\begin{align*}
\widetilde{e}_{I\udot} - \mathsf{E}_{J\udot}^{\mathrm{3D}} = &\frac{t_1^{-1}}{1-t_1^{-1}} \Big((1-yt_1t_2t_3) t_3^{-d} C \Big) - \\
&\frac{1}{1-t_1^{-1}} \Big( \big( (1-yt_1t_2t_3) t_3^{-d} C \big) \Big|_{(t_1^{-1},t_2t_1^{-m},t_3t_1^{-m'})} \Big).
\end{align*}
In both cases, we see that $[-\widetilde{e}_{I\udot}]|_{y=t_4}$ is well-defined and equal to $[-\mathsf{E}_{J\udot}^{\mathrm{3D}}]$.
\end{proof}

The specialization $y = t_4$ provides a dimensional reduction, which we think of as a \emph{local} virtual Lefschetz principle as opposed to the global virtual Lefschetz principles studied in Sections \ref{sec:Lefglob} and \ref{sec:Leflocal}.

\subsection{Local projective plane} \label{sec:localP2P3}

Consider $X = \mathrm{Tot}_{\PP^2}(\O_{\PP^2}(-2) \oplus \O_{\PP^2}(-1))$. Note that
\[
v = (0,0,\gamma,\beta,n - \gamma \cdot \td_2(X)) \in H_c^*(X,\Q)
\]
has the property that $\gamma = d [\PP^2]$ and $\beta = m [\PP^1]$. We therefore write $\gamma = d$ and $\beta = m$. Then, for $q \in \{-1,0,1\}$, the virtual dimension of $\PTqvX$ equals
\[
\vd = n - d^2.
\]
For $q=1$, and fixed $d,m$, there are only finitely many $n$ such that $\vd \geq 0$ \cite[Prop.~2.14]{BKP}.

For a fixed $d$, we are interested in the generating series
\begin{align*}
\mathsf{G}^{\PT_q}_d(q,Q,y) := \sum_{m,n} \langle\!\langle \O_X \rangle\!\rangle_{d,m,n}^{\PT_q} Q^m q^n,
\end{align*}
for $q \in \{0,1\}$. By Proposition \ref{prop:localFanos}, the moduli spaces $\PTqvX$ are projective. Hence the specialization $t_1=t_2=t_3=t_4 = 1$ is well-defined and we will always make this specialization.

We use the description of fixed loci (Section \ref{sec:fixloc}) and the vertex formalism (Section \ref{sec:redist}) to calculate part of these generating series. We focus mostly on the calculation of $\PT_1$ invariants. Since $X$ is the canonical bundle over a smooth quasi-projective toric 3-fold, we can in fact use Section \ref{sec:local3fold}, i.e., Proposition \ref{prop:DTsumsqrtKY}.

\begin{remark}
Since $X$ is a local surface, Theorem \ref{prop:fixlocmain1}(3) implies that, for $q=0$, $\PTqvX^{T_X} = \PTqvX^T$ is 0-dimensional and reduced. For $q=1$, the calculations in this section depend on Conjecture \ref{conj:0dimTXfix}.
\end{remark}

We used an implementation into a Maple program to show that there are signs such that the following formula holds:
\begin{align*}
[y]^{-1} \mathsf{G}^{\PT_0}_1 = &\, q Q^{\frac{3}{2}} + \Big( q^3 + (y^{\frac{3}{2}}+y^{-\frac{3}{2}}) q^4 +(y^{2}+y^{-2}) q^5 + (y^{\frac{5}{2}}+y^{-\frac{5}{2}}) q^6 + \cdots \Big)Q^{\frac{5}{2}} + \\
&\Big( 0 \cdot q^4 + 0 \cdot q^5 + \cdots \Big) Q^\frac{7}{2} +  \cdots\,,
\end{align*}
where $\cdots$ means higher order terms in $q$ or $Q$. The term with $Q^{\frac{3}{2}}$ has no embedded curve, i.e., $\PT_0 = \PT_1$. The terms with $Q^{\frac{5}{2}}$ have one embedded $\PP^1$. The terms with $Q^{\frac{7}{2}}$ have two embedded $\PP^1$s. Based on this, it is reasonable to expect that for $m = \frac{5}{2}$ and any $n > 3$, we have
\begin{equation} \label{eqn:localPT0eqn}
[y]^{-1} \langle\!\langle \O_X \rangle\!\rangle_{1,\frac{5}{2},n}^{\PT_0} = y^{\frac{n-1}{2}}+y^{-\frac{n-1}{2}}.
\end{equation}
Equation \eqref{eqn:localPT0eqn} will be proved using different methods a sequel to this paper. We note that the two zeroes at the end are non-trivial. E.g., the term $0 \cdot q^5 Q^{\frac{7}{2}}$ involves the summation of 15 fixed points which, individually, have non-zero contributions. 

Next, we used the vertex formalism to show that there are signs such that:
\begin{align} 
\begin{split} \label{eqn:localPT1eqn}
\mathsf{G}^{\PT_1}_1 = &[y]qQ^{\frac{3}{2}}+\Big( 2[y]q+[y^2]q^2+[y^3]q^3 \Big)Q^{\frac{5}{2}} + \\
&\Big(...q+...q^2+...q^3+2[y^4]q^4+[y^5]q^5+[y^6]q^6 \Big) Q^\frac{7}{2} + \cdots,
\end{split}
\end{align}
where $...q^i$ means that we did not determine the coefficient. In a sequel to this paper, we determine $\mathsf{G}^{\PT_1}_1$ up to all orders (essentially by using the nested description of Proposition \ref{cor:PT1smsupp}). We present the answer. Let $S$ be a smooth projective surface and $\cZ \subset S \times \Hilb^m(S)$ the universal subscheme. Denote the corresponding ideal sheaf by $I_{\cZ}$ and let $\pi : S \times \Hilb^m(S) \to \Hilb^m(S)$ be the projection. For any $L \in \Pic(S)$, consider the complex
\[
T(L) := R\hom_{\pi}(I_{\cZ},I_{\cZ} \otimes L)_0[1],
\]
which is in fact a vector bundle concentrated in degree 0.
Then Carlsson-Okounkov \cite{CO} showed 
\[
\sum_{m=0}^{\infty} \int_{\Hilb^m(S)} c_{2m}(T(L)) q^m = \Bigg( \prod_{m=1}^{\infty} \frac{1}{1-q^m} \Bigg)^{e(S) + L(L-K_S)},
\]
where $e(S)$ is the topological Euler characteristic of $S$. Denote the coefficient of $q^m$ on the left hand side by $\mathsf{CO}(S,L)_m$. Then in a sequel we show that 
\begin{align*} 
\langle\!\langle \O_X \rangle\!\rangle_{1,\frac{3}{2}+\delta,n}^{\PT_1} &= \mathsf{CO}_{\chi(\O_{\PP^2}(\delta))-n}(\PP^2,\O(-1)) \cdot [y^n] \\
&= \mathsf{CO}_{\chi(\O_{\PP^2}(\delta))-n}(\PP^2,\O(-1)) [y] (y^{\frac{n-1}{2}} + \cdots + y^{-\frac{n-1}{2}}),
\end{align*}
for all $\delta \geq 0$ and $1 \leq n \leq \chi(\O_{\PP^2}(\delta))$. This formula is consistent with the toric calculations.

The proofs of \eqref{eqn:localPT0eqn} and \eqref{eqn:localPT1eqn} in the sequels strongly depend on the fact that, so far, we considered the \emph{irreducible} class $\gamma = [\PP^2]$. The main strength of the vertex formalism is that it allows us to consider $\gamma = d [\PP^2]$ with $d>1$. 

For $d=2$, we have shown that 
\begin{align*}
[y]^{-4} \mathsf{G}^{\PT_0}_2 = q^4 Q^4 + \Big( 2(y^{\frac{1}{2}} + y^{-\frac{1}{2}})^3 q^7 + \cdots \Big) Q^5 + \cdots.
\end{align*}
The term $q^7 Q^5$ involves contributions of the zero section $\PP^2 \subset X$ thickened in the direction of $\O(-1)$ as well as in the direction of $\O(-2)$. In the former case, it involves contributions with an embedded curve thickened in the direction of $\O(-1)$ as well as in the direction of $\O(-2)$. Moreover, the signs in Proposition \ref{prop:DTsumsqrtKY} for fixed points scheme theoretically inside the zero section $\mathrm{Tot}(\O_{\PP^2}(-1)) \subset X$ and those which are not are \emph{opposite}. Next, we determined the following terms of the $\PT_1$ generating series for $d=2$
\begin{align*}
[y]^{-4} \mathsf{G}^{\PT_1}_2 = &\, q^4 Q^{4}+\Big( 4q^4 +2(y^{\frac{1}{2}}+y^{-\frac{1}{2}})q^5+2(y+y^{-1})q^6 + \\
&2(y^{\frac{1}{2}} + y^{-\frac{1}{2}})(2y  + 1 + 2 y^{-1} ) q^7 \Big)Q^5 + \\
&\Big( ...q^4 +\cdots + 2(y^{\frac{1}{2}}+y^{-\frac{1}{2}})(y+y^{-1})(5y+1+5y^{-1}) q^{9}  +  \\
&(7 y^{3}  +12 y^2 + 15 y + 16 + 15 y^{-1} + 12 y^{-2} + 7 y^{-3}) q^{10}  + \\
&2(y^{\frac{1}{2}} + y^{-\frac{1}{2}})(4 y^{3}  + 3 y^2 + 6 y + 4 + 6 y^{-1} + 3 y^{-2} + 4 y^{-3}) q^{11}   \Big) Q^6 + \cdots.
\end{align*}
The calculation involves thickenings of the zero section $\PP^2$ in the $\O(-1)$ as well as the $\O(-2)$ direction. The computation of the coefficient of $q^9Q^6$ involves 48 fixed points. For $d=3$ we found
\begin{align*}
[y]^{-9} \mathsf{G}^{\PT_1}_3 = &\, \Big(0q^{9} + 3(y^{\frac{1}{2}}+y^{-\frac{1}{2}})q^{10} \Big)Q^{\frac{15}{2}}+ \Big(2q^{9} +18(y^{\frac{1}{2}}+y^{-\frac{1}{2}})q^{10}+ \\
&9(y^{\frac{1}{2}}+y^{-\frac{1}{2}})^2 q^{11}+9(y^{\frac{1}{2}}+y^{-\frac{1}{2}})(y+y^{-1})q^{12}+ \\
&3(3y^2+2y-3+2y^{-1}+3y^{-2}) q^{13}+ \\
&3(y^{\frac{1}{2}}+y^{-\frac{1}{2}})(7y^2+8y+12+8y^{-1}+7y^{-2}) q^{14} \Big)Q^{\frac{17}{2}} + \cdots.
\end{align*}
The calculation of the coefficient of $q^{9} Q^{\frac{17}{2}}$ involves 51 fixed points. For $d=4$, we obtained:
\begin{align*}
[y]^{-16} \mathsf{G}^{\PT_1}_4 = &\, q^{16}  Q^{11}+\Big(8q^{16}+4(y^{\frac{1}{2}}+y^{-\frac{1}{2}}) q^{17}+4(y+y^{-1}) q^{18}+ \\
&4(y^{\frac{3}{2}}+y^{-\frac{3}{2}}) q^{19} +10(y^2+8y+12+8y^{-1}+y^{-2})q^{20} \Big)Q^{12} + \cdots.
\end{align*}
Generally, the poles in the equivariant parameters $t_i$ in the calculations only cancel when the contributions of the fixed points are summed with very specific signs, which we expect to be induced from global canonical orientations. We observe that all terms calculated in this section are invariant under $y \to 1/y$ and, more generally, are consistent with the structure described in Proposition \ref{prop:palin}. This is evidence for the expectation that the sign choices found in this section are induced from global orientations. 

We note that all coefficients appearing in the calculations of this section are \emph{integer}. This is no longer true when one does not include the $K$-theoretic tautological insertion (Definition \ref{def:Ktheorinv}). E.g.~for $d=1$, $m = \frac{7}{2}$, $n=5$, we have $\chi(\curly P_v^{(1)}(X), \widehat{\O}^{\vir}) = \frac{3}{128}$ (for the correct sign choice).

\noindent {\tt{y.bae@uu.nl, m.kool1@uu.nl, hyeonjunpark@kias.re.kr}}
\end{document}

%% file: fig0a.tex
% Plane partition
% Author: Jang Soo Kim
%\documentclass{minimal}
%\usepackage{tikz}
% Three counters
\newcounter{x}
\newcounter{y}
\newcounter{z}

% The angles of x,y,z-axes
\newcommand\xaxis{210}
\newcommand\yaxis{-30}
\newcommand\zaxis{90}

% The top side of a cube
\newcommand\topside[3]{
  \fill[fill=cyan, draw=black,shift={(\xaxis:#1)},shift={(\yaxis:#2)},
  shift={(\zaxis:#3)}] (0,0) -- (30:1) -- (0,1) --(150:1)--(0,0);
}

% The left side of a cube
\newcommand\leftside[3]{
  \fill[fill=cyan, draw=black,shift={(\xaxis:#1)},shift={(\yaxis:#2)},
  shift={(\zaxis:#3)}] (0,0) -- (0,-1) -- (210:1) --(150:1)--(0,0);
}

% The right side of a cube
\newcommand\rightside[3]{
  \fill[fill=cyan, draw=black,shift={(\xaxis:#1)},shift={(\yaxis:#2)},
  shift={(\zaxis:#3)}] (0,0) -- (30:1) -- (-30:1) --(0,-1)--(0,0);
}

% The cube 
\newcommand\cube[3]{
  \topside{#1}{#2}{#3} \leftside{#1}{#2}{#3} \rightside{#1}{#2}{#3}
}

% Definition of \planepartition
% To draw the following plane partition, just write \planepartition{ {a, b, c}, {d,e} }.
%  a b c
%  d e
\newcommand\planepartition[1]{
 \setcounter{x}{-1}
  \foreach \a in {#1} {
    \addtocounter{x}{1}
    \setcounter{y}{-1}
    \foreach \b in \a {
      \addtocounter{y}{1}
      \setcounter{z}{-1}
      \foreach \c in {1,...,\b} {
        \addtocounter{z}{1}
        \cube{\value{x}}{\value{y}}{\value{z}}
      }
    }
  }
}

%\begin{document} 

\begin{tikzpicture}[scale=0.25]

\newcommand{\boxt}[3]{ 
\draw[fill=violet, draw=black,shift={(\xaxis:#1)},shift={(\yaxis:#2)},
  shift={(\zaxis:#3)}] (0,0) -- (30:1) -- (0,1) --(150:1)--(0,0);
}

\newcommand{\boxl}[3]{ 
\draw[fill=violet, draw=black,shift={(\xaxis:#1)},shift={(\yaxis:#2)},
  shift={(\zaxis:#3)}] (0,0) -- (0,-1) -- (210:1) --(150:1)--(0,0);
}

\newcommand{\boxr}[3]{ 
\draw[fill=violet, draw=black,shift={(\xaxis:#1)},shift={(\yaxis:#2)},
  shift={(\zaxis:#3)}] (0,0) -- (30:1) -- (-30:1) --(0,-1)--(0,0);
}

\newcommand{\boxxt}[3]{ 
\draw[fill=cyan, draw=black,shift={(\xaxis:#1)},shift={(\yaxis:#2)},
  shift={(\zaxis:#3)}] (0,0) -- (30:1) -- (0,1) --(150:1)--(0,0);
}

\newcommand{\boxxl}[3]{ 
\draw[fill=cyan, draw=black,shift={(\xaxis:#1)},shift={(\yaxis:#2)},
  shift={(\zaxis:#3)}] (0,0) -- (0,-1) -- (210:1) --(150:1)--(0,0);
}

\newcommand{\boxxr}[3]{ 
\draw[fill=cyan, draw=black,shift={(\xaxis:#1)},shift={(\yaxis:#2)},
  shift={(\zaxis:#3)}] (0,0) -- (30:1) -- (-30:1) --(0,-1)--(0,0);
}

\planepartition{{10,1,1,1,1,1,1,1,1,1},{1,1,1,1,1,1,1,1,1,1},{1,1,1,1,1,1,1,1,1,1},{1,1,1,1,1,1,1,1,1,1},{1,1,1,1,1,1,1,1,1,1},{1,1,1,1,1,1,1,1,1,1},{1,1,1,1,1,1,1,1,1,1},{1,1,1,1,1,1,1,1,1,1},{1,1,1,1,1,1,1,1,1,1},{1,1,1,1,1,1,1,1,1,1}}

\boxt{0}{0}{9}
\boxl{0}{0}{9}
\boxr{0}{0}{9}
\boxl{0}{0}{8}
\boxr{0}{0}{8}
\boxl{0}{0}{7}
\boxr{0}{0}{7}
\boxl{0}{0}{6}
\boxr{0}{0}{6}
\boxl{0}{0}{5}
\boxr{0}{0}{5}
\boxl{0}{0}{4}
\boxr{0}{0}{4}
\boxl{0}{0}{3}
\boxr{0}{0}{3}
\boxl{0}{0}{2}
\boxr{0}{0}{2}
\boxl{0}{0}{1}
\boxr{0}{0}{1}

\draw [thick,->] (0,10)-- (0,12);
\draw [thick,->] (5*1.74,-5*1.74*0.57735)-- (10,-10*0.57735);
\draw [thick,->]  (-8.6,-8.6*0.57735)-- (-10,-10*0.57735);

\node [right] at (0,12.5) {$x_3$};
\node [right] at (9.7,-9.7*0.57735) {$x_2$};
\node [left] at  (-9.7,-9.7*0.57735) {$x_1$};
\node [above] at  (-9,-7*0.57735) {$ $};
\node [right] at (2,10) {$ $};
\node  at (6,6*0.57735) {$ $};
\node  at (-6,6*0.57735) {$ $};
\node  at (0,-6) {$ $};

%\shade [ball color=gray] (-2.6,-.6) circle [radius=0.4cm];
%\shade [ball color=gray] (1.8,-1.2) circle [radius=0.4cm];
%\shade [ball color=gray] (-0.85,-1.7) circle [radius=0.4cm];

\end{tikzpicture}

%\end{document} 

%% file: fig0b.tex
% Plane partition
% Author: Jang Soo Kim
%\documentclass{minimal}
%\usepackage{tikz}
% Three counters
%\newcounter{x}
%\newcounter{y}
%\newcounter{z}

% The angles of x,y,z-axes
\newcommand\xaxis{210}
\newcommand\yaxis{-30}
\newcommand\zaxis{90}

% The top side of a cube
\newcommand\topside[3]{
  \fill[fill=cyan, draw=black,shift={(\xaxis:#1)},shift={(\yaxis:#2)},
  shift={(\zaxis:#3)}] (0,0) -- (30:1) -- (0,1) --(150:1)--(0,0);
}

% The left side of a cube
\newcommand\leftside[3]{
  \fill[fill=cyan, draw=black,shift={(\xaxis:#1)},shift={(\yaxis:#2)},
  shift={(\zaxis:#3)}] (0,0) -- (0,-1) -- (210:1) --(150:1)--(0,0);
}

% The right side of a cube
\newcommand\rightside[3]{
  \fill[fill=cyan, draw=black,shift={(\xaxis:#1)},shift={(\yaxis:#2)},
  shift={(\zaxis:#3)}] (0,0) -- (30:1) -- (-30:1) --(0,-1)--(0,0);
}

% The cube 
\newcommand\cube[3]{
  \topside{#1}{#2}{#3} \leftside{#1}{#2}{#3} \rightside{#1}{#2}{#3}
}

% Definition of \planepartition
% To draw the following plane partition, just write \planepartition{ {a, b, c}, {d,e} }.
%  a b c
%  d e
\newcommand\planepartition[1]{
 \setcounter{x}{-1}
  \foreach \a in {#1} {
    \addtocounter{x}{1}
    \setcounter{y}{-1}
    \foreach \b in \a {
      \addtocounter{y}{1}
      \setcounter{z}{-1}
      \foreach \c in {1,...,\b} {
        \addtocounter{z}{1}
        \cube{\value{x}}{\value{y}}{\value{z}}
      }
    }
  }
}

%\begin{document} 

\begin{tikzpicture}[scale=0.25]

\newcommand{\boxt}[3]{ 
\draw[fill=violet, draw=black,shift={(\xaxis:#1)},shift={(\yaxis:#2)},
  shift={(\zaxis:#3)}] (0,0) -- (30:1) -- (0,1) --(150:1)--(0,0);
}

\newcommand{\boxl}[3]{ 
\draw[fill=violet, draw=black,shift={(\xaxis:#1)},shift={(\yaxis:#2)},
  shift={(\zaxis:#3)}] (0,0) -- (0,-1) -- (210:1) --(150:1)--(0,0);
}

\newcommand{\boxr}[3]{ 
\draw[fill=violet, draw=black,shift={(\xaxis:#1)},shift={(\yaxis:#2)},
  shift={(\zaxis:#3)}] (0,0) -- (30:1) -- (-30:1) --(0,-1)--(0,0);
}

\newcommand{\boxxt}[3]{ 
\draw[fill=cyan, draw=black,shift={(\xaxis:#1)},shift={(\yaxis:#2)},
  shift={(\zaxis:#3)}] (0,0) -- (30:1) -- (0,1) --(150:1)--(0,0);
}

\newcommand{\boxxl}[3]{ 
\draw[fill=cyan, draw=black,shift={(\xaxis:#1)},shift={(\yaxis:#2)},
  shift={(\zaxis:#3)}] (0,0) -- (0,-1) -- (210:1) --(150:1)--(0,0);
}

\newcommand{\boxxr}[3]{ 
\draw[fill=cyan, draw=black,shift={(\xaxis:#1)},shift={(\yaxis:#2)},
  shift={(\zaxis:#3)}] (0,0) -- (30:1) -- (-30:1) --(0,-1)--(0,0);
}

\planepartition{{10,1,1,1,1,1,1,1,1,1},{1,1,1,1,1,1,1,1,1,1},{1,1,1,1,1,1,1,1,1,1},{1,1,1,1,1,1,1,1,1,1},{1,1,1,1,1,1,1,1,1,1},{1,1,1,1,1,1,1,1,1,1},{1,1,1,1,1,1,1,1,1,1},{1,1,1,1,1,1,1,1,1,1},{1,1,1,1,1,1,1,1,1,1},{1,1,1,1,1,1,1,1,1,1}}

\boxt{0}{0}{9}
\boxl{0}{0}{9}
\boxr{0}{0}{9}
\boxl{0}{0}{8}
\boxr{0}{0}{8}
\boxl{0}{0}{7}
\boxr{0}{0}{7}
\boxl{0}{0}{6}
\boxr{0}{0}{6}
\boxl{0}{0}{5}
\boxr{0}{0}{5}
\boxl{0}{0}{4}
\boxr{0}{0}{4}
\boxl{0}{0}{3}
\boxr{0}{0}{3}
\boxl{0}{0}{2}
\boxr{0}{0}{2}
\boxl{0}{0}{1}
\boxr{0}{0}{1}

\boxxr{1}{0}{1}
\boxxl{1}{0}{2}
\boxxr{1}{0}{2}
\boxxt{1}{0}{3}
\boxxl{1}{0}{3}
\boxxr{1}{0}{3}
\boxxt{2}{0}{1}
\boxxl{2}{0}{1}
\boxxr{2}{0}{1}

\boxxt{0}{1}{1}
\boxxl{0}{1}{1}
\boxxr{0}{1}{1}

\draw [thick,->] (0,10)-- (0,12);
\draw [thick,->] (5*1.74,-5*1.74*0.57735)-- (10,-10*0.57735);
\draw [thick,->]  (-8.6,-8.6*0.57735)-- (-10,-10*0.57735);

\node [right] at (0,12.5) {$x_3$};
\node [right] at (9.7,-9.7*0.57735) {$x_2$};
\node [left] at  (-9.7,-9.7*0.57735) {$x_1$};
\node [above] at  (-9,-7*0.57735) {$ $};
\node [right] at (2,10) {$ $};
\node  at (6,6*0.57735) {$ $};
\node  at (-6,6*0.57735) {$ $};
\node  at (0,-6) {$ $};

%\shade [ball color=gray] (-2.6,-.6) circle [radius=0.4cm];
%\shade [ball color=gray] (1.8,-1.2) circle [radius=0.4cm];
%\shade [ball color=gray] (-0.85,-1.7) circle [radius=0.4cm];

\end{tikzpicture}

%\end{document} 

%% file: fig1.tex
% Plane partition
% Author: Jang Soo Kim
%\documentclass{minimal}
%\usepackage{tikz}
% Three counters
%\newcounter{x}
%\newcounter{y}
%\newcounter{z}

% The angles of x,y,z-axes
\newcommand\xaxis{210}
\newcommand\yaxis{-30}
\newcommand\zaxis{90}

% The top side of a cube
\newcommand\topside[3]{
  \fill[fill=cyan, draw=black,shift={(\xaxis:#1)},shift={(\yaxis:#2)},
  shift={(\zaxis:#3)}] (0,0) -- (30:1) -- (0,1) --(150:1)--(0,0);
}

% The left side of a cube
\newcommand\leftside[3]{
  \fill[fill=cyan, draw=black,shift={(\xaxis:#1)},shift={(\yaxis:#2)},
  shift={(\zaxis:#3)}] (0,0) -- (0,-1) -- (210:1) --(150:1)--(0,0);
}

% The right side of a cube
\newcommand\rightside[3]{
  \fill[fill=cyan, draw=black,shift={(\xaxis:#1)},shift={(\yaxis:#2)},
  shift={(\zaxis:#3)}] (0,0) -- (30:1) -- (-30:1) --(0,-1)--(0,0);
}

% The cube 
\newcommand\cube[3]{
  \topside{#1}{#2}{#3} \leftside{#1}{#2}{#3} \rightside{#1}{#2}{#3}
}

% Definition of \planepartition
% To draw the following plane partition, just write \planepartition{ {a, b, c}, {d,e} }.
%  a b c
%  d e
\newcommand\planepartition[1]{
 \setcounter{x}{-1}
  \foreach \a in {#1} {
    \addtocounter{x}{1}
    \setcounter{y}{-1}
    \foreach \b in \a {
      \addtocounter{y}{1}
      \setcounter{z}{-1}
      \foreach \c in {1,...,\b} {
        \addtocounter{z}{1}
        \cube{\value{x}}{\value{y}}{\value{z}}
      }
    }
  }
}

%\begin{document} 

\begin{tikzpicture}[scale=0.25]

\newcommand{\boxt}[3]{ 
\draw[fill=violet, draw=black,shift={(\xaxis:#1)},shift={(\yaxis:#2)},
  shift={(\zaxis:#3)}] (0,0) -- (30:1) -- (0,1) --(150:1)--(0,0);
}

\newcommand{\boxl}[3]{ 
\draw[fill=violet, draw=black,shift={(\xaxis:#1)},shift={(\yaxis:#2)},
  shift={(\zaxis:#3)}] (0,0) -- (0,-1) -- (210:1) --(150:1)--(0,0);
}

\newcommand{\boxr}[3]{ 
\draw[fill=violet, draw=black,shift={(\xaxis:#1)},shift={(\yaxis:#2)},
  shift={(\zaxis:#3)}] (0,0) -- (30:1) -- (-30:1) --(0,-1)--(0,0);
}

\planepartition{{11,11,2,2,2,2,2,2,2,2,2},{11,11,2,2,2,2,2,2,2,2,2},{2,2,2,2,2,2,2,2,2,2,2},{2,2,2,2,2,2,2,2,2,2,2},{2,2,2,2,2,2,2,2,2,2,2},{2,2,2,2,2,2,2,2,2,2,2},{2,2,2,2,2,2,2,2,2,2,2},{2,2,2,2,2,2,2,2,2,2,2},{2,2,2,2,2,2,2,2,2,2,2},{2,2,2,2,2,2,2,2,2,2,2},{2,2,2,2,2,2,2,2,2,2,2}}

\boxt{0}{0}{10}
\boxt{0}{1}{10}
\boxt{1}{0}{10}
\boxt{1}{1}{10}
\boxl{1}{0}{2}
\boxl{1}{1}{2}
\boxl{1}{0}{3}
\boxl{1}{1}{3}
\boxl{1}{1}{4}
\boxl{1}{0}{4}
\boxl{1}{1}{5}
\boxl{1}{0}{5}
\boxl{1}{1}{6}
\boxl{1}{0}{6}
\boxl{1}{1}{7}
\boxl{1}{0}{7}
\boxl{1}{1}{8}
\boxl{1}{0}{8}
\boxl{1}{1}{9}
\boxl{1}{0}{9}
\boxl{1}{1}{10}
\boxl{1}{0}{10}
\boxr{0}{1}{2}
\boxr{1}{1}{2}
\boxr{1}{1}{3}
\boxr{0}{1}{3}
\boxr{1}{1}{4}
\boxr{0}{1}{4}
\boxr{1}{1}{5}
\boxr{0}{1}{5}
\boxr{1}{1}{6}
\boxr{0}{1}{6}
\boxr{1}{1}{7}
\boxr{0}{1}{7}
\boxr{1}{1}{8}
\boxr{0}{1}{8}
\boxr{1}{1}{9}
\boxr{0}{1}{9}
\boxr{1}{1}{10}
\boxr{0}{1}{10}

\draw [thick,->] (0,11)-- (0,13);
\draw [thick,->] (5.5*1.74,-5.5*1.74*0.57735)-- (11,-11*0.57735);
\draw [thick,->]  (-9.5,-9.5*0.57735)-- (-11,-11*0.57735);

\node [right] at (0,13.5) {$x_3$};
\node [right] at (11,-11*0.57735) {$x_2$};
\node [left] at  (-11,-11*0.57735) {$x_1$};
\node [above] at  (-9,-7*0.57735) {$ $};
\node [right] at (2,10) {$ $};
\node  at (6,6*0.57735) {$ $};
\node  at (-6,6*0.57735) {$ $};
\node  at (0,-6) {$ $};

%\shade [ball color=gray] (-2.6,-.6) circle [radius=0.4cm];
%\shade [ball color=gray] (1.8,-1.2) circle [radius=0.4cm];
%\shade [ball color=gray] (-0.85,-1.7) circle [radius=0.4cm];

\end{tikzpicture}

%\end{document} 

%% file: fig1a.tex
% Plane partition
% Author: Jang Soo Kim
%\documentclass{minimal}
%\usepackage{tikz}
% Three counters
%\newcounter{x}
%\newcounter{y}
%\newcounter{z}

% The angles of x,y,z-axes
\newcommand\xaxis{210}
\newcommand\yaxis{-30}
\newcommand\zaxis{90}

% The top side of a cube
\newcommand\topside[3]{
  \fill[fill=cyan, draw=black,shift={(\xaxis:#1)},shift={(\yaxis:#2)},
  shift={(\zaxis:#3)}] (0,0) -- (30:1) -- (0,1) --(150:1)--(0,0);
}

% The left side of a cube
\newcommand\leftside[3]{
  \fill[fill=cyan, draw=black,shift={(\xaxis:#1)},shift={(\yaxis:#2)},
  shift={(\zaxis:#3)}] (0,0) -- (0,-1) -- (210:1) --(150:1)--(0,0);
}

% The right side of a cube
\newcommand\rightside[3]{
  \fill[fill=cyan, draw=black,shift={(\xaxis:#1)},shift={(\yaxis:#2)},
  shift={(\zaxis:#3)}] (0,0) -- (30:1) -- (-30:1) --(0,-1)--(0,0);
}

% The cube 
\newcommand\cube[3]{
  \topside{#1}{#2}{#3} \leftside{#1}{#2}{#3} \rightside{#1}{#2}{#3}
}

% Definition of \planepartition
% To draw the following plane partition, just write \planepartition{ {a, b, c}, {d,e} }.
%  a b c
%  d e
\newcommand\planepartition[1]{
 \setcounter{x}{-1}
  \foreach \a in {#1} {
    \addtocounter{x}{1}
    \setcounter{y}{-1}
    \foreach \b in \a {
      \addtocounter{y}{1}
      \setcounter{z}{-1}
      \foreach \c in {1,...,\b} {
        \addtocounter{z}{1}
        \cube{\value{x}}{\value{y}}{\value{z}}
      }
    }
  }
}

%\begin{document} 

\begin{tikzpicture}[scale=0.25]

\newcommand{\boxt}[3]{ 
\draw[fill=violet, draw=black,shift={(\xaxis:#1)},shift={(\yaxis:#2)},
  shift={(\zaxis:#3)}] (0,0) -- (30:1) -- (0,1) --(150:1)--(0,0);
}

\newcommand{\boxl}[3]{ 
\draw[fill=violet, draw=black,shift={(\xaxis:#1)},shift={(\yaxis:#2)},
  shift={(\zaxis:#3)}] (0,0) -- (0,-1) -- (210:1) --(150:1)--(0,0);
}

\newcommand{\boxr}[3]{ 
\draw[fill=violet, draw=black,shift={(\xaxis:#1)},shift={(\yaxis:#2)},
  shift={(\zaxis:#3)}] (0,0) -- (30:1) -- (-30:1) --(0,-1)--(0,0);
}

\newcommand{\boxxt}[3]{ 
\draw[fill=cyan, draw=black,shift={(\xaxis:#1)},shift={(\yaxis:#2)},
  shift={(\zaxis:#3)}] (0,0) -- (30:1) -- (0,1) --(150:1)--(0,0);
}

\newcommand{\boxxl}[3]{ 
\draw[fill=cyan, draw=black,shift={(\xaxis:#1)},shift={(\yaxis:#2)},
  shift={(\zaxis:#3)}] (0,0) -- (0,-1) -- (210:1) --(150:1)--(0,0);
}

\newcommand{\boxxr}[3]{ 
\draw[fill=cyan, draw=black,shift={(\xaxis:#1)},shift={(\yaxis:#2)},
  shift={(\zaxis:#3)}] (0,0) -- (30:1) -- (-30:1) --(0,-1)--(0,0);
}

\planepartition{{11,11,2,2,2,2,2,2,2,2,2},{11,11,2,2,2,2,2,2,2,2,2},{2,2,2,2,2,2,2,2,2,2,2},{2,2,2,2,2,2,2,2,2,2,2},{2,2,2,2,2,2,2,2,2,2,2},{2,2,2,2,2,2,2,2,2,2,2},{2,2,2,2,2,2,2,2,2,2,2},{2,2,2,2,2,2,2,2,2,2,2},{2,2,2,2,2,2,2,2,2,2,2},{2,2,2,2,2,2,2,2,2,2,2},{2,2,2,2,2,2,2,2,2,2,2}}
\boxt{0}{0}{10}
\boxt{0}{1}{10}
\boxt{1}{0}{10}
\boxt{1}{1}{10}
\boxl{1}{0}{2}
\boxl{1}{1}{2}
\boxl{1}{0}{3}
\boxl{1}{1}{3}
\boxl{1}{1}{4}
\boxl{1}{0}{4}
\boxl{1}{1}{5}
\boxl{1}{0}{5}
\boxl{1}{1}{6}
\boxl{1}{0}{6}
\boxl{1}{1}{7}
\boxl{1}{0}{7}
\boxl{1}{1}{8}
\boxl{1}{0}{8}
\boxl{1}{1}{9}
\boxl{1}{0}{9}
\boxl{1}{1}{10}
\boxl{1}{0}{10}
\boxr{0}{1}{2}
\boxr{1}{1}{2}
\boxr{1}{1}{3}
\boxr{0}{1}{3}
\boxr{1}{1}{4}
\boxr{0}{1}{4}
\boxr{1}{1}{5}
\boxr{0}{1}{5}
\boxr{1}{1}{6}
\boxr{0}{1}{6}
\boxr{1}{1}{7}
\boxr{0}{1}{7}
\boxr{1}{1}{8}
\boxr{0}{1}{8}
\boxr{1}{1}{9}
\boxr{0}{1}{9}
\boxr{1}{1}{10}
\boxr{0}{1}{10}

\boxxl{2}{0}{3}
\boxxr{2}{0}{3}
\boxxr{2}{0}{2}

\boxxt{3}{0}{2}
\boxxr{3}{0}{2}

\boxxl{4}{0}{2}
\boxxt{4}{0}{2}
\boxxr{4}{0}{2}

\boxxl{2}{0}{4}
\boxxr{2}{0}{4}

\boxxl{2}{0}{5}
\boxxr{2}{0}{5}

\boxxl{2}{0}{6}
\boxxr{2}{0}{6}

\boxxl{2}{0}{7}
\boxxr{2}{0}{7}

\boxxl{2}{0}{8}
\boxxr{2}{0}{8}

\boxxl{2}{0}{9}
\boxxr{2}{0}{9}

\boxxl{2}{0}{10}
\boxxt{2}{0}{10}
\boxxr{2}{0}{10}

\boxxl{0}{2}{2}
\boxxt{0}{2}{2}
\boxxr{0}{2}{2}

\draw [thick,->] (0,11)-- (0,13);
\draw [thick,->] (5.5*1.74,-5.5*1.74*0.57735)-- (11,-11*0.57735);
\draw [thick,->]  (-9.5,-9.5*0.57735)-- (-11,-11*0.57735);

\node [right] at (0,13.5) {$x_3$};
\node [right] at (11,-11*0.57735) {$x_2$};
\node [left] at  (-11,-11*0.57735) {$x_1$};
\node [above] at  (-9,-7*0.57735) {$ $};
\node [right] at (2,10) {$ $};
\node  at (6,6*0.57735) {$ $};
\node  at (-6,6*0.57735) {$ $};
\node  at (0,-6) {$ $};

%\shade [ball color=gray] (-2.6,-.6) circle [radius=0.4cm];
%\shade [ball color=gray] (1.8,-1.2) circle [radius=0.4cm];
%\shade [ball color=gray] (-0.85,-1.7) circle [radius=0.4cm];

\end{tikzpicture}

%\end{document} 

%% file: fig1b.tex
% Plane partition
% Author: Jang Soo Kim
%\documentclass{minimal}
%\usepackage{tikz}
% Three counters
%\newcounter{x}
%\newcounter{y}
%\newcounter{z}

% The angles of x,y,z-axes
\newcommand\xaxis{210}
\newcommand\yaxis{-30}
\newcommand\zaxis{90}

% The top side of a cube
\newcommand\topside[3]{
  \fill[fill=cyan, draw=black,shift={(\xaxis:#1)},shift={(\yaxis:#2)},
  shift={(\zaxis:#3)}] (0,0) -- (30:1) -- (0,1) --(150:1)--(0,0);
}

% The left side of a cube
\newcommand\leftside[3]{
  \fill[fill=cyan, draw=black,shift={(\xaxis:#1)},shift={(\yaxis:#2)},
  shift={(\zaxis:#3)}] (0,0) -- (0,-1) -- (210:1) --(150:1)--(0,0);
}

% The right side of a cube
\newcommand\rightside[3]{
  \fill[fill=cyan, draw=black,shift={(\xaxis:#1)},shift={(\yaxis:#2)},
  shift={(\zaxis:#3)}] (0,0) -- (30:1) -- (-30:1) --(0,-1)--(0,0);
}

% The cube 
\newcommand\cube[3]{
  \topside{#1}{#2}{#3} \leftside{#1}{#2}{#3} \rightside{#1}{#2}{#3}
}

% Definition of \planepartition
% To draw the following plane partition, just write \planepartition{ {a, b, c}, {d,e} }.
%  a b c
%  d e
\newcommand\planepartition[1]{
 \setcounter{x}{-1}
  \foreach \a in {#1} {
    \addtocounter{x}{1}
    \setcounter{y}{-1}
    \foreach \b in \a {
      \addtocounter{y}{1}
      \setcounter{z}{-1}
      \foreach \c in {1,...,\b} {
        \addtocounter{z}{1}
        \cube{\value{x}}{\value{y}}{\value{z}}
      }
    }
  }
}

%\begin{document} 

\begin{tikzpicture}[scale=0.25]

\planepartition{{2,2},{2,2}}

\draw [thick,->] (0,2)-- (0,4);
\draw [thick,->] (1*1.74,-1*1.74*0.57735)-- (3,-3*0.57735);
\draw [thick,->]  (-1.7,-1.7*0.57735)-- (-3,-3*0.57735);

\node [right] at (0,4.5) {$x_3$};
\node [right] at (3,-3*0.57735) {$x_2$};
\node [left] at  (-3,-3*0.57735) {$x_1$};
\node [above] at  (-9,-7*0.57735) {$ $};
\node [right] at (2,10) {$ $};
\node  at (6,6*0.57735) {$ $};
\node  at (-6,6*0.57735) {$ $};
\node  at (0,-6) {$ $};

%\shade [ball color=gray] (-2.6,-.6) circle [radius=0.4cm];
%\shade [ball color=gray] (1.8,-1.2) circle [radius=0.4cm];
%\shade [ball color=gray] (-0.85,-1.7) circle [radius=0.4cm];

\end{tikzpicture}

%\end{document} 

%% file: fig3.tex
% Plane partition
% Author: Jang Soo Kim
%\documentclass{minimal}
%\usepackage{tikz}
% Three counters
%\newcounter{x}
%\newcounter{y}
%\newcounter{z}

% The angles of x,y,z-axes
\newcommand\xaxis{210}
\newcommand\yaxis{-30}
\newcommand\zaxis{90}

% The top side of a cube
\newcommand\topside[3]{
  \fill[fill=cyan, draw=black,shift={(\xaxis:#1)},shift={(\yaxis:#2)},
  shift={(\zaxis:#3)}] (0,0) -- (30:1) -- (0,1) --(150:1)--(0,0);
}

% The left side of a cube
\newcommand\leftside[3]{
  \fill[fill=cyan, draw=black,shift={(\xaxis:#1)},shift={(\yaxis:#2)},
  shift={(\zaxis:#3)}] (0,0) -- (0,-1) -- (210:1) --(150:1)--(0,0);
}

% The right side of a cube
\newcommand\rightside[3]{
  \fill[fill=cyan, draw=black,shift={(\xaxis:#1)},shift={(\yaxis:#2)},
  shift={(\zaxis:#3)}] (0,0) -- (30:1) -- (-30:1) --(0,-1)--(0,0);
}

% The cube 
\newcommand\cube[3]{
  \topside{#1}{#2}{#3} \leftside{#1}{#2}{#3} \rightside{#1}{#2}{#3}
}

% Definition of \planepartition
% To draw the following plane partition, just write \planepartition{ {a, b, c}, {d,e} }.
%  a b c
%  d e
\newcommand\planepartition[1]{
 \setcounter{x}{-1}
  \foreach \a in {#1} {
    \addtocounter{x}{1}
    \setcounter{y}{-1}
    \foreach \b in \a {
      \addtocounter{y}{1}
      \setcounter{z}{-1}
      \foreach \c in {1,...,\b} {
        \addtocounter{z}{1}
        \cube{\value{x}}{\value{y}}{\value{z}}
      }
    }
  }
}

%\begin{document} 

\begin{tikzpicture}[scale=0.25]

\newcommand{\boxt}[3]{ 
\draw[fill=blue, draw=black,shift={(\xaxis:#1)},shift={(\yaxis:#2)},
  shift={(\zaxis:#3)}] (0,0) -- (30:1) -- (0,1) --(150:1)--(0,0);
}

\newcommand{\boxl}[3]{ 
\draw[fill=blue, draw=black,shift={(\xaxis:#1)},shift={(\yaxis:#2)},
  shift={(\zaxis:#3)}] (0,0) -- (0,-1) -- (210:1) --(150:1)--(0,0);
}

\newcommand{\boxr}[3]{ 
\draw[fill=blue, draw=black,shift={(\xaxis:#1)},shift={(\yaxis:#2)},
  shift={(\zaxis:#3)}] (0,0) -- (30:1) -- (-30:1) --(0,-1)--(0,0);
}

\newcommand{\boxxt}[3]{ 
\draw[fill=violet, draw=black,shift={(\xaxis:#1)},shift={(\yaxis:#2)},
  shift={(\zaxis:#3)}] (0,0) -- (30:1) -- (0,1) --(150:1)--(0,0);
}

\newcommand{\boxxl}[3]{ 
\draw[fill=violet, draw=black,shift={(\xaxis:#1)},shift={(\yaxis:#2)},
  shift={(\zaxis:#3)}] (0,0) -- (0,-1) -- (210:1) --(150:1)--(0,0);
}

\newcommand{\boxxr}[3]{ 
\draw[fill=violet, draw=black,shift={(\xaxis:#1)},shift={(\yaxis:#2)},
  shift={(\zaxis:#3)}] (0,0) -- (30:1) -- (-30:1) --(0,-1)--(0,0);
}

\planepartition{{11,2,2,2,2,2,2,2,2,2,2},{11,2,2,2,2,2,2,2,2,2,2},{11,2,2,2,2,2,2,2,2,2,2},{11,2,2,2,2,2,2,2,2,2,2},{11,2,2,2,2,2,2,2,2,2,2},{11,2,2,2,2,2,2,2,2,2,2},{11,2,2,2,2,2,2,2,2,2,2},{11,2,2,2,2,2,2,2,2,2,2},{11,2,2,2,2,2,2,2,2,2,2},{11,2,2,2,2,2,2,2,2,2,2},{11,2,2,2,2,2,2,2,2,2,2}}

\boxt{0}{1}{1}
\boxt{0}{2}{1}
\boxt{0}{3}{1}
\boxt{0}{4}{1}
\boxt{0}{5}{1}
\boxt{0}{6}{1}
\boxt{0}{7}{1}
\boxt{0}{8}{1}
\boxt{0}{9}{1}
\boxt{0}{10}{1}

\boxt{1}{1}{1}
\boxt{1}{2}{1}
\boxt{1}{3}{1}
\boxt{1}{4}{1}
\boxt{1}{5}{1}
\boxt{1}{6}{1}
\boxt{1}{7}{1}
\boxt{1}{8}{1}
\boxt{1}{9}{1}
\boxt{1}{10}{1}

\boxt{2}{1}{1}
\boxt{2}{2}{1}
\boxt{2}{3}{1}
\boxt{2}{4}{1}
\boxt{2}{5}{1}
\boxt{2}{6}{1}
\boxt{2}{7}{1}
\boxt{2}{8}{1}
\boxt{2}{9}{1}
\boxt{2}{10}{1}

\boxt{3}{1}{1}
\boxt{3}{2}{1}
\boxt{3}{3}{1}
\boxt{3}{4}{1}
\boxt{3}{5}{1}
\boxt{3}{6}{1}
\boxt{3}{7}{1}
\boxt{3}{8}{1}
\boxt{3}{9}{1}
\boxt{3}{10}{1}

\boxt{4}{1}{1}
\boxt{4}{2}{1}
\boxt{4}{3}{1}
\boxt{4}{4}{1}
\boxt{4}{5}{1}
\boxt{4}{6}{1}
\boxt{4}{7}{1}
\boxt{4}{8}{1}
\boxt{4}{9}{1}
\boxt{4}{10}{1}

\boxt{5}{1}{1}
\boxt{5}{2}{1}
\boxt{5}{3}{1}
\boxt{5}{4}{1}
\boxt{5}{5}{1}
\boxt{5}{6}{1}
\boxt{5}{7}{1}
\boxt{5}{8}{1}
\boxt{5}{9}{1}
\boxt{5}{10}{1}

\boxt{6}{1}{1}
\boxt{6}{2}{1}
\boxt{6}{3}{1}
\boxt{6}{4}{1}
\boxt{6}{5}{1}
\boxt{6}{6}{1}
\boxt{6}{7}{1}
\boxt{6}{8}{1}
\boxt{6}{9}{1}
\boxt{6}{10}{1}

\boxt{7}{1}{1}
\boxt{7}{2}{1}
\boxt{7}{3}{1}
\boxt{7}{4}{1}
\boxt{7}{5}{1}
\boxt{7}{6}{1}
\boxt{7}{7}{1}
\boxt{7}{8}{1}
\boxt{7}{9}{1}
\boxt{7}{10}{1}

\boxt{8}{1}{1}
\boxt{8}{2}{1}
\boxt{8}{3}{1}
\boxt{8}{4}{1}
\boxt{8}{5}{1}
\boxt{8}{6}{1}
\boxt{8}{7}{1}
\boxt{8}{8}{1}
\boxt{8}{9}{1}
\boxt{8}{10}{1}

\boxt{9}{1}{1}
\boxt{9}{2}{1}
\boxt{9}{3}{1}
\boxt{9}{4}{1}
\boxt{9}{5}{1}
\boxt{9}{6}{1}
\boxt{9}{7}{1}
\boxt{9}{8}{1}
\boxt{9}{9}{1}
\boxt{9}{10}{1}

\boxt{10}{1}{1}
\boxt{10}{2}{1}
\boxt{10}{3}{1}
\boxt{10}{4}{1}
\boxt{10}{5}{1}
\boxt{10}{6}{1}
\boxt{10}{7}{1}
\boxt{10}{8}{1}
\boxt{10}{9}{1}
\boxt{10}{10}{1}

\boxl{10}{0}{0}
\boxl{10}{0}{1}
\boxl{10}{1}{0}
\boxl{10}{1}{1}
\boxl{10}{2}{0}
\boxl{10}{2}{1}
\boxl{10}{3}{0}
\boxl{10}{3}{1}
\boxl{10}{4}{0}
\boxl{10}{4}{1}
\boxl{10}{5}{0}
\boxl{10}{5}{1}
\boxl{10}{6}{0}
\boxl{10}{6}{1}
\boxl{10}{7}{0}
\boxl{10}{7}{1}
\boxl{10}{8}{0}
\boxl{10}{8}{1}
\boxl{10}{9}{0}
\boxl{10}{9}{1}
\boxl{10}{10}{0}
\boxl{10}{10}{1}

\boxr{0}{10}{0}
\boxr{0}{10}{1}
\boxr{1}{10}{0}
\boxr{1}{10}{1}
\boxr{2}{10}{0}
\boxr{2}{10}{1}
\boxr{3}{10}{0}
\boxr{3}{10}{1}
\boxr{4}{10}{0}
\boxr{4}{10}{1}
\boxr{5}{10}{0}
\boxr{5}{10}{1}
\boxr{6}{10}{0}
\boxr{6}{10}{1}
\boxr{7}{10}{0}
\boxr{7}{10}{1}
\boxr{8}{10}{0}
\boxr{8}{10}{1}
\boxr{9}{10}{0}
\boxr{9}{10}{1}
\boxr{10}{10}{0}
\boxr{10}{10}{1}

\boxxr{0}{0}{2}
\boxxr{0}{0}{3}
\boxxr{0}{0}{4}
\boxxr{0}{0}{5}
\boxxr{0}{0}{6}
\boxxr{0}{0}{7}
\boxxr{0}{0}{8}
\boxxr{0}{0}{9}
\boxxr{0}{0}{10}
\boxxt{0}{0}{10}

\boxxr{1}{0}{2}
\boxxr{1}{0}{3}
\boxxr{1}{0}{4}
\boxxr{1}{0}{5}
\boxxr{1}{0}{6}
\boxxr{1}{0}{7}
\boxxr{1}{0}{8}
\boxxr{1}{0}{9}
\boxxr{1}{0}{10}
\boxxt{1}{0}{10}

\draw [thick,->] (0,11)-- (0,13);
\draw [thick,->] (5.5*1.74,-5.5*1.74*0.57735)-- (11,-11*0.57735);
\draw [thick,->]  (-9.5,-9.5*0.57735)-- (-11,-11*0.57735);

\node [right] at (0,13.5) {$x_3$};
\node [right] at (11,-11*0.57735) {$x_2$};
\node [left] at  (-11,-11*0.57735) {$x_1$};
\node [above] at  (-9,-7*0.57735) {$ $};
\node [right] at (2,10) {$ $};
\node  at (6,6*0.57735) {$ $};
\node  at (-6,6*0.57735) {$ $};
\node  at (0,-6) {$ $};

%\shade [ball color=gray] (-2.6,-.6) circle [radius=0.4cm];
%\shade [ball color=gray] (1.8,-1.2) circle [radius=0.4cm];
%\shade [ball color=gray] (-0.85,-1.7) circle [radius=0.4cm];

\end{tikzpicture}

%\end{document} 

%% file: fig3a.tex
% Plane partition
% Author: Jang Soo Kim
%\documentclass{minimal}
%\usepackage{tikz}
% Three counters
%\newcounter{x}
%\newcounter{y}
%\newcounter{z}

% The angles of x,y,z-axes
\newcommand\xaxis{210}
\newcommand\yaxis{-30}
\newcommand\zaxis{90}

% The top side of a cube
\newcommand\topside[3]{
  \fill[fill=cyan, draw=black,shift={(\xaxis:#1)},shift={(\yaxis:#2)},
  shift={(\zaxis:#3)}] (0,0) -- (30:1) -- (0,1) --(150:1)--(0,0);
}

% The left side of a cube
\newcommand\leftside[3]{
  \fill[fill=cyan, draw=black,shift={(\xaxis:#1)},shift={(\yaxis:#2)},
  shift={(\zaxis:#3)}] (0,0) -- (0,-1) -- (210:1) --(150:1)--(0,0);
}

% The right side of a cube
\newcommand\rightside[3]{
  \fill[fill=cyan, draw=black,shift={(\xaxis:#1)},shift={(\yaxis:#2)},
  shift={(\zaxis:#3)}] (0,0) -- (30:1) -- (-30:1) --(0,-1)--(0,0);
}

% The cube 
\newcommand\cube[3]{
  \topside{#1}{#2}{#3} \leftside{#1}{#2}{#3} \rightside{#1}{#2}{#3}
}

% Definition of \planepartition
% To draw the following plane partition, just write \planepartition{ {a, b, c}, {d,e} }.
%  a b c
%  d e
\newcommand\planepartition[1]{
 \setcounter{x}{-1}
  \foreach \a in {#1} {
    \addtocounter{x}{1}
    \setcounter{y}{-1}
    \foreach \b in \a {
      \addtocounter{y}{1}
      \setcounter{z}{-1}
      \foreach \c in {1,...,\b} {
        \addtocounter{z}{1}
        \cube{\value{x}}{\value{y}}{\value{z}}
      }
    }
  }
}

%\begin{document} 

\begin{tikzpicture}[scale=0.25]

\newcommand{\boxt}[3]{ 
\draw[fill=blue, draw=black,shift={(\xaxis:#1)},shift={(\yaxis:#2)},
  shift={(\zaxis:#3)}] (0,0) -- (30:1) -- (0,1) --(150:1)--(0,0);
}

\newcommand{\boxl}[3]{ 
\draw[fill=blue, draw=black,shift={(\xaxis:#1)},shift={(\yaxis:#2)},
  shift={(\zaxis:#3)}] (0,0) -- (0,-1) -- (210:1) --(150:1)--(0,0);
}

\newcommand{\boxr}[3]{ 
\draw[fill=blue, draw=black,shift={(\xaxis:#1)},shift={(\yaxis:#2)},
  shift={(\zaxis:#3)}] (0,0) -- (30:1) -- (-30:1) --(0,-1)--(0,0);
}

\newcommand{\boxxt}[3]{ 
\draw[fill=purple, draw=black,shift={(\xaxis:#1)},shift={(\yaxis:#2)},
  shift={(\zaxis:#3)}] (0,0) -- (30:1) -- (0,1) --(150:1)--(0,0);
}

\newcommand{\boxxl}[3]{ 
\draw[fill=purple, draw=black,shift={(\xaxis:#1)},shift={(\yaxis:#2)},
  shift={(\zaxis:#3)}] (0,0) -- (0,-1) -- (210:1) --(150:1)--(0,0);
}

\newcommand{\boxxr}[3]{ 
\draw[fill=purple, draw=black,shift={(\xaxis:#1)},shift={(\yaxis:#2)},
  shift={(\zaxis:#3)}] (0,0) -- (30:1) -- (-30:1) --(0,-1)--(0,0);
}

\planepartition{{11,2},{11,2}}

\draw [thick,->] (0,11)-- (0,13);
\draw [thick,->] (1*1.74,-1*1.74*0.57735)-- (3,-3*0.57735);
\draw [thick,->]  (-1.7,-1.7*0.57735)-- (-3,-3*0.57735);

\node [right] at (0,13.5) {$x_3$};
\node [right] at (3,-3*0.57735) {$x_4$};
\node [left] at  (-3,-3*0.57735) {$x_1$};
\node [above] at  (-9,-7*0.57735) {$ $};
\node [right] at (2,10) {$ $};
\node  at (6,6*0.57735) {$ $};
\node  at (-6,6*0.57735) {$ $};
\node  at (0,-6) {$ $};

%\shade [ball color=gray] (-2.6,-.6) circle [radius=0.4cm];
%\shade [ball color=gray] (1.8,-1.2) circle [radius=0.4cm];
%\shade [ball color=gray] (-0.85,-1.7) circle [radius=0.4cm];

\end{tikzpicture}

%\end{document} 

%% file: fig2.tex
% Plane partition
% Author: Jang Soo Kim
%\documentclass{minimal}
%\usepackage{tikz}
% Three counters
%\newcounter{x}
%\newcounter{y}
%\newcounter{z}

% The angles of x,y,z-axes
\newcommand\xaxis{210}
\newcommand\yaxis{-30}
\newcommand\zaxis{90}

% The top side of a cube
\newcommand\topside[3]{
  \fill[fill=cyan, draw=black,shift={(\xaxis:#1)},shift={(\yaxis:#2)},
  shift={(\zaxis:#3)}] (0,0) -- (30:1) -- (0,1) --(150:1)--(0,0);
}

% The left side of a cube
\newcommand\leftside[3]{
  \fill[fill=cyan, draw=black,shift={(\xaxis:#1)},shift={(\yaxis:#2)},
  shift={(\zaxis:#3)}] (0,0) -- (0,-1) -- (210:1) --(150:1)--(0,0);
}

% The right side of a cube
\newcommand\rightside[3]{
  \fill[fill=cyan, draw=black,shift={(\xaxis:#1)},shift={(\yaxis:#2)},
  shift={(\zaxis:#3)}] (0,0) -- (30:1) -- (-30:1) --(0,-1)--(0,0);
}

% The cube 
\newcommand\cube[3]{
  \topside{#1}{#2}{#3} \leftside{#1}{#2}{#3} \rightside{#1}{#2}{#3}
}

% Definition of \planepartition
% To draw the following plane partition, just write \planepartition{ {a, b, c}, {d,e} }.
%  a b c
%  d e
\newcommand\planepartition[1]{
 \setcounter{x}{-1}
  \foreach \a in {#1} {
    \addtocounter{x}{1}
    \setcounter{y}{-1}
    \foreach \b in \a {
      \addtocounter{y}{1}
      \setcounter{z}{-1}
      \foreach \c in {1,...,\b} {
        \addtocounter{z}{1}
        \cube{\value{x}}{\value{y}}{\value{z}}
      }
    }
  }
}

%\begin{document} 

\begin{tikzpicture}[scale=0.25]

\planepartition{{4,3,3,3,3,3,3,3,3,3,3},{4,3,3,3,3,3,3,3,3,3,3},{4,2,2,2,2,2,2,2,2,2,2},{4,2,2,2,2,2,2,2,2,2,2},{4,2,2,2,2,2,2,2,2,2,2},{4,2,2,2,2,2,2,2,2,2,2},{4,2,2,2,2,2,2,2,2,2,2},{4,2,2,2,2,2,2,2,2,2,2},{4,2,2,2,2,2,2,2,2,2,2},{4,2,2,2,2,2,2,2,2,2,2},{4,2,2,2,2,2,2,2,2,2,2}}

\draw [thick,->] (0,4)-- (0,6);
\draw [thick,->] (5.5*1.74,-5.5*1.74*0.57735)-- (11,-11*0.57735);
\draw [thick,->]  (-9.5,-9.5*0.57735)-- (-11,-11*0.57735);

\node [right] at (0,7) {$x_3$};
\node [right] at (11,-11*0.57735) {$x_2$};
\node [left] at  (-11,-11*0.57735) {$x_1$};
\node [above] at  (-9,-7*0.57735) {$ $};
\node [right] at (2,10) {$ $};
\node  at (6,6*0.57735) {$ $};
\node  at (-6,6*0.57735) {$ $};
\node  at (0,-6) {$ $};

%\shade [ball color=gray] (-2.6,-.6) circle [radius=0.4cm];
%\shade [ball color=gray] (1.8,-1.2) circle [radius=0.4cm];
%\shade [ball color=gray] (-0.85,-1.7) circle [radius=0.4cm];

\end{tikzpicture}

%\end{document} 

%% file: fig2a.tex
% Plane partition
% Author: Jang Soo Kim
%\documentclass{minimal}
%\usepackage{tikz}
% Three counters
%\newcounter{x}
%\newcounter{y}
%\newcounter{z}

% The angles of x,y,z-axes
\newcommand\xaxis{210}
\newcommand\yaxis{-30}
\newcommand\zaxis{90}

% The top side of a cube
\newcommand\topside[3]{
  \fill[fill=cyan, draw=black,shift={(\xaxis:#1)},shift={(\yaxis:#2)},
  shift={(\zaxis:#3)}] (0,0) -- (30:1) -- (0,1) --(150:1)--(0,0);
}

% The left side of a cube
\newcommand\leftside[3]{
  \fill[fill=cyan, draw=black,shift={(\xaxis:#1)},shift={(\yaxis:#2)},
  shift={(\zaxis:#3)}] (0,0) -- (0,-1) -- (210:1) --(150:1)--(0,0);
}

% The right side of a cube
\newcommand\rightside[3]{
  \fill[fill=cyan, draw=black,shift={(\xaxis:#1)},shift={(\yaxis:#2)},
  shift={(\zaxis:#3)}] (0,0) -- (30:1) -- (-30:1) --(0,-1)--(0,0);
}

% The cube 
\newcommand\cube[3]{
  \topside{#1}{#2}{#3} \leftside{#1}{#2}{#3} \rightside{#1}{#2}{#3}
}

% Definition of \planepartition
% To draw the following plane partition, just write \planepartition{ {a, b, c}, {d,e} }.
%  a b c
%  d e
\newcommand\planepartition[1]{
 \setcounter{x}{-1}
  \foreach \a in {#1} {
    \addtocounter{x}{1}
    \setcounter{y}{-1}
    \foreach \b in \a {
      \addtocounter{y}{1}
      \setcounter{z}{-1}
      \foreach \c in {1,...,\b} {
        \addtocounter{z}{1}
        \cube{\value{x}}{\value{y}}{\value{z}}
      }
    }
  }
}

%\begin{document} 

\begin{tikzpicture}[scale=0.25]

\newcommand{\boxt}[3]{ 
\draw[fill=cyan, draw=black,shift={(\xaxis:#1)},shift={(\yaxis:#2)},
  shift={(\zaxis:#3)}] (0,0) -- (30:1) -- (0,1) --(150:1)--(0,0);
}
\newcommand{\boxl}[3]{ 
\draw[fill=cyan, draw=black,shift={(\xaxis:#1)},shift={(\yaxis:#2)},
  shift={(\zaxis:#3)}] (0,0) -- (0,-1) -- (210:1) --(150:1)--(0,0);
}
\newcommand{\boxr}[3]{ 
\draw[fill=cyan, draw=black,shift={(\xaxis:#1)},shift={(\yaxis:#2)},
  shift={(\zaxis:#3)}] (0,0) -- (30:1) -- (-30:1) --(0,-1)--(0,0);
}

%\planepartition{{3,2,2,2,2,2,2,2,2,2,2},{3,2,2,2,2,2,2,2,2,2,2},{3,1,1,1,1,1,1,1,1,1,1},{3,1,1,1,1,1,1,1,1,1,1},{3,1,1,1,1,1,1,1,1,1,1},{3,1,1,1,1,1,1,1,1,1,1},{3,1,1,1,1,1,1,1,1,1,1},{3,1,1,1,1,1,1,1,1,1,1},{3,1,1,1,1,1,1,1,1,1,1},{3,1,1,1,1,1,1,1,1,1,1},{3,1,1,1,1,1,1,1,1,1,1}}

\boxt{0}{0}{2}
\boxl{0}{0}{2}
\boxr{0}{0}{2}
\boxt{0}{-1}{2}
%\boxl{1}{0}{2}
%\boxr{1}{0}{2}
\boxt{0}{-2}{2}
%\boxl{1}{0}{2}
%\boxr{1}{0}{2}
\boxt{0}{-3}{2}
%\boxl{1}{0}{2}
%\boxr{1}{0}{2}
\boxt{0}{-4}{2}
%\boxl{1}{0}{2}
%\boxr{1}{0}{2}
\boxt{0}{-5}{2}
%\boxl{1}{0}{2}
%\boxr{1}{0}{2}
\boxt{0}{-6}{2}
%\boxl{1}{0}{2}
%\boxr{1}{0}{2}
\boxt{0}{-7}{2}
%\boxl{1}{0}{2}
%\boxr{1}{0}{2}
\boxt{0}{-8}{2}
%\boxl{1}{0}{2}
%\boxr{1}{0}{2}
\boxt{0}{-9}{2}
%\boxl{1}{0}{2}
%\boxr{1}{0}{2}
\boxt{0}{-10}{2}
%\boxl{1}{0}{2}
%\boxr{1}{0}{2}

\boxt{1}{0}{2}
\boxl{1}{0}{2}
\boxr{1}{0}{2}
\boxt{1}{-1}{2}
\boxl{1}{-1}{2}
%\boxr{1}{-1}{2}
\boxt{1}{-2}{2}
\boxl{1}{-2}{2}
%\boxr{1}{-1}{2}
\boxt{1}{-3}{2}
\boxl{1}{-3}{2}
%\boxr{1}{-1}{2}
\boxt{1}{-4}{2}
\boxl{1}{-4}{2}
%\boxr{1}{-1}{2}
\boxt{1}{-5}{2}
\boxl{1}{-5}{2}
%\boxr{1}{-1}{2}
\boxt{1}{-6}{2}
\boxl{1}{-6}{2}
%\boxr{1}{-1}{2}
\boxt{1}{-7}{2}
\boxl{1}{-7}{2}
%\boxr{1}{-1}{2}
\boxt{1}{-8}{2}
\boxl{1}{-8}{2}
%\boxr{1}{-1}{2}
\boxt{1}{-9}{2}
\boxl{1}{-9}{2}
%\boxr{1}{-1}{2}
\boxt{1}{-10}{2}
\boxl{1}{-10}{2}
%\boxr{1}{-1}{2}

%\boxt{-1}{0}{2}
%\boxl{-1}{0}{2}
\boxr{-1}{0}{2}
%\boxt{-2}{0}{2}
%\boxl{-1}{0}{2}
\boxr{-2}{0}{2}
%\boxt{-3}{0}{2}
%\boxl{-1}{0}{2}
\boxr{-3}{0}{2}
%\boxt{-4}{0}{2}
%\boxl{-1}{0}{2}
\boxr{-4}{0}{2}
%\boxt{-5}{0}{2}
%\boxl{-1}{0}{2}
\boxr{-5}{0}{2}
%\boxt{-6}{0}{2}
%\boxl{-1}{0}{2}
\boxr{-6}{0}{2}
%\boxt{-7}{0}{2}
%\boxl{-1}{0}{2}
\boxr{-7}{0}{2}
%\boxt{-8}{0}{2}
%\boxl{-1}{0}{2}
\boxr{-8}{0}{2}
%\boxt{-9}{0}{2}
%\boxl{-1}{0}{2}
\boxr{-9}{0}{2}
%\boxt{-10}{0}{2}
%\boxl{-1}{0}{2}
\boxr{-10}{0}{2}
%\boxt{-11}{0}{2}
%\boxl{-1}{0}{2}
\boxr{-11}{0}{2}

\boxt{-1}{0}{3}
\boxl{-1}{0}{3}
\boxr{-1}{0}{3}
\boxt{-2}{0}{3}
%\boxl{-2}{0}{3}
\boxr{-2}{0}{3}
\boxt{-3}{0}{3}
%\boxl{-2}{0}{3}
\boxr{-3}{0}{3}
\boxt{-4}{0}{3}
%\boxl{-2}{0}{3}
\boxr{-4}{0}{3}
\boxt{-5}{0}{3}
%\boxl{-2}{0}{3}
\boxr{-5}{0}{3}
\boxt{-6}{0}{3}
%\boxl{-2}{0}{3}
\boxr{-6}{0}{3}
\boxt{-7}{0}{3}
%\boxl{-2}{0}{3}
\boxr{-7}{0}{3}
\boxt{-8}{0}{3}
%\boxl{-2}{0}{3}
\boxr{-8}{0}{3}
\boxt{-9}{0}{3}
%\boxl{-2}{0}{3}
\boxr{-9}{0}{3}
\boxt{-10}{0}{3}
%\boxl{-2}{0}{3}
\boxr{-10}{0}{3}
\boxt{-11}{0}{3}
%\boxl{-2}{0}{3}
\boxr{-11}{0}{3}

\draw [thick,->] (0,3)-- (0,6);
\draw [thick,-] (0,0)-- (0,1);
\draw [thick,->] (0*1.74,-0*1.74*0.57735)-- (6,-6*0.57735);
\draw [thick,->]  (0,-0*0.57735)-- (-6,-6*0.57735);

\node [right] at (0,7) {$x_3$};
\node [right] at (6,-6*0.57735) {$x_2$};
\node [left] at  (-6,-6*0.57735) {$x_1$};
\node [above] at  (-9,-7*0.57735) {$ $};
\node [right] at (2,10) {$ $};
\node  at (6,6*0.57735) {$ $};
\node  at (-6,6*0.57735) {$ $};
\node  at (0,-6) {$ $};

%\shade [ball color=gray] (-2.6,-.6) circle [radius=0.4cm];
%\shade [ball color=gray] (1.8,-1.2) circle [radius=0.4cm];
%\shade [ball color=gray] (-0.85,-1.7) circle [radius=0.4cm];

\end{tikzpicture}

%\end{document} 

%% file: fig1c.tex
% Plane partition
% Author: Jang Soo Kim
%\documentclass{minimal}
%\usepackage{tikz}
% Three counters
%\newcounter{x}
%\newcounter{y}
%\newcounter{z}

% The angles of x,y,z-axes
\newcommand\xaxis{210}
\newcommand\yaxis{-30}
\newcommand\zaxis{90}

% The top side of a cube
\newcommand\topside[3]{
  \fill[fill=cyan, draw=black,shift={(\xaxis:#1)},shift={(\yaxis:#2)},
  shift={(\zaxis:#3)}] (0,0) -- (30:1) -- (0,1) --(150:1)--(0,0);
}

% The left side of a cube
\newcommand\leftside[3]{
  \fill[fill=cyan, draw=black,shift={(\xaxis:#1)},shift={(\yaxis:#2)},
  shift={(\zaxis:#3)}] (0,0) -- (0,-1) -- (210:1) --(150:1)--(0,0);
}

% The right side of a cube
\newcommand\rightside[3]{
  \fill[fill=cyan, draw=black,shift={(\xaxis:#1)},shift={(\yaxis:#2)},
  shift={(\zaxis:#3)}] (0,0) -- (30:1) -- (-30:1) --(0,-1)--(0,0);
}

% The cube 
\newcommand\cube[3]{
  \topside{#1}{#2}{#3} \leftside{#1}{#2}{#3} \rightside{#1}{#2}{#3}
}

% Definition of \planepartition
% To draw the following plane partition, just write \planepartition{ {a, b, c}, {d,e} }.
%  a b c
%  d e
\newcommand\planepartition[1]{
 \setcounter{x}{-1}
  \foreach \a in {#1} {
    \addtocounter{x}{1}
    \setcounter{y}{-1}
    \foreach \b in \a {
      \addtocounter{y}{1}
      \setcounter{z}{-1}
      \foreach \c in {1,...,\b} {
        \addtocounter{z}{1}
        \cube{\value{x}}{\value{y}}{\value{z}}
      }
    }
  }
}

%\begin{document} 

\begin{tikzpicture}[scale=0.25]

\newcommand{\boxt}[3]{ 
\draw[fill=red, draw=black,shift={(\xaxis:#1)},shift={(\yaxis:#2)},
  shift={(\zaxis:#3)}] (0,0) -- (30:1) -- (0,1) --(150:1)--(0,0);
}

\newcommand{\boxl}[3]{ 
\draw[fill=red, draw=black,shift={(\xaxis:#1)},shift={(\yaxis:#2)},
  shift={(\zaxis:#3)}] (0,0) -- (0,-1) -- (210:1) --(150:1)--(0,0);
}

\newcommand{\boxr}[3]{ 
\draw[fill=red, draw=black,shift={(\xaxis:#1)},shift={(\yaxis:#2)},
  shift={(\zaxis:#3)}] (0,0) -- (30:1) -- (-30:1) --(0,-1)--(0,0);
}

\planepartition{{2,2},{2,2}}

\boxt{1}{1}{1}
\boxt{0}{1}{1}
\boxt{1}{0}{1}
\boxl{1}{1}{1}
\boxl{1}{1}{0}
\boxl{1}{0}{1}
\boxr{1}{1}{1}
\boxr{1}{1}{0}
\boxr{0}{1}{1}

\draw [thick,->] (0,2)-- (0,4);
\draw [thick,->] (1*1.74,-1*1.74*0.57735)-- (3,-3*0.57735);
\draw [thick,->]  (-1.7,-1.7*0.57735)-- (-3,-3*0.57735);

\node [right] at (0,4.5) {$x_3$};
\node [right] at (3,-3*0.57735) {$x_2$};
\node [left] at  (-3,-3*0.57735) {$x_1$};
\node [above] at  (-9,-7*0.57735) {$ $};
\node [right] at (2,10) {$ $};
\node  at (6,6*0.57735) {$ $};
\node  at (-6,6*0.57735) {$ $};
\node  at (0,-6) {$ $};

%\shade [ball color=gray] (-2.6,-.6) circle [radius=0.4cm];
%\shade [ball color=gray] (1.8,-1.2) circle [radius=0.4cm];
%\shade [ball color=gray] (-0.85,-1.7) circle [radius=0.4cm];

\end{tikzpicture}

%\end{document} 

%% file: fig3b.tex
% Plane partition
% Author: Jang Soo Kim
%\documentclass{minimal}
%\usepackage{tikz}
% Three counters
%\newcounter{x}
%\newcounter{y}
%\newcounter{z}

% The angles of x,y,z-axes
\newcommand\xaxis{210}
\newcommand\yaxis{-30}
\newcommand\zaxis{90}

% The top side of a cube
\newcommand\topside[3]{
  \fill[fill=cyan, draw=black,shift={(\xaxis:#1)},shift={(\yaxis:#2)},
  shift={(\zaxis:#3)}] (0,0) -- (30:1) -- (0,1) --(150:1)--(0,0);
}

% The left side of a cube
\newcommand\leftside[3]{
  \fill[fill=cyan, draw=black,shift={(\xaxis:#1)},shift={(\yaxis:#2)},
  shift={(\zaxis:#3)}] (0,0) -- (0,-1) -- (210:1) --(150:1)--(0,0);
}

% The right side of a cube
\newcommand\rightside[3]{
  \fill[fill=cyan, draw=black,shift={(\xaxis:#1)},shift={(\yaxis:#2)},
  shift={(\zaxis:#3)}] (0,0) -- (30:1) -- (-30:1) --(0,-1)--(0,0);
}

% The cube 
\newcommand\cube[3]{
  \topside{#1}{#2}{#3} \leftside{#1}{#2}{#3} \rightside{#1}{#2}{#3}
}

% Definition of \planepartition
% To draw the following plane partition, just write \planepartition{ {a, b, c}, {d,e} }.
%  a b c
%  d e
\newcommand\planepartition[1]{
 \setcounter{x}{-1}
  \foreach \a in {#1} {
    \addtocounter{x}{1}
    \setcounter{y}{-1}
    \foreach \b in \a {
      \addtocounter{y}{1}
      \setcounter{z}{-1}
      \foreach \c in {1,...,\b} {
        \addtocounter{z}{1}
        \cube{\value{x}}{\value{y}}{\value{z}}
      }
    }
  }
}

%\begin{document} 

\begin{tikzpicture}[scale=0.25]

\newcommand{\boxt}[3]{ 
\draw[fill=cyan, draw=black,shift={(\xaxis:#1)},shift={(\yaxis:#2)},
  shift={(\zaxis:#3)}] (0,0) -- (30:1) -- (0,1) --(150:1)--(0,0);
}

\newcommand{\boxl}[3]{ 
\draw[fill=cyan, draw=black,shift={(\xaxis:#1)},shift={(\yaxis:#2)},
  shift={(\zaxis:#3)}] (0,0) -- (0,-1) -- (210:1) --(150:1)--(0,0);
}

\newcommand{\boxr}[3]{ 
\draw[fill=cyan, draw=black,shift={(\xaxis:#1)},shift={(\yaxis:#2)},
  shift={(\zaxis:#3)}] (0,0) -- (30:1) -- (-30:1) --(0,-1)--(0,0);
}

\newcommand{\boxxt}[3]{ 
\draw[fill=red, draw=black,shift={(\xaxis:#1)},shift={(\yaxis:#2)},
  shift={(\zaxis:#3)}] (0,0) -- (30:1) -- (0,1) --(150:1)--(0,0);
}

\newcommand{\boxxl}[3]{ 
\draw[fill=red, draw=black,shift={(\xaxis:#1)},shift={(\yaxis:#2)},
  shift={(\zaxis:#3)}] (0,0) -- (0,-1) -- (210:1) --(150:1)--(0,0);
}

\newcommand{\boxxr}[3]{ 
\draw[fill=red, draw=black,shift={(\xaxis:#1)},shift={(\yaxis:#2)},
  shift={(\zaxis:#3)}] (0,0) -- (30:1) -- (-30:1) --(0,-1)--(0,0);
}

%\planepartition{{11,2},{11,2}}

\boxt{0}{1}{1}
\boxl{1}{1}{0}
\boxr{0}{1}{0}
\boxr{1}{1}{0}
\boxr{0}{1}{1}

\boxxt{1}{1}{1}
\boxxl{1}{1}{1}
\boxxr{1}{1}{1}

\draw [thick,->] (0,1)-- (0,4);
\draw [thick,->] (1*1.74,-1*1.74*0.57735)-- (3,-3*0.57735);
\draw [thick,->]  (-0.9,-1.3*0.57735)-- (-3,-3*0.57735);

\node [right] at (0,4.5) {$x_3$};
\node [right] at (3,-3*0.57735) {$x_4$};
\node [left] at  (-3,-3*0.57735) {$x_1$};
\node [above] at  (-9,-7*0.57735) {$ $};
\node [right] at (2,10) {$ $};
\node  at (6,6*0.57735) {$ $};
\node  at (-6,6*0.57735) {$ $};
\node  at (0,-6) {$ $};

%\shade [ball color=gray] (-2.6,-.6) circle [radius=0.4cm];
%\shade [ball color=gray] (1.8,-1.2) circle [radius=0.4cm];
%\shade [ball color=gray] (-0.85,-1.7) circle [radius=0.4cm];

\end{tikzpicture}

%\end{document} 

%% file: fig2b.tex
% Plane partition
% Author: Jang Soo Kim
%\documentclass{minimal}
%\usepackage{tikz}
% Three counters
%\newcounter{x}
%\newcounter{y}
%\newcounter{z}

% The angles of x,y,z-axes
\newcommand\xaxis{210}
\newcommand\yaxis{-30}
\newcommand\zaxis{90}

% The top side of a cube
\newcommand\topside[3]{
  \fill[fill=cyan, draw=black,shift={(\xaxis:#1)},shift={(\yaxis:#2)},
  shift={(\zaxis:#3)}] (0,0) -- (30:1) -- (0,1) --(150:1)--(0,0);
}

% The left side of a cube
\newcommand\leftside[3]{
  \fill[fill=cyan, draw=black,shift={(\xaxis:#1)},shift={(\yaxis:#2)},
  shift={(\zaxis:#3)}] (0,0) -- (0,-1) -- (210:1) --(150:1)--(0,0);
}

% The right side of a cube
\newcommand\rightside[3]{
  \fill[fill=cyan, draw=black,shift={(\xaxis:#1)},shift={(\yaxis:#2)},
  shift={(\zaxis:#3)}] (0,0) -- (30:1) -- (-30:1) --(0,-1)--(0,0);
}

% The cube 
\newcommand\cube[3]{
  \topside{#1}{#2}{#3} \leftside{#1}{#2}{#3} \rightside{#1}{#2}{#3}
}

% Definition of \planepartition
% To draw the following plane partition, just write \planepartition{ {a, b, c}, {d,e} }.
%  a b c
%  d e
\newcommand\planepartition[1]{
 \setcounter{x}{-1}
  \foreach \a in {#1} {
    \addtocounter{x}{1}
    \setcounter{y}{-1}
    \foreach \b in \a {
      \addtocounter{y}{1}
      \setcounter{z}{-1}
      \foreach \c in {1,...,\b} {
        \addtocounter{z}{1}
        \cube{\value{x}}{\value{y}}{\value{z}}
      }
    }
  }
}

%\begin{document} 

\begin{tikzpicture}[scale=0.25]

\newcommand{\boxt}[3]{ 
\draw[fill=cyan, draw=black,shift={(\xaxis:#1)},shift={(\yaxis:#2)},
  shift={(\zaxis:#3)}] (0,0) -- (30:1) -- (0,1) --(150:1)--(0,0);
}
\newcommand{\boxl}[3]{ 
\draw[fill=cyan, draw=black,shift={(\xaxis:#1)},shift={(\yaxis:#2)},
  shift={(\zaxis:#3)}] (0,0) -- (0,-1) -- (210:1) --(150:1)--(0,0);
}
\newcommand{\boxr}[3]{ 
\draw[fill=cyan, draw=black,shift={(\xaxis:#1)},shift={(\yaxis:#2)},
  shift={(\zaxis:#3)}] (0,0) -- (30:1) -- (-30:1) --(0,-1)--(0,0);
}

\newcommand{\boxxt}[3]{ 
\draw[fill=red, draw=black,shift={(\xaxis:#1)},shift={(\yaxis:#2)},
  shift={(\zaxis:#3)}] (0,0) -- (30:1) -- (0,1) --(150:1)--(0,0);
}
\newcommand{\boxxl}[3]{ 
\draw[fill=red, draw=black,shift={(\xaxis:#1)},shift={(\yaxis:#2)},
  shift={(\zaxis:#3)}] (0,0) -- (0,-1) -- (210:1) --(150:1)--(0,0);
}
\newcommand{\boxxr}[3]{ 
\draw[fill=red, draw=black,shift={(\xaxis:#1)},shift={(\yaxis:#2)},
  shift={(\zaxis:#3)}] (0,0) -- (30:1) -- (-30:1) --(0,-1)--(0,0);
}

\boxt{0}{0}{2}
\boxl{0}{0}{2}
\boxr{0}{0}{2}
\boxt{0}{-1}{2}
%\boxl{1}{0}{2}
%\boxr{1}{0}{2}
\boxt{0}{-2}{2}
%\boxl{1}{0}{2}
%\boxr{1}{0}{2}
\boxt{0}{-3}{2}
%\boxl{1}{0}{2}
%\boxr{1}{0}{2}
\boxt{0}{-4}{2}
%\boxl{1}{0}{2}
%\boxr{1}{0}{2}
\boxt{0}{-5}{2}
%\boxl{1}{0}{2}
%\boxr{1}{0}{2}
\boxt{0}{-6}{2}
%\boxl{1}{0}{2}
%\boxr{1}{0}{2}
\boxt{0}{-7}{2}
%\boxl{1}{0}{2}
%\boxr{1}{0}{2}
\boxt{0}{-8}{2}
%\boxl{1}{0}{2}
%\boxr{1}{0}{2}
\boxt{0}{-9}{2}
%\boxl{1}{0}{2}
%\boxr{1}{0}{2}
\boxt{0}{-10}{2}
%\boxl{1}{0}{2}
%\boxr{1}{0}{2}

\boxt{1}{0}{2}
\boxl{1}{0}{2}
\boxr{1}{0}{2}
\boxt{1}{-1}{2}
\boxl{1}{-1}{2}
%\boxr{1}{-1}{2}
\boxt{1}{-2}{2}
\boxl{1}{-2}{2}
%\boxr{1}{-1}{2}
\boxt{1}{-3}{2}
\boxl{1}{-3}{2}
%\boxr{1}{-1}{2}
\boxt{1}{-4}{2}
\boxl{1}{-4}{2}
%\boxr{1}{-1}{2}
\boxt{1}{-5}{2}
\boxl{1}{-5}{2}
%\boxr{1}{-1}{2}
\boxt{1}{-6}{2}
\boxl{1}{-6}{2}
%\boxr{1}{-1}{2}
\boxt{1}{-7}{2}
\boxl{1}{-7}{2}
%\boxr{1}{-1}{2}
\boxt{1}{-8}{2}
\boxl{1}{-8}{2}
%\boxr{1}{-1}{2}
\boxt{1}{-9}{2}
\boxl{1}{-9}{2}
%\boxr{1}{-1}{2}
\boxt{1}{-10}{2}
\boxl{1}{-10}{2}
%\boxr{1}{-1}{2}

%\boxt{-1}{0}{2}
%\boxl{-1}{0}{2}
\boxr{-1}{0}{2}
%\boxt{-2}{0}{2}
%\boxl{-1}{0}{2}
\boxr{-2}{0}{2}
%\boxt{-3}{0}{2}
%\boxl{-1}{0}{2}
\boxr{-3}{0}{2}
%\boxt{-4}{0}{2}
%\boxl{-1}{0}{2}
\boxr{-4}{0}{2}
%\boxt{-5}{0}{2}
%\boxl{-1}{0}{2}
\boxr{-5}{0}{2}
%\boxt{-6}{0}{2}
%\boxl{-1}{0}{2}
\boxr{-6}{0}{2}
%\boxt{-7}{0}{2}
%\boxl{-1}{0}{2}
\boxr{-7}{0}{2}
%\boxt{-8}{0}{2}
%\boxl{-1}{0}{2}
\boxr{-8}{0}{2}
%\boxt{-9}{0}{2}
%\boxl{-1}{0}{2}
\boxr{-9}{0}{2}
%\boxt{-10}{0}{2}
%\boxl{-1}{0}{2}
\boxr{-10}{0}{2}
%\boxt{-11}{0}{2}
%\boxl{-1}{0}{2}
\boxr{-11}{0}{2}

\boxt{-1}{0}{3}
\boxl{-1}{0}{3}
\boxr{-1}{0}{3}
\boxt{-2}{0}{3}
%\boxl{-2}{0}{3}
\boxr{-2}{0}{3}
\boxt{-3}{0}{3}
%\boxl{-2}{0}{3}
\boxr{-3}{0}{3}
\boxt{-4}{0}{3}
%\boxl{-2}{0}{3}
\boxr{-4}{0}{3}
\boxt{-5}{0}{3}
%\boxl{-2}{0}{3}
\boxr{-5}{0}{3}
\boxt{-6}{0}{3}
%\boxl{-2}{0}{3}
\boxr{-6}{0}{3}
\boxt{-7}{0}{3}
%\boxl{-2}{0}{3}
\boxr{-7}{0}{3}
\boxt{-8}{0}{3}
%\boxl{-2}{0}{3}
\boxr{-8}{0}{3}
\boxt{-9}{0}{3}
%\boxl{-2}{0}{3}
\boxr{-9}{0}{3}
\boxt{-10}{0}{3}
%\boxl{-2}{0}{3}
\boxr{-10}{0}{3}
\boxt{-11}{0}{3}
%\boxl{-2}{0}{3}
\boxr{-11}{0}{3}

\boxxt{0}{0}{2}
\boxxl{0}{0}{2}
\boxxr{0}{0}{2}

\boxxt{1}{0}{2}
\boxxl{1}{0}{2}
\boxxr{1}{0}{2}
\boxxt{1}{-1}{2}
\boxxl{1}{-1}{2}
%\boxxr{1}{-1}{2}
\boxxt{1}{-2}{2}
\boxxl{1}{-2}{2}
%\boxxr{1}{-1}{2}
\boxxt{1}{-3}{2}
\boxxl{1}{-3}{2}
%\boxxr{1}{-1}{2}
\boxxt{1}{-4}{2}
\boxxl{1}{-4}{2}
%\boxxr{1}{-1}{2}

\boxxt{-1}{0}{3}
\boxxl{-1}{0}{3}
\boxxr{-1}{0}{3}
\boxxt{-2}{0}{3}
%\boxxl{-2}{0}{3}
\boxxr{-2}{0}{3}
\boxxt{-3}{0}{3}
%\boxxl{-2}{0}{3}
\boxxr{-3}{0}{3}

\draw [thick,->] (0,4)-- (0,6);
\draw [thick,-] (0,0)-- (0,1);
\draw [thick,->] (0*1.74,-0*1.74*0.57735)-- (6,-6*0.57735);
\draw [thick,->]  (0,-0*0.57735)-- (-6,-6*0.57735);

\node [right] at (0,7) {$x_3$};
\node [right] at (6,-6*0.57735) {$x_2$};
\node [left] at  (-6,-6*0.57735) {$x_1$};
\node [above] at  (-9,-7*0.57735) {$ $};
\node [right] at (2,10) {$ $};
\node  at (6,6*0.57735) {$ $};
\node  at (-6,6*0.57735) {$ $};
\node  at (0,-6) {$ $};

%\shade [ball color=gray] (-2.6,-.6) circle [radius=0.4cm];
%\shade [ball color=gray] (1.8,-1.2) circle [radius=0.4cm];
%\shade [ball color=gray] (-0.85,-1.7) circle [radius=0.4cm];

\end{tikzpicture}

%\end{document} 

%% file: fig4.tex
% Plane partition
% Author: Jang Soo Kim
%\documentclass{minimal}
%\usepackage{tikz}
% Three counters
%\newcounter{x}
%\newcounter{y}
%\newcounter{z}

% The angles of x,y,z-axes
\newcommand\xaxis{210}
\newcommand\yaxis{-30}
\newcommand\zaxis{90}

% The top side of a cube
\newcommand\topside[3]{
  \fill[fill=cyan, draw=black,shift={(\xaxis:#1)},shift={(\yaxis:#2)},
  shift={(\zaxis:#3)}] (0,0) -- (30:1) -- (0,1) --(150:1)--(0,0);
}

% The left side of a cube
\newcommand\leftside[3]{
  \fill[fill=cyan, draw=black,shift={(\xaxis:#1)},shift={(\yaxis:#2)},
  shift={(\zaxis:#3)}] (0,0) -- (0,-1) -- (210:1) --(150:1)--(0,0);
}

% The right side of a cube
\newcommand\rightside[3]{
  \fill[fill=cyan, draw=black,shift={(\xaxis:#1)},shift={(\yaxis:#2)},
  shift={(\zaxis:#3)}] (0,0) -- (30:1) -- (-30:1) --(0,-1)--(0,0);
}

% The cube 
\newcommand\cube[3]{
  \topside{#1}{#2}{#3} \leftside{#1}{#2}{#3} \rightside{#1}{#2}{#3}
}

% Definition of \planepartition
% To draw the following plane partition, just write \planepartition{ {a, b, c}, {d,e} }.
%  a b c
%  d e
\newcommand\planepartition[1]{
 \setcounter{x}{-1}
  \foreach \a in {#1} {
    \addtocounter{x}{1}
    \setcounter{y}{-1}
    \foreach \b in \a {
      \addtocounter{y}{1}
      \setcounter{z}{-1}
      \foreach \c in {1,...,\b} {
        \addtocounter{z}{1}
        \cube{\value{x}}{\value{y}}{\value{z}}
      }
    }
  }
}

%\begin{document} 

\begin{tikzpicture}[scale=0.25]

\newcommand{\boxt}[3]{ 
\draw[fill=cyan, draw=black,shift={(\xaxis:#1)},shift={(\yaxis:#2)},
  shift={(\zaxis:#3)}] (0,0) -- (30:1) -- (0,1) --(150:1)--(0,0);
}

\newcommand{\boxl}[3]{ 
\draw[fill=cyan, draw=black,shift={(\xaxis:#1)},shift={(\yaxis:#2)},
  shift={(\zaxis:#3)}] (0,0) -- (0,-1) -- (210:1) --(150:1)--(0,0);
}

\newcommand{\boxr}[3]{ 
\draw[fill=cyan, draw=black,shift={(\xaxis:#1)},shift={(\yaxis:#2)},
  shift={(\zaxis:#3)}] (0,0) -- (30:1) -- (-30:1) --(0,-1)--(0,0);
}

\boxl{8}{-1}{0}
\boxl{8}{-1}{1}
\boxt{8}{-1}{1}
\boxt{7}{-1}{1}
\boxt{6}{-1}{1}
\boxt{5}{-1}{1}
\boxt{4}{-1}{1}
\boxt{3}{-1}{1}
\boxt{2}{-1}{1}
\boxt{1}{-1}{1}
\boxt{0}{-1}{1}
\boxt{-1}{-1}{1}
\boxt{-2}{-1}{1}
\boxt{-3}{-1}{1}
\boxt{-4}{-1}{1}
\boxt{-5}{-1}{1}
\boxt{-6}{-1}{1}
\boxt{-7}{-1}{1}
\boxt{-8}{-1}{1}

\boxl{8}{-2}{0}
\boxl{8}{-2}{1}
\boxt{8}{-2}{1}
\boxt{7}{-2}{1}
\boxt{6}{-2}{1}
\boxt{5}{-2}{1}
\boxt{4}{-2}{1}
\boxt{3}{-2}{1}
\boxt{2}{-2}{1}
\boxt{1}{-2}{1}
\boxt{0}{-2}{1}
\boxt{-1}{-2}{1}
\boxt{-2}{-2}{1}
\boxt{-3}{-2}{1}
\boxt{-4}{-2}{1}
\boxt{-5}{-2}{1}
\boxt{-6}{-2}{1}
\boxt{-7}{-2}{1}
\boxt{-8}{-2}{1}

\boxl{8}{-3}{0}
\boxl{8}{-3}{1}
\boxt{8}{-3}{1}
\boxt{7}{-3}{1}
\boxt{6}{-3}{1}
\boxt{5}{-3}{1}
\boxt{4}{-3}{1}
\boxt{3}{-3}{1}
\boxt{2}{-3}{1}
\boxt{1}{-3}{1}
\boxt{0}{-3}{1}
\boxt{-1}{-3}{1}
\boxt{-2}{-3}{1}
\boxt{-3}{-3}{1}
\boxt{-4}{-3}{1}
\boxt{-5}{-3}{1}
\boxt{-6}{-3}{1}
\boxt{-7}{-3}{1}
\boxt{-8}{-3}{1}

\boxl{8}{-4}{0}
\boxl{8}{-4}{1}
\boxt{8}{-4}{1}
\boxt{7}{-4}{1}
\boxt{6}{-4}{1}
\boxt{5}{-4}{1}
\boxt{4}{-4}{1}
\boxt{3}{-4}{1}
\boxt{2}{-4}{1}
\boxt{1}{-4}{1}
\boxt{0}{-4}{1}
\boxt{-1}{-4}{1}
\boxt{-2}{-4}{1}
\boxt{-3}{-4}{1}
\boxt{-4}{-4}{1}
\boxt{-5}{-4}{1}
\boxt{-6}{-4}{1}
\boxt{-7}{-4}{1}
\boxt{-8}{-4}{1}

\boxl{8}{-5}{0}
\boxl{8}{-5}{1}
\boxt{8}{-5}{1}
\boxt{7}{-5}{1}
\boxt{6}{-5}{1}
\boxt{5}{-5}{1}
\boxt{4}{-5}{1}
\boxt{3}{-5}{1}
\boxt{2}{-5}{1}
\boxt{1}{-5}{1}
\boxt{0}{-5}{1}
\boxt{-1}{-5}{1}
\boxt{-2}{-5}{1}
\boxt{-3}{-5}{1}
\boxt{-4}{-5}{1}
\boxt{-5}{-5}{1}
\boxt{-6}{-5}{1}
\boxt{-7}{-5}{1}
\boxt{-8}{-5}{1}

\boxl{8}{-6}{0}
\boxl{8}{-6}{1}
\boxt{8}{-6}{1}
\boxt{7}{-6}{1}
\boxt{6}{-6}{1}
\boxt{5}{-6}{1}
\boxt{4}{-6}{1}
\boxt{3}{-6}{1}
\boxt{2}{-6}{1}
\boxt{1}{-6}{1}
\boxt{0}{-6}{1}
\boxt{-1}{-6}{1}
\boxt{-2}{-6}{1}
\boxt{-3}{-6}{1}
\boxt{-4}{-6}{1}
\boxt{-5}{-6}{1}
\boxt{-6}{-6}{1}
\boxt{-7}{-6}{1}
\boxt{-8}{-6}{1}

\boxl{8}{-7}{0}
\boxl{8}{-7}{1}
\boxt{8}{-7}{1}
\boxt{7}{-7}{1}
\boxt{6}{-7}{1}
\boxt{5}{-7}{1}
\boxt{4}{-7}{1}
\boxt{3}{-7}{1}
\boxt{2}{-7}{1}
\boxt{1}{-7}{1}
\boxt{0}{-7}{1}
\boxt{-1}{-7}{1}
\boxt{-2}{-7}{1}
\boxt{-3}{-7}{1}
\boxt{-4}{-7}{1}
\boxt{-5}{-7}{1}
\boxt{-6}{-7}{1}
\boxt{-7}{-7}{1}
\boxt{-8}{-7}{1}

\boxl{8}{-8}{0}
\boxl{8}{-8}{1}
\boxt{8}{-8}{1}
\boxt{7}{-8}{1}
\boxt{6}{-8}{1}
\boxt{5}{-8}{1}
\boxt{4}{-8}{1}
\boxt{3}{-8}{1}
\boxt{2}{-8}{1}
\boxt{1}{-8}{1}
\boxt{0}{-8}{1}
\boxt{-1}{-8}{1}
\boxt{-2}{-8}{1}
\boxt{-3}{-8}{1}
\boxt{-4}{-8}{1}
\boxt{-5}{-8}{1}
\boxt{-6}{-8}{1}
\boxt{-7}{-8}{1}
\boxt{-8}{-8}{1}

\boxl{8}{-9}{0}
\boxl{8}{-9}{1}
\boxt{8}{-9}{1}
\boxt{7}{-9}{1}
\boxt{6}{-9}{1}
\boxt{5}{-9}{1}
\boxt{4}{-9}{1}
\boxt{3}{-9}{1}
\boxt{2}{-9}{1}
\boxt{1}{-9}{1}
\boxt{0}{-9}{1}
\boxt{-1}{-9}{1}
\boxt{-2}{-9}{1}
\boxt{-3}{-9}{1}
\boxt{-4}{-9}{1}
\boxt{-5}{-9}{1}
\boxt{-6}{-9}{1}
\boxt{-7}{-9}{1}
\boxt{-8}{-9}{1}

\boxl{-1}{0}{0}
\boxl{-1}{0}{1}
\boxt{-1}{0}{1}
\boxt{-2}{0}{1}
\boxt{-2}{0}{1}
\boxt{-3}{0}{1}
\boxt{-4}{0}{1}
\boxt{-5}{0}{1}
\boxt{-6}{0}{1}
\boxt{-7}{0}{1}
\boxt{-8}{0}{1}

\boxl{-1}{1}{0}
\boxl{-1}{1}{1}
\boxt{-1}{1}{1}
\boxt{-2}{1}{1}
\boxt{-2}{1}{1}
\boxt{-3}{1}{1}
\boxt{-4}{1}{1}
\boxt{-5}{1}{1}
\boxt{-6}{1}{1}
\boxt{-7}{1}{1}
\boxt{-8}{1}{1}

\boxl{-1}{2}{0}
\boxl{-1}{2}{1}
\boxt{-1}{2}{1}
\boxt{-2}{2}{1}
\boxt{-2}{2}{1}
\boxt{-3}{2}{1}
\boxt{-4}{2}{1}
\boxt{-5}{2}{1}
\boxt{-6}{2}{1}
\boxt{-7}{2}{1}
\boxt{-8}{2}{1}

\boxl{-1}{3}{0}
\boxl{-1}{3}{1}
\boxt{-1}{3}{1}
\boxt{-2}{3}{1}
\boxt{-2}{3}{1}
\boxt{-3}{3}{1}
\boxt{-4}{3}{1}
\boxt{-5}{3}{1}
\boxt{-6}{3}{1}
\boxt{-7}{3}{1}
\boxt{-8}{3}{1}

\boxl{-1}{4}{0}
\boxl{-1}{4}{1}
\boxt{-1}{4}{1}
\boxt{-2}{4}{1}
\boxt{-2}{4}{1}
\boxt{-3}{4}{1}
\boxt{-4}{4}{1}
\boxt{-5}{4}{1}
\boxt{-6}{4}{1}
\boxt{-7}{4}{1}
\boxt{-8}{4}{1}

\boxl{-1}{5}{0}
\boxl{-1}{5}{1}
\boxt{-1}{5}{1}
\boxt{-2}{5}{1}
\boxt{-2}{5}{1}
\boxt{-3}{5}{1}
\boxt{-4}{5}{1}
\boxt{-5}{5}{1}
\boxt{-6}{5}{1}
\boxt{-7}{5}{1}
\boxt{-8}{5}{1}

\boxl{-1}{6}{0}
\boxl{-1}{6}{1}
\boxt{-1}{6}{1}
\boxt{-2}{6}{1}
\boxt{-2}{6}{1}
\boxt{-3}{6}{1}
\boxt{-4}{6}{1}
\boxt{-5}{6}{1}
\boxt{-6}{6}{1}
\boxt{-7}{6}{1}
\boxt{-8}{6}{1}

\boxl{-1}{7}{0}
\boxl{-1}{7}{1}
\boxt{-1}{7}{1}
\boxt{-2}{7}{1}
\boxt{-2}{7}{1}
\boxt{-3}{7}{1}
\boxt{-4}{7}{1}
\boxt{-5}{7}{1}
\boxt{-6}{7}{1}
\boxt{-7}{7}{1}
\boxt{-8}{7}{1}

\boxl{-1}{8}{0}
\boxl{-1}{8}{1}
\boxt{-1}{8}{1}
\boxt{-2}{8}{1}
\boxt{-2}{8}{1}
\boxt{-3}{8}{1}
\boxt{-4}{8}{1}
\boxt{-5}{8}{1}
\boxt{-6}{8}{1}
\boxt{-7}{8}{1}
\boxt{-8}{8}{1}

\boxr{0}{-1}{0}
\boxr{0}{-1}{1}
\boxr{1}{-1}{0}
\boxr{1}{-1}{1}
\boxr{2}{-1}{0}
\boxr{2}{-1}{1}
\boxr{3}{-1}{0}
\boxr{3}{-1}{1}
\boxr{4}{-1}{0}
\boxr{4}{-1}{1}
\boxr{5}{-1}{0}
\boxr{5}{-1}{1}
\boxr{6}{-1}{0}
\boxr{6}{-1}{1}
\boxr{7}{-1}{0}
\boxr{7}{-1}{1}
\boxr{8}{-1}{0}
\boxr{8}{-1}{1}

\boxr{-1}{8}{0}
\boxr{-1}{8}{1}
\boxr{-2}{8}{0}
\boxr{-2}{8}{1}
\boxr{-3}{8}{0}
\boxr{-3}{8}{1}
\boxr{-4}{8}{0}
\boxr{-4}{8}{1}
\boxr{-5}{8}{0}
\boxr{-5}{8}{1}
\boxr{-6}{8}{0}
\boxr{-6}{8}{1}
\boxr{-7}{8}{0}
\boxr{-7}{8}{1}
\boxr{-8}{8}{0}
\boxr{-8}{8}{1}

\draw [thick,->] (0,0)-- (0,13);
\draw [thick,->] (0*1.74,-0*1.74*0.57735)-- (11,-11*0.57735);
\draw [thick,->]  (-0,-0*0.57735)-- (-11,-11*0.57735);

\node [right] at (0,13) {$x_3$};
\node [right] at (11,-11*0.57735) {$x_2$};
\node [left] at  (-11,-11*0.57735) {$x_1$};
\node [above] at  (-9,-7*0.57735) {$ $};
\node [right] at (2,10) {$ $};
\node  at (6,6*0.57735) {$ $};
\node  at (-6,6*0.57735) {$ $};
\node  at (0,-6) {$ $};

%\shade [ball color=gray] (-2.6,-.6) circle [radius=0.4cm];
%\shade [ball color=gray] (1.8,-1.2) circle [radius=0.4cm];
%\shade [ball color=gray] (-0.85,-1.7) circle [radius=0.4cm];

\end{tikzpicture}

%\end{document} 

%% file: fig4a.tex
% Plane partition
% Author: Jang Soo Kim
%\documentclass{minimal}
%\usepackage{tikz}
% Three counters
%\newcounter{x}
%\newcounter{y}
%\newcounter{z}

% The angles of x,y,z-axes
\newcommand\xaxis{210}
\newcommand\yaxis{-30}
\newcommand\zaxis{90}

% The top side of a cube
\newcommand\topside[3]{
  \fill[fill=cyan, draw=black,shift={(\xaxis:#1)},shift={(\yaxis:#2)},
  shift={(\zaxis:#3)}] (0,0) -- (30:1) -- (0,1) --(150:1)--(0,0);
}

% The left side of a cube
\newcommand\leftside[3]{
  \fill[fill=cyan, draw=black,shift={(\xaxis:#1)},shift={(\yaxis:#2)},
  shift={(\zaxis:#3)}] (0,0) -- (0,-1) -- (210:1) --(150:1)--(0,0);
}

% The right side of a cube
\newcommand\rightside[3]{
  \fill[fill=cyan, draw=black,shift={(\xaxis:#1)},shift={(\yaxis:#2)},
  shift={(\zaxis:#3)}] (0,0) -- (30:1) -- (-30:1) --(0,-1)--(0,0);
}

% The cube 
\newcommand\cube[3]{
  \topside{#1}{#2}{#3} \leftside{#1}{#2}{#3} \rightside{#1}{#2}{#3}
}

% Definition of \planepartition
% To draw the following plane partition, just write \planepartition{ {a, b, c}, {d,e} }.
%  a b c
%  d e
\newcommand\planepartition[1]{
 \setcounter{x}{-1}
  \foreach \a in {#1} {
    \addtocounter{x}{1}
    \setcounter{y}{-1}
    \foreach \b in \a {
      \addtocounter{y}{1}
      \setcounter{z}{-1}
      \foreach \c in {1,...,\b} {
        \addtocounter{z}{1}
        \cube{\value{x}}{\value{y}}{\value{z}}
      }
    }
  }
}

%\begin{document} 

\begin{tikzpicture}[scale=0.25]

\newcommand{\boxt}[3]{ 
\draw[fill=cyan, draw=black,shift={(\xaxis:#1)},shift={(\yaxis:#2)},
  shift={(\zaxis:#3)}] (0,0) -- (30:1) -- (0,1) --(150:1)--(0,0);
}

\newcommand{\boxl}[3]{ 
\draw[fill=cyan, draw=black,shift={(\xaxis:#1)},shift={(\yaxis:#2)},
  shift={(\zaxis:#3)}] (0,0) -- (0,-1) -- (210:1) --(150:1)--(0,0);
}

\newcommand{\boxr}[3]{ 
\draw[fill=cyan, draw=black,shift={(\xaxis:#1)},shift={(\yaxis:#2)},
  shift={(\zaxis:#3)}] (0,0) -- (30:1) -- (-30:1) --(0,-1)--(0,0);
}

\newcommand{\boxxt}[3]{ 
\draw[fill=red, draw=black,shift={(\xaxis:#1)},shift={(\yaxis:#2)},
  shift={(\zaxis:#3)}] (0,0) -- (30:1) -- (0,1) --(150:1)--(0,0);
}

\newcommand{\boxxl}[3]{ 
\draw[fill=red, draw=black,shift={(\xaxis:#1)},shift={(\yaxis:#2)},
  shift={(\zaxis:#3)}] (0,0) -- (0,-1) -- (210:1) --(150:1)--(0,0);
}

\newcommand{\boxxr}[3]{ 
\draw[fill=red, draw=black,shift={(\xaxis:#1)},shift={(\yaxis:#2)},
  shift={(\zaxis:#3)}] (0,0) -- (30:1) -- (-30:1) --(0,-1)--(0,0);
}

\boxl{8}{-1}{0}
\boxl{8}{-1}{1}
\boxt{8}{-1}{1}
\boxt{7}{-1}{1}
\boxt{6}{-1}{1}
\boxt{5}{-1}{1}
\boxt{4}{-1}{1}
\boxt{3}{-1}{1}
\boxt{2}{-1}{1}
\boxt{1}{-1}{1}
\boxt{0}{-1}{1}
\boxt{-1}{-1}{1}
\boxt{-2}{-1}{1}
\boxt{-3}{-1}{1}
\boxt{-4}{-1}{1}
\boxt{-5}{-1}{1}
\boxt{-6}{-1}{1}
\boxt{-7}{-1}{1}
\boxt{-8}{-1}{1}

\boxl{8}{-2}{0}
\boxl{8}{-2}{1}
\boxt{8}{-2}{1}
\boxt{7}{-2}{1}
\boxt{6}{-2}{1}
\boxt{5}{-2}{1}
\boxt{4}{-2}{1}
\boxt{3}{-2}{1}
\boxt{2}{-2}{1}
\boxt{1}{-2}{1}
\boxt{0}{-2}{1}
\boxt{-1}{-2}{1}
\boxt{-2}{-2}{1}
\boxt{-3}{-2}{1}
\boxt{-4}{-2}{1}
\boxt{-5}{-2}{1}
\boxt{-6}{-2}{1}
\boxt{-7}{-2}{1}
\boxt{-8}{-2}{1}

\boxl{8}{-3}{0}
\boxl{8}{-3}{1}
\boxt{8}{-3}{1}
\boxt{7}{-3}{1}
\boxt{6}{-3}{1}
\boxt{5}{-3}{1}
\boxt{4}{-3}{1}
\boxt{3}{-3}{1}
\boxt{2}{-3}{1}
\boxt{1}{-3}{1}
\boxt{0}{-3}{1}
\boxt{-1}{-3}{1}
\boxt{-2}{-3}{1}
\boxt{-3}{-3}{1}
\boxt{-4}{-3}{1}
\boxt{-5}{-3}{1}
\boxt{-6}{-3}{1}
\boxt{-7}{-3}{1}
\boxt{-8}{-3}{1}

\boxl{8}{-4}{0}
\boxl{8}{-4}{1}
\boxt{8}{-4}{1}
\boxt{7}{-4}{1}
\boxt{6}{-4}{1}
\boxt{5}{-4}{1}
\boxt{4}{-4}{1}
\boxt{3}{-4}{1}
\boxt{2}{-4}{1}
\boxt{1}{-4}{1}
\boxt{0}{-4}{1}
\boxt{-1}{-4}{1}
\boxt{-2}{-4}{1}
\boxt{-3}{-4}{1}
\boxt{-4}{-4}{1}
\boxt{-5}{-4}{1}
\boxt{-6}{-4}{1}
\boxt{-7}{-4}{1}
\boxt{-8}{-4}{1}

\boxl{8}{-5}{0}
\boxl{8}{-5}{1}
\boxt{8}{-5}{1}
\boxt{7}{-5}{1}
\boxt{6}{-5}{1}
\boxt{5}{-5}{1}
\boxt{4}{-5}{1}
\boxt{3}{-5}{1}
\boxt{2}{-5}{1}
\boxt{1}{-5}{1}
\boxt{0}{-5}{1}
\boxt{-1}{-5}{1}
\boxt{-2}{-5}{1}
\boxt{-3}{-5}{1}
\boxt{-4}{-5}{1}
\boxt{-5}{-5}{1}
\boxt{-6}{-5}{1}
\boxt{-7}{-5}{1}
\boxt{-8}{-5}{1}

\boxl{8}{-6}{0}
\boxl{8}{-6}{1}
\boxt{8}{-6}{1}
\boxt{7}{-6}{1}
\boxt{6}{-6}{1}
\boxt{5}{-6}{1}
\boxt{4}{-6}{1}
\boxt{3}{-6}{1}
\boxt{2}{-6}{1}
\boxt{1}{-6}{1}
\boxt{0}{-6}{1}
\boxt{-1}{-6}{1}
\boxt{-2}{-6}{1}
\boxt{-3}{-6}{1}
\boxt{-4}{-6}{1}
\boxt{-5}{-6}{1}
\boxt{-6}{-6}{1}
\boxt{-7}{-6}{1}
\boxt{-8}{-6}{1}

\boxl{8}{-7}{0}
\boxl{8}{-7}{1}
\boxt{8}{-7}{1}
\boxt{7}{-7}{1}
\boxt{6}{-7}{1}
\boxt{5}{-7}{1}
\boxt{4}{-7}{1}
\boxt{3}{-7}{1}
\boxt{2}{-7}{1}
\boxt{1}{-7}{1}
\boxt{0}{-7}{1}
\boxt{-1}{-7}{1}
\boxt{-2}{-7}{1}
\boxt{-3}{-7}{1}
\boxt{-4}{-7}{1}
\boxt{-5}{-7}{1}
\boxt{-6}{-7}{1}
\boxt{-7}{-7}{1}
\boxt{-8}{-7}{1}

\boxl{8}{-8}{0}
\boxl{8}{-8}{1}
\boxt{8}{-8}{1}
\boxt{7}{-8}{1}
\boxt{6}{-8}{1}
\boxt{5}{-8}{1}
\boxt{4}{-8}{1}
\boxt{3}{-8}{1}
\boxt{2}{-8}{1}
\boxt{1}{-8}{1}
\boxt{0}{-8}{1}
\boxt{-1}{-8}{1}
\boxt{-2}{-8}{1}
\boxt{-3}{-8}{1}
\boxt{-4}{-8}{1}
\boxt{-5}{-8}{1}
\boxt{-6}{-8}{1}
\boxt{-7}{-8}{1}
\boxt{-8}{-8}{1}

\boxl{8}{-9}{0}
\boxl{8}{-9}{1}
\boxt{8}{-9}{1}
\boxt{7}{-9}{1}
\boxt{6}{-9}{1}
\boxt{5}{-9}{1}
\boxt{4}{-9}{1}
\boxt{3}{-9}{1}
\boxt{2}{-9}{1}
\boxt{1}{-9}{1}
\boxt{0}{-9}{1}
\boxt{-1}{-9}{1}
\boxt{-2}{-9}{1}
\boxt{-3}{-9}{1}
\boxt{-4}{-9}{1}
\boxt{-5}{-9}{1}
\boxt{-6}{-9}{1}
\boxt{-7}{-9}{1}
\boxt{-8}{-9}{1}

\boxl{-1}{0}{0}
\boxl{-1}{0}{1}
\boxt{-1}{0}{1}
\boxt{-2}{0}{1}
\boxt{-2}{0}{1}
\boxt{-3}{0}{1}
\boxt{-4}{0}{1}
\boxt{-5}{0}{1}
\boxt{-6}{0}{1}
\boxt{-7}{0}{1}
\boxt{-8}{0}{1}

\boxl{-1}{1}{0}
\boxl{-1}{1}{1}
\boxt{-1}{1}{1}
\boxt{-2}{1}{1}
\boxt{-2}{1}{1}
\boxt{-3}{1}{1}
\boxt{-4}{1}{1}
\boxt{-5}{1}{1}
\boxt{-6}{1}{1}
\boxt{-7}{1}{1}
\boxt{-8}{1}{1}

\boxl{-1}{2}{0}
\boxl{-1}{2}{1}
\boxt{-1}{2}{1}
\boxt{-2}{2}{1}
\boxt{-2}{2}{1}
\boxt{-3}{2}{1}
\boxt{-4}{2}{1}
\boxt{-5}{2}{1}
\boxt{-6}{2}{1}
\boxt{-7}{2}{1}
\boxt{-8}{2}{1}

\boxl{-1}{3}{0}
\boxl{-1}{3}{1}
\boxt{-1}{3}{1}
\boxt{-2}{3}{1}
\boxt{-2}{3}{1}
\boxt{-3}{3}{1}
\boxt{-4}{3}{1}
\boxt{-5}{3}{1}
\boxt{-6}{3}{1}
\boxt{-7}{3}{1}
\boxt{-8}{3}{1}

\boxl{-1}{4}{0}
\boxl{-1}{4}{1}
\boxt{-1}{4}{1}
\boxt{-2}{4}{1}
\boxt{-2}{4}{1}
\boxt{-3}{4}{1}
\boxt{-4}{4}{1}
\boxt{-5}{4}{1}
\boxt{-6}{4}{1}
\boxt{-7}{4}{1}
\boxt{-8}{4}{1}

\boxl{-1}{5}{0}
\boxl{-1}{5}{1}
\boxt{-1}{5}{1}
\boxt{-2}{5}{1}
\boxt{-2}{5}{1}
\boxt{-3}{5}{1}
\boxt{-4}{5}{1}
\boxt{-5}{5}{1}
\boxt{-6}{5}{1}
\boxt{-7}{5}{1}
\boxt{-8}{5}{1}

\boxl{-1}{6}{0}
\boxl{-1}{6}{1}
\boxt{-1}{6}{1}
\boxt{-2}{6}{1}
\boxt{-2}{6}{1}
\boxt{-3}{6}{1}
\boxt{-4}{6}{1}
\boxt{-5}{6}{1}
\boxt{-6}{6}{1}
\boxt{-7}{6}{1}
\boxt{-8}{6}{1}

\boxl{-1}{7}{0}
\boxl{-1}{7}{1}
\boxt{-1}{7}{1}
\boxt{-2}{7}{1}
\boxt{-2}{7}{1}
\boxt{-3}{7}{1}
\boxt{-4}{7}{1}
\boxt{-5}{7}{1}
\boxt{-6}{7}{1}
\boxt{-7}{7}{1}
\boxt{-8}{7}{1}

\boxl{-1}{8}{0}
\boxl{-1}{8}{1}
\boxt{-1}{8}{1}
\boxt{-2}{8}{1}
\boxt{-2}{8}{1}
\boxt{-3}{8}{1}
\boxt{-4}{8}{1}
\boxt{-5}{8}{1}
\boxt{-6}{8}{1}
\boxt{-7}{8}{1}
\boxt{-8}{8}{1}

\boxr{0}{-1}{0}
\boxr{0}{-1}{1}
\boxr{1}{-1}{0}
\boxr{1}{-1}{1}
\boxr{2}{-1}{0}
\boxr{2}{-1}{1}
\boxr{3}{-1}{0}
\boxr{3}{-1}{1}
\boxr{4}{-1}{0}
\boxr{4}{-1}{1}
\boxr{5}{-1}{0}
\boxr{5}{-1}{1}
\boxr{6}{-1}{0}
\boxr{6}{-1}{1}
\boxr{7}{-1}{0}
\boxr{7}{-1}{1}
\boxr{8}{-1}{0}
\boxr{8}{-1}{1}

\boxr{-1}{8}{0}
\boxr{-1}{8}{1}
\boxr{-2}{8}{0}
\boxr{-2}{8}{1}
\boxr{-3}{8}{0}
\boxr{-3}{8}{1}
\boxr{-4}{8}{0}
\boxr{-4}{8}{1}
\boxr{-5}{8}{0}
\boxr{-5}{8}{1}
\boxr{-6}{8}{0}
\boxr{-6}{8}{1}
\boxr{-7}{8}{0}
\boxr{-7}{8}{1}
\boxr{-8}{8}{0}
\boxr{-8}{8}{1}

\boxxt{-1}{-1}{1}
\boxxt{0}{-1}{1}
\boxxt{1}{-1}{1}
\boxxt{2}{-1}{1}
\boxxt{3}{-1}{1}
\boxxt{4}{-1}{1}
\boxxt{5}{-1}{1}
\boxxt{6}{-1}{1}
\boxxt{7}{-1}{1}
\boxxt{8}{-1}{1}

\boxxt{-1}{0}{1}
\boxxt{-1}{1}{1}
\boxxt{-1}{2}{1}
\boxxt{-1}{3}{1}
\boxxt{-1}{4}{1}
\boxxt{-1}{5}{1}
\boxxt{-1}{6}{1}
\boxxt{-1}{7}{1}
\boxxt{-1}{8}{1}

\boxxr{0}{-1}{1}
\boxxr{1}{-1}{1}
\boxxr{2}{-1}{1}
\boxxr{3}{-1}{1}
\boxxr{4}{-1}{1}
\boxxr{5}{-1}{1}
\boxxr{6}{-1}{1}
\boxxr{7}{-1}{1}
\boxxr{8}{-1}{1}

\boxxr{1}{-1}{0}
\boxxr{2}{-1}{0}
\boxxr{3}{-1}{0}
\boxxr{4}{-1}{0}
\boxxr{5}{-1}{0}
\boxxr{6}{-1}{0}
\boxxr{7}{-1}{0}
\boxxr{8}{-1}{0}

\boxxr{-1}{8}{1}

\boxxl{8}{-1}{0}
\boxxl{8}{-1}{1}

\boxxl{-1}{0}{1}
\boxxl{-1}{1}{1}
\boxxl{-1}{2}{1}
\boxxl{-1}{3}{1}
\boxxl{-1}{4}{1}
\boxxl{-1}{5}{1}
\boxxl{-1}{6}{1}
\boxxl{-1}{7}{1}
\boxxl{-1}{8}{1}

\draw [thick,->] (0,0)-- (0,13);
\draw [thick,->] (0*1.74,-0*1.74*0.57735)-- (11,-11*0.57735);
\draw [thick,->]  (-0,-0*0.57735)-- (-11,-11*0.57735);

\node [right] at (0,13) {$x_3$};
\node [right] at (11,-11*0.57735) {$x_2$};
\node [left] at  (-11,-11*0.57735) {$x_1$};
\node [above] at  (-9,-7*0.57735) {$ $};
\node [right] at (2,10) {$ $};
\node  at (6,6*0.57735) {$ $};
\node  at (-6,6*0.57735) {$ $};
\node  at (0,-6) {$ $};

%\shade [ball color=gray] (-2.6,-.6) circle [radius=0.4cm];
%\shade [ball color=gray] (1.8,-1.2) circle [radius=0.4cm];
%\shade [ball color=gray] (-0.85,-1.7) circle [radius=0.4cm];

\end{tikzpicture}

%\end{document} 

%% file: fig5a.tex
% Plane partition
% Author: Jang Soo Kim
%\documentclass{minimal}
%\usepackage{tikz}
% Three counters
%\newcounter{x}
%\newcounter{y}
%\newcounter{z}

% The angles of x,y,z-axes
\newcommand\xaxis{210}
\newcommand\yaxis{-30}
\newcommand\zaxis{90}

% The top side of a cube
\newcommand\topside[3]{
  \fill[fill=cyan, draw=black,shift={(\xaxis:#1)},shift={(\yaxis:#2)},
  shift={(\zaxis:#3)}] (0,0) -- (30:1) -- (0,1) --(150:1)--(0,0);
}

% The left side of a cube
\newcommand\leftside[3]{
  \fill[fill=cyan, draw=black,shift={(\xaxis:#1)},shift={(\yaxis:#2)},
  shift={(\zaxis:#3)}] (0,0) -- (0,-1) -- (210:1) --(150:1)--(0,0);
}

% The right side of a cube
\newcommand\rightside[3]{
  \fill[fill=cyan, draw=black,shift={(\xaxis:#1)},shift={(\yaxis:#2)},
  shift={(\zaxis:#3)}] (0,0) -- (30:1) -- (-30:1) --(0,-1)--(0,0);
}

% The cube 
\newcommand\cube[3]{
  \topside{#1}{#2}{#3} \leftside{#1}{#2}{#3} \rightside{#1}{#2}{#3}
}

% Definition of \planepartition
% To draw the following plane partition, just write \planepartition{ {a, b, c}, {d,e} }.
%  a b c
%  d e
\newcommand\planepartition[1]{
 \setcounter{x}{-1}
  \foreach \a in {#1} {
    \addtocounter{x}{1}
    \setcounter{y}{-1}
    \foreach \b in \a {
      \addtocounter{y}{1}
      \setcounter{z}{-1}
      \foreach \c in {1,...,\b} {
        \addtocounter{z}{1}
        \cube{\value{x}}{\value{y}}{\value{z}}
      }
    }
  }
}

%\begin{document} 

\begin{tikzpicture}[scale=0.22]

\newcommand{\boxt}[3]{ 
\draw[fill=cyan, draw=black,shift={(\xaxis:#1)},shift={(\yaxis:#2)},
  shift={(\zaxis:#3)}] (0,0) -- (30:1) -- (0,1) --(150:1)--(0,0);
}

\newcommand{\boxl}[3]{ 
\draw[fill=cyan, draw=black,shift={(\xaxis:#1)},shift={(\yaxis:#2)},
  shift={(\zaxis:#3)}] (0,0) -- (0,-1) -- (210:1) --(150:1)--(0,0);
}

\newcommand{\boxr}[3]{ 
\draw[fill=cyan, draw=black,shift={(\xaxis:#1)},shift={(\yaxis:#2)},
  shift={(\zaxis:#3)}] (0,0) -- (30:1) -- (-30:1) --(0,-1)--(0,0);
}

\newcommand{\boxxt}[3]{ 
\draw[fill=red, draw=black,shift={(\xaxis:#1)},shift={(\yaxis:#2)},
  shift={(\zaxis:#3)}] (0,0) -- (30:1) -- (0,1) --(150:1)--(0,0);
}

\newcommand{\boxxl}[3]{ 
\draw[fill=red, draw=black,shift={(\xaxis:#1)},shift={(\yaxis:#2)},
  shift={(\zaxis:#3)}] (0,0) -- (0,-1) -- (210:1) --(150:1)--(0,0);
}

\newcommand{\boxxr}[3]{ 
\draw[fill=red, draw=black,shift={(\xaxis:#1)},shift={(\yaxis:#2)},
  shift={(\zaxis:#3)}] (0,0) -- (30:1) -- (-30:1) --(0,-1)--(0,0);
}

\newcommand{\boxxxt}[3]{ 
\draw[fill=cyan, draw=black,shift={(\xaxis:#1)},shift={(\yaxis:#2)},
  shift={(\zaxis:#3)}] (0,0) -- (30:1) -- (0,1) --(150:1)--(0,0);
}

\newcommand{\boxxxl}[3]{ 
\draw[fill=cyan, draw=black,shift={(\xaxis:#1)},shift={(\yaxis:#2)},
  shift={(\zaxis:#3)}] (0,0) -- (0,-1) -- (210:1) --(150:1)--(0,0);
}

\newcommand{\boxxxr}[3]{ 
\draw[fill=cyan, draw=black,shift={(\xaxis:#1)},shift={(\yaxis:#2)},
  shift={(\zaxis:#3)}] (0,0) -- (30:1) -- (-30:1) --(0,-1)--(0,0);
}

\boxl{8}{-1}{0}
\boxt{8}{-1}{0}
\boxt{7}{-1}{0}
\boxt{6}{-1}{0}
\boxt{5}{-1}{0}
\boxt{4}{-1}{0}
\boxt{3}{-1}{0}
\boxt{2}{-1}{0}
\boxt{1}{-1}{0}
\boxt{0}{-1}{0}
\boxt{-1}{-1}{0}
\boxt{-2}{-1}{0}
\boxt{-3}{-1}{0}
\boxt{-4}{-1}{0}
\boxt{-5}{-1}{0}
\boxt{-6}{-1}{0}
\boxt{-7}{-1}{0}
\boxt{-8}{-1}{0}

\boxl{8}{-2}{0}
\boxt{8}{-2}{0}
\boxt{7}{-2}{0}
\boxt{6}{-2}{0}
\boxt{5}{-2}{0}
\boxt{4}{-2}{0}
\boxt{3}{-2}{0}
\boxt{2}{-2}{0}
\boxt{1}{-2}{0}
\boxt{0}{-2}{0}
\boxt{-1}{-2}{0}
\boxt{-2}{-2}{0}
\boxt{-3}{-2}{0}
\boxt{-4}{-2}{0}
\boxt{-5}{-2}{0}
\boxt{-6}{-2}{0}
\boxt{-7}{-2}{0}
\boxt{-8}{-2}{0}

\boxl{8}{-3}{0}
\boxt{8}{-3}{0}
\boxt{7}{-3}{0}
\boxt{6}{-3}{0}
\boxt{5}{-3}{0}
\boxt{4}{-3}{0}
\boxt{3}{-3}{0}
\boxt{2}{-3}{0}
\boxt{1}{-3}{0}
\boxt{0}{-3}{0}
\boxt{-1}{-3}{0}
\boxt{-2}{-3}{0}
\boxt{-3}{-3}{0}
\boxt{-4}{-3}{0}
\boxt{-5}{-3}{0}
\boxt{-6}{-3}{0}
\boxt{-7}{-3}{0}
\boxt{-8}{-3}{0}

\boxl{8}{-4}{0}
\boxt{8}{-4}{0}
\boxt{7}{-4}{0}
\boxt{6}{-4}{0}
\boxt{5}{-4}{0}
\boxt{4}{-4}{0}
\boxt{3}{-4}{0}
\boxt{2}{-4}{0}
\boxt{1}{-4}{0}
\boxt{0}{-4}{0}
\boxt{-1}{-4}{0}
\boxt{-2}{-4}{0}
\boxt{-3}{-4}{0}
\boxt{-4}{-4}{0}
\boxt{-5}{-4}{0}
\boxt{-6}{-4}{0}
\boxt{-7}{-4}{0}
\boxt{-8}{-4}{0}

\boxl{8}{-5}{0}
\boxt{8}{-5}{0}
\boxt{7}{-5}{0}
\boxt{6}{-5}{0}
\boxt{5}{-5}{0}
\boxt{4}{-5}{0}
\boxt{3}{-5}{0}
\boxt{2}{-5}{0}
\boxt{1}{-5}{0}
\boxt{0}{-5}{0}
\boxt{-1}{-5}{0}
\boxt{-2}{-5}{0}
\boxt{-3}{-5}{0}
\boxt{-4}{-5}{0}
\boxt{-5}{-5}{0}
\boxt{-6}{-5}{0}
\boxt{-7}{-5}{0}
\boxt{-8}{-5}{0}

\boxl{8}{-6}{0}
\boxt{8}{-6}{0}
\boxt{7}{-6}{0}
\boxt{6}{-6}{0}
\boxt{5}{-6}{0}
\boxt{4}{-6}{0}
\boxt{3}{-6}{0}
\boxt{2}{-6}{0}
\boxt{1}{-6}{0}
\boxt{0}{-6}{0}
\boxt{-1}{-6}{0}
\boxt{-2}{-6}{0}
\boxt{-3}{-6}{0}
\boxt{-4}{-6}{0}
\boxt{-5}{-6}{0}
\boxt{-6}{-6}{0}
\boxt{-7}{-6}{0}
\boxt{-8}{-6}{0}

\boxl{8}{-7}{0}
\boxt{8}{-7}{0}
\boxt{7}{-7}{0}
\boxt{6}{-7}{0}
\boxt{5}{-7}{0}
\boxt{4}{-7}{0}
\boxt{3}{-7}{0}
\boxt{2}{-7}{0}
\boxt{1}{-7}{0}
\boxt{0}{-7}{0}
\boxt{-1}{-7}{0}
\boxt{-2}{-7}{0}
\boxt{-3}{-7}{0}
\boxt{-4}{-7}{0}
\boxt{-5}{-7}{0}
\boxt{-6}{-7}{0}
\boxt{-7}{-7}{0}
\boxt{-8}{-7}{0}

\boxl{8}{-8}{0}
\boxt{8}{-8}{0}
\boxt{7}{-8}{0}
\boxt{6}{-8}{0}
\boxt{5}{-8}{0}
\boxt{4}{-8}{0}
\boxt{3}{-8}{0}
\boxt{2}{-8}{0}
\boxt{1}{-8}{0}
\boxt{0}{-8}{0}
\boxt{-1}{-8}{0}
\boxt{-2}{-8}{0}
\boxt{-3}{-8}{0}
\boxt{-4}{-8}{0}
\boxt{-5}{-8}{0}
\boxt{-6}{-8}{0}
\boxt{-7}{-8}{0}
\boxt{-8}{-8}{0}

\boxl{8}{-9}{0}
\boxt{8}{-9}{0}
\boxt{7}{-9}{0}
\boxt{6}{-9}{0}
\boxt{5}{-9}{0}
\boxt{4}{-9}{0}
\boxt{3}{-9}{0}
\boxt{2}{-9}{0}
\boxt{1}{-9}{0}
\boxt{0}{-9}{0}
\boxt{-1}{-9}{0}
\boxt{-2}{-9}{0}
\boxt{-3}{-9}{0}
\boxt{-4}{-9}{0}
\boxt{-5}{-9}{0}
\boxt{-6}{-9}{0}
\boxt{-7}{-9}{0}
\boxt{-8}{-9}{0}

\boxl{-1}{0}{0}
\boxt{-1}{0}{0}
\boxt{-2}{0}{0}
\boxt{-2}{0}{0}
\boxt{-3}{0}{0}
\boxt{-4}{0}{0}
\boxt{-5}{0}{0}
\boxt{-6}{0}{0}
\boxt{-7}{0}{0}
\boxt{-8}{0}{0}

\boxl{-1}{1}{0}
\boxt{-1}{1}{0}
\boxt{-2}{1}{0}
\boxt{-2}{1}{0}
\boxt{-3}{1}{0}
\boxt{-4}{1}{0}
\boxt{-5}{1}{0}
\boxt{-6}{1}{0}
\boxt{-7}{1}{0}
\boxt{-8}{1}{0}

\boxl{-1}{2}{0}
\boxt{-1}{2}{0}
\boxt{-2}{2}{0}
\boxt{-2}{2}{0}
\boxt{-3}{2}{0}
\boxt{-4}{2}{0}
\boxt{-5}{2}{0}
\boxt{-6}{2}{0}
\boxt{-7}{2}{0}
\boxt{-8}{2}{0}

\boxl{-1}{3}{0}
\boxt{-1}{3}{0}
\boxt{-2}{3}{0}
\boxt{-2}{3}{0}
\boxt{-3}{3}{0}
\boxt{-4}{3}{0}
\boxt{-5}{3}{0}
\boxt{-6}{3}{0}
\boxt{-7}{3}{0}
\boxt{-8}{3}{0}

\boxl{-1}{4}{0}
\boxt{-1}{4}{0}
\boxt{-2}{4}{0}
\boxt{-2}{4}{0}
\boxt{-3}{4}{0}
\boxt{-4}{4}{0}
\boxt{-5}{4}{0}
\boxt{-6}{4}{0}
\boxt{-7}{4}{0}
\boxt{-8}{4}{0}

\boxl{-1}{5}{0}
\boxt{-1}{5}{0}
\boxt{-2}{5}{0}
\boxt{-2}{5}{0}
\boxt{-3}{5}{0}
\boxt{-4}{5}{0}
\boxt{-5}{5}{0}
\boxt{-6}{5}{0}
\boxt{-7}{5}{0}
\boxt{-8}{5}{0}

\boxl{-1}{6}{0}
\boxt{-1}{6}{0}
\boxt{-2}{6}{0}
\boxt{-2}{6}{0}
\boxt{-3}{6}{0}
\boxt{-4}{6}{0}
\boxt{-5}{6}{0}
\boxt{-6}{6}{0}
\boxt{-7}{6}{0}
\boxt{-8}{6}{0}

\boxl{-1}{7}{0}
\boxt{-1}{7}{0}
\boxt{-2}{7}{0}
\boxt{-2}{7}{0}
\boxt{-3}{7}{0}
\boxt{-4}{7}{0}
\boxt{-5}{7}{0}
\boxt{-6}{7}{0}
\boxt{-7}{7}{0}
\boxt{-8}{7}{0}

\boxl{-1}{8}{0}
\boxt{-1}{8}{0}
\boxt{-2}{8}{0}
\boxt{-2}{8}{0}
\boxt{-3}{8}{0}
\boxt{-4}{8}{0}
\boxt{-5}{8}{0}
\boxt{-6}{8}{0}
\boxt{-7}{8}{0}
\boxt{-8}{8}{0}

\boxr{0}{-1}{0}
\boxr{1}{-1}{0}
\boxr{2}{-1}{0}
\boxr{3}{-1}{0}
\boxr{4}{-1}{0}
\boxr{5}{-1}{0}
\boxr{6}{-1}{0}
\boxr{7}{-1}{0}
\boxr{8}{-1}{0}

\boxr{-1}{8}{0}
\boxr{-2}{8}{0}
\boxr{-3}{8}{0}
\boxr{-4}{8}{0}
\boxr{-5}{8}{0}
\boxr{-6}{8}{0}
\boxr{-7}{8}{0}
\boxr{-8}{8}{0}

\boxxt{-1}{-1}{0}
\boxxt{0}{-1}{0}
\boxxt{1}{-1}{0}
\boxxt{2}{-1}{0}
\boxxt{3}{-1}{0}
\boxxt{4}{-1}{0}
\boxxt{5}{-1}{0}
\boxxt{6}{-1}{0}
\boxxt{7}{-1}{0}
\boxxt{8}{-1}{0}

\boxxt{-1}{-2}{0}
\boxxt{0}{-2}{0}
\boxxt{1}{-2}{0}
\boxxt{2}{-2}{0}
\boxxt{3}{-2}{0}
\boxxt{4}{-2}{0}
\boxxt{5}{-2}{0}
\boxxt{6}{-2}{0}
\boxxt{7}{-2}{0}
\boxxt{8}{-2}{0}

\boxxt{-1}{-3}{0}
\boxxt{0}{-3}{0}
\boxxt{1}{-3}{0}
\boxxt{2}{-3}{0}
\boxxt{3}{-3}{0}
\boxxt{4}{-3}{0}
\boxxt{5}{-3}{0}
\boxxt{6}{-3}{0}
\boxxt{7}{-3}{0}
\boxxt{8}{-3}{0}

\boxxt{-1}{-4}{0}
\boxxt{0}{-4}{0}
\boxxt{1}{-4}{0}
\boxxt{2}{-4}{0}
\boxxt{3}{-4}{0}
\boxxt{4}{-4}{0}
\boxxt{5}{-4}{0}
\boxxt{6}{-4}{0}
\boxxt{7}{-4}{0}
\boxxt{8}{-4}{0}

\boxxt{-1}{0}{0}
\boxxt{-1}{1}{0}
\boxxt{-1}{2}{0}
\boxxt{-1}{3}{0}
\boxxt{-1}{4}{0}
\boxxt{-1}{5}{0}
\boxxt{-1}{6}{0}
\boxxt{-1}{7}{0}
\boxxt{-1}{8}{0}

\boxxt{-2}{-4}{0}
\boxxt{-2}{-3}{0}
\boxxt{-2}{-2}{0}
\boxxt{-2}{-1}{0}
\boxxt{-2}{0}{0}
\boxxt{-2}{1}{0}
\boxxt{-2}{2}{0}
\boxxt{-2}{3}{0}
\boxxt{-2}{4}{0}
\boxxt{-2}{5}{0}
\boxxt{-2}{6}{0}
\boxxt{-2}{7}{0}
\boxxt{-2}{8}{0}

\boxxr{0}{-1}{0}
\boxxr{1}{-1}{0}
\boxxr{2}{-1}{0}
\boxxr{3}{-1}{0}
\boxxr{4}{-1}{0}
\boxxr{5}{-1}{0}
\boxxr{6}{-1}{0}
\boxxr{7}{-1}{0}
\boxxr{8}{-1}{0}

\boxxr{-1}{8}{0}
\boxxr{-2}{8}{0}

\boxxl{8}{-1}{0}
\boxxl{8}{-2}{0}
\boxxl{8}{-3}{0}
\boxxl{8}{-4}{0}

\boxxl{-1}{0}{0}
\boxxl{-1}{1}{0}
\boxxl{-1}{2}{0}
\boxxl{-1}{3}{0}
\boxxl{-1}{4}{0}
\boxxl{-1}{5}{0}
\boxxl{-1}{6}{0}
\boxxl{-1}{7}{0}
\boxxl{-1}{8}{0}

\boxxxt{-2}{-4}{0}
\boxxxt{-2}{-3}{0}
\boxxxt{-2}{-2}{0}

\boxxxt{-1}{-4}{0}
\boxxxt{-1}{-3}{0}

\boxxxt{0}{-4}{0}
\boxxxt{1}{-4}{0}

\draw [thick,->] (0,0)-- (0,13);
\draw [thick,->] (0*1.74,-0*1.74*0.57735)-- (11,-11*0.57735);
\draw [thick,->]  (-0,-0*0.57735)-- (-11,-11*0.57735);

\node [right] at (0,13) {$x_3$};
\node [right] at (11,-11*0.57735) {$x_2$};
\node [left] at  (-11,-11*0.57735) {$x_1$};
\node [above] at  (-9,-7*0.57735) {$ $};
\node [right] at (2,10) {$ $};
\node  at (6,6*0.57735) {$ $};
\node  at (-6,6*0.57735) {$ $};
\node  at (0,-6) {$ $};

%\shade [ball color=gray] (-2.6,-.6) circle [radius=0.4cm];
%\shade [ball color=gray] (1.8,-1.2) circle [radius=0.4cm];
%\shade [ball color=gray] (-0.85,-1.7) circle [radius=0.4cm];

\end{tikzpicture}

%\end{document} 

%% file: fig5b.tex
% Plane partition
% Author: Jang Soo Kim
%\documentclass{minimal}
%\usepackage{tikz}
% Three counters
%\newcounter{x}
%\newcounter{y}
%\newcounter{z}

% The angles of x,y,z-axes
\newcommand\xaxis{210}
\newcommand\yaxis{-30}
\newcommand\zaxis{90}

% The top side of a cube
\newcommand\topside[3]{
  \fill[fill=cyan, draw=black,shift={(\xaxis:#1)},shift={(\yaxis:#2)},
  shift={(\zaxis:#3)}] (0,0) -- (30:1) -- (0,1) --(150:1)--(0,0);
}

% The left side of a cube
\newcommand\leftside[3]{
  \fill[fill=cyan, draw=black,shift={(\xaxis:#1)},shift={(\yaxis:#2)},
  shift={(\zaxis:#3)}] (0,0) -- (0,-1) -- (210:1) --(150:1)--(0,0);
}

% The right side of a cube
\newcommand\rightside[3]{
  \fill[fill=cyan, draw=black,shift={(\xaxis:#1)},shift={(\yaxis:#2)},
  shift={(\zaxis:#3)}] (0,0) -- (30:1) -- (-30:1) --(0,-1)--(0,0);
}

% The cube 
\newcommand\cube[3]{
  \topside{#1}{#2}{#3} \leftside{#1}{#2}{#3} \rightside{#1}{#2}{#3}
}

% Definition of \planepartition
% To draw the following plane partition, just write \planepartition{ {a, b, c}, {d,e} }.
%  a b c
%  d e
\newcommand\planepartition[1]{
 \setcounter{x}{-1}
  \foreach \a in {#1} {
    \addtocounter{x}{1}
    \setcounter{y}{-1}
    \foreach \b in \a {
      \addtocounter{y}{1}
      \setcounter{z}{-1}
      \foreach \c in {1,...,\b} {
        \addtocounter{z}{1}
        \cube{\value{x}}{\value{y}}{\value{z}}
      }
    }
  }
}

%\begin{document} 

\begin{tikzpicture}[scale=0.22]

\newcommand{\boxt}[3]{ 
\draw[fill=cyan, draw=black,shift={(\xaxis:#1)},shift={(\yaxis:#2)},
  shift={(\zaxis:#3)}] (0,0) -- (30:1) -- (0,1) --(150:1)--(0,0);
}

\newcommand{\boxl}[3]{ 
\draw[fill=cyan, draw=black,shift={(\xaxis:#1)},shift={(\yaxis:#2)},
  shift={(\zaxis:#3)}] (0,0) -- (0,-1) -- (210:1) --(150:1)--(0,0);
}

\newcommand{\boxr}[3]{ 
\draw[fill=cyan, draw=black,shift={(\xaxis:#1)},shift={(\yaxis:#2)},
  shift={(\zaxis:#3)}] (0,0) -- (30:1) -- (-30:1) --(0,-1)--(0,0);
}

\newcommand{\boxxt}[3]{ 
\draw[fill=red, draw=black,shift={(\xaxis:#1)},shift={(\yaxis:#2)},
  shift={(\zaxis:#3)}] (0,0) -- (30:1) -- (0,1) --(150:1)--(0,0);
}

\newcommand{\boxxl}[3]{ 
\draw[fill=red, draw=black,shift={(\xaxis:#1)},shift={(\yaxis:#2)},
  shift={(\zaxis:#3)}] (0,0) -- (0,-1) -- (210:1) --(150:1)--(0,0);
}

\newcommand{\boxxr}[3]{ 
\draw[fill=red, draw=black,shift={(\xaxis:#1)},shift={(\yaxis:#2)},
  shift={(\zaxis:#3)}] (0,0) -- (30:1) -- (-30:1) --(0,-1)--(0,0);
}

\newcommand{\boxxxt}[3]{ 
\draw[fill=green, draw=black,shift={(\xaxis:#1)},shift={(\yaxis:#2)},
  shift={(\zaxis:#3)}] (0,0) -- (30:1) -- (0,1) --(150:1)--(0,0);
}

\newcommand{\boxxxl}[3]{ 
\draw[fill=green, draw=black,shift={(\xaxis:#1)},shift={(\yaxis:#2)},
  shift={(\zaxis:#3)}] (0,0) -- (0,-1) -- (210:1) --(150:1)--(0,0);
}

\newcommand{\boxxxr}[3]{ 
\draw[fill=green, draw=black,shift={(\xaxis:#1)},shift={(\yaxis:#2)},
  shift={(\zaxis:#3)}] (0,0) -- (30:1) -- (-30:1) --(0,-1)--(0,0);
}

\boxl{8}{-1}{0}
\boxt{8}{-1}{0}
\boxt{7}{-1}{0}
\boxt{6}{-1}{0}
\boxt{5}{-1}{0}
\boxt{4}{-1}{0}
\boxt{3}{-1}{0}
\boxt{2}{-1}{0}
\boxt{1}{-1}{0}
\boxt{0}{-1}{0}
\boxt{-1}{-1}{0}
\boxt{-2}{-1}{0}
\boxt{-3}{-1}{0}
\boxt{-4}{-1}{0}
\boxt{-5}{-1}{0}
\boxt{-6}{-1}{0}
\boxt{-7}{-1}{0}
\boxt{-8}{-1}{0}

\boxl{8}{-2}{0}
\boxt{8}{-2}{0}
\boxt{7}{-2}{0}
\boxt{6}{-2}{0}
\boxt{5}{-2}{0}
\boxt{4}{-2}{0}
\boxt{3}{-2}{0}
\boxt{2}{-2}{0}
\boxt{1}{-2}{0}
\boxt{0}{-2}{0}
\boxt{-1}{-2}{0}
\boxt{-2}{-2}{0}
\boxt{-3}{-2}{0}
\boxt{-4}{-2}{0}
\boxt{-5}{-2}{0}
\boxt{-6}{-2}{0}
\boxt{-7}{-2}{0}
\boxt{-8}{-2}{0}

\boxl{8}{-3}{0}
\boxt{8}{-3}{0}
\boxt{7}{-3}{0}
\boxt{6}{-3}{0}
\boxt{5}{-3}{0}
\boxt{4}{-3}{0}
\boxt{3}{-3}{0}
\boxt{2}{-3}{0}
\boxt{1}{-3}{0}
\boxt{0}{-3}{0}
\boxt{-1}{-3}{0}
\boxt{-2}{-3}{0}
\boxt{-3}{-3}{0}
\boxt{-4}{-3}{0}
\boxt{-5}{-3}{0}
\boxt{-6}{-3}{0}
\boxt{-7}{-3}{0}
\boxt{-8}{-3}{0}

\boxl{8}{-4}{0}
\boxt{8}{-4}{0}
\boxt{7}{-4}{0}
\boxt{6}{-4}{0}
\boxt{5}{-4}{0}
\boxt{4}{-4}{0}
\boxt{3}{-4}{0}
\boxt{2}{-4}{0}
\boxt{1}{-4}{0}
\boxt{0}{-4}{0}
\boxt{-1}{-4}{0}
\boxt{-2}{-4}{0}
\boxt{-3}{-4}{0}
\boxt{-4}{-4}{0}
\boxt{-5}{-4}{0}
\boxt{-6}{-4}{0}
\boxt{-7}{-4}{0}
\boxt{-8}{-4}{0}

\boxl{8}{-5}{0}
\boxt{8}{-5}{0}
\boxt{7}{-5}{0}
\boxt{6}{-5}{0}
\boxt{5}{-5}{0}
\boxt{4}{-5}{0}
\boxt{3}{-5}{0}
\boxt{2}{-5}{0}
\boxt{1}{-5}{0}
\boxt{0}{-5}{0}
\boxt{-1}{-5}{0}
\boxt{-2}{-5}{0}
\boxt{-3}{-5}{0}
\boxt{-4}{-5}{0}
\boxt{-5}{-5}{0}
\boxt{-6}{-5}{0}
\boxt{-7}{-5}{0}
\boxt{-8}{-5}{0}

\boxl{8}{-6}{0}
\boxt{8}{-6}{0}
\boxt{7}{-6}{0}
\boxt{6}{-6}{0}
\boxt{5}{-6}{0}
\boxt{4}{-6}{0}
\boxt{3}{-6}{0}
\boxt{2}{-6}{0}
\boxt{1}{-6}{0}
\boxt{0}{-6}{0}
\boxt{-1}{-6}{0}
\boxt{-2}{-6}{0}
\boxt{-3}{-6}{0}
\boxt{-4}{-6}{0}
\boxt{-5}{-6}{0}
\boxt{-6}{-6}{0}
\boxt{-7}{-6}{0}
\boxt{-8}{-6}{0}

\boxl{8}{-7}{0}
\boxt{8}{-7}{0}
\boxt{7}{-7}{0}
\boxt{6}{-7}{0}
\boxt{5}{-7}{0}
\boxt{4}{-7}{0}
\boxt{3}{-7}{0}
\boxt{2}{-7}{0}
\boxt{1}{-7}{0}
\boxt{0}{-7}{0}
\boxt{-1}{-7}{0}
\boxt{-2}{-7}{0}
\boxt{-3}{-7}{0}
\boxt{-4}{-7}{0}
\boxt{-5}{-7}{0}
\boxt{-6}{-7}{0}
\boxt{-7}{-7}{0}
\boxt{-8}{-7}{0}

\boxl{8}{-8}{0}
\boxt{8}{-8}{0}
\boxt{7}{-8}{0}
\boxt{6}{-8}{0}
\boxt{5}{-8}{0}
\boxt{4}{-8}{0}
\boxt{3}{-8}{0}
\boxt{2}{-8}{0}
\boxt{1}{-8}{0}
\boxt{0}{-8}{0}
\boxt{-1}{-8}{0}
\boxt{-2}{-8}{0}
\boxt{-3}{-8}{0}
\boxt{-4}{-8}{0}
\boxt{-5}{-8}{0}
\boxt{-6}{-8}{0}
\boxt{-7}{-8}{0}
\boxt{-8}{-8}{0}

\boxl{8}{-9}{0}
\boxt{8}{-9}{0}
\boxt{7}{-9}{0}
\boxt{6}{-9}{0}
\boxt{5}{-9}{0}
\boxt{4}{-9}{0}
\boxt{3}{-9}{0}
\boxt{2}{-9}{0}
\boxt{1}{-9}{0}
\boxt{0}{-9}{0}
\boxt{-1}{-9}{0}
\boxt{-2}{-9}{0}
\boxt{-3}{-9}{0}
\boxt{-4}{-9}{0}
\boxt{-5}{-9}{0}
\boxt{-6}{-9}{0}
\boxt{-7}{-9}{0}
\boxt{-8}{-9}{0}

\boxl{-1}{0}{0}
\boxt{-1}{0}{0}
\boxt{-2}{0}{0}
\boxt{-2}{0}{0}
\boxt{-3}{0}{0}
\boxt{-4}{0}{0}
\boxt{-5}{0}{0}
\boxt{-6}{0}{0}
\boxt{-7}{0}{0}
\boxt{-8}{0}{0}

\boxl{-1}{1}{0}
\boxt{-1}{1}{0}
\boxt{-2}{1}{0}
\boxt{-2}{1}{0}
\boxt{-3}{1}{0}
\boxt{-4}{1}{0}
\boxt{-5}{1}{0}
\boxt{-6}{1}{0}
\boxt{-7}{1}{0}
\boxt{-8}{1}{0}

\boxl{-1}{2}{0}
\boxt{-1}{2}{0}
\boxt{-2}{2}{0}
\boxt{-2}{2}{0}
\boxt{-3}{2}{0}
\boxt{-4}{2}{0}
\boxt{-5}{2}{0}
\boxt{-6}{2}{0}
\boxt{-7}{2}{0}
\boxt{-8}{2}{0}

\boxl{-1}{3}{0}
\boxt{-1}{3}{0}
\boxt{-2}{3}{0}
\boxt{-2}{3}{0}
\boxt{-3}{3}{0}
\boxt{-4}{3}{0}
\boxt{-5}{3}{0}
\boxt{-6}{3}{0}
\boxt{-7}{3}{0}
\boxt{-8}{3}{0}

\boxl{-1}{4}{0}
\boxt{-1}{4}{0}
\boxt{-2}{4}{0}
\boxt{-2}{4}{0}
\boxt{-3}{4}{0}
\boxt{-4}{4}{0}
\boxt{-5}{4}{0}
\boxt{-6}{4}{0}
\boxt{-7}{4}{0}
\boxt{-8}{4}{0}

\boxl{-1}{5}{0}
\boxt{-1}{5}{0}
\boxt{-2}{5}{0}
\boxt{-2}{5}{0}
\boxt{-3}{5}{0}
\boxt{-4}{5}{0}
\boxt{-5}{5}{0}
\boxt{-6}{5}{0}
\boxt{-7}{5}{0}
\boxt{-8}{5}{0}

\boxl{-1}{6}{0}
\boxt{-1}{6}{0}
\boxt{-2}{6}{0}
\boxt{-2}{6}{0}
\boxt{-3}{6}{0}
\boxt{-4}{6}{0}
\boxt{-5}{6}{0}
\boxt{-6}{6}{0}
\boxt{-7}{6}{0}
\boxt{-8}{6}{0}

\boxl{-1}{7}{0}
\boxt{-1}{7}{0}
\boxt{-2}{7}{0}
\boxt{-2}{7}{0}
\boxt{-3}{7}{0}
\boxt{-4}{7}{0}
\boxt{-5}{7}{0}
\boxt{-6}{7}{0}
\boxt{-7}{7}{0}
\boxt{-8}{7}{0}

\boxl{-1}{8}{0}
\boxt{-1}{8}{0}
\boxt{-2}{8}{0}
\boxt{-2}{8}{0}
\boxt{-3}{8}{0}
\boxt{-4}{8}{0}
\boxt{-5}{8}{0}
\boxt{-6}{8}{0}
\boxt{-7}{8}{0}
\boxt{-8}{8}{0}

\boxr{0}{-1}{0}
\boxr{1}{-1}{0}
\boxr{2}{-1}{0}
\boxr{3}{-1}{0}
\boxr{4}{-1}{0}
\boxr{5}{-1}{0}
\boxr{6}{-1}{0}
\boxr{7}{-1}{0}
\boxr{8}{-1}{0}

\boxr{-1}{8}{0}
\boxr{-2}{8}{0}
\boxr{-3}{8}{0}
\boxr{-4}{8}{0}
\boxr{-5}{8}{0}
\boxr{-6}{8}{0}
\boxr{-7}{8}{0}
\boxr{-8}{8}{0}

\boxxt{-1}{-1}{0}
\boxxt{0}{-1}{0}
\boxxt{1}{-1}{0}
\boxxt{2}{-1}{0}
\boxxt{3}{-1}{0}
\boxxt{4}{-1}{0}
\boxxt{5}{-1}{0}
\boxxt{6}{-1}{0}
\boxxt{7}{-1}{0}
\boxxt{8}{-1}{0}

\boxxt{-1}{-2}{0}
\boxxt{0}{-2}{0}
\boxxt{1}{-2}{0}
\boxxt{2}{-2}{0}
\boxxt{3}{-2}{0}
\boxxt{4}{-2}{0}
\boxxt{5}{-2}{0}
\boxxt{6}{-2}{0}
\boxxt{7}{-2}{0}
\boxxt{8}{-2}{0}

\boxxt{-1}{-3}{0}
\boxxt{0}{-3}{0}
\boxxt{1}{-3}{0}
\boxxt{2}{-3}{0}
\boxxt{3}{-3}{0}
\boxxt{4}{-3}{0}
\boxxt{5}{-3}{0}
\boxxt{6}{-3}{0}
\boxxt{7}{-3}{0}
\boxxt{8}{-3}{0}

\boxxt{-1}{-4}{0}
\boxxt{0}{-4}{0}
\boxxt{1}{-4}{0}
\boxxt{2}{-4}{0}
\boxxt{3}{-4}{0}
\boxxt{4}{-4}{0}
\boxxt{5}{-4}{0}
\boxxt{6}{-4}{0}
\boxxt{7}{-4}{0}
\boxxt{8}{-4}{0}

\boxxt{-1}{0}{0}
\boxxt{-1}{1}{0}
\boxxt{-1}{2}{0}
\boxxt{-1}{3}{0}
\boxxt{-1}{4}{0}
\boxxt{-1}{5}{0}
\boxxt{-1}{6}{0}
\boxxt{-1}{7}{0}
\boxxt{-1}{8}{0}

\boxxt{-2}{-4}{0}
\boxxt{-2}{-3}{0}
\boxxt{-2}{-2}{0}
\boxxt{-2}{-1}{0}
\boxxt{-2}{0}{0}
\boxxt{-2}{1}{0}
\boxxt{-2}{2}{0}
\boxxt{-2}{3}{0}
\boxxt{-2}{4}{0}
\boxxt{-2}{5}{0}
\boxxt{-2}{6}{0}
\boxxt{-2}{7}{0}
\boxxt{-2}{8}{0}

\boxxr{0}{-1}{0}
\boxxr{1}{-1}{0}
\boxxr{2}{-1}{0}
\boxxr{3}{-1}{0}
\boxxr{4}{-1}{0}
\boxxr{5}{-1}{0}
\boxxr{6}{-1}{0}
\boxxr{7}{-1}{0}
\boxxr{8}{-1}{0}

\boxxr{-1}{8}{0}
\boxxr{-2}{8}{0}

\boxxl{8}{-1}{0}
\boxxl{8}{-2}{0}
\boxxl{8}{-3}{0}
\boxxl{8}{-4}{0}

\boxxl{-1}{0}{0}
\boxxl{-1}{1}{0}
\boxxl{-1}{2}{0}
\boxxl{-1}{3}{0}
\boxxl{-1}{4}{0}
\boxxl{-1}{5}{0}
\boxxl{-1}{6}{0}
\boxxl{-1}{7}{0}
\boxxl{-1}{8}{0}

\boxxxt{-2}{-4}{0}
\boxxxt{-2}{-3}{0}
\boxxxt{-2}{-2}{0}

\boxxxt{-1}{-4}{0}
\boxxxt{-1}{-3}{0}

\boxxxt{0}{-4}{0}
\boxxxt{1}{-4}{0}

\draw [thick,->] (0,0)-- (0,13);
\draw [thick,->] (0*1.74,-0*1.74*0.57735)-- (11,-11*0.57735);
\draw [thick,->]  (-0,-0*0.57735)-- (-11,-11*0.57735);

\node [right] at (0,13) {$x_3$};
\node [right] at (11,-11*0.57735) {$x_2$};
\node [left] at  (-11,-11*0.57735) {$x_1$};
\node [above] at  (-9,-7*0.57735) {$ $};
\node [right] at (2,10) {$ $};
\node  at (6,6*0.57735) {$ $};
\node  at (-6,6*0.57735) {$ $};
\node  at (0,-6) {$ $};

%\shade [ball color=gray] (-2.6,-.6) circle [radius=0.4cm];
%\shade [ball color=gray] (1.8,-1.2) circle [radius=0.4cm];
%\shade [ball color=gray] (-0.85,-1.7) circle [radius=0.4cm];

\end{tikzpicture}

%\end{document} 